\documentclass{amsart}

\usepackage{amsmath,amsthm}
\usepackage {latexsym}
\usepackage{amssymb}

\newcommand{\lambdainf}{\lambda_{\infty}}
\newcommand{\tilxi}{\tilde{\xi}}
\newcommand{\tileta}{\tilde{\eta}}
\newcommand{\calT}{\mathcal{T}}
\newcommand{\calC}{\mathcal{C}}

\newcommand{\beeq}{\begin{equation}}
\newcommand{\eneq}{\end{equation}}

\newcommand{\bear}{\begin{eqnarray}}
\newcommand{\eear}{\end{eqnarray}}
\newcommand{\beq}{\begin{equation}}
\newcommand{\eeq}{\end{equation}}

\newcommand{\half}{\frac{1}{2}}

\newcommand{\R}{{\mathbb R}}

\newcommand{\calA}{{\mathcal A}}
\newcommand{\calg}{\,{\mathfrak g}}

\newcommand{\calF}{{\mathcal F}}
\newcommand{\calL}{{\mathcal L}}
\newcommand{\calS}{{\mathcal S}}

\newcommand{\calM}{{\mathcal M}}

\newcommand{\calH}{{\mathcal H}}

\newcommand{\calN}{{\mathcal N}}

\newcommand{\la}{\langle}
\newcommand{\ra}{\rangle}

\def\nn{\nonumber}
\def\calge1{\calg_{\vec{e_1}}}

\def\bm{\left( \begin{array}{cc}}
\def\endm{\end{array}\right)}

\def\Hil{{\mathcal H}}

\newtheorem{theorem}{Theorem}[section]
\newtheorem{lemma}[theorem]{Lemma}
\newtheorem{defi}[theorem]{Definition}

\newtheorem{proposition}[theorem]{Proposition}
\newtheorem{conj}[theorem]{Conjecture}

\theoremstyle{remark}

\def\norm[#1][#2]{\Vert #1 \Vert_{#2}}

\def\Util3{\tilde{U}^{(3)}}

\def\xitil{\tilde{\xi}}

\def\calW{{\mathcal W}}

\numberwithin{equation}{section}

\begin{document}

\title{Non-generic blow-up solutions for the critical focusing NLS in
1-d.}
\author{J.\ Krieger, W. Schlag}
\thanks{The first author was partially supported by NSF grant DMS-0401177 and the second author by the NSF grant
DMS-0300081 and a Sloan fellowship.}
\address{Harvard University, Dept. of Mathematics, Science Center,
1 Oxford Street, Cambridge, MA 02138, USA} \email{
jkrieger@math.harvard.edu}
\address{Department of Mathematics, University of Chicago, 5734
South University Avenue, Chicago, Il 60637, USA}
\email{schlag@math.uchicago.edu}
\date{}
\maketitle

\section{Introduction}

We consider the critical focusing NLS in 1-d of the form
\begin{equation}\label{TheEquation}
i\partial_{t}\psi+\partial_{x}^{2}\psi=-|\psi|^{4}\psi,\,i=\sqrt{-1},\,\psi=\psi(t,x),
\end{equation}
and $\psi$ complex valued. It is well-known that this equation
permits standing wave solutions of the form
\begin{equation}\nonumber
\phi(t,x)=e^{i\alpha t} \phi_{0}(x,\alpha),\,\alpha>0
\end{equation}
Indeed, requiring positivity and evenness in $x$ for
$\phi_{0}(x,\alpha)$ implies for example
\begin{equation}\nonumber
\phi_{0}(x,\alpha)=\frac{\alpha^{\frac{1}{2}}(\frac{3}{2})^{\frac{1}{4}}}{\cosh^{\frac{1}{2}}(\frac{\alpha}{2}
x)}
\end{equation}
Another remarkable feature of the equation \eqref{TheEquation} is
the large symmetry group carrying solutions into solutions: this
is generated by
\\

{\bf{Galilei transformations:}}
\begin{equation}\nonumber
\psi(t,x)\longrightarrow
e^{i(\gamma+vx-v^{2}t)}e^{-i(2tv+\mu)p}\psi(t,x)
=e^{i(\gamma+vx-v^{2}t)}\psi(t,x-2tv-\mu),\,p=-i\frac{d}{dx}
\end{equation}
\\

{\bf{$SL(2,\mathbf{R})$-transformations:}}
\begin{equation}\nonumber
\psi(t,x)\longrightarrow
(a+bt)^{-\frac{1}{2}}e^{\frac{ibx^{2}}{4(a+bt)}}\psi(\frac{c+dt}{a+bt},\frac{x}{a+bt}),\,\bm
a&b\\c&d\endm\in SL(2,{\mathbf{R}})
\end{equation}

Observe that the latter subsume re-scalings $\psi(t,x)\rightarrow
a^{\frac{1}{2}}\psi(a^{2}t,ax)$ while the former subsume
phase-shifts $\psi(t,x)\rightarrow e^{i\gamma}\psi(t,x)$ as well
as translations. We usually identify a matrix $\bm
a&b\\c&d\endm\in SL(2,{\mathbf{R}})$ with the corresponding
transformation. It is the $SL(2,{\mathbf{R}})$-transformations
that distinguish the critical NLS from the sub- and supercritical
NLS, and allows us to exhibit explicit blow-up solutions: indeed,
fixing $\bm a& b\\c&d\endm\in SL(2,{\mathbf{R}})$, we have the
explicit solution
\begin{equation}\label{explicit}
f(t,x)=(a+bt)^{-\frac{1}{2}}e^{i\frac{c+dt}{a+bt}}e^{\frac{ibx^{2}}{4(a+bt)}}\phi_{0}(\frac{x}{a+bt},1),
\end{equation}
which blows up for $t=-\frac{a}{b}$. Fixing $a\sim 1$, $b\sim -1$,
it is then a natural question to ask whether one may perturb the
initial data of \eqref{explicit} at time $t=0$ such that the
corresponding solution exhibits the same type of blow-up behavior.
More precisely, the solution should asymptotically behave like
$\sqrt{\frac{1}{T-t}}e^{i\Psi(t,x)}\phi(\frac{x-\mu(t)}{T-t})$ for
a bounded function $\mu(t)$ and suitable Schwartz function $\phi$,
with blow up time $T$.  The recent work of Merle-Raphael [MeRa]
has demonstrated that this is generically impossible, i. e. there
are open sets of initial data containing $f(0,x)$ in their
closure\footnote{With respect to any reasonable norm.} and such
that their blow-up behavior is of the following type, which we
henceforth refer to as 'generic':
\begin{equation}\nonumber
\psi(t,x)\sim
e^{i\mu(t)}\lambda^{\frac{1}{2}}(t)\phi(\lambda(t)x),\,\lambda(t)\sim
\sqrt{\frac{\log|\log(T-t)|}{T-t}}
\end{equation}
Blow-up solutions of this type were first constructed in a
remarkable paper by G. Perelman [Per2], but for non-generic
initial data sets. This blow-up rate was shown to be stable in
[Ra]. Moreover, in [MeRa] the authors showed that for initial data
in a sufficiently small neighborhood of $f(0,x)$ the only possible
blow-up speeds are the generic speed or else at least as fast as
the explicit speed;  we now refer to the latter as 'non-generic'.
The issue remains as to whether perturbations of the initial data
$f(0,x)$ in certain directions would result in the non-generic
blow-up type. The first and to our knowledge only result of this
type was established by Bourgain-Wang [BW]\footnote{The authors
state this Theorem for the case of $d=1,2$ dimensions.}, and
asserts the following:

\begin{theorem}\label{B-W}[Bourgain-Wang] Let $z_{\phi}$ be the local-in-time solution of
\begin{equation}\nonumber
i\psi_{t}+\triangle\psi+|\psi|^{4}\psi=0,\,\psi(0)=\phi,
\end{equation}
which for smooth $\phi$ exists on an interval $[-\delta,\delta]$
for $\delta=\delta(\phi)$ small enough. Then provided $\phi$ is
smooth and vanishes sufficiently fast at $0$, i.e.
$|\phi(x)|\lesssim |x|^{A}$ for $A$ large enough\footnote{The
numerology in [BW] appears to imply $A\geq 16$}, there exists
smooth $w(t,x)$ in a suitable function space with $w(0,x)=0$ and
such that
\begin{equation}\label{staticcoupling}
\psi(t,x)=t^{-\frac{1}{2}}e^{\frac{x^{2}-4}{4it}}\phi_{0}(\frac{x}{t},1)+z_{\phi}(t,x)+w(t,x)
\end{equation}
solves \eqref{TheEquation} on $[-\delta,0]$. One may let
$\delta\rightarrow \infty$ by letting $\phi\rightarrow 0$ .
\end{theorem}
The key behind this result is to first undo the blow-up by
applying a pseudo-conformal transformation $C^{-1}$ where
$C\psi(t,x)=t^{-\frac{1}{2}}e^{\frac{x^{2}}{4it}}\psi(\frac{x}{t},-\frac{1}{t})$
and then employ the properties of the linear evolution associated
with the linearization around the standing wave
$e^{it}\phi_{0}(x)$. More precisely, one passes to the vector
valued function $\bm \psi(t,x)\\\overline{\psi(t,x)}\endm$ and
observes that if $\psi(t,x)=e^{it}(\phi_{0}(x,1)+u(t,x))$ solves
\eqref{TheEquation} then we have
\begin{equation}\nonumber
(i\partial_{t}+\calH )\bm u(t,x)\\ \overline{u(t,x)}\endm = N(u)
\end{equation}
where we put
\begin{equation}\label{H}
\calH=\bm \partial_{x}^{2}-1+3\phi_{0}^{4}(x,1)& 2\phi_{0}^{4}(x,1)\\
-2\phi_{0}^{4}(x,1)&-\partial_{x}^{2}+1-3\phi_{0}^{4}(x,1)\endm
\end{equation}
and $N(u)$ is of order $\geq 2$ in $u$. The spectral properties of
the operator $\calH$ are well-known after the pioneering work of
Weinstein [Wei1] as well as Buslaev-Perelman [BusPer] and Perelman
[Per2]. In
particular, the linear equation $(i\partial_{t}+\calH)\bm u\\
\bar{u}\endm =0$ only displays algebraic instabilities. More
precisely, the spectrum of $\calH$ has essential part
$(-\infty,-1]\cup [1,\infty)$ and discrete spectrum $\{0\}$ of
geometric multiplicity $2$ and algebraic multiplicity $6$. A
solution $\Phi(x)=\bm \phi(x)\\ \bar{\phi}(x)\endm$ in the
generalized root space satisfies for example
$||e^{it\calH}\Phi(x)||_{L_{x}^{2}}\lesssim
(1+t^{3})\int_{-\infty}^{\infty}e^{-c|x|}|\Phi(x)|dx$. In order to
counteract this growth behavior at infinity, Bourgain and Wang use
the ansatz $u(t,x)=C^{-1}z_{\phi}e^{-it}+\tilde{w}(t,x)$, see
\eqref{staticcoupling}, and then observe that the
non-linearity of the resulting equation for $\bm \tilde{w}\\
\overline{\tilde{w}}\endm$ decays sufficiently rapidly at infinity
(due to the local decay of $C^{-1}z_{\phi}e^{-it}$) that it
overwhelms any losses due to the algebraic instability of $\calH$.
The fact that the 'static coupling' \eqref{staticcoupling} barely
exploits the symmetries of the equation and in particular doesn't
allow the standing wave to 'drift' certainly implies the
sub-optimality\footnote{In the sense that the set of initial data
resulting in the non-generic blow-up should be significantly
larger than indicated there.} of
Theorem~\ref{B-W}.\\
Indeed, a careful analysis of the root-space of $\calH$ as in
[Wei1] shows that $5$ of the generalized root modes (the 'good
modes') 'are due to' the internal symmetries of
\eqref{TheEquation}, while there is one 'exotic mode', see the
ensuing discussion. This intimates that upon applying suitable
internal symmetries to the standing wave $e^{it}\phi_{0}(x,1)$ in
time-dependent fashion (i. e. using a modulation-theoretic
approach), one should be able to control the root part of the
radiation corresponding to the good modes, and indeed obtain a
co-dimension 1 stable manifold of initial data (due to the 'exotic
mode' which cannot be so controlled) resulting in the non-generic
blow-up profile:
\\

\begin{conj}\label{MainConjecture} [Galina Perelman] There
exists a co-dimension 1 manifold of initial data resulting in the
non-generic blow-up behavior.
\end{conj}
We note that this is also implicitly mentioned although in less
precise form in [B].\\
This falls in neatly with recent results in [Sch] and [KriSch1],
the latter closely following the former, which in the context of
the $L^{2}$-super-critical NLS (the cubic in 3-d in [Sch] and the
full super-critical range in 1-d in [KriSch1]) established
existence of co-dimension 1 manifolds of initial data resulting in
globally (for $t\rightarrow+\infty$) smooth solutions. The
co-dimension 1 here has to do with one exponentially unstable mode
(in the forward time direction); the generalized root space has
only dimension 4, in one-one correspondence with the internal
symmetries.
\\

In this paper we attempt to make some progress toward
Conjecture~\ref{MainConjecture}: let\footnote{We identify the
matrix with its associated transformation.}
\begin{equation}\nonumber
\calT_{\infty}=e^{-i(v_{\infty}^{2}s+\gamma_{\infty}+v_{\infty}y)}e^{i(2v_{\infty}s+y_{\infty})p}\bm
a_{\infty}& b_{\infty}\\0&a_{\infty}^{-1}\endm
\end{equation}
Also, let the generalized root space of $\calH$ be generated by
the (vector-valued) Schwartz functions $\eta_{i,\text{proper}}=\bm
\eta_{i,\text{proper}}^{1}\\\overline{\eta_{i,\text{proper}}^{1}}\endm$,
$i=1,\ldots,6$, and that of $\calH^{*}$ be generated by the
Schwartz functions $\xi_{i,\text{proper}}(x),\,i=1,\ldots,6$, viz.
the ensuing discussion.
\\

Our main result is the following
\begin{theorem}\label{Main} Fix real parameters $\lambda\sim 1, \beta\sim 1,
\omega\lesssim 1,\gamma\lesssim 1,\mu\lesssim 1$. Given a vector
valued function $x\longrightarrow \bm U\\\bar{U}\endm(0,x)$
satisfying $\la \bm U\\\bar{U}\endm(0,.), \xi_{i,\text{proper}}\ra
=0 \forall i$, as well as the smallness condition
$|||U(0,.)|||<\delta$ for a suitable norm\footnote{see
definition~\ref{Abounds}} $|||.|||$ and sufficiently small
$\delta>0$, there exist numbers
$\tilde{\lambda}_{i}\in{\mathbf{R}}$ with
$|\tilde{\lambda}_{i}|\lesssim |||U(0,.)|||^{2}$ and parameters
$\{a_{\infty},b_{\infty},v_{\infty},y_{\infty},\gamma_{\infty}\}$
with $|a_{\infty}-\lambda|\lesssim |||U(0,.)|||^{2}$,
$|b_{\infty}-\beta\lambda|\lesssim |||U(0,.)|||^{2}$,
$|v_{\infty}-\frac{\beta\lambda\mu}{2}-\omega|\lesssim
|||U(0,.)|||^{2}$, $|y_{\infty}-\lambda\mu|\lesssim
|||U(0,.)|||^{2}$, $\gamma_{\infty}=\gamma_{\infty}(,\gamma,
\lambda,\beta,\omega,\mu)+O(|||U(0,x)|||^{2})$, such that the
initial data
\begin{equation}\nonumber
\psi(0,x):=W(0,x)+\calT_{\infty}^{-1}[U(0,x)+\sum_{i=1}^{6}\tilde{\lambda}_{i}\eta_{i,\text{proper}}^{1}]
\end{equation}
lead to solutions of \eqref{TheEquation} blowing up in finite time
according to the non-generic profile, where
\begin{equation}\nonumber
W(0,x)=e^{i(\gamma+\omega(x-\mu))}e^{-i\frac{\beta}{4}\lambda^{2}(x-\mu)^{2}}\sqrt{\lambda}\phi_{0}(\lambda(x-\mu),1)
\end{equation}
More precisely, the solution decouples as
\begin{equation}\nonumber
\psi(t,x)=e^{i(\gamma(t)+\omega(t)(x-\mu(t)))}e^{-i\frac{\beta(t)}{4}\lambda^{2}(t)(x-\mu(t))^{2}}\sqrt{\lambda(t)}\phi_{0}(\lambda(t)(x-\mu(t)),1)+
R(t,x),
\end{equation}
where $\lambda(t)\sim \frac{1}{T-t}$ for a suitable $T>0$, and
$\mu(t)$ is bounded, while we have the bounds
\begin{equation}\nonumber
\sup_{0\leq t<T}||R(t,x)||_{L_{t}^{\infty}}\lesssim
\delta,\,||R(t,x)||_{H^{1}}\lesssim \delta(T-t)^{-1},\,0\leq t<T
\end{equation}
In particular, $||\psi(t,x)||_{H^{1}}\sim (T-t)^{-1}$ for  $0\leq
t<T$.
\end{theorem}

{\bf{Remark:}} We observe that this result would imply
Conjecture~\ref{MainConjecture} if we could show Lipschitz
continuous dependence of the $\tilde{\lambda}_{i}$, $a_{\infty}$
etc on $U(0,x)$, along the lines in [KriSch1]. However, we cannot
establish this. Indeed, even demonstrating the possibility or
impossibility of {\it{choosing these parameters in continuous
fashion}} appears extremely difficult.
\\

We now outline the strategy used to prove this Theorem: there are
the following four stages:
\\

{\bf{Stage A: Setting up the equations for radiation part and
modulation parameters.}}
\\

Instead of the static coupling \eqref{staticcoupling}, we make the
ansatz
\begin{equation}\label{dynamiccoupling}
\psi(t,x)=W(t,x)+R(t,x),\,W(t,x)=e^{i(\theta(t,x))}e^{-i\frac{\beta}{4}(t)\lambda^{2}(x-\mu(t))^{2}}\sqrt{\lambda(t)}
\phi_{0}(\lambda(t)(x-\mu(t)),1)
\end{equation}
In order to ensure that this solution behaves like a non-generic
blow-up solution, we impose the condition $\lambda(t)\sim
\frac{1}{t_{*}-t}$ for suitable $t_{*}\in{\mathbf{R}}_{>0}$. We
shall similarly have to carefully specify the 'asymptotic
behavior' of the remaining parameters as we approach blow-up time.
In order to specify the evolution of these parameters, we impose
suitable orthogonality conditions: letting
$\xi_{i,\text{proper}}$,\,$i=1,\ldots,6$ denote a certain basis of
the generalized root space of $\calH^{*}$ (recall \eqref{H}) to be
specified below, we impose\footnote{The root functions
$\xi_{i,\text{proper}}$, $i=1,\ldots,5$ here are chosen to be the
'good modes' in one-one relation with the internal symmetries,
while the root function $\xi_{6,\text{proper}}$ is the 'exotic
mode' due to the degeneracy in the critical case.}
\begin{equation}\label{orthogonality}
\la \bm R\\ \bar{R}\endm, \xi_{i}\ra=0,\,\xi_{i}=\bm
e^{i\Psi(t,x)}&0\\0&e^{-i\Psi(t,x)}\endm
\sqrt{\lambda(t)}\xi_{i,\text{proper}}(\lambda(t)(x-\mu(t))),i=2,\ldots,6
\end{equation}
where we have introduced the notation
\begin{equation}\nonumber
\Psi(t,x)=\theta(t,x)-\frac{\beta(t)}{4}\lambda^{2}(t)(x-\mu(t))^{2}
\end{equation}
For later reference , we define $\eta_{i}$ correspondingly, with
$\xi_{i,\text{proper}}$ replaced by $\eta_{i,\text{proper}}$. This
is analogous to the procedure in [Sch], [KriSch1], where the
generalized root space is only 4-dimensional. The above
orthogonality condition then implies that at time $t$ the
radiation part {\it{when projected onto the generalized root space
of the instantaneous linearization around the drifting soliton}}
gives zero. Note that the fact the we no longer work with a static
standing wave forces us to work with modifications of the operator
$\calH$.\\

Instead of working with the formulation \eqref{dynamiccoupling},
we then revert to a 'different Gauge' as in [BW]. Specifically, we
apply a suitable transformation $\calT_{\infty}$,
\begin{equation}\label{Tinfty}
\calT_{\infty}=e^{-i(v_{\infty}^{2}s+\gamma_{\infty}+v_{\infty}x)}e^{i(2v_{\infty}s+y_{\infty})p}\bm
a_{\infty}& b_{\infty}\\ 0&
a_{\infty}^{-1}\endm,\,p=-i\frac{d}{dx},
\end{equation}
to $\psi(t,x)$ which is to undo the singular behavior and should
map the blow-up time $t_{*}$ to $s=\infty$. To see how the
'coefficients at infinity' $a_{\infty}$ etc. should be chosen, we
observe that
\begin{equation}\nonumber
(\calT_{\infty}F)(s,y)=e^{-i\tilde{\Psi}(s,y)}\lambda_{\infty}^{-\frac{1}{2}}(s)F(\int_{0}^{s}\lambda_{\infty}^{-2}(\sigma)d\sigma,\lambda_{\infty}^{-1}(s)y+\mu_{\infty}(s))
\end{equation}
where
\begin{equation}\nonumber
\lambda_{\infty}(s)=a_{\infty}+b_{\infty}s,\,\mu_{\infty}(s)=\frac{2v_{\infty}s+y_{\infty}}{a_{\infty}+b_{\infty}s},
\tilde{\Psi}_{\infty}(s,y)=v_{\infty}^{2}s+\gamma_{\infty}+v_{\infty}y-\frac{b_{\infty}(y+2v_{\infty}s+y_{\infty})^{2}}{4(a_{\infty}+b_{\infty}s)},
\end{equation}
and therefore
\begin{equation}\nonumber\begin{split}
&e^{-is}(\calT_{\infty}W)(s,y)\\&=e^{-i(s+\tilde{\Psi}_{\infty}(s,y))+i\Psi(t(s),\mu_{\infty}(s)+\lambda_{\infty}^{-1}(s)y)}\lambda_{\infty}^{-\frac{1}{2}}(s)\lambda^{\frac{1}{2}}(t(s))\phi_{0}(\lambda(t(s))(\mu_{\infty}(s)-\mu(t(s))+\lambda_{\infty}^{-1}(s)y),1),
\end{split}\end{equation}
where we have put
$t(s)=\int_{0}^{s}\lambda_{\infty}^{-2}(\sigma)d\sigma=\frac{a_{\infty}^{-1}s}{a_{\infty}+b_{\infty}s}$.
The above suggests that we should impose $\lambda(t(s))\sim
\lambda_{\infty}(s)$, $\mu_{\infty}(s)\sim \mu(t(s))$ as
$s\rightarrow \infty$ in a precise sense to be specified. In
particular, we have $t_{*}=\frac{1}{a_{\infty}b_{\infty}}$ for the
blow-up time.
\\

We shall now work with the vector valued function
\begin{equation}\nonumber
\bm U\\ \bar{U}\endm:=\calM \calT_{\infty}\bm R\\
\bar{R}\endm,\,\calM=\bm e^{-is}&0\\0&e^{is}\endm
\end{equation}
Then introduce the functions
$\tilde{\eta}_{i}=\calM\calT_{\infty}\eta_{i}$,
$\tilde{\xi}_{i}=\calM\calT_{\infty}\xi_{i}$. One deduces the
following equation for $\bm U\\ \bar{U}\endm$:
\begin{align}\label{Ueqn}
i\partial_{s}\bm U\\ \bar{U}\endm+\calH(s)\bm U\\ \bar{U}\endm& =-i(\dot\lambda\lambda^{-1}-\beta\nu^2)(\tileta_2-\beta\tileta_5/2+\omega\tileta_4)\\
&\quad +\frac{i}{4}(\dot\beta+\beta^2\nu^2)\tileta_5 + i(\nu^2-\dot\gamma+\nu^2\omega^2)\tileta_1 \\
&\quad -i(\dot\omega+\beta\omega\nu^2)\tileta_4 -
i\nu(\dot\mu\lambda_\infty-2\nu\omega)(-\omega\tileta_1-\tileta_3+\beta\tileta_4/2)
+ N(U,\pi),
\end{align}
where we use $\nu(s)=\frac{\lambda(t(s))}{\lambda_{\infty}(s)}$,
$\lambda(t(s))=\lambda(s)$,
$\dot{\lambda}=\frac{\partial}{\partial s}\lambda(s)$, and
$N(U,\pi)$ is quadratic in $U$ but also depends on the modulation
parameters $\lambda(s)$ etc., as well as the parameters at
infinity $a_{\infty}$ etc. We denote the latter collectively as
$\pi$, following the notation in [Sch], [KriSch1]. The operator
$\calH(s)$ in the preceding is given by
\begin{align*}
\calH(s):=\bm
\partial_{y}^{2}-1+3\nu^{2}(s)\phi_{0}^{4}(\lambda(\mu_{\infty}-\mu+\lambda_{\infty}^{-1}y))&
2\nu^{2}\phi_{0}^{4}(\lambda(\mu_{\infty}-\mu+\lambda_{\infty}^{-1}y))e^{2i(\Psi-\Psi_{\infty})}\\
-2\nu^{2}\phi_{0}^{4}(\lambda(\mu_{\infty}-\mu+\lambda_{\infty}^{-1}y))e^{-2i(\Psi-\Psi_{\infty})}&
-\partial_{y}^{2}+1-3\nu^{2}(s)\phi_{0}^{4}(\lambda(\mu_{\infty}-\mu+\lambda_{\infty}^{-1}y))\endm,
\end{align*}
where we use the notation
$\Psi_{\infty}(s,y):=\tilde{\Psi}_{\infty}(s,y)+s$. The
orthogonality relations \eqref{orthogonality} become the
following:
\begin{equation}\label{orthogonality'}
\la \bm U\\ \bar{U}\endm, \tilde{\xi}_{i}\ra =0,\,i=2,\ldots,6,
\end{equation}
and upon differentiating with respect to $s$ imply a set of ODE's
for the parameters $\lambda(s)$ etc. The crux now is to deduce a
priori estimates for the transformed radiation part $\bm U\\
\bar{U}\endm$ as well as for the modulation parameters; the latter
need to satisfy the required asymptotic estimates for
$s\rightarrow\infty$. In order to control the radiation part, one
essentially\footnote{For technical reasons, one uses such a
decomposition for a slightly transformed function $\bm \tilde{U}\\
\overline{\tilde{U}}\endm$.} invokes a decomposition
\begin{equation}\nonumber
\bm U\\
\bar{U}\endm(s,y)=\bm U\\
\bar{U}\endm_{dis}(s,y)+\sum_{i=1}^{6}\lambda_{i}(s)\eta_{i,\text{proper}}(y),
\end{equation}
where the coefficients $\lambda_{i}(s)$, $i=1,\ldots,5$, are
determined by the orthogonality relations \eqref{orthogonality'}.
The coefficient $\lambda_{6}(s)$ is determined by means of the
requirement $\lim_{s\rightarrow\infty}\lambda_{6}(s)=0$, which
forces an initial condition $\lambda_{6}(0)$, similarly to the
super-critical case treated in [Sch], [KriSch1]. By comparison to
the latter, though, controlling $\lambda_{6}(s)$ appears more
difficult, and requires the development of rather new technology.
Specifically, a careful analysis of the modulation equations
reveals that one needs to control quantities of the form
$\int_{T}^{\infty}t\lambda_{6}(t)dt$, which upon substituting the
solution for the ODE satisfied by $\lambda_{6}(t)$ results in
quadratic expressions of at worst the form\footnote{We are again
careless here; the expression should really involve
$\tilde{U}^{2}-\overline{\tilde{U}}^{2}$.}
$\int_{T}^{\infty}t\int_{t}^{\infty}\la U^{2}(s)-\bar{U}^{2}(s),
\phi\ra ds dt$ etc., where $\phi$ stands for a suitable Schwartz
function. This shows that one should aim for a local decay of the
radiation part of at least the strength $|\la U(t,.),\phi\ra
|\lesssim \la t\ra ^{-\frac{3}{2}}$ in order to be able to
estimate this expression; indeed, this local decay rate is in
accordance with the linear estimates derived in [KriSch1].
However, we are dealing with a nonlinear problem here. This is the
first significant difficulty to be overcome:
\\

{\bf{Stage B: Deducing the strong local dispersive\footnote{More
precisely, we establish the strong local dispersive estimate up to
an arbitrarily small error.} estimate for the radiation part.}}
\\

Schematically, the equation \eqref{Ueqn} etc can be recast as
\begin{equation}\nonumber
(i\partial_{s}+\calH)\bm U\\ \bar{U}\endm = VU +\bm |U|^{4}U \\
-|U|^{4}\bar{U}\endm,
\end{equation}
where $\calH=\bm \partial_{y}^{2}-1+3\phi_{0}^{4}(.,1)& 2\phi_{0}^{4}(.,1)\\
-2\phi_{0}^{4}(.,1)& -\partial_{y}^{2}+1-3\phi_{0}^{4}(.,1)\endm$,
and $V$, a Schwartz function, depends on $U$ as well as the
modulation parameters etc. The local\footnote{We refer to
expressions which are Schwartz functions alternatively as
'local'.} expression $VU$ is due to interactions of the drifting
soliton with itself as well as to interactions of the radiation
with the drifting soliton, while the non-local quintilinear
expressions $|U|^{4}U$ come from interactions of the radiation
part with itself. While the root part of $\bm U\\ \bar{U}\endm $
is controlled in terms of the coefficients $\lambda_{i}(s)$,
$i=1,\ldots,6$, whose estimation is relegated to the third stage,
the dispersive part\footnote{Again,
one should really use $\bm \tilde{U}\\
\overline{\tilde{U}}\endm_{dis}$.}
(viz. the next section for the linear background) $\bm U\\
\bar{U}\endm_{dis}$ satisfies
\begin{equation}\nonumber
(i\partial_{s}+\calH)\bm U\\ \bar{U}\endm_{dis} = [VU +\bm |U|^{4}U \\
-|U|^{4}\bar{U}\endm]_{dis},
\end{equation}
The really difficult contribution on the right comes from the
non-local quintilinear term: note that the standard way to deduce
the local estimate is to combine the linear estimate\footnote{We
abuse notation here and use letters $\phi,\psi$ to denote vector
valued functions. Also, we let $\la x\ra=|x|+1$.}
\begin{equation}\nonumber
|\la e^{it\calH}\phi_{dis},\psi\ra| \lesssim t^{-\frac{3}{2}}||\la
x\ra\phi||_{L_{x}^{1}}||\la x\ra\psi||_{L_{x}^{1}}
\end{equation}
with Duhamel's formula, which then forces us to estimate the
expression
\begin{equation}\label{1}
\int_{0}^{t}\la t-s\ra^{-\frac{3}{2}} ||\la x\ra
|U|^{4}(s,x)U(s,x)||_{L_{x}^{1}}ds
\end{equation}
On the other hand, again from the linear theory summarized in the
next section we expect the estimate
\begin{equation}\nonumber
 ||\la x\ra
|U|^{4}(s,x)U(s,x)||_{L_{x}^{1}}\lesssim
||U(s,.)||_{L_{x}^{\infty}}^{3}||xU(s,x)||_{L_{x}^{2}}||U(s,x)||_{L_{x}^{2}}\lesssim
\la s\ra^{-\frac{3}{2}}s=s^{-\frac{1}{2}},
\end{equation}
which only gives the decay $t^{-\frac{1}{2}}$ when substituted
into \eqref{1}. One can modify this argument to eke out a local
dispersive decay of $\la t\ra ^{-1+}$, which however is
insufficient for controlling the root part and modulation
parameters.
\\

The way out of this is to observe that the quintilinear expression
exhibits a {\it{special algebraic cancellation structure}}, which
in combination with the linear theory of $\calH$(and in particular
the absence of resonances at the edges of the essential spectrum)
allows one to significantly improve on the preceding. To explain
the use of this algebraic structure heuristically, note that one
expects the
small-frequency part of $\bm |U|^{4}U\\
-|U|^{4}\bar{U}\endm_{dis}$ to contribute less due to the absence
of resonances at the end of the essential spectrum of $\calH$.
Another reason is that the small frequency part propagates more
slowly, and hence when hit with a weight $\la x\ra $ should cost
less than the $s$ used in the above calculation. On the other
hand, assume that we localize one of the factors $|U|^{2}$ in
$|U|^{4}=|U|^{2}|U|^{2}$ to relatively large frequency\footnote{By
this we mean here frequency in the Littlewood-Paley sense. One has
to be a bit careful to keep this separate from frequency in the
sense of $\calH$. The relation of the two will become clear thanks
to the explicit distorted Fourier basis explained in the next
section; the general heuristic is that a function with small
frequency with respect to $\calH$ is the sum of a (negligible)
Schwartz function plus a function of small frequency in the
Littlewood-Paley sense.}. In that case the key is to use the
following simple identity:
\begin{equation}\label{null-form1}
2is\partial_{x}[|U|^{2}](s,.)=(x+2is\partial_{x})U(s,.)
\bar{U}(s,.)-U(s,.)\overline{(x+2is\partial_{x})U(s,.)}
\end{equation}
The operator $C=(x+2is\partial_{x})$ is the standard
pseudo-conformal operator, and one expects an estimate
$\sup_{s\geq 0}||(x+2is\partial_{x})U(s,.)||_{L_{x}^{2}}\lesssim
1$. Indeed, this turns out to be true (although establishing it
requires a couple of tricks, as we don't deal with the free
evolution $e^{i t\triangle}$ here.) Thus provided we restrict the
frequency of $|U|^{2}(s,.)$ sufficiently far away from $0$, we
expect to be able to score an extra gain here, and one can play
these two considerations against each other to almost obtain the
optimal estimate. This very crudely summarizes the strategy for
stage {\bf{B}}.
\\

{\bf{Stage C: Controlling the root part of $\bm U\\ \bar{U}\endm$
and the modulation parameters.}}

We now return to controlling $\lambda_{6}(t)$ as well as the
modulation parameters, which we recall involved estimating
expressions such as $\int_{T}^{\infty}t\int_{t}^{\infty}\la
U^{2}(s)-\bar{U}^{2}(s), \phi\ra ds dt$, as well as similar ones.
Clearly even the strong local dispersive estimate isn't good
enough for this purpose, and we have to resort to more refined
considerations. In case of the displayed expression, this involves
identifying another instance of an algebraic cancellation
structure, in this case a {\it{symplectic structure}}. Again this
shall rely on the spectral properties of $\calH$.
\\

{\bf{Stage D: Locating a fixed point.}}

The a priori estimates suggest running a Banach iteration;
unfortunately, the presence of the phase
$e^{i(\Psi-\Psi_{\infty})(s,y)}$ with $(\Psi-\Psi_{\infty})(s,y)$
growing like $s^{\frac{1}{2}+}$ doesn't allow one to deduce good
estimates for the differences of iterates. This is the fundamental
obstacle to proving Conjecture~\ref{MainConjecture}.  We thus have
to resort to an abstract fixed point Theorem (Schauder-Tychonoff
Theorem) to prove Theorem~\ref{Main}. It is to be hoped that the
techniques developed in this paper help to further elucidate the
nature of the non-generic blow up solutions.
\\

{\bf{Acknowledgements:}} The authors would like to thank Jean
Bourgain and Carlos Kenig for helpful discussions and their
interest in our work. The first author would also like to express
his gratitude to the California Institute of Technology and the
University of Chicago for generously hosting him during the Summer
2004 and Summer 2005, respectively.

\section{Background material on $\calH$.}

The material in this section quickly summarizes certain facts
established in the last section of [KriSch1], much of which was
based on the work of Buslaev-Perelman and Perelman as well as
earlier work by the 2nd author. We refer to [KriSch1] as well as
[Per2] for proofs. Consider the operator
\begin{equation}\nonumber
\calH=\bm \partial_{x}^{2}-1+3\phi_{0}^{4}(x,1)& 2\phi_{0}^{4}(x,1)\\
-2\phi_{0}^{4}(x,1)& -\partial_{x}^{2}+1-3\phi_{0}^{4}(x,1)\endm
\end{equation}
The spectrum consists of $(-\infty,-1]\cup [1,\infty)\cup\{0\}$,
with essential spectrum $(-\infty,-1]\cup [1,\infty)$ and discrete
spectrum $\{0\}$ of geometric multiplicity $2$ and algebraic
multiplicity $6$.  The generalized root space $\calN$ is generated
by the following vector valued functions: from now on, we adhere
to the convention $\phi_{0}:=\phi_{0}(.,1)$, see the preceding
section.
\begin{align*}
\eta_{1,\text{proper}}(z) &:= \binom{i \phi_0(z)}{-i\phi_0(z)},
\quad \eta_{2,\text{proper}}(z) := \binom{(z \phi'_0(z)+\phi_0(z)/2)}{(z \phi'_0(z)+\phi_0(z)/2)}\\
\eta_{3,\text{proper}}(z) &:= \binom{\phi'_0(z)}{ \phi'_0(z)},
\quad \eta_{4,\text{proper}}(z):= \binom{i z \phi_0(z)}{-iz \phi_0(z)}\\
\eta_{5,\text{proper}}(z)&:= \binom{i z^2\phi_0(z)}{-i z^2
\phi_0(z)}, \quad \eta_{6,\text{proper}}(z) := \binom{
\rho(z)}{\rho(z)}
\end{align*}
The first five are in one-one correspondence with internal
symmetries('good modes'), while the last is the 'exotic mode',
characterized by
\begin{equation}\nonumber
L_{+}\rho=z^{2}\phi_{0}(z),\,L_{+}=-\partial_{x}^{2}+1-5\phi_{0}^{4}
\end{equation}
The root space is generated by $\eta_{1,\text{proper}}$,
$\eta_{3,\text{proper}}$. \\
As for the essential spectrum, its edges $\pm 1$ {\it{are no
resonances}}. This means that there are no solutions
$f_{\pm}(z)\in L^{\infty}$ satisfying
\begin{equation}\nonumber
\calH f_{\pm}=\pm f_{\pm}
\end{equation}
This is in marked contrast to the operator $\calH_{0}:=\bm
\partial_{x}^{2}-1&0\\0&-\partial_{x}^{2}+1\endm$, and responsible
for much improved local decay estimates. \\
Identical observations apply to the operator $\calH^{*}=\bm \partial_{x}^{2}-1+3\phi_{0}^{4}(x,1)& -2\phi_{0}^{4}(x,1)\\
+2\phi_{0}^{4}(x,1)& -\partial_{x}^{2}+1-3\phi_{0}^{4}(x,1)\endm$,
its generalized root space $\calN^{*}$ being generated by
\begin{align}\label{xiproper}
\xi_{1,\text{proper}}(z) &:= \binom{\phi_0(z)}{\phi_0(z)},
\quad \xi_{2,\text{proper}}(z) := \binom{i(z \phi'_0(z)+\phi_0(z)/2)}{-i(z \phi'_0(z)+\phi_0(z)/2)}\\
\xi_{3,\text{proper}}(z) &:= \binom{i\phi'_0(z)}{ -i\phi'_0(z)},
\quad \xi_{4,\text{proper}}(z):= \binom{z \phi_0(z)}{z \phi_0(z)}\\
\xi_{5,\text{proper}}(z)&:= \binom{z^2\phi_0(z)}{z^2 \phi_0(z)},
\quad \xi_{6,\text{proper}}(z) := \binom{ i\rho(z)}{-i\rho(z)}
\end{align}
Then we have the direct sum decomposition
\begin{equation}\nonumber
L^{2}({\mathbf{R}})\times L^{2}({\mathbf{R}})=\calN
+(\calN^{*})^{\perp}
\end{equation}
This means that every vector valued function $\bm U\\
\bar{U}\endm(x) $ with $U(.)\in L^{2}({\mathbf{R}})$ can be
uniquely represented as
\begin{equation}\nonumber
\bm U\\
\bar{U}\endm=\sum_{i=1}^{6}\lambda_{i}\eta_{i,\text{proper}}+\bm
U\\ \bar{U}\endm_{dis},\,\bm U\\ \bar{U}\endm_{dis}\in
(\calN^{*})^{\perp}
\end{equation}
In order to determine the $\lambda_{i}$, one uses the following
table of orthogonality relations\footnote{We use the convention
$\la \bm U_{1}\\ U_{2}\endm, \bm V_{1}\\ V_{2}\endm\ra =\la U_{1},
V_{1}\ra +\la U_{2}, V_{2}\ra$.}:
\begin{align}\label{orthorelations}
\langle \eta_{j, \text{proper}}, \xi_{1,\text{proper}} \rangle &= 0, \; 1\le j\le 5, \;\langle \eta_{6,\text{proper}}, \xi_{1, \text{proper}} \rangle = 2\kappa_2 \\
\langle \eta_{j,\text{proper}},\xi_{2,\text{proper}} \rangle &= 0,\; j=1,2,3,4,6,\; \quad \langle \eta_{5,\text{proper}},\xi_{2,\text{proper}}\rangle = -4\kappa_2 \\
\langle \eta_{j,\text{proper}},\xi_{3,\text{proper}} \rangle &= 0,\; j=1,2,3,5,6,\; \quad \langle \eta_{4,\text{proper}},\xi_{3,\text{proper}}\rangle = -\kappa_1 \\
\langle \eta_{j,\text{proper}},\xi_{4,\text{proper}} \rangle &= 0,\; j=1,2,4,5,6,\; \quad \langle \eta_{3,\text{proper}},\xi_{4,\text{proper}}\rangle = -\kappa_1 \\
\langle \eta_{j,\text{proper}},\xi_{5,\text{proper}} \rangle &=
0,\; j=1,3,4,5,\; \quad \langle
\eta_{2,\text{proper}},\xi_{5,\text{proper}}\rangle = -4\kappa_2,
\; \quad \langle \eta_{6,\text{proper}},\xi_{5,\text{proper}}\rangle = 2\kappa_3 \\
\langle \eta_{j,\text{proper}},\xi_{6,\text{proper}} \rangle &=
0,\; j=2,3,4,6,\; \quad \langle
\eta_{1,\text{proper}},\xi_{6,\text{proper}}\rangle = 2\kappa_2,\;
\quad \langle \eta_{5,\text{proper}},\xi_{6,\text{proper}}\rangle
= 2\kappa_3,
\end{align}
where we use the notation
\[ \langle \phi_0,\phi_0\rangle = \kappa_1>0, \qquad \langle \rho,\phi_0\rangle = \half\int x^2\phi_0^2(x)\, dx=:\kappa_2>0,
\qquad \langle x^2\phi_0,\rho\rangle=:\kappa_3,
\]
We also write
\begin{equation}\nonumber
\sum_{i=1}^{6}\lambda_{i}\eta_{i,\text{proper}}=P_{root}\bm U\\
\bar{U}\endm,\,\bm U\\ \bar{U}\endm_{dis}=P_{s}\bm U\\
\bar{U}\endm
\end{equation}
We have the following important linear estimates:
\begin{theorem}\label{linearestimates} The following estimates hold for vector valued
functions $f(.)\in  L^{1}({\mathbf{R}})\cap  L^{2}({\mathbf{R}})$
and $0\leq\theta\leq 1$:
\begin{equation}\nonumber
||e^{it\calH}P_{s}f||_{L^{2}}\lesssim
||f||_{L^{2}},\,||e^{it\calH}P_{s}f||_{L^{\infty}}\lesssim
|t|^{-\frac{1}{2}}||f||_{L^{1}},\, ||\la
x\ra^{-\theta}e^{it\calH}P_{s}f||_{L^{\infty}}\lesssim
|t|^{-\frac{1}{2}-\theta|}||\la x\ra f||_{L_{x}^{1}}
\end{equation}
\end{theorem}

The first two of these are just as for $e^{it\triangle}$, while
the last is not true for the latter and due to the absence of
resonances at the edges of the essential spectrum of $\calH$.
\\

By analogy to Fourier transformation associated with $\triangle$,
there is a {\it{distorted Fourier transform}} associated with
$\calH$:
\begin{theorem}\label{DistortedFourier}
There exist Lipschitz continuous vector valued functions
$e_{\pm}(x,\xi)$ with the property
\begin{equation}\nonumber
P_{s}\bm U\\ \bar{U}\endm(x) =
\sum_{\pm}\int_{-\infty}^{\infty}e_{\pm}(x,\xi)\la \bm U\\
\bar{U}\endm, \sigma_{3}e_{\pm}(x,\xi)\ra d\xi,\,\sigma_{3}=\bm
1&0\\0&-1\endm
\end{equation}
for every rapidly\footnote{We are being overly restrictive in the formulation here; all facts about the distorted Fourier transform
apply in the same degree of generality as for the ordinary Fourier transform.} decaying function $x\rightarrow \bm U\\
\bar{U}\endm(x)$. Moreover, there exist smooth functions $s(\xi)$,
$r(\xi)$ satisfying $s(0)=0$, $r(0)=-1$, as well as suitable
numbers $\gamma>0,\mu>0$, such that
\begin{equation}\nonumber
e_{+}(x,\xi)=s(\xi)[e^{ix\xi}\bm 1\\0\endm
+O((1+|\xi|)^{-1}e^{-\gamma x})]+O(\xi(1+|\xi|)^{-2}e^{-\mu
x}),\,(x,\xi)\in {\mathbf{R}}_{\geq 0}\times  {\mathbf{R}}_{\geq
0}
\end{equation}
\begin{equation}\nonumber
e_{+}(x,\xi)=[e^{ix\xi}+r(\xi)e^{-ix\xi}]\bm 1\\0\endm+O(\xi
(1+|\xi|)^{-2}e^{\gamma x}),\,(x,\xi)\in {\mathbf{R}}_{<0}\times
{\mathbf{R}}_{\geq 0}
\end{equation}
\begin{equation}\nonumber
e_{+}(x,\xi)=s(-\xi)[e^{ix\xi}\bm 1\\0\endm
+O((1+|\xi|)^{-1}e^{+\gamma x})]+O(\xi(1+|\xi|)^{-2}e^{+\mu
x}),\,(x,\xi)\in {\mathbf{R}}_{< 0}\times  {\mathbf{R}}_{<0}
\end{equation}
\begin{equation}\nonumber
e_{+}(x,\xi)=[e^{ix\xi}+r(-\xi)e^{-ix\xi}]\bm 1\\0\endm+O(\xi
(1+|\xi|)^{-2}e^{-\gamma x}),\,(x,\xi)\in {\mathbf{R}}_{\geq
0}\times {\mathbf{R}}_{<0}
\end{equation}
Also, we have the relation
\begin{equation}\nonumber
e_{-}(x,\xi)=\sigma_{1}e_{+}(x,\xi),\,\sigma_{1}=\bm 0&1\\1&0\endm
\end{equation}
\end{theorem}

In analogy to the usual Fourier transform, there is a
{\it{distorted Plancherel's Theorem}}:
\begin{theorem}\label{DistortedPlancherel}
Let $\phi,\psi\in {\mathcal{S}}({\mathbf{R}})$ be vector valued
functions. Then we have
\begin{equation}\nonumber
\la P_{s}\phi, \psi\ra
=\sum_{\pm}\frac{1}{2\pi}\int_{-\infty}^{\infty}\la \phi,
\sigma_{3}e_{\pm}(.,\xi)\ra \overline{\la \psi,
e_{\pm}(.,\xi)\ra}d\xi
\end{equation}
\end{theorem}
We shall use the notation
\begin{equation}\nonumber
\calF_{\pm}(\phi)(\xi):=\la \phi,
\sigma_{3}e_{\pm}(.,\xi)\ra,\,\tilde{\calF}_{\pm}(\phi)(\xi):=\la
\phi, e_{\pm}(x,\xi)\ra
\end{equation}

When working with $e_{\pm}(x,\xi)$, we shall for example write
\begin{equation}\nonumber
e_{+}(x,\xi)=s(\xi)e^{ix\xi}\underline{e}+\phi(x,\xi),\,(x,\xi)\in
{\mathbf{R}}_{\geq 0}\times {\mathbf{R}}_{\geq 0},
\end{equation}
and similarly for the other values of $(x,\xi)$. The functions
$\phi(x,\xi)$, which are Schwartz with respect to $x$, are
understood to vary accordingly, but all vanish uniformly in $x$ at
$\xi=0$ and decay like $|\xi|^{-1}$ for $|\xi|\rightarrow \infty$
uniformly in $x$.

\section{Setting up the equations}
\subsection{Algebraic manipulations I; analysis of the modulation parameters.}

We now flesh out the discussion of the first section. In this
section as well as the next, we shall use formal algebraic
manipulations to derive the equations which will serve to define
the iterative step. We shall also mention the required estimates.
In the final sections of the paper, we shall then show that the
iterative step indeed makes sense when performed on suitable
function spaces. Thus consider now a solution $\psi(t,x)=
W(t,x)+R(t,x)$ of \eqref{TheEquation}, where
\begin{equation}\nonumber
W(t,x)=e^{i(\theta(t,x))}e^{-i\frac{\beta}{4}(t)\lambda^{2}(x-\mu(t))^{2}}\sqrt{\lambda(t)}
\phi_{0}(\lambda(t)(x-\mu(t)))
\end{equation}
Use the notation
$\Psi(t,x)=\theta(t,x)-\frac{\beta}{4}\lambda^{2}(x-\mu)^{2}$,
$z=\lambda(x-\mu)$. Also, write
$\tilde{\theta}(t,z):=\theta(t,x)$. An elementary calculation then
shows that we have
\begin{equation}\nonumber\begin{split}
&i\partial_{t}W+\partial_{x}^{2}W+|W|^{4}W\\&=(\lambda_{t}\lambda^{-1}-\beta\lambda^{2})[ie^{i\Psi}z\lambda^{\frac{1}{2}}\phi_{0}'(z)+\frac{\beta}{2}z^{2}W+\frac{i}{2}W
-\tilde{\theta}_{z}zW]\\&+\frac{z^{2}}{4}(\beta_{t}+\lambda^{2}\beta^{2})W+i\lambda^{2}\tilde{\theta}_{zz}W+(\lambda^{2}-\tilde{\theta}_{t}+\lambda^{2}\tilde{\theta}_{z}^{2}-\beta\lambda^{2}z\tilde{\theta}_{z})W
\\&+(\mu_{t}-2\lambda\tilde{\theta}_{z})[\lambda\tilde{\theta}_{z}-i\lambda^{\frac{3}{2}}e^{i\Psi}\phi_{0}'(z)-\frac{\beta}{2}\lambda
zW]
\end{split}\end{equation}
Introduce the following notation:
\begin{align*}
\eta_1 &:= \binom{ie^{i\Psi} \sqrt{\lambda}\phi_0(z)}{-ie^{-i\Psi}
\sqrt{\lambda}\phi_0(z)},
\quad \eta_2 := \binom{e^{i\Psi}(z \sqrt{\lambda}\phi'_0(z)+\sqrt{\lambda}\phi_0(z)/2)}{e^{-i\Psi}(z \sqrt{\lambda}\phi'_0(z)+\sqrt{\lambda}\phi_0(z)/2)}\\
\eta_3 &:= \binom{e^{i\Psi} \sqrt{\lambda}\phi'_0(z)}{e^{-i\Psi}
\sqrt{\lambda}\phi'_0(z)},
\quad \eta_4 := \binom{ie^{i\Psi} z \sqrt{\lambda}\phi_0(z)}{-ie^{-i\Psi}z \sqrt{\lambda}\phi_0(z)}\\
\eta_5 &:= \binom{ie^{i\Psi}
z^2\sqrt{\lambda}\phi_0(z)}{-ie^{-i\Psi} z^2
\sqrt{\lambda}\phi_0(z)}, \quad \eta_6 := \binom{e^{i\Psi}
\sqrt{\lambda}\rho(z)}{e^{-i\Psi}\sqrt{\lambda}\rho(z)}
\end{align*}
Now impose the relation $\tilde{\theta}_{zz}=0$, whence
$\tilde{\theta}=\gamma(t)+\omega(t)z$. The function $\rho$ here is
defined via
\begin{align*}
L_- &:= -\partial_{xx} +1 - \phi_0^4, \qquad L_+ := -\partial_{xx} +1 - 5\phi_0^4 \\
L_-\phi_0 &= 0, \quad L_+ (\frac12 \phi_0 + x\phi_0')=-2\phi_0 \\
L_-(x^2\phi_0) &= -4(\frac12\phi_0+x\phi_0'), \quad L_+\rho =
x^2\phi_0
\end{align*}
Then the vector-function $\calW(t,x) :=
\binom{W(t,x)}{\bar{W}(t,x)}$ satisfies
\begin{align*}
& i\partial_t \calW + \left[\begin{array}{ll} \partial_{xx} & 0 \\
0 & -\partial_{xx}
\end{array}\right] \calW +\binom{|W|^4W}{-|W|^4\bar{W}} \\
& = i(\dot{\lambda}\lambda^{-1}-\beta\lambda^2)
(\eta_2-\beta\eta_5/2+\omega\eta_4) -
\frac{i}{4}(\dot{\beta}+\lambda^2\beta^2)\eta_5
 -i(\lambda^2-\dot{\gamma} + \lambda^2\omega^2)\eta_1 \\
&\qquad +i(\dot{\omega}+
\beta\lambda^2\omega)\eta_4+\lambda(\dot{\mu}-2\lambda\omega)(-i\omega
\eta_1 -i\eta_3 + i\beta \eta_4/2)
\end{align*}
Write $\bm \psi(t,x)\\ \overline{\psi(t,x)}\endm
=\calW(t,x)+Z(t,x)$ where $Z(t,x)=\bm R(t,x)\\
\overline{R(t,x)}\endm$. We deduce the following equation for $Z$:
\begin{align}\label{Zeqn}
i\partial_t Z + \Hil(t)Z &=
-i(\dot{\lambda}\lambda^{-1}-\beta\lambda^2)
(\eta_2-{\beta}\eta_5/2+\omega\eta_4) +
\frac{i}{4}(\dot{\beta}+\lambda^2\beta^2)\eta_5
 +i(\lambda^2-\dot{\gamma} + \lambda^2\omega^2)\eta_1 \\
&\qquad -i(\dot{\omega}+ \beta\lambda^2\omega)\eta_4-\lambda(\dot{\mu}-2\lambda\omega)(-i\omega \eta_1 -i\eta_3 + i{\beta} \eta_4/2) + N(Z)\\
\Hil(t) &:= \left[\begin{array}{ll} \partial_{xx} + 3|W|^4 & 2|W|^2 W^2 \\ -2|W|^2\bar{W}^2 & -\partial_{xx} - 3|W|^4 \\
\end{array} \right] \\
N(Z) &:= \binom{-|R+W|^4(R+W)+|W|^4W+3|W|^4R+2|W|^2W^2\bar{R}}{|R+W|^4(\bar R+\bar W)-|W|^4\bar W-3|W|^4\bar{R}-2|W|^2\bar{W}^2R} \\
&=\binom{-3R^2|W|^2\bar{W}-6|R|^2|W|^2W-W^3\bar{R}^2+O(|R|^3|W|^2+|R|^5)}{3\bar{R}^2|W|^2{W}+6|R|^2|W|^2\bar{W}+\bar{W}^3{R}^2
+O(|R|^3|W|^2+|R|^5)}
\end{align}
In order to formulate the modulation equations, it will be
convenient to introduce the following family of auxiliary
functions, which are in some sense dual to the $\eta_{i}$:
\begin{align*}
\xi_1 &:= \binom{e^{i\Psi} \sqrt{\lambda}\phi_0(z)}{e^{-i\Psi}
\sqrt{\lambda}\phi_0(z)},
\quad \xi_2 := \binom{ie^{i\Psi}(z \sqrt{\lambda}\phi'_0(z)+\sqrt{\lambda}\phi_0(z)/2)}{-ie^{-i\Psi}(z \sqrt{\lambda}\phi'_0(z)+\sqrt{\lambda}\phi_0(z)/2)}\\
\xi_3 &:= \binom{ie^{i\Psi} \sqrt{\lambda}\phi'_0(z)}{-ie^{-i\Psi}
\sqrt{\lambda}\phi'_0(z)},
\quad \xi_4 := \binom{e^{i\Psi} z \sqrt{\lambda}\phi_0(z)}{e^{-i\Psi}z \sqrt{\lambda}\phi_0(z)}\\
\xi_5 &:= \binom{e^{i\Psi} z^2\sqrt{\lambda}\phi_0(z)}{e^{-i\Psi}
z^2 \sqrt{\lambda}\phi_0(z)}, \quad \xi_6 := \binom{ie^{i\Psi}
\sqrt{\lambda}\rho(z)}{-ie^{-i\Psi}\sqrt{\lambda}\rho(z)}
\end{align*}
Then, analogously to \eqref{orthorelations}, we have the following
(using the same notation)
\begin{align}\label{orthorelations'}
\langle \eta_j, \xi_1 \rangle &= 0, \; 1\le j\le 5, \;\langle \eta_6, \xi_1 \rangle = 2\kappa_2 \\
\langle \eta_j,\xi_2 \rangle &= 0,\; j=1,2,3,4,6,\; \quad \langle \eta_5,\xi_2\rangle = -4\kappa_2 \\
\langle \eta_j,\xi_3 \rangle &= 0,\; j=1,2,3,5,6,\; \quad \langle \eta_4,\xi_3\rangle = -\kappa_1 \\
\langle \eta_j,\xi_4 \rangle &= 0,\; j=1,2,4,5,6,\; \quad \langle \eta_3,\xi_4\rangle = -\kappa_1 \\
\langle \eta_j,\xi_5 \rangle &= 0,\; j=1,3,4,5,\; \quad \langle
\eta_2,\xi_5\rangle = -4\kappa_2,
\; \quad \langle \eta_6,\xi_5\rangle = 2\kappa_3 \\
\langle \eta_j,\xi_6 \rangle &= 0,\; j=2,3,4,6,\; \quad \langle
\eta_1,\xi_6\rangle = 2\kappa_2,\; \quad \langle
\eta_5,\xi_6\rangle = 2\kappa_3
\end{align}
Recall from the discussion in the first section that we impose the
orthogonality relations $\la Z, \xi_{i}\ra =0$, $i=2,\ldots,6$.
This allows us to control the 'good component' of the root part of
the radiation. Using Leibnitz' rule we get
\begin{equation}\nonumber
\la i\partial_{t}Z+\calH(t)Z, \xi_{j}(t)\ra =\la Z,
(i\partial_{t}+\calH(t)^{*})\xi_{j}(t)\ra =: \la Z,
\calL^{*}\xi_{j}\ra
\end{equation}
Explicitly, using \eqref{Zeqn}, these read as follows:
\begin{align*}
 -2i\beta (\dot{\lambda}\lambda^{-1} - \beta\lambda^2) \kappa_2 &- i(\dot{\beta}+\lambda^2\beta^2) \kappa_2
 = \langle Z,\calL^* \xi_2 \rangle - \langle N(Z),\xi_2 \rangle \\
 \frac{i}{2}\lambda\beta\kappa_1 (\dot{\mu} - 2\lambda \omega) &+ i\omega\kappa_1(\dot{\lambda}\lambda^{-1} - \beta\lambda^2) + i\kappa_1 (\dot{\omega}+\beta\lambda^2\omega) = \langle Z,\calL^* \xi_3\rangle -
\langle N(Z), \xi_3\rangle \\
 -i\lambda\kappa_1 (\dot{\mu}-2\lambda\omega) &= \langle Z,\calL^* \xi_4\rangle -
\langle N(Z), \xi_4\rangle \\
 4i\kappa_2 (\dot{\lambda}\lambda^{-1}-\beta\lambda^2) &= \langle Z,\calL^* \xi_5\rangle -
\langle N(Z), \xi_5\rangle \\
2i\kappa_2 (\lambda^2-\dot{\gamma}+\lambda^2\omega^2) &+ 2i\lambda\omega\kappa_2 (\dot{\mu}-2\lambda\omega) + i\beta\kappa_3 (\dot{\lambda}\lambda^{-1}-\beta\lambda^2)+\frac{i}{2}\kappa_3(\dot{\beta}+\lambda^2\beta^2)\\
& = \langle Z,\calL^* \xi_6\rangle - \langle N(Z), \xi_6\rangle
\end{align*}
Before proceeding, let's carry out a consistency check: we know
that the case $Z=0$ corresponds to a transformed standing wave. In
this case, the above relations simplify to
\begin{align*}
\dot{\lambda}\lambda^{-1}-\beta\lambda^2 &=0, \quad
\dot{\beta}+\lambda^2\beta^2 = 0, \quad\frac{d}{dt} (\lambda\beta)=0, \; \beta=-b\lambda^{-1} \\
\dot{\lambda} &= -b\lambda^2,\quad \lambda(t) = (a+bt)^{-1}, \; \beta(t)= -b(a+bt),\; \\
\dot{\mu}-2\lambda\omega &= 0, \quad \dot{\omega}+\beta\lambda^2\omega =0,\quad \dot{\omega}-b\lambda\omega=0\\
\omega(t) &=v(a+bt), \quad \dot{\mu} = 2\lambda \omega = 2v,\quad \mu(t)=2tv+\mu_0 \\
\lambda^2-\dot{\gamma}+\lambda^2\omega^2 &= 0,\; \dot{\gamma} = \lambda^2 + v^2 = -\frac{1}{b}\dot{\lambda} + v^2\\
\gamma(t) &= -\frac{\lambda}{b} + tv^2 + \gamma_0 =
-\frac{1}{b(a+bt)}+v^2t+\gamma_0
\end{align*}
So the exact solution looks as follows (with
$\gamma=\gamma_0-v\mu_0$):
\begin{align*}
{\theta}(t,z) &= -\frac{1}{b(a+bt)}-v^2t+\gamma+vx \\
\lambda(t) &= (a+bt)^{-1} \\
\beta(t) &= -b(a+bt) \\
\mu(t) &= 2tv+ \mu_0
\end{align*}
and the transformed standing wave $W$ is
\begin{align} W(t,x)& = \exp\Big(i\Big[-\frac{1}{b(a+bt)}-v^2 t+\gamma + vx +\frac{b}{4(a+bt)}(x-2tv-\mu_0)^2\Big]\Big)\cdot \nn\\
& \qquad (a+bt)^{-\half}\phi_0\big((a+bt)^{-1}(x-2tv-\mu_0)).
\label{eq:Wexact}
\end{align}
Now recall the pseudo-conformal transformation
\[ \calC_{M}:\quad \psi(t,x) \to (a+bt)^{-\half} \exp\Big(i\frac{bx^2}{4(a+bt)}\Big) \psi\Big(\frac{c+dt}{a+bt},\frac{x}{a+bt}\Big)\]
where $M=\bm a&b\\c&d \endm \in SL(2,\R)$. Starting from the
standing wave $e^{it}\phi_0(x)$, apply the pseudo-conformal
transformation with matrix $\bm a & b \\ -b^{-1} & 0 \endm$:
\[ \exp\big(-i\frac{1}{b(a+bt)}+i\frac{b}{4(a+bt)}x^2\big)(a+bt)^{-\half}\phi_0((a+bt)^{-1}x)\]
and then the Galilei transform
\[ \calg_{\gamma,v,\mu_0}(t)=e^{i(\gamma+vx-tv^2)} e^{-i(2tv+\mu_0)p}.\]
This leads to the exact same expression as in~\eqref{eq:Wexact}.
\\

We now intend to translate the above equations from the $(t,x)$
coordinates to a new coordinate system $(s,y)$, in which we
'de-singularize' the equations. The blow-up time $t_{*}$ shall be
transformed into $s=+\infty$, and the $t$-interval
$(-\infty,t_{*}]$ shall correspond to $(-c,\infty]$ for suitable
$c>0$. Thus make the ansatz\footnote{We use the notation $U(s,y)$,
$U(s)$, $U$ {\it{all for the same function of two variables}} in
order to streamline the notation in some places.}
\begin{equation}\nonumber
\bm U\\ \bar{U}\endm:= \calM\calT_{\infty}\bm R\\ \bar{R}\endm,
\end{equation}
where
\begin{equation}\nonumber
\calT_{\infty}=e^{-i(v_{\infty}^{2}s+\gamma_{\infty}+v_{\infty}y)}e^{i(2v_{\infty}s+y_{\infty})p}\bm
a_{\infty}&b_{\infty}\\0&a_{\infty}^{-1}\endm,\,
\calM=\calM(s)=\bm e^{-is}&0\\0&e^{is}\endm, p=-i\frac{d}{dy}
\end{equation}
Then we have the succinct identity
\begin{equation}\nonumber
(\calT_{\infty}F)(s,y)=e^{-i\tilde{\Psi}_{\infty}(s,y)}\lambda_{\infty}^{-\frac{1}{2}}(s)F(\int_{0}^{s}\lambda_{\infty}^{-2}(\sigma)d\sigma,
\lambda_{\infty}^{-1}(s)y+\mu_{\infty}(s)),
\end{equation}
where
\begin{equation}\label{tilPsi}
\tilde{\Psi}_{\infty}(s,y)=v_{\infty}^{2}s+\gamma_{\infty}+v_{\infty}y-\frac{b_{\infty}(y+2v_{\infty}s+y_{\infty})^{2}}{4(a_{\infty}+b_{\infty}s)},\,\mu_{\infty}(s)=\frac{2v_{\infty}s+y_{\infty}}{a_{\infty}+b_{\infty}s},\,\lambda_{\infty}=a_{\infty}+b_{\infty}s
\end{equation}
Thus we have
\begin{equation}\nonumber
e^{-is}(\calT_{\infty}W)(s,y)=e^{-i(s+\tilde{\Psi}_{\infty})+i\Psi(t(s),\mu_{\infty}+\lambda_{\infty}^{-1}(s)y)}\lambda_{\infty}^{-\frac{1}{2}}(s)\lambda(t(s))^{\frac{1}{2}}\phi_{0}(\lambda(t(s))(\mu_{\infty}(s)-\mu(t(s))+\lambda_{\infty}^{-1}(s)y)),
\end{equation}
where we put
$t(s):=\int_{0}^{s}\lambda_{\infty}^{-2}(\sigma)d\sigma=\frac{a_{\infty}^{-1}s}{a_{\infty}+b_{\infty}s}$.
We shall now impose the asymptotic conditions
$\nu(s):=\frac{\lambda(s)}{\lambda_{\infty}(s)}\rightarrow 1$,
$\lambda(t(s))(\mu_{\infty}(s)-\mu(t(s)))\rightarrow 0$ as
$s\rightarrow +\infty$. Unfortunately, it appears that no such
requirement can be applied to
$\Psi(t(s),\mu_{\infty}+\lambda_{\infty}^{-1}(s)y)-s-\tilde{\Psi}_{\infty}(s,y)$,
as will follow from the ensuing discussion.
\\

Now introduce the Schwartz functions
$\tilde{\eta}_{i}:=\calM\calT_{\infty}\eta_{i}$,
$\tilde{\xi}_{i}:=\calM\calT_{\infty}\xi_{i}$.  Then we can deduce
the following equation for $\bm U\\ \bar{U}\endm$:
\begin{align}\label{Ueqn'}
i\partial_{s}\bm U\\ \bar{U}\endm+\calH(s)\bm U\\ \bar{U}\endm =&-i(\dot\lambda\lambda^{-1}-\beta\nu^2)(\tileta_2-\beta\tileta_5/2+\omega\tileta_4)\\
&\quad +\frac{i}{4}(\dot\beta+\beta^2\nu^2)\tileta_5 + i(\nu^2-\dot\gamma+\nu^2\omega^2)\tileta_1 \\
&\quad -i(\dot\omega+\beta\omega\nu^2)\tileta_4 -
i\nu(\dot\mu\lambda_\infty-2\nu\omega)(-\omega\tileta_1-\tileta_3+\beta\tileta_4/2)
+ N(U,\pi),
\end{align}
In this equation $\dot{\lambda}=\partial_{s}[\lambda(t(s))]$. We
use the
abbreviations\footnote{Also, recall that $\Psi_{\infty}(s,y)=\tilde{\Psi}_{\infty}(s,y)+s$.} (with $\bm U\\ \bar{U}\endm =\bm U_{1}\\
U_{2}\endm$)
\begin{align*}
N(U,\pi) &:= \binom{-3U_1^2
\tilde{\phi}_0^3e^{i(\Psi_\infty-\Psi)}\nu^{\frac32}-6|U_1|^2
e^{i(\Psi-\Psi_\infty)} \nu^{\frac32} \tilde{\phi}_0^3 - U_2^2
e^{3i(\Psi-\Psi_\infty)} \nu^{\frac{3}{2}}\tilde{\phi}_0^3 +
O(|U|^3+|U|^5)}{3U_1^2
\tilde{\phi}_0^3e^{-i(\Psi_\infty-\Psi)}\nu^{\frac32}+6|U_1|^2
e^{-i(\Psi-\Psi_\infty)} \nu^{\frac32} \tilde{\phi}_0^3 + U_2^2
e^{-3i(\Psi-\Psi_\infty)} \nu^{\frac32}\tilde{\phi}_0^3 +
O(|U|^3+|U|^5)}
\end{align*}
where
$\tilde{\phi}_{0}(y)=\phi_{0}(\lambda(\mu_{\infty}-\mu+\lambda_{\infty}^{-1}y))$.
\begin{align*}
\calH(s):=\bm
\partial_{y}^{2}-1+3\nu^{2}(s)\tilde{\phi}_{0}^{4}&
2\nu^{2}\tilde{\phi}_{0}^{4}e^{2i(\Psi-\Psi_{\infty})}\\
-2\nu^{2}\tilde{\phi}_{0}^{4}e^{-2i(\Psi-\Psi_{\infty})}&
-\partial_{y}^{2}+1-3\nu^{2}(s)\tilde{\phi}_{0}^{4}\endm
\end{align*}
The orthogonality conditions $\la Z, \xi_{i}\ra=0$,
$i=2,\ldots,6$, translate to $\la \bm U\\ \bar{U}\endm,
\tilde{\xi}_{i}\ra=0$, $i=2,\ldots,6$. If one differentiates this
relation with respect to $s$ and uses the Leibnitz rule as well as
\eqref{Ueqn'}, this leads to the following system of 'ODE's', the
modulation equations:
\begin{align}\label{modulation}
-2\kappa_2i\beta(\dot\lambda\lambda^{-1}-\beta\nu^2) &- i\kappa_2(\dot\beta+\beta^2\nu^2) = -\langle N,\tilxi_2\rangle + \la U,(i\partial_s+\Hil(s)^*)\xitil_2\ra\\
i\omega\kappa_1 (\dot\lambda\lambda^{-1}-\beta\nu^2) &+
i\kappa_1(\dot\omega+\beta\omega\nu^2)
+\frac{i}{2}\beta\nu(\dot\mu\lambdainf-2\omega\nu)\kappa_1 = -\langle N,\tilxi_3\rangle + \la U,(i\partial_s+\Hil(s)^*)\xitil_3\ra\\
-i\kappa_1\nu(\dot\mu\lambdainf-2\omega\nu) &= -\langle N,\tilxi_4\rangle + \la U,(i\partial_s+\Hil(s)^*)\xitil_4\ra\\
4i\kappa_2 (\dot\lambda\lambda^{-1}-\beta\nu^2) &= -\langle N,\tilxi_5\rangle + \la U,(i\partial_s+\Hil(s)^*)\xitil_5\ra\\
2i\kappa_2(\nu^2-\dot\gamma+\nu^2\omega^2)
&+2i\omega\nu\kappa_2(\dot\mu\lambdainf-2\nu\omega)
+\frac{i}{2}\kappa_3(\dot\beta+\beta^2\nu^2) + i\beta\kappa_3(\dot\lambda\lambda^{-1}-\beta\nu^2)\\
&\nonumber\quad = -\langle N,\tilxi_6\rangle + \la
U,(i\partial_s+\Hil(s)^*)\xitil_6\ra,
\end{align}
Of course, for all this to make sense we need to specify the
'parameters at infinity'\\ $\{a_{\infty}, b_{\infty},
v_{\infty},y_{\infty},\gamma_{\infty}\}$. We shall soon see that
their value is forced by the asymptotic conditions on the
modulation parameters.
\\

In accordance with the statement of Theorem~\ref{Main}, we now fix
the values of $\lambda(s)$,$\beta(s)$,$\mu(s)$, $\omega(s)$,
$\gamma(s)$ at time $s=0$, where we require $\lambda(0)\sim 1$,
$\beta(0)\sim 1$. Instead of working with these parameters,
though, we shall work with
$\nu(s)=\frac{\lambda(s)}{\lambda_{\infty}(s)}$,
$\beta(s)\nu(s)-\frac{b_{\infty}}{\lambda_{\infty}(s)}$, $\mu(s)$,
$\omega(s)$, $\gamma(s)$. Start with the fourth modulation
equation. Formulate this as follows:
\begin{equation}\nonumber
\dot{\nu}\nu^{-1}-\beta\nu^{2}=(4i\kappa_{2})^{-1}[-\la
N,\tilde{\xi}_{5}\ra +\la
U,(i\partial_{s}+\calH(s)^{*})\tilde{\xi}_{5}\ra]-b_{\infty}\lambda_{\infty}^{-1},
\end{equation}
From the fourth and 2nd equation, we get
\begin{equation}\nonumber
\dot{\beta}+\beta^{2}\nu^{2}=-(i\kappa_{2})^{-1}E_{2}-\frac{\beta}{2i\kappa_{2}}E_{5},
\end{equation}
where we use the notation $E_{j}:=-\la N, \tilde{\xi}_{j}\ra+\la
U,(i\partial_{s}+\calH^{*}(s))\tilde{\xi}_{j}\ra$. Noting the
simple identity
$(b_{\infty}\lambda_{\infty}^{-1})_{s}+(b_{\infty}\lambda_{\infty}^{-1})^{2}=0$,
we get
\begin{equation}\nonumber
\dot{\nu}-b_{\infty}\lambda_{\infty}^{-1}(\nu-1)=\nu(4i\kappa_{2})^{-1}E_{5}+\beta\nu(\nu-1)^{2}+(2\nu-1)(\beta\nu-b_{\infty}\lambda_{\infty}^{-1})
\end{equation}
\begin{equation}\nonumber
\frac{d}{ds}(\beta\nu-b_{\infty}\lambda_{\infty}^{-1})+(\beta\nu-b_{\infty}\lambda_{\infty}^{-1})b_{\infty}\lambda_{\infty}^{-1}
=-\nu(i\kappa_{2})^{-1}E_{2}+[-\beta\frac{\nu}{2i\kappa_{2}}+\nu\beta(4i\kappa_{2})^{-1}]E_{5}
\end{equation}
We can further reformulate these equations as follows:
\begin{equation}\nonumber
\frac{d}{ds}[(\nu-1)\lambda_{\infty}^{-1}](s)=\lambda_{\infty}^{-1}[\nu(4i\kappa_{2})^{-1}E_{5}+\beta\nu(\nu-1)^{2}+(2\nu-1)(\beta\nu-b_{\infty}\lambda_{\infty}^{-1})]
\end{equation}
\begin{equation}\nonumber
\frac{d}{ds}([\beta\nu-b_{\infty}\lambda_{\infty}^{-1}]\lambda_{\infty})=\lambda_{\infty}[-\nu(i\kappa_{2})^{-1}E_{2}-\frac{\beta\nu}{4i\kappa_{2}}E_{5}](s)
\end{equation}
The condition that $\nu(s)\rightarrow 1$ as $s\rightarrow+\infty$,
as well as the condition
$\beta\nu(s)-b_{\infty}\lambda_{\infty}^{-1}(s)\rightarrow 0$
imply the following identities:
\begin{equation}\nonumber
0=(\beta\nu-b_{\infty}\lambda_{\infty}^{-1})\lambda_{\infty}(0)+\int_{0}^{\infty}\lambda_{\infty}(s)[-\nu(i\kappa_{2})^{-1}E_{2}-\frac{\beta\nu}{4i\kappa_{2}}E_{5}](s)ds
\end{equation}
\begin{equation}\nonumber
0=(\nu-1)\lambda_{\infty}^{-1}(0)+\int_{0}^{\infty}\lambda_{\infty}(s)^{-1}[\nu(4i\kappa_{2})^{-1}E_{5}+\beta\nu(\nu-1)^{2}+(2\nu-1)(\beta\nu-b_{\infty}\lambda_{\infty}^{-1})](s)ds
\end{equation}
whence
\begin{equation}\label{binfty}
0=(\beta(0)\lambda(0)-b_{\infty})+\int_{0}^{\infty}\lambda_{\infty}(s)[-\nu(i\kappa_{2})^{-1}E_{2}+\frac{3\beta\nu}{4i\kappa_{2}}E_{5}](s)ds
\end{equation}
\begin{equation}\label{ainfty}
0=\lambda(0)-a_{\infty}+a_{\infty}^{2}\int_{0}^{\infty}\lambda_{\infty}(s)^{-1}[\nu(4i\kappa_{2})^{-1}E_{5}+\beta\nu(\nu-1)^{2}+(2\nu-1)(\beta\nu-b_{\infty}\lambda_{\infty}^{-1})](s)ds
\end{equation}
Assuming the integral expressions known, this allows for solving
for the coefficients $a_{\infty}, b_{\infty}$, using the Implicit
function Theorem. Moreover, we get the formulae
\begin{equation}\label{nu}
\nu(s)-1=-\lambda_{\infty}(s)\int_{s}^{\infty}\lambda_{\infty}(\sigma)^{-1}[\nu(4i\kappa_{2})^{-1}E_{5}+\beta\nu(\nu-1)^{2}+(2\nu-1)(\beta\nu-b_{\infty}\lambda_{\infty}^{-1})](\sigma)d\sigma
\end{equation}
\begin{equation}\label{beta}
(\beta\nu-b_{\infty}\lambda_{\infty}^{-1})(s)=-\lambda_{\infty}(s)^{-1}\int_{s}^{\infty}\lambda_{\infty}(\sigma)[-\nu(i\kappa_{2})^{-1}E_{2}-\frac{\beta\nu}{4i\kappa_{2}}E_{5}](\sigma)d\sigma
\end{equation}
\\

Next, from the 2nd, 3rd and 4th modulation equation we gather
\begin{equation}\nonumber
\dot{\omega}+\beta\nu^{2}\omega=(i\kappa_{1})^{-1}E_{3}-\omega(4i\kappa_{2})^{-1}E_{5}-\frac{\beta}{2i\kappa_{1}}E_{4}
\end{equation}
Introduce the quantity
\begin{equation}\nonumber
B(s)=\exp(\int_{0}^{s}[\beta\nu^{2}+\frac{1}{4i\kappa_{2}}E_{5}](\sigma)d\sigma)
\end{equation}
We can then write
\begin{equation}\label{omega}
\omega(s)=B(s)^{-1}\omega(0)+\int_{0}^{s}\frac{B(\sigma)}{B(s)}[(i\kappa_{1})^{-1}E_{3}+\frac{\beta}{2i\kappa_{1}}E_{4}](\sigma)d\sigma
\end{equation}
Decompose
\begin{equation}\nonumber
\beta\nu^{2}=(\beta\nu-b_{\infty}\lambda_{\infty}^{-1})\nu+(\nu-1)b_{\infty}\lambda_{\infty}^{-1}+b_{\infty}\lambda_{\infty}^{-1}
\end{equation}
The stipulations $\lim_{s\rightarrow+\infty}\nu(s)=1$,
$\lim_{s\rightarrow+\infty}\beta\nu-b_{\infty}\lambda_{\infty}^{-1}=0$
then yield\footnote{We shall soon specify the precise decay
rates.}
\begin{equation}\label{omegaasympto}
B^{-1}(s)=c\lambda_{\infty}^{-1}(s)+o(\frac{1}{s})
\end{equation}
We then reformulate \eqref{omega} as follows:
\begin{equation}\label{omega1}\begin{split}
&\omega(s)=c\lambda_{\infty}(s)^{-1}[\omega(0)+\int_{0}^{\infty}B(\sigma)[(i\kappa_{1})^{-1}E_{3}+\frac{\beta}{2i\kappa_{1}}E_{4}](\sigma)d\sigma]
\\&\hspace{3cm}-c\lambda_{\infty}(s)^{-1}\int_{s}^{\infty}B(\sigma)[(i\kappa_{1})^{-1}E_{3}+\frac{\beta}{2i\kappa_{1}}E_{4}](\sigma)d\sigma]+o(s^{-1}),
\end{split}\end{equation}
from which we obtain for suitable $c_{\infty}$ the asymptotic
relation $\omega(s)=c_{\infty}\lambda_{\infty}^{-1}+o(s^{-1})$,
provided we can control all the integrals. From the 3rd modulation
equation we obtain
\begin{equation}\label{mu}
\mu(s)=\mu(0)+\int_{0}^{s}\lambda_{\infty}(\sigma)^{-1}[2\omega\nu+(i\kappa_{1}\nu(\sigma))^{-1}E_{4}(\sigma)]d\sigma
\end{equation}
If we feed in the relation \eqref{omegaasympto}, we infer the
existence of parameters $v_{\infty}$, $y_{\infty}$ with the
property
\begin{equation}\label{musasympto}
\mu(s)=\frac{2v_{\infty}s+y_{\infty}}{a_{\infty}+b_{\infty}s}+o(s^{-1})
\end{equation}
Finally, the 5th modulation equation gives
\begin{equation}\label{gamma}\begin{split}
&\gamma(s)=\gamma(0)+\int_{0}^{s}[\nu^{2}(\sigma)-(2i\kappa_{2})^{-1}E_{6}(\sigma)+\nu^{2}(\sigma)\omega^{2}(\sigma)-\frac{1}{i\kappa_{1}}\omega(\sigma)E_{4}(\sigma)
-\frac{1}{2\kappa_{2}}(i\kappa_{2})^{-1}E_{2}(\sigma)]d\sigma
\end{split}\end{equation}
Last but not least, we choose $\gamma_{\infty}$ such that
$(\Psi-\Psi_{\infty})_{1}(0)=0$, where
$\Psi_{\infty}(s,y)=\tilde{\Psi}_{\infty}(s,y)+s$, and we define
$(\Psi-\Psi_{\infty})_{1}(s)$ to be that part of
$\Psi-\Psi_{\infty}$ which only depends on $s$, see the ensuing
subsection.
\\

We now state the {\bf{precise estimates for the modulation
parameters}}: first, choose small positive numbers
$\delta_{i},\,i=1,2,3$, and $\delta>0$ with the property
$\delta<<\delta_{2}<<\delta_{3}<<\delta_{1}$. These shall be fixed
throughout. The number $\delta$ will control the
size\footnote{With respect to suitable norms.} of radiation part
as well as modulation parameters, while the parameters
$\delta_{i}$, measure parameters in certain norms. Then we need
for a sufficiently large\footnote{We shall need
$N<<N_{1}(\delta_{2}, \delta_{1})$ and $N>N_{2}(\delta_{3})$. The
parameter $\delta_{2}$ will appear in the estimates
\eqref{disasympto}, where we specify $N_{1}(\delta_{2},
\delta_{1})\sim |\log( \frac{\delta_{1}}{\delta_{2}})|$. The bound
$N>N_{2}(\delta_{3})$ is needed in bootstrapping the strong local
dispersive estimate.} $N=N(\delta_{2}, \delta_{3}, \delta_{1})$
and very large $M>M(\delta_{2})$ held fixed throughout
\begin{equation}\label{modulationasympto}\begin{split}
&|\nu(s)-1|\lesssim \delta^{2}\la
s\ra^{-\frac{1}{2}+\delta_{1}},\,\sup_{1\leq i\leq
[\frac{N}{2}]}||\la
s\ra^{\frac{3}{2}-2\delta_{1}}\frac{d^{i}}{ds^{i}}\nu(s)||_{L^{M}}\lesssim
\delta^{2},\,|\beta(s)\nu(s)-b_{\infty}\lambda_{\infty}^{-1}(s)|\lesssim
\la s\ra^{\frac{3}{2}-\delta_{1}}\delta^{2},\\&\sup_{1\leq i\leq
[\frac{N}{2}]}||\la
s\ra^{2-2\delta_{1}}\frac{d^{i}}{ds^{i}}[\beta(s)\nu(s)-b_{\infty}\lambda_{\infty}^{-1}(s)]||_{L^{M}}\lesssim
\delta^{2},\,|\omega(s)-c_{\infty}\lambda_{\infty}^{-1}(s)|\lesssim
\delta^{2}\la s\ra^{-\frac{3}{2}+\delta_{1}},\\&\sup_{1\leq i\leq
[\frac{N}{2}]}||\la
s\ra^{2-2\delta_{1}}\frac{d^{i}}{ds^{i}}[\omega(s)-c_{\infty}\lambda_{\infty}^{-1}(s)]||_{L^{M}}\lesssim
\delta^{2},|\partial_{s}(\gamma(s)-s-c_{\infty}^{2}\lambda_{\infty}^{-2}(s))|\lesssim
\delta^{2}\la s\ra^{-\frac{1}{2}+\delta_{1}},\\&\sup_{2\leq i\leq
[\frac{N}{2}]}||\la
s\ra^{\frac{3}{2}-2\delta_{1}}\frac{d^{i}}{ds^{i}}(\gamma(s)-s-c_{\infty}^{2}\lambda_{\infty}^{-2}(s))||_{L^{M}}\lesssim
\delta^{2},
|\mu(s)-\frac{2v_{\infty}s+y_{\infty}}{a_{\infty}+b_{\infty}s}|\lesssim
\delta^{2}\la s\ra^{-\frac{3}{2}+\delta_{1}},\\&\sup_{1\leq i\leq
[\frac{N}{2}]}||\la
s\ra^{\frac{5}{2}-\delta_{1}}\frac{d^{i}}{ds^{i}}[\mu(s)-\frac{2v_{\infty}s+y_{\infty}}{a_{\infty}+b_{\infty}s}]||_{L^{M}}\lesssim
\delta^{2}
\end{split}\end{equation}
The fact that we work with $L^{M}$ instead of $L^{\infty}$ for the
derivatives is a technical complication due to the fact that we
need a {\it{compactness property}} for the fixed point Theorem to
apply, see below.

\subsection{Algebraic manipulations II; analysis of the radiation
part.}

We now look at $\bm U\\ \bar{U}\endm$. As mentioned in the first
section, we essentially break this into its root and dispersive
part; more precisely, we first tweak this function a bit, after a
careful analysis of the phase $(\Psi-\Psi_{\infty})(s,y)=(\Psi-\tilde{\Psi}_{\infty})(s,y)-s$. \\
From \eqref{tilPsi} we infer the relation\footnote{Recall that we
also defined
$\mu_{\infty}(s)=\frac{2v_{\infty}s+y_{\infty}}{a_{\infty}+b_{\infty}s}$}
\begin{align*}
&\Psi-\Psi_{\infty}(s,y)\\&=\gamma(s)-s+y[\omega(s)\nu(s)-\frac{\beta(s)}{2}\nu(s)\lambda(s)(\mu_{\infty}-\mu)-\frac{a_{\infty}v_{\infty}-\frac{b_{\infty}y_{\infty}}{2}}{a_{\infty}+b_{\infty}s}]+\omega(s)\lambda(s)(\mu_{\infty}-\mu)(s)
\\&+\frac{b_{\infty}y^{2}}{4(a_{\infty}+b_{\infty}s)}-\frac{\beta}{4}[\nu
y]^{2}-\frac{\beta}{4}[\lambda[\mu_{\infty}-\mu]]^{2}-[\frac{v_{\infty}^{2}s
a_{\infty}}{a_{\infty}+b_{\infty}s}-\frac{b_{\infty}v_{\infty}s
y_{\infty}}{a_{\infty}+b_{\infty}s}+\gamma_{\infty}-\frac{b_{\infty}y_{\infty}^{2}}{4(a_{\infty}+s
b_{\infty})}]
\end{align*}
We decompose this into two parts,
$(\Psi-\Psi_{\infty})(s,y)=(\Psi-\Psi_{\infty})_{1}(s)+(\Psi-\Psi_{\infty})_{2}(s,y)$,
where
\begin{equation}\label{Psi-Psiinfty}\begin{split}
&(\Psi-\Psi_{\infty})_{1}(s)\\&=\gamma(s)-s+\omega(s)\lambda(s)(\mu_{\infty}-\mu)(s)-\frac{\beta}{4}[\lambda[\mu_{\infty}-\mu]]^{2}-[\frac{v_{\infty}^{2}s
a_{\infty}}{a_{\infty}+b_{\infty}s}-\frac{b_{\infty}v_{\infty}s
y_{\infty}}{a_{\infty}+b_{\infty}s}+\gamma_{\infty}-\frac{b_{\infty}y_{\infty}^{2}}{4(a_{\infty}+sb_{\infty})}],
\end{split}\end{equation}
i. e. this is the part of $\Psi-\Psi_{\infty}$ which only depends
on $s$ and not on $y$. Then we define
\begin{equation}\nonumber
\bm \tilde{U}(s,y)\\ \bar{\tilde{U}}(s,y)\endm:=\bm
e^{-i(\Psi-\Psi_{\infty})_{1}(s)}U(s,y+\lambda_{\infty}(\mu-\mu_{\infty})(s))\\
e^{+i(\Psi-\Psi_{\infty})_{1}(s)}\bar{U}(s,y+\lambda_{\infty}(\mu-\mu_{\infty})(s))\endm
\end{equation}
We decompose
\begin{equation}\nonumber
\bm \tilde{U}(s,y)\\
\bar{\tilde{U}}(s,y)\endm=\sum_{i=1}^{6}\lambda_{i}\eta_{i,\text{proper}}+\bm
\tilde{U}(s,y)\\ \bar{\tilde{U}}(s,y)\endm_{dis}
\end{equation}
We can then infer the parameters $\lambda_{i},\,i=1,\ldots,5$ from
the orthogonality condition \eqref{orthogonality'}, while the
parameter $\lambda_{6}$ is governed by a suitable ODE. We now
carefully analyze these equations. First, for $j=2,\ldots,6$, we
have explicitly (recall \eqref{xiproper} as well as
\eqref{orthogonality})
\begin{equation}\nonumber
\la \bm
e^{i(\Psi_{\infty}-\Psi)(t)}U(t,y+\lambda_{\infty}(\mu-\mu_{\infty})(t))\\e^{-i(\Psi_{\infty}-\Psi)(t)}\bar{U}(t,y+\lambda_{\infty}(\mu-\mu_{\infty})(t))\endm,
\xi_{j,\text{proper}}(\nu(t)y)\ra=0
\end{equation}
This may be recast as\footnote{We also use the notation
$(\Psi_{\infty}-\Psi)_{1}(s):=-(\Psi-\Psi_{\infty})_{1}(s)$}
\begin{equation}\nonumber\begin{split}
&\la\bm\tilde{U}\\ \bar{\tilde{U}}\endm,
\xi_{j,\text{proper}}(\nu(t)y)\ra=\la \bm
e^{i(\Psi_{\infty}-\Psi)_{1}(t)}U(t,y+\lambda_{\infty}(\mu-\mu_{\infty})(t))\\e^{-i(\Psi_{\infty}-\Psi)_{1}(t)}\bar{U}(t,y+\lambda_{\infty}(\mu-\mu_{\infty})(t))\endm,
\xi_{j,\text{proper}}(\nu(t)y)\ra\\&=\la\bm
e^{i(\Psi_{\infty}-\Psi)_{1}(t)}(1-e^{i(\Psi_{\infty}-\Psi)_{2}(t)})&0\\0&e^{-i(\Psi_{\infty}-\Psi)_{1}(t)}(1-e^{-i(\Psi_{\infty}-\Psi)_{2}(t)})\endm\\&\hspace{5cm}\bm
U(t,y+\lambda_{\infty}(\mu-\mu_{\infty})(t))\\
\bar{U}(t,y+\lambda_{\infty}(\mu-\mu_{\infty})(t))\endm,
\xi_{j,\text{proper}}(\nu(t)y) \ra
\end{split}\end{equation}
The intuition here, to be made precise below, is that
$(\Psi-\Psi_{\infty})_{2}(s,y)$, when localized in $y$, decays
quite rapidly in $s$.  Our first task is filtering out the
$\lambda_{i},\,i=1,\ldots,5$ from this relation, while avoiding
$\lambda_{6}$ if possible. From the above we have
\begin{equation}\nonumber\begin{split}
&\sum_{i=1}^{5}\lambda_{i}\la
\eta_{i,\text{proper}},\xi_{l,\text{proper}}(\nu(t)y)\ra
+\lambda_{6}\la
\eta_{6,\text{proper}},\xi_{l,\text{proper}}(\nu(t)y)\ra+\la \bm
\tilde{U}\\
\bar{\tilde{U}}\endm_{dis},\xi_{l,\text{proper}}(\nu(t)y)\ra\\&
=\la\bm
e^{i(\Psi_{\infty}-\Psi)_{1}(t)}(1-e^{i(\Psi_{\infty}-\Psi)_{2}(t)})&0\\0&e^{-i(\Psi_{\infty}-\Psi)_{1}(t)}(1-e^{-i(\Psi_{\infty}-\Psi)_{2}(t)})\endm\\&\hspace{8cm}\bm
U(t,y+\lambda_{\infty}(\mu-\mu_{\infty})(t))\\
\bar{U}(t,y+\lambda_{\infty}(\mu-\mu_{\infty})(t))\endm,
\xi_{l,\text{proper}}(\nu(t)y) \ra,
\end{split}\end{equation}
where $l=2,\ldots,6$. Commence with the case $l=2$. Observe that
\begin{equation}\nonumber
\la
\eta_{6,\text{proper}},\xi_{2,\text{proper}}(\nu(t)y)\ra=0=\la\eta_{2,\text{proper}},\xi_{2,\text{proper}}(\nu(t)y)\ra
\end{equation}
Hence
\begin{equation}\nonumber\begin{split}
&\sum_{i=1}^{5}\lambda_{i}\la
\eta_{i,\text{proper}},\xi_{2,\text{proper}}(\nu(t)y)\ra
+\lambda_{6}\la
\eta_{6,\text{proper}},\xi_{2,\text{proper}}(\nu(t)y)\ra=
\sum_{i\neq 2,6}\lambda_{i}\la
\eta_{i,\text{proper}},\xi_{2,\text{proper}}(\nu(t)y)\ra
\end{split}\end{equation}
Next, we observe that
\begin{equation}\nonumber
\la
\eta_{6,\text{proper}},\xi_{3,\text{proper}}(\nu(t)y)\ra=0=\la\eta_{2,\text{proper}},\xi_{3,\text{proper}}(\nu(t)y)\ra
\end{equation}
Thus we have
\begin{equation}\nonumber
\sum_{i=1}^{5}\lambda_{i}\la
\eta_{i,\text{proper}},\xi_{3,\text{proper}}(\nu(t)y)\ra
+\lambda_{6}\la
\eta_{6,\text{proper}},\xi_{3,\text{proper}}(\nu(t)y)\ra=
\sum_{i\neq 2,6}\lambda_{i}\la
\eta_{i,\text{proper}},\xi_{3,\text{proper}}(\nu(t)y)\ra
\end{equation}
Further, observe that for reasons of parity, we have
\begin{equation}\nonumber
\la
\eta_{6,\text{proper}},\xi_{4,\text{proper}}(\nu(t)y)\ra=0=\la\eta_{2,\text{proper}},\xi_{4,\text{proper}}(\nu(t)y)\ra
\end{equation}
The conclusion is that
\begin{equation}\nonumber
\sum_{i=1}^{5}\lambda_{i}\la
\eta_{i,\text{proper}},\xi_{4,\text{proper}}(\nu(t)y)\ra
+\lambda_{6}\la
\eta_{6,\text{proper}},\xi_{4,\text{proper}}(\nu(t)y)\ra=
\sum_{i\neq 2,6}\lambda_{i}\la
\eta_{i,\text{proper}},\xi_{4,\text{proper}}(\nu(t)y)\ra
\end{equation}
One concludes similarly for the inner product with
$\xi_{6}(\nu(t)y)$.  Next we consider the inner product with
$\xi_{5,\text{proper}}(\nu(t)y)$. We note the following inner
product relations:
\begin{equation}\nonumber
\la \eta_{1,\text{proper}},\xi_{5,\text{proper}}(\nu(t)y)\ra
=0,\,\la \eta_{2,\text{proper}}, \xi_{5,\text{proper}}(\nu(t)y)\ra
= a(t),\,\la
\eta_{3,\text{proper}},\xi_{5,\text{proper}}(\nu(t)y)\ra =0
\end{equation}
\begin{equation}\nonumber
\la \eta_{4,\text{proper}},\xi_{5,\text{proper}}(\nu(t)y)\ra =0,
\la \eta_{5,\text{proper}},\xi_{5,\text{proper}}(\nu(t)y)\ra=0,\la
\eta_{6,\text{proper}},\xi_{5,\text{proper}}(\nu(t)y)\ra=b(t)
\end{equation}
In the immediately preceding the function $a(t)$ can be forced to
vanish nowhere upon choosing $\delta$ small enough. We don't need
this information concerning $b(t)$. We can now infer the following
relations: first
\begin{equation}\label{lambda2}
\lambda_{2}(t)=(\nu(t)-1)\la \bm \tilde{U}\\
\overline{\tilde{U}}\endm_{dis}, \phi(t)\ra + \sum _{i\neq
2}a_{2i}(t)\lambda_{i}(t),
\end{equation}
where $\phi(t)$ denotes a certain time dependent vector-valued
Schwartz function (with uniform decay estimates for all its
derivatives, including its time derivatives), while the parameters
$a_{2i}(t),\,i\neq 2,6$ decay at the same rate as
$\phi(t,x)(\Psi-\Psi_{\infty})_{2}(t,x))$, for another Schwartz
function $\phi(t,x)$ (we shall henceforth denote (vector valued)
Schwartz functions (with respect to the 2nd argument) in this
manner, without distinguishing between them, it being understood
that they satisfy uniform-in-time decay estimates,
including all their derivatives.)\\
In the same vein, the preceding calculations allow us to infer
that the coefficients $\lambda_{i}(t),\,i\neq 2,6$ satisfy the
relations
\begin{equation}\label{lambdai}
\lambda_{i}(t)= (\nu(t)-1)\la \bm
\tilde{U}\\\overline{\tilde{U}}\endm_{dis}, \phi(t)\ra
+\lambda_{2}(t)a_{i2}(t)+\lambda_{6}(t)a_{i6}(t),
\end{equation}
where the coefficients $a_{i2}(t), a_{i6}(t)$ satisfy the same
estimates as $a_{2i}(t)$(with $i\neq 2,6$) above. Of course if we
substitute \eqref{lambda2} here we can get rid of the 2nd term on
the right (choosing $\delta$ small enough). We have used the fact
that
\begin{equation}\nonumber
\la \bm\tilde{U}\\
\overline{\tilde{U}}\endm_{dis},
\xi_{l,\text{proper}}(\nu(t).)\ra =\la \bm\tilde{U}\\
\overline{\tilde{U}}\endm_{dis},
\xi_{l,\text{proper}}(\nu(t)y)-\xi_{l,\text{proper}}(.)\ra,\,l\neq
1
\end{equation}
In order to complete the control of the root part, we thus need to
finally consider $\lambda_{6}(t)$, which controls the contribution
of the 'exotic mode'. This we filter out by means of
\begin{equation}\nonumber
2\kappa_{2}\lambda_{6}(t)=\la \bm\tilde{U}\\
\overline{\tilde{U}}\endm, \xi_{1,\text{proper}}\ra
\end{equation}
Upon differentiation, this relation implies the following:
\begin{equation}\label{lambda6'}\begin{split}
&i2\kappa_{2}\dot{\lambda_{6}}(t)=\la i\partial_{t}\bm \tilde{U}\\
\bar{\tilde{U}}\endm, \xi_{1,\text{proper}}\ra\\& = \la \bm \tilde{U}\\
\bar{\tilde{U}}\endm, \bm
\partial_{t}[\Psi-\Psi_{\infty}]_{1}&0\\0&-\partial_{t}[\Psi-\Psi_{\infty}]_{1}\endm\xi_{1,\text{proper}}\ra+i\la
\partial_{t}[\lambda_{\infty}(\mu-\mu_{\infty})]\bm
\partial_{x}\tilde{U}\\ \overline{\partial_{x}\tilde{U}}\endm,
\xi_{1,\text{proper}}\ra\\& +\la
\bm e^{-i(\Psi-\Psi_{\infty})_{1}(t)}i\partial_{t}U(t,y+\lambda_{\infty}(\mu-\mu_{\infty})(t))\\
\overline{
-e^{i(\Psi-\Psi_{\infty})_{1}(t)}i\partial_{t}U(t,y+\lambda_{\infty}(\mu-\mu_{\infty})(t))}\endm,\xi_{1,\text{proper}}\ra
\end{split}\end{equation}
We now carefully analyze each of the three expressions on the
right. The key is to show that no quantity
morally\footnote{Observe that for example the quantity
$\partial_{t}(\Psi-\Psi_{\infty})(t,y)$ decays like $\nu(t)-1$.}
of the form $(\nu(t)-1)^{a}\lambda_{6}(t)$,
$(\nu(t)-1)^{a}\lambda_{2}(t)$, $a=1,2$, occurs, as this would
sabotage any attempt at controlling $\lambda_{6}$ by means of ODE
techniques, on account of the estimates \eqref{modulationasympto}.
This appears to require a lot of careful bookkeeping: start with
the first expression on the right. We have
\begin{equation}\nonumber\begin{split}
&\la \bm \tilde{U}\\
\bar{\tilde{U}}\endm, \bm
\partial_{t}[\Psi-\Psi_{\infty}]_{1}&0\\0&-\partial_{t}[\Psi-\Psi_{\infty}]_{1}\endm\xi_{1,\text{proper}}\ra
\\&=\sum_{i\neq 2,6}\lambda_{j}(t)\la \eta_{j,\text{proper}}, \bm
\partial_{t}[\Psi-\Psi_{\infty}]_{1}&0\\0&-\partial_{t}[\Psi-\Psi_{\infty}]_{1}\endm\xi_{1,\text{proper}}\ra
\\&^+\la \bm \tilde{U}\\
\bar{\tilde{U}}\endm_{dis}, \bm
\partial_{t}[\Psi-\Psi_{\infty}]_{1}&0\\0&-\partial_{t}[\Psi-\Psi_{\infty}]_{1}\endm\xi_{1,\text{proper}}\ra
\end{split}\end{equation}
Next, write
\begin{equation}\nonumber\begin{split}
&\la
\partial_{t}[\lambda_{\infty}(\mu-\mu_{\infty})]\bm
\partial_{x}\tilde{U}\\ \overline{\partial_{x}\tilde{U}}\endm,
\xi_{1,\text{proper}}\ra=\partial_{t}[\lambda_{\infty}(\mu-\mu_{\infty})]\sum_{j=1}^{6}\lambda_{j}\la
\partial_{x}\eta_{j,\text{proper}},\xi_{1,\text{proper}}\ra\\&+\partial_{t}[\lambda_{\infty}(\mu-\mu_{\infty})]\la
\bm
\partial_{x}\tilde{U}_{dis}\\
\overline{\partial_{x}\tilde{U}}_{dis}\endm,
\xi_{1,\text{proper}}\ra
\end{split}\end{equation}
Carefully observe from \eqref{modulationasympto} that we get
$|\partial_{t}[\lambda_{\infty}(\mu-\mu_{\infty})](t)|\lesssim \la
t\ra^{-\frac{3}{2}+\delta_{1}}$. Finally, consider the
contribution of
\begin{equation}\label{2}
\la
\bm e^{-i(\Psi-\Psi_{\infty})_{1}(t)}i\partial_{t}U(t,y+\lambda_{\infty}(\mu-\mu_{\infty})(t))\\
\overline{
-e^{i(\Psi-\Psi_{\infty})_{1}(t)}i\partial_{t}U(t,y+\lambda_{\infty}(\mu-\mu_{\infty})(t))}\endm,\xi_{1,\text{proper}}\ra
\end{equation}
This we reformulate using \eqref{Ueqn'}. Observe that we have
\begin{equation}\nonumber\begin{split}
&\bm
\partial_{y}^{2}-1+3\nu^{2}\phi_{0}^{4}(\nu(t)y)&2\nu^{2}\phi_{0}^{4}(\nu(t)y)\\-2\nu^{2}\phi_{0}^{4}(\nu(t)y)&-\partial_{y}^{2}+1-3\nu^{2}\phi_{0}^{4}(\nu(t)y)\endm
\bm
e^{i(\Psi_{\infty}-\Psi)_{1}(t)}U(y+\lambda_{\infty}[\mu-\mu_{\infty}](t),t)\\
e^{-i(\Psi_{\infty}-\Psi)_{1}(t)}\bar{U}(y+\lambda_{\infty}[\mu-\mu_{\infty}](t),t)\endm\\
&=\bm e^{i(\Psi_{\infty}-\Psi)_{1}(t)}&0\\0&
e^{-i(\Psi_{\infty}-\Psi)_{1}(t)}\endm\bm
\partial_{y}^{2}-1+3\nu^{2}\phi_{0}^{4}(\nu(t)y)&2\nu^{2}e^{2i(\Psi-\Psi_{\infty})_{1}}\phi_{0}^{4}(\nu(t)y)\\-2\nu^{2}e^{-2i(\Psi-\Psi_{\infty})_{1}}\phi_{0}^{4}(\nu(t)y)&-\partial_{y}^{2}+1-3\nu^{2}\phi_{0}^{4}(\nu(t)y)\endm\bm
U\\\bar{U}\endm(.)
\end{split}\end{equation}
This shows that we can reformulate \eqref{2} as follows:
\begin{equation}\nonumber\begin{split}
&\la\bm e^{i(\Psi_{\infty}-\Psi)_{1}(t)}&0\\0&
e^{-i(\Psi_{\infty}-\Psi)_{1}(t)}\endm[i\partial_{t}+\calH_{1}(t)](\bm U\\
\bar{U}\endm(y+\lambda_{\infty}[\mu-\mu_{\infty}](t),t)),\xi_{1,\text{proper}}\ra\\&-\la
\bm
e^{i(\Psi_{\infty}-\Psi)_{1}(t)}U(y+\lambda_{\infty}[\mu-\mu_{\infty}](t),t)\\
e^{-i(\Psi_{\infty}-\Psi)_{1}(t)}\bar{U}(y+\lambda_{\infty}[\mu-\mu_{\infty}](t),t)\endm,
\tilde{\calH}(t)^{*}\xi_{1,\text{proper}}\ra,
\end{split}\end{equation}
where
\begin{equation}\nonumber
\tilde{\calH}(t)^{*}=\bm
\partial_{y}^{2}-1+3\nu^{2}\phi_{0}^{4}(\nu(t)y)&-2\nu^{2}\phi_{0}^{4}(\nu(t)y)\\+2\nu^{2}\phi_{0}^{4}(\nu(t)y)&-\partial_{y}^{2}+1-3\nu^{2}\phi_{0}^{4}(\nu(t)y)\endm
\end{equation}
\begin{equation}\nonumber
\calH_{1}(t)=\bm\partial_{y}^{2}-1+3\nu^{2}\phi_{0}^{4}(\nu(t)y)&2\nu^{2}e^{2i(\Psi-\Psi_{\infty})_{1}}\phi_{0}^{4}(\nu(t)y)\\-2\nu^{2}e^{-2i(\Psi-\Psi_{\infty})_{1}}\phi_{0}^{4}(\nu(t)y)&-\partial_{y}^{2}+1-3\nu^{2}\phi_{0}^{4}(\nu(t)y)\endm
\end{equation}
Moreover, we have
\begin{equation}\nonumber
\tilde{\calH}(t)^{*}\xi_{1,\text{proper}}=\bm
3\nu^{2}\phi_{0}^{4}(\nu(t)y)-3\phi_{0}^{4}(y)&-2\nu^{2}\phi_{0}^{4}(\nu(t)y)+2\phi_{0}^{4}\\-[-2\nu^{2}\phi_{0}^{4}(\nu(t)y)+2\phi_{0}^{4}]&
-[3\nu^{2}\phi_{0}^{4}(\nu(t)y)-3\phi_{0}^{4}(y)]\endm\xi_{1,\text{proper}},
\end{equation}
which is of the form $\bm \alpha\\ -\alpha\endm$. This reveals
that
\begin{equation}\nonumber\begin{split}
&\la \bm
e^{i(\Psi_{\infty}-\Psi)_{1}(t)}U(y+\lambda_{\infty}[\mu-\mu_{\infty}](t),t)\\
e^{-i(\Psi_{\infty}-\Psi)_{1}(t)}\bar{U}(y+\lambda_{\infty}[\mu-\mu_{\infty}](t),t)\endm,
\tilde{\calH}(t)^{*}\xi_{1,\text{proper}}\ra\\&=\sum_{j\neq
2,6}\lambda_{j}(t)\la\eta_{j,\text{proper}},\tilde{\calH}(t)^{*}\xi_{1,\text{proper}}\ra+\la
\bm \tilde{U}\\
\bar{\tilde{U}}\endm_{dis},\tilde{\calH}(t)^{*}\xi_{1,\text{proper}}\ra
\end{split}\end{equation}
Also, note that
$\tilde{\calH}(t)^{*}\xi_{1,\text{proper}}=(\nu(t)-1)\phi(t,.)$
for a suitable (vector-valued) Schwartz function $\phi(t,.)$. We
now need to carefully analyze the expression
\begin{equation}\nonumber
\la\bm e^{i(\Psi_{\infty}-\Psi)_{1}(t)}&0\\0&
e^{-i(\Psi_{\infty}-\Psi)_{1}(t)}\endm([i\partial_{t}+\calH_{1}(t)]\bm U\\
\bar{U}\endm)(y+\lambda_{\infty}[\mu-\mu_{\infty}](t),t),\xi_{1,\text{proper}}\ra
\end{equation}
We reformulate this as\footnote{Here we define
$\tilde{\calH}_{1}(s)$ like $\calH(s)$ but with
$\Psi-\Psi_{\infty}$ replaced by $(\Psi-\Psi_{\infty})_{1}$.}
\begin{equation}\nonumber
\la [i\partial_{t}+\tilde{\calH}_{1}(t)](\bm U\\
\bar{U}\endm(y,t)),\bm e^{-i(\Psi_{\infty}-\Psi)_{1}(t)}&0\\0&
e^{+i(\Psi_{\infty}-\Psi)_{1}(t)}\endm\xi_{1,\text{proper}}(y-\lambda_{\infty}[\mu-\mu_{\infty}](t))\ra
\end{equation}
and use \eqref{Ueqn'}, in which we schematically write the
righthand side as $\dot{\pi}\partial_{\pi}W+N(U,\pi)$, where the
first summand refers to those expressions which only involve the
modulation parameters and their derivatives, and not (explicitly)
the radiation. Then we can schematically rewrite the above as
\begin{equation}\label{3}
\la \dot{\pi}\partial_{\pi}W,\tilde{\xi}_{1,\text{proper}}\ra+\la
N(U,\pi),\tilde{\xi}_{1,\text{proper}}\ra+\la(e^{i(\Psi-\Psi_{\infty})_{2}(t)}-1)
\bm U\\ \bar{U}\endm, \phi(t,.)\ra,
\end{equation}
where we have introduced the notation
\begin{equation}\nonumber
\bm e^{-i(\Psi_{\infty}-\Psi)_{1}(t)}&0\\0&
e^{-i(\Psi-\Psi_{\infty})_{1}(t)}\endm
\xi_{1,\text{proper}}(y-\lambda_{\infty}(\mu-\mu_{\infty})(t))=\tilde{\xi}_{1,\text{proper}}
\end{equation}
We now carefully analyze the first two expressions in \eqref{3},
again in order to check that these don't implicitly contain
expressions of the form $(\nu(t)-1)^{a}\lambda_{2}(t)$,
$(\nu(t)-1)^{a}\lambda_{6}(t)$, $a=1,2$. The third expression in
\eqref{3}
turns out to be small, as we'll see later on. \\
Now expand the schematic expression $\la
\dot{\pi}\partial_{\pi}W,\tilde{\xi}_{1,\text{proper}}\ra$,
invoking \eqref{Ueqn'}. First one obtains
\begin{equation}\nonumber
\la
(\dot{\lambda}\lambda^{-1}-\beta\nu^{2})(\tileta_{2}-\frac{\beta}{2}\tileta_{5}+\omega\tileta_{4}),\tilde{\xi}_{1,\text{proper}}\ra
\end{equation}
We note that the vectors $\tilde{\eta}_{j}$ appearing here carry
the phases $e^{\pm i(\Psi-\Psi_{\infty})}$. Thus by modifying them
by errors of size $O((e^{i(\Psi_{\infty}-\Psi)_{2}}-1)\phi(t,.))$,
we can replace these phases by $e^{\pm
i(\Psi-\Psi_{\infty})_{1}}$. By abuse of notation we shall refer
to these vectors again as $\tilde{\eta}_{j}$. Then we potentially
have $\la \tileta_{2},\tilde{\xi}_{1,\text{proper}}\ra\neq 0$.
Using the 5th modulation equation, we recall that
\begin{equation}\nonumber
4i\kappa_{2}(\dot{\lambda}\lambda^{-1}-\beta\nu^{2})=-\la
N(U,\pi),\tilde{\xi}_{5}\ra+\la \bm U\\
\bar{U}\endm,(i\partial_{t}+\calH(t)^{*})\tilde{\xi}_{5}\ra,
\end{equation}
where
\begin{equation}\nonumber
\tilde{\xi}_{5}(t,y)= \bm
e^{i(\Psi-\Psi_{\infty})(t,y)}\sqrt{\nu}(t)\phi_{0}(\nu(t)[y-\lambda_{\infty}(\mu-\mu_{\infty})(t)]\\e^{-i(\Psi-\Psi_{\infty})(t)}\sqrt{\nu}(t)\phi_{0}(\nu(t)[y-\lambda_{\infty}(\mu-\mu_{\infty})(t)]\endm
\end{equation}
Calculate
\begin{equation}\nonumber\begin{split}
&i\partial_{t}\bm
e^{i(\Psi-\Psi_{\infty})(t,y)}\sqrt{\nu}(t)\phi_{0}(\nu(t)[y-\lambda_{\infty}(\mu-\mu_{\infty})(t)]\\e^{-i(\Psi-\Psi_{\infty})(t,y)}\sqrt{\nu}(t)\phi_{0}(\nu(t)[y-\lambda_{\infty}(\mu-\mu_{\infty})(t)]\endm
\\&=i\bm i\partial_{t}(\Psi-\Psi_{\infty})(t,y)
e^{i(\Psi-\Psi_{\infty})(t,y)}\sqrt{\nu}(t)\phi_{0}(\nu(t)[y-\lambda_{\infty}(\mu-\mu_{\infty})(t)]\\-i\partial_{t}(\Psi-\Psi_{\infty})(t,y)e^{-i(\Psi-\Psi_{\infty})(t,y)}\sqrt{\nu}(t)\phi_{0}(\nu(t)[y-\lambda_{\infty}(\mu-\mu_{\infty})(t)]\endm
\\&+i\sqrt{\nu}(t)[-\partial_{t}[\lambda(t)(\mu-\mu_{\infty})(t)]+y\dot{\nu}(t)]\bm
e^{i(\Psi-\Psi_{\infty})(t,y)}\nabla\phi_{0}(\nu(t)[y-\lambda_{\infty}(\mu-\mu_{\infty})(t)]\\e^{-i(\Psi-\Psi_{\infty})(t,y)}\nabla\phi_{0}(\nu(t)[y-\lambda_{\infty}(\mu-\mu_{\infty})(t)]\endm
\\&+\frac{i\dot{\nu}(t)}{2\sqrt{\nu}(t)}\bm
e^{i(\Psi-\Psi_{\infty})(t,y)}\sqrt{\nu}(t)\phi_{0}(\nu(t)[y-\lambda_{\infty}(\mu-\mu_{\infty})(t)]\\e^{-i(\Psi-\Psi_{\infty})(t,y)}\sqrt{\nu}(t)\phi_{0}(\nu(t)[y-\lambda_{\infty}(\mu-\mu_{\infty})(t)]\endm
\end{split}\end{equation}
Further note that for a certain $k\neq 1$
\begin{equation}\nonumber\begin{split}
&\calH(t)^{*}\tilde{\xi}_{5}\\&= \bm
\partial_{y}^{2}-1+3\nu^{2}\phi_{0}^{4}(\nu(t)[y+\lambda_{\infty}(\mu_{\infty}-\mu)])&-2\nu^{2}\phi_{0}^{4}(\nu(t)[y+\lambda_{\infty}(\mu_{\infty}-\mu)])e^{2i(\Psi-\Psi_{\infty})}\\
+2\nu^{2}\phi_{0}^{4}(\nu(t)[y+\lambda_{\infty}(\mu_{\infty}-\mu)])e^{2i(\Psi-\Psi_{\infty})}&-[\partial_{y}^{2}-1+3\nu^{2}\phi_{0}^{4}(\nu(t)[y+\lambda_{\infty}(\mu_{\infty}-\mu)])]\endm
\tilde{\xi}_{5}\\
&=\bm
\nu^{2}-1&0\\0&-(\nu^{2}-1)\endm\tilde{\xi}_{5}+\nu^{2}\tilde{\xi}_{k}+O(\partial_{y}^{1,2}[\Psi-\Psi_{\infty}]\phi(t,.))
\end{split}\end{equation}
Then we observe that
\begin{equation}\nonumber\begin{split}
&\la \bm U\\
\bar{U}\endm, \bm
\nu^{2}-1&0\\0&-(\nu^{2}-1)\endm\tilde{\xi}_{5}+\nu^{2}\tilde{\xi}_{k}\ra
=\sum_{j\neq 2,6}\lambda_{j}\la \eta_{j,\text{proper}}, \bm
\nu^{2}-1&0\\0&-(\nu^{2}-1)\endm\xi_{5,\text{proper}}(\nu(t)y)\ra\\&+\la \bm \tilde{U}\\
\bar{\tilde{U}}\endm_{dis}, \bm
\nu^{2}-1&0\\0&-(\nu^{2}-1)\endm\xi_{5,\text{proper}}(\nu(t)y)\ra+\la(e^{i(\Psi-\Psi_{\infty})_{2}(t,.)}-1)\bm
U\\ \bar{U}\endm, \phi(t,.)\ra
\end{split}\end{equation}
Continue by observing that $\la
\tilde{\eta}_{5,4,1},\tilde{\xi}_{1,\text{proper}}\ra =0$,
provided we abuse notation and change the phase in
$\tilde{\eta}_{i}$ to $(\Psi-\Psi_{\infty})_{1}$, which generates
errors of the type
$O((e^{i(\Psi-\Psi_{\infty})_{2}(t)}-1)\phi(t,.))$. Continuing in
this fashion, we note that $\la \tilde{\eta}_{3},
\xi_{1,\text{proper}}\ra \neq 0$ generically (but this function
will only be of size $O((\nu(t)-1))$. From the 3rd modulation
equation we get
\begin{equation}\nonumber
-i\kappa_{1}\nu(\dot{\mu}\lambda_{\infty}-2\omega\nu)=-\la
N(U,\pi), \tilde{\xi}_{4}\ra +\la\bm U\\ \bar{U}\endm,
(i\partial_{s}+\calH(s)^{*})\tilde{\xi}_{4}\ra
\end{equation}
Write as before
\begin{equation}\nonumber\begin{split}
&\partial_{s}\tilde{\xi}_{4}=\bm
i\partial_{s}[\Psi-\Psi_{\infty}](s)&0\\0&-i\partial_{s}[\Psi-\Psi_{\infty}](s)\endm
\tilde{\xi}_{4}+[\dot{\nu}y+\partial_{s}[\lambda(\mu_{\infty}-\mu)]]\partial_{x}\tilde{\xi}_{4}\\
&+O(\dot{\nu}\phi(t,.))
\end{split}\end{equation}
Moreover, we have (for a suitable $j\neq 1$)
\begin{equation}\nonumber\begin{split}
&\calH(s)^{*}\tilde{\xi}_{4}\\&=\bm
\partial_{y}^{2}-1+3\nu^{2}(s)\phi_{0}^{4}(\nu(s)(y+\lambda_{\infty}(\mu_{\infty}-\mu)))&-2\nu^{2}\phi_{0}^{4}(\nu(s)(y+\lambda_{\infty}(\mu_{\infty}-\mu)))e^{2i(\Psi-\Psi_{\infty})(s)}\\
+2\nu^{2}\phi_{0}^{4}(\nu(s)(y+\lambda_{\infty}(\mu_{\infty}-\mu)))e^{-2i(\Psi-\Psi_{\infty})(s)}&-(\partial_{y}^{2}-1+3\nu^{2}(s)\phi_{0}^{4}(\nu(s)(y+\lambda_{\infty}(\mu_{\infty}-\mu)))\endm\tilde{\xi}_{4}\\
&=\bm
\nu^{2}-1&0\\0&-(\nu^{2}-1)\endm\tilde{\xi}_{4}+\nu^{2}\tilde{\xi}_{j}+O(\partial_{y}^{1,2}[\Psi-\Psi_{\infty}]_{2}\phi(t,.))
\end{split}\end{equation}
Thus we obtain
\begin{equation}\nonumber\begin{split}
&\la \calH(s)^{*}\tilde{\xi}_{4}, \bm U\\
\bar{U}\endm\ra = \sum_{j\neq 2,6}\lambda_{j}\la
\eta_{j,\text{proper}},\bm
\nu^{2}-1&0\\0&-(\nu^{2}-1)\endm\xi_{4,\text{proper}}(\nu(s)y)\ra\\&+\la
\bm \tilde{U}\\ \bar{\tilde{U}}\endm_{dis},\bm
\nu^{2}-1&0\\0&-(\nu^{2}-1)\endm\xi_{4,\text{proper}}(\nu(s)y)\ra+O(\partial_{y}^{1,2}[\Psi-\Psi_{\infty}]_{2}\phi(t,.))\\
&+O((e^{i(\Psi-\Psi_{\infty})_{2}}-1)\phi(t,.))
\end{split}\end{equation}
The preceding observations allow us to control the expression $\la
\dot{\pi}\partial_{\pi}W,\tilde{\xi}_{1,\text{proper}}\ra$ in
\eqref{3}. Observe that in the preceding we also generated the
(schematic) terms $(\nu(t)-1)\la N(U,\pi), \tilde{\xi}_{4}\ra$,
$(\nu(t)-1)\la N(U,\pi), \tilde{\xi}_{5}\ra$. Now consider the
terms at least quadratic with respect to the radiation in
\eqref{3}, i. e. the expression $\la N(U,\pi),
\tilde{\xi}_{1,\text{proper}}\ra$. We are predominantly concerned
with the quadratic contribution, which we spell out explicitly:
\begin{equation}\label{quadratic}
\la\bm
-3U^{2}\tilde{\phi}_{0}^{3}e^{i(\Psi_{\infty}-\Psi)}\nu^{\frac{3}{2}}-6|U|^{2}e^{i(\Psi-\Psi_{\infty})}\nu^{\frac{3}{2}}\tilde{\phi}_{0}^{3}-\bar{U}^{2}e^{3i(\Psi-\Psi_{\infty})}\nu^{\frac{3}{2}}\tilde{\phi}_{0}^{3}\\
3\bar{U}^{2}\tilde{\phi}_{0}^{3}e^{-i(\Psi_{\infty}-\Psi)}\nu^{\frac{3}{2}}+6|U|^{2}e^{-i(\Psi-\Psi_{\infty})}\nu^{\frac{3}{2}}\tilde{\phi}_{0}^{3}-U^{2}e^{-3i(\Psi-\Psi_{\infty})}\nu^{\frac{3}{2}}\tilde{\phi}_{0}^{3}\endm,
\tilde{\xi}_{1,\text{proper}}\ra
\end{equation}
where we recall $\tilde{\phi}_{0}$ is $\phi_{0}$  evaluated at
$\nu(s)y+\lambda(\mu_{\infty}-\mu)$. Now we substitute
\begin{equation}\nonumber
\bm \tilde{U}\\ \bar{\tilde{U}}\endm
=\sum_{i=1}^{5}\lambda_{i}\eta_{i,\text{proper}}+\lambda_{6}\eta_{6,\text{proper}}+
\bm \tilde{U}\\ \bar{\tilde{U}}\endm_{dis}
\end{equation}
Use \eqref{lambda2}, to reformulate this as\footnote{As usual,
$\phi(t,.)$ represents various Schwartz functions, which in
addition to all their derivatives, both with respect to $t$ and
$x$, satisfy uniform decay estimates. Also, $\partial_{t}\phi$ is
of size at most $\dot{\nu}$.}
\begin{equation}\nonumber\begin{split}
\bm \tilde{U}\\ \bar{\tilde{U}}\endm(t,.)= \sum_{i\neq
2,6}\lambda_{i}(t)(\eta_{i,\text{proper}}+&a_{2i}(t)\eta_{2,\text{proper}})+\lambda_{6}(t)(\eta_{6,\text{proper}}+a_{26}(t)\eta_{2,\text{proper}})\\
&+\bm \tilde{U}\\
\bar{\tilde{U}}\endm_{dis}(t,.)+(\nu(t)-1)\phi(t,.)\la
\bm \tilde{U}\\ \bar{\tilde{U}}\endm_{dis}(t,.), \phi(t,.)\ra \\
\end{split}\end{equation}
Back to \eqref{quadratic}, we first rewrite
\begin{equation}\nonumber
\la\bm
-U^{2}\tilde{\phi}_{0}^{3}e^{i(\Psi_{\infty}-\Psi)(t,.)}\nu^{\frac{3}{2}}\\
\overline{U^{2}}\tilde{\phi}_{0}^{3}e^{i(\Psi_{\infty}-\Psi)(t,.)}\nu^{\frac{3}{2}}\endm,
\tilde{\xi}_{1,\text{proper}}\ra = \la \bm
-\tilde{U}^{2}\phi_{0}^{3}\nu^{\frac{3}{2}}\\
\overline{\tilde{U}}^{2}\phi_{0}^{3}\nu^{\frac{3}{2}}\endm,
\xi_{1,\text{proper}}\ra+O((e^{i(\Psi-\Psi_{\infty})_{2}(t,.)}-1)\phi(t,.))
\end{equation}
Then note that schematically
\begin{equation}\nonumber\begin{split}
\frac{1}{2}\la \bm
-\tilde{U}^{2}\phi_{0}^{3}\nu^{\frac{3}{2}}\\
\overline{\tilde{U}}^{2}\phi_{0}^{3}\nu^{\frac{3}{2}}\endm,
\xi_{1,\text{proper}}\ra=&\nu^{\frac{3}{2}}\Im \la
\tilde{U}_{dis}^{2}\phi_{0}^{3}(\nu .), \phi_{0}\ra \\
&+2\lambda_{6}[\sum_{i\neq 2,6}\lambda_{i}\nu^{\frac{3}{2}}\Im \la
\eta_{i,\text{proper}}^{1}+a_{2i}\eta_{2,\text{proper}},\phi_{0}^{3}(\nu.)(\eta_{6,\text{proper}}+a_{26}(t)\eta_{2,\text{proper}})^{1}\ra\\
&+\Im\la \phi(t,.),
\tilde{U}_{dis}\ra]+\sum_{i,j}a_{ij}(t)\lambda_{i}\lambda_{j}(t)+\sum_{j=1,2}(\nu(t)-1)^{j}\la
\tilde{U}_{dis},\phi_{1}(t,.)\ra\la
\tilde{U}_{dis},\phi_{2}(t,.)\ra
\end{split}\end{equation}
The expression
\begin{equation}\nonumber
\la \bm -\overline{\tilde{U}}^{2}e^{3i(\Psi-\Psi_{\infty})}\phi_{0}^{3}(\nu(t)y)\nu^{\frac{3}{2}}\\
\tilde{U}^{2}e^{-3i(\Psi-\Psi_{\infty})}\phi_{0}^{3}(\nu(t)y)\nu^{\frac{3}{2}}\endm,\tilde{\xi}_{1,\text{proper}}\ra
\end{equation}
is handled similarly. Moreover, it is easily seen that
\begin{equation}\nonumber
\la \bm-|\tilde{U}|^{2}\nu^{\frac{3}{2}}\phi_{0}^{3}(\nu y)\\
|\tilde{U}|^{2}\nu^{\frac{3}{2}}\phi_{0}^{3}(\nu
y)\endm,\xi_{1,\text{proper}}\ra =0
\end{equation}
Finally, we can summarize the discussion following
\eqref{lambda6'} in the following schematic equality\footnote{The
first instance of $\la N(\tilde{U}_{dis}^{2},\pi), \phi_{1}\ra$
refers to the symplectic form $\la
\tilde{U}_{dis}^{2}-\overline{\tilde{U}_{dis}}^{2},
\phi_{1}\ra$.}:
\begin{equation}\label{lambda6}\begin{split}
&\dot{\lambda}_{6}=\lambda_{6}[\la \bm\tilde{U}\\
\bar{\tilde{U}}\endm_{dis},\phi\ra+(\nu-1)\la \bm\tilde{U}\\
\bar{\tilde{U}}\endm_{dis},\phi\ra
+\dot{\nu}+\lambda_{6}(\nu-1)+\la\phi(e^{i(\Psi-\Psi_{\infty})_{2}}-1),\psi\ra]\\&+\la N(\tilde{U}_{dis}^{2},\pi), \phi_{1}\ra +(\nu-1)\la N(\tilde{U}_{dis}^{2},\pi), \phi_{2}\ra+\sum_{a=1,2}(\nu-1)^{a}\la \bm \tilde{U}\\
\bar{\tilde{U}}\endm_{dis},\phi\ra +\dot{\nu}\la \bm \tilde{U}\\
\bar{\tilde{U}}\endm_{dis},\phi\ra\\&+(\nu-1)\la \bm \tilde{U}\\
\bar{\tilde{U}}\endm_{dis},\phi_{1}\ra\la \bm \tilde{U}\\
\bar{\tilde{U}}\endm_{dis},\phi_{2}\ra+O(\la
U^{3}+U^{5}+[(\nu-1)^{2}+\dot{\nu}]U^{2},\phi\ra)
\end{split}\end{equation}
As usual the functions $\phi, \psi$ etc represent Schwartz
functions (with respect to the spatial variable) with uniform
decay estimates in time. One easily checks that all these
functions have time derivatives decaying like $\dot{\nu}$. In the
arguments below, we shall omit the time dependence, as one easily
checks that any additional terms generated by this additional time
dependence of the $\phi$ etc (for example when performing
integrations by parts in $t$) can be handled by exactly the same
methods or are much simpler to estimate. We now impose the
condition $\lim_{s\rightarrow+\infty}\lambda_{6}(s)=0$.
Introducing the integrating factor
\begin{equation}\nonumber\begin{split}
\Lambda(t)=&\int_{t}^{\infty}\la \bm\tilde{U}(s,.)\\
\overline{\tilde{U}(s,.)}\endm_{dis},\phi\ra+(\nu-1)(s)\la \bm\tilde{U}(s,.)\\
\overline{\tilde{U}(s,.)}\endm_{dis},\phi\ra
+\dot{\nu}(s)+\lambda_{6}(s)(\nu(s)-1)\\&\hspace{9cm}+\la\phi(e^{i(\Psi-\Psi_{\infty})_{2}(s,.)}-1),\psi\ra
ds,
\end{split}\end{equation}
this leads to the following relation:
\begin{equation}\label{lambda6''}
\lambda_{6}(s)=-e^{-\Lambda(t)}\int_{t}^{\infty}e^{\Lambda(s)}[...]ds,
\end{equation}
where $[...]$ stands for the part on the righthand side of
\eqref{lambda6} without the expression $\lambda_{6}[...]$. The
equations \eqref{lambda2}, \eqref{lambdai}, \eqref{lambda6''}
completely govern the evolution of the root part $\bm \tilde{U}\\
\overline{\tilde{U}}\endm_{root}$. Thus, to conclude the
discussion of the radiation part, we need to describe the
evolution of $\bm \tilde{U}\\ \overline{\tilde{U}}\endm_{dis}$.
This is straightforward from Duhamel's principle. Recall from
Theorem~\ref{Main} that we need to match the initial data $\bm
U(0,.)\\ \overline{U(0,.)}\endm_{dis}=\bm A(.)\\ \bar{A}(.)\endm$.
To this end write
\begin{equation}\label{initialdecomp}
\bm U\\
\bar{U}\endm(0,.)=\sum_{i=1}^{6}\alpha_{i}\eta_{i,\text{proper}}+\bm
A\\ \bar{A}\endm
\end{equation}
The coefficients $\alpha_{i}$ here can be inferred from the
orthogonality relations \eqref{orthorelations}. Thus
schematically\footnote{We really get a linear combination of
expressions of the indicated form.}
we get $\alpha_{i}=\la \bm U\\
\bar{U}\endm(0,.), \xi_{k(i),\text{proper}}\ra=\la \bm \tilde{U}\\
\overline{\tilde{U}}\endm(0,.),
\tilde{\xi}_{k(i),\text{proper}}\ra$. Using our standard
decomposition we now get
\begin{equation}\nonumber\begin{split}
&\alpha_{i}=\la \bm \tilde{U}\\
\overline{\tilde{U}}\endm(0,.),
\tilde{\xi}_{k(i),\text{proper}}(0,.)\ra\\&
=\sum_{j}\lambda_{j}(0)\la \eta_{j,\text{proper}},
\tilde{\xi}_{k(i),\text{proper}}(0,.)\ra + \la \bm \tilde{U}\\
\overline{\tilde{U}}\endm_{dis}(0,.),
\tilde{\xi}_{k(i),\text{proper}}(0,.)-\xi_{k(i),\text{proper}}\ra\\
&=\sum_{j}\lambda_{j}(0)\la \eta_{j,\text{proper}},
\tilde{\xi}_{k(i),\text{proper}}(0,.)\ra+(\nu(0)-1)\la\bm \tilde{U}\\
\overline{\tilde{U}}\endm_{dis}(0,.),\phi\ra
\end{split}\end{equation}
Now, to prescribe the evolution of $\bm \tilde{U}\\
\overline{\tilde{U}}\endm_{dis}$, we refer to \eqref{Ueqn}.
Introduce the function
\begin{equation}\label{Uts}
\bm \tilde{U}^{(t)}(s,y)\\ \bar{\tilde{U}}^{(t)}(s,y)\endm:=\bm
e^{-i(\Psi-\Psi_{\infty})_{1}(t)}U(s,y+\lambda_{\infty}(\mu-\mu_{\infty})(t))\\
e^{+i(\Psi-\Psi_{\infty})_{1}(t)}\bar{U}(s,y+\lambda_{\infty}(\mu-\mu_{\infty})(t))\endm,\,\bm
\tilde{U}^{(t)}(t,y)\\ \bar{\tilde{U}}^{(t)}(t,y)\endm=\bm
\tilde{U}(t,y)\\ \bar{\tilde{U}}(t,y)\endm,
\end{equation}
Then we deduce the following equation:
\begin{equation}\label{4}\begin{split} &[i\partial_{s}+\bm
\partial_{y}^{2}-1+3\phi_{0}^{4}&
2\phi_{0}^{4}\\-\phi_{0}^{4}&-\partial_{y}^{2}+1-3\phi_{0}^{4}\endm]
\bm \tilde{U}^{(t)}\\ \overline{\tilde{U}^{(t)}}\endm(s,y)=\bm
e^{-i(\Psi-\Psi_{\infty})_{1}(t)}&0
\\ 0&
e^{i(\Psi-\Psi_{\infty})_{1}(t)}\endm [...]\\&+2\bm
0&-e^{2i(\Psi-\Psi_{\infty})_{1}(s)-2i(\Psi-\Psi_{\infty})_{1}(t)}+1
\\ e^{-2i(\Psi-\Psi_{\infty})_{1}(s)+2i(\Psi-\Psi_{\infty})_{1}(t)}-1&0\endm
\phi_{0}^{4}\bm \tilde{U}^{(t)}(s)\\
\overline{\tilde{U}^{(t)}(s)}\endm\\&+O([\nu^{2}-1]\phi_{0}^{4}U)+O([\mu_{\infty}-\mu](s)\lambda(s)\phi_{0}^{4}U)+O([e^{i(\Psi-\Psi_{\infty})_{2}(s,.)-1}]\phi_{0}^{4}U)
\end{split}\end{equation}
The quantity $[...]$ on the righthand side refers to the righthand
expression in \eqref{Ueqn} translated by the amount
$+\lambda_{\infty}(\mu-\mu_{\infty})(t)$ in the spatial variable,
but one uses the identifications
\begin{equation}\nonumber
4i\kappa_{2}(\dot{\lambda}\lambda^{-1}-\beta\nu^{2})=-\la N,
\tilde{\xi}_{5}\ra +\la
U,(i\partial_{s}+\calH(s)^{*})\tilde{\xi}_{5}\ra\,\text{etc},
\end{equation}
coming from the modulation equations, within\footnote{Recall that
we use the schematic notation $(i\partial_{s}+\calH(s))\bm U\\
\bar{U}\endm(s,.)=\dot{\pi}\partial_{\pi}W+N(U,\pi)$.}
$\dot{\pi}\partial_{\pi}W$; thus we replace the left hand
expressions by the ones on the right. We then project the
preceding equation onto the dispersive part, and invoke Duhamel's
principle, which results in the following equation
governing the evolution of $\bm \tilde{U}\\
\overline{\tilde{U}}\endm_{dis}$:
\begin{equation}\label{Udis}\begin{split}
\bm \tilde{U}\\
\overline{\tilde{U}}\endm_{dis}(t,.)=&e^{it\calH}P_{s}[\bm
e^{i(\Psi_{\infty}-\Psi)_{1}(t)}&0\\0&e^{-i(\Psi_{\infty}-\Psi)_{1}(t)}\endm[\bm
A(.+\lambda_{\infty}(\mu-\mu_{\infty})(t))\\\bar{A}(.+\lambda_{\infty}(\mu-\mu_{\infty})(t))\endm
\\&+\sum_{j=1}^{6}\alpha_{j}
[\eta_{j,\text{proper}}(.+\lambda_{\infty}(\mu-\mu_{\infty})(t))]]]-i\int_{0}^{t}e^{i(t-s)\calH}[...]_{dis}(s,.)ds,
\end{split}\end{equation}
in which $[...]$ refers to the righthand side of \eqref{4}. Also,
the coefficients $\alpha_{i}$ are given by the formula detailed
further above, i. e.
\begin{equation}\label{alphai}
\alpha_{i}=\sum_{j}\lambda_{j}(0)\la \eta_{j,\text{proper}},
\tilde{\xi}_{k(i),\text{proper}}(0,.)\ra+(\nu(0)-1)\la\bm \tilde{U}\\
\overline{\tilde{U}}\endm_{dis}(0,.),\phi\ra
\end{equation}
\\

{\bf{Summary:}} {\it{The coefficients $\lambda_{i}(t),\,i\neq 6$,
are given by the relation \eqref{lambda2}, \eqref{lambdai}, while
the coefficient $\lambda_{6}(t)$ which controls the 'exotic mode'
is given by the formula \eqref{lambda6''}. In particular, the
latter forces a value for $\lambda_{6}(0)$. The dispersive part of
the
'tweaked radiation part', namely $\bm \tilde{U}\\
\overline{\tilde{U}}\endm_{dis}$, is governed by \eqref{Udis},
upon fixing the condition $\bm U(0,.)\\
\overline{U(0,.)}\endm_{dis}=\bm A(.)\\ \overline{A(.)}\endm$.}}
\\

To conclude this subsection, we still need to {\bf{specify the
estimates to be satisfied by the radiation part}}; as for the root
part, we shall need
\begin{equation}\label{rootasympto}
\sup_{0\leq k\leq [\frac{N}{2}]}||\la t\ra
^{2-2\delta_{1}}\partial_{t}^{k}\lambda_{i}(t)||_{L^{M}}\lesssim
\delta^{2},
\end{equation}
where $\delta_{i},\,\delta, N, M$ are as in
\eqref{modulationasympto}. The reason why we don't work with
$L^{\infty}$ is again the compactness property\footnote{Of course
we can recover $L^{\infty}$ bounds for all derivatives with
exception of the top derivative from this information. We avoid
this distinction in order to simplify matters.}. As for the
dispersive part, let $C_{k}$ be sufficiently rapidly\footnote{As
usual the necessary rate of growth can be inferred from the
proof.} growing numbers, $1\leq k\leq N$; we shall impose
$25^{N}\delta_{2}<<\delta_{1}$. Then we need
\begin{equation}\label{disasympto}\begin{split}
&\sup_{0\leq k\leq N}\sup_{2i+j\leq k}\sup_{s\geq
0}C_{k}^{-1}||\la s
\ra^{\frac{1}{2}-25^{k}\delta_{2}}\partial_{s}^{i}\partial_{y}^{j}U(s,y)||_{L_{s}^{M}L_{y}^{M}}\lesssim
\delta,\\&\sup_{\phi\in\calA}\sup_{0\leq k\leq N}\sup_{2i+j\leq
k}\sup_{s\geq 0}C_{k}^{-1}||\la s
\ra^{1-20^{k}\delta_{2}}\phi\partial_{s}^{i}\partial_{y}^{j}U(s,y)||_{L_{s}^{M}L_{y}^{M}}\lesssim
\delta\\
&\sup_{\phi\in\calA}\sup_{s\geq 0}||\la
s\ra^{\frac{3}{2}-\delta_{3}}\phi
U(s,.)||_{L_{y}^{\infty}}\lesssim \delta,\,\sup_{s\geq
0}||CU(s,y)||_{L_{y}^{2}}\lesssim \delta,\\&\sup_{0\leq k\leq
N}\sup_{2i+j\leq k}\sup_{s\geq 0}C_{k}^{-1}||\la s
\ra^{-10^{k}\delta_{2}}\partial_{s}^{i}\partial_{y}^{j}U(s,y)||_{L_{s}^{M}L_{y}^{2}}\lesssim
\delta,
\end{split}\end{equation}
where $C$ refers to the standard {\bf{pseudo-conformal operator}}
$C=y+2is\partial_{y}$. Also, we denote by $\calA$ the set of all
Schwartz functions satisfying
\begin{equation}\nonumber
\sup_{i\leq 100N}\sup_{x\in{\mathbf{R}}}|\la
x\ra^{100}\partial_{x}^{i}\phi(x)|<1
\end{equation}

\section{The iterative step.}
\subsection{Deducing Theorem~\ref{Main} from a fixed point Theorem.}

It is now straightforward with our setup to formulate the
iterative\footnote{However, we shall not be able to construct the
solution by iteration alone.} step: We commence with a tuple of
functions and parameters, as follows:
\begin{equation}\nonumber
\{\bm \tilde{U}\\
\overline{\tilde{U}}\endm_{dis}(.,.),\lambda_{1}(.),\ldots,\lambda_{6}(.),\nu_{1}(.),\beta_{1}(.),\omega_{1}(.),\mu_{1}(.),\gamma_{1}(.),a_{\infty}^{(1)},b_{\infty}^{(1)},v_{\infty}^{(1)},
y_{\infty}^{(1)}\},
\end{equation}
We define here\footnote{Keep in mind that the parameters
$\lambda(0),\beta(0)$ etc are fixed throughout.}
$a_{\infty}^{(1)}:=a_{\infty}-\lambda(0)$,
$b_{\infty}^{(1)}=b_{\infty}-(\beta\lambda)(0)$,
$y_{\infty}^{(1)}=y_{\infty}-(\lambda\mu)(0)$,
$v_{\infty}^{(1)}=v_{\infty}-\frac{(\beta\lambda\mu)(0)}{2}-\omega(0)$.
The functions $\nu_{1}(s)$, $\beta_{1}(s)$ stand for $\nu(s)-1$,
$(\beta\nu-b_{\infty}\lambda_{\infty}^{-1})(s)$, respectively,
while the functions $\mu_{1}(s)$, $\omega_{1}(s)$, $\gamma_{1}(s)$
stand for
$(\mu(s)-\frac{2v_{\infty}s+y_{\infty}}{a_{\infty}+b_{\infty}s})$,
$(\omega-\frac{c_{\infty}}{\lambda_{\infty}})(s)$,
$\gamma(s)-s-\frac{c_{\infty}^{2}}{\lambda_{\infty}^{2}(s)}$,
respectively, see the discussion in the last section; we let
$c_{\infty}=v_{\infty}a_{\infty}-\frac{b_{\infty}y_{\infty}}{2}$.
Of course we
require the orthogonality conditions $\la \bm \tilde{U}\\
\overline{\tilde{U}}\endm_{dis}, \xi_{i,\text{proper}}\ra=0 $,
$i=1,\ldots,6$. Moreover, the function $\bm \tilde{U}\\
\overline{\tilde{U}}\endm_{dis}$ is to satisfy the estimates
\eqref{disasympto}, the functions $\lambda_{i}(t)$ are to satisfy
the inequalities \eqref{rootasympto}, while the functions
$\nu_{1}$, $\beta_{1}$ etc. are to satisfy the corresponding
estimates in \eqref{modulationasympto}. Finally, we require the
bounds $|a_{\infty}^{(1)}|+|b_{\infty}^{(1)}|+|y_{\infty}|^{1}+|v_{\infty}|^{1}\lesssim \delta^{2}$.
Upon fixing the 'initial condition' $\bm A\\
\bar{A}\endm$ as in Theorem~\ref{Main}, we then construct a map
$T_{A}$ which associates another tuple, characterized by primes,
i. e. we get a tuple $\{\bm U'\\ \bar{U'}\endm_{dis},\ldots\}$, as
follows: first, re-construct the original quantities $\lambda$,
$\beta$ etc. from the tuple in the obvious fashion. This then also
defines $\Psi-\Psi_{\infty}$ etc, see the beginning of subsection
3.2,
provided we choose $\gamma_{\infty}$ in such fashion that $(\Psi-\Psi_{\infty})_{1}(0)=0$. In particular, we can reconstruct $\bm U\\
\bar{U}\endm(t,.)=\bm
e^{i(\Psi-\Psi_{\infty})_{1}(t)}&0\\0&e^{-i(\Psi-\Psi_{\infty})_{1}(t)}\endm
\bm
\tilde{U}\\
\overline{\tilde{U}}\endm(t,.-\lambda_{\infty}(\mu-\mu_{\infty})(t))$.
Then we use \eqref{Udis} to construct $\bm \tilde{U'}\\
\overline{\tilde{U'}}\endm_{dis}$; simply use the un-primed
quantities for the right-hand side. Next, define $\lambda_{2}'$
via the righthand side of \eqref{lambda2}, $\lambda_{i}'$ via the
righthand side of \eqref{lambdai}, and $\lambda_{6}'$ via
\eqref{lambda6''}. We next turn to the modulation parameters:
define $\nu_{1}'$ as the righthand side of \eqref{nu}, and
$\beta_{1}'$ via \eqref{beta}. We can then also define
$b_{\infty}'$ and $a_{\infty}'$ by means of \eqref{binfty},
\eqref{ainfty}, respectively\footnote{In these relations use the
unprimed quantities inside the integrals.}, which in turn defines
$\lambda_{\infty}'(.)$. Further, put\footnote{The quantity $B'(s)$
{\it{is not the derivative }}, but the new $B(s)$.}
\begin{equation}\nonumber
B'(s)=\exp(\int_{0}^{s}[\beta'\nu'^{2}+\frac{1}{4i\kappa_{2}}E_{5}](\sigma)d\sigma),
\end{equation}
where $E_{5}$ is defined with respect to the un-primed quantities.
We shall show later that under suitable assumptions on the
tuple $\{\bm \tilde{U}\\
\overline{\tilde{U}}\endm_{dis},\ldots\}$, we have
$B'(s)^{-1}=c'\lambda_{\infty}'^{-1}(s)+o(\frac{1}{s})$, for
suitable $c'$. Then define $\omega'$ via the formula
\begin{equation}\nonumber\begin{split}
&\omega'(s)=c'\lambda_{\infty}'(s)^{-1}[\omega(0)+\int_{0}^{\infty}B'(\sigma)[(i\kappa_{1})^{-1}E_{3}+\frac{\beta}{2i\kappa_{1}}E_{4}](\sigma)d\sigma]
\\&\hspace{3cm}-c'\lambda_{\infty}'(s)^{-1}\int_{s}^{\infty}B'(\sigma)[(i\kappa_{1})^{-1}E_{3}+\frac{\beta}{2i\kappa_{1}}E_{4}](\sigma)d\sigma]+o(s^{-1}),
\end{split}\end{equation}
see the discussion preceding \eqref{omega1}. We can infer from
this a number $c_{\infty}'$ with the property\footnote{We shall
verify later that this definition is indeed meaningful.}
\begin{equation}\nonumber
\omega(s)'=\frac{c_{\infty}'}{\lambda_{\infty}'(s)}+o(\frac{1}{s}),
\end{equation}
whence we can define
$\omega_{1}(s)'=\omega(s)'-c_{\infty}'\lambda_{\infty}'^{-1}(s)$.
Continue by setting
\begin{equation}\label{mu'}
\mu(s)'=\mu(0)+\int_{0}^{s}\lambda_{\infty}'(\sigma)^{-1}[2\omega'\nu'+(i\kappa_{1}\nu(\sigma))^{-1}E_{4}(\sigma)]d\sigma
\end{equation}
Again, under suitable assumptions on the original tuple we shall
be able to infer the existence of numbers $v_{\infty}'$,
$y_{\infty}'$ with the property\footnote{See last footnote}
\begin{equation}\nonumber
\mu(s)'=\frac{2v_{\infty}'s+y_{\infty}'}{a_{\infty}'+b_{\infty}'s}+o(\frac{1}{s}),
\end{equation}
whence we can define
$\mu_{1}(s)'=\mu(s)'-\frac{2v_{\infty}'s+y_{\infty}'}{a_{\infty}'+b_{\infty}'s}$.
Also, we shall have
$c_{\infty}'=v_{\infty}'a_{\infty}'-\frac{b_{\infty}'y_{\infty}'}{2}$.
Finally, put
\begin{equation}\nonumber\begin{split}
&\gamma(s)'=\gamma(0)+\int_{0}^{s}[\nu'^{2}(\sigma)-(2i\kappa_{2})^{-1}E_{6}(\sigma)+\nu'^{2}(\sigma)\omega'^{2}(\sigma)-\frac{1}{i\kappa_{1}}\omega(\sigma)E_{4}(\sigma)
-\frac{1}{2\kappa_{2}}(i\kappa_{2})^{-1}E_{2}(\sigma)]d\sigma,
\end{split}\end{equation}
whence we can define
$\gamma_{1}(s)'=\gamma(s)'-s-\frac{c_{\infty}'^{2}}{\lambda_{\infty}'^{2}}$.
\\

We can now reduce the proof of Theorem~\ref{Main} to locating a
fixed point of the map $T_{A}$. This follows from the next
Proposition; let the norm $|||.|||$ be defined as in
definition~\ref{Abounds}:

\begin{proposition}\label{fixedpointreduction} Let $A: {\mathbf{R}}\rightarrow{\mathbf{C}}$ be a smooth function satisfying $|||A|||<\delta$ for suitably small $\delta>0$,
as well as the orthogonality conditions $\la \bm A\\ \bar{A}\endm, \xi_{i,\text{proper}}\ra=0$, $i=1,\ldots,6$.  Let $\{\bm\tilde{U}\\
\overline{\tilde{U}}\endm_{dis}(.,.),\ldots\}$ be a tuple as
above\footnote{In particular satisfying all the above-specified
estimates.} satisfying the fixed point property\footnote{In
particular, we assume that the operation of $T_{A}$ is
well-defined on this tuple. We shall soon analyze where $T_{A}$ is
well-defined.}
\begin{equation}\nonumber
T_{A}\{\bm\tilde{U}\\
\overline{\tilde{U}}\endm_{dis}(.,.),\ldots\}=\{\bm\tilde{U}\\
\overline{\tilde{U}}\endm_{dis}(.,.),\ldots\}
\end{equation}
Also, assume that $\sup_{s\geq 0}||U(s,.)||_{L_{x}^{2}}\lesssim
\delta$, $\sup_{s\geq
0}s^{\frac{1}{2}}||U(s,.)||_{L_{x}^{\infty}}\lesssim \delta$.
Define\footnote{The function $(\Psi-\Psi_{\infty})_{1}(s)$ is
given by \eqref{Psi-Psiinfty}.}
\begin{equation}\nonumber\begin{split}
&\bm U\\ \bar{U}\endm(s,.)\\&=\bm
e^{i(\Psi-\Psi_{\infty})_{1}(s)}&0\\0&e^{-i(\Psi-\Psi_{\infty})_{1}(s)}\endm
[\bm \tilde{U}\\
\overline{\tilde{U}}\endm_{dis}+\sum_{j=1}^{6}\lambda_{j}\eta_{j,\text{proper}}](s,.-\lambda_{\infty}(\mu-\mu_{\infty})(s)),
\end{split}\end{equation}
Then the function (with $\calT_{\infty}$ as in \eqref{Tinfty})
\begin{equation}\nonumber
\psi(t,x)=W(t,x)+\calT_{\infty}^{-1}[e^{is}U(s,.)](t,x)
\end{equation}
is a non-generic blow-up solution of \eqref{TheEquation} exploding
at time $t_{*}=\frac{1}{a_{\infty}b_{\infty}}$ (recall
$a_{\infty}, b_{\infty}\sim 1$). We have
\begin{equation}\nonumber
W(t,x)=e^{i(\gamma(s(t))+[\omega(x-\mu)](s(t))}e^{-i\frac{\beta}{4}[\lambda(x-\mu)(s(t))]^{2}}\sqrt{\lambda(s(t))}\phi_{0}([\lambda(x-\mu)](s(t)),1),
\end{equation}
where $s(t)=\frac{a_{\infty}t}{a_{\infty}^{-1}-b_{\infty}t}$.
Finally, we have
\begin{equation}\nonumber
\bm U\\ \bar{U}\endm(0,x)=\bm A\\
\bar{A}\endm(x)+\sum_{i=1}^{6}\tilde{\lambda}_{i}\eta_{i,\text{proper}}
\end{equation}
for certain numbers $\tilde{\lambda_{i}}$ with
$|\tilde{\lambda_{i}}|\lesssim \delta^{2}$.
\end{proposition}

\begin{proof} To begin with, note that the modulation equations \eqref{modulation} etc are satisfied. We continue by verifying that $\bm U\\
\bar{U}\endm(s,y)$ satisfies \eqref{Ueqn}. Thus define a function
$\bm \Omega\\ \bar{\Omega}\endm(s,y)$ by means of the
inhomogeneous linear equation
\begin{align}\label{Omega}
i\partial_{s}\bm \Omega\\ \bar{\Omega}\endm+\calH(s)\bm \Omega\\ \bar{\Omega}\endm& =-i(\dot\lambda\lambda^{-1}-\beta\nu^2)(\tileta_2-\beta\tileta_5/2+\omega\tileta_4)\\
&\quad +\frac{i}{4}(\dot\beta+\beta^2\nu^2)\tileta_5 + i(\nu^2-\dot\gamma+\nu^2\omega^2)\tileta_1 \\
&\quad -i(\dot\omega+\beta\omega\nu^2)\tileta_4 -
i\nu(\dot\mu\lambda_\infty-2\nu\omega)(-\omega\tileta_1-\tileta_3+\beta\tileta_4/2)
+ N(U,\pi),\\
\bm \Omega\\ \bar{\Omega}\endm(0,.)=\bm U\\ \bar{U}\endm(0,.)
\end{align}
Define $\bm \tilde{\Omega}\\
\overline{\tilde{\Omega}}\endm(t,.):=\bm
e^{-i(\Psi-\Psi_{\infty})_{1}(t)}&0\\0&e^{+i(\Psi-\Psi_{\infty})_{1}(t)}\endm
\bm
\Omega\\
\overline{\Omega}\endm(t,.+\lambda_{\infty}(\mu-\mu_{\infty})(t))$,
and use a decomposition
\begin{equation}\label{Omegadecomp}
\bm \tilde{\Omega}\\
\overline{\tilde{\Omega}}\endm(t,.)=\bm \tilde{\Omega}\\
\overline{\tilde{\Omega}}\endm_{dis}(t,.)+\sum_{j=1}^{6}\mu_{j}(t)\eta_{j,\text{proper}}
\end{equation}
Now we deduce the equation
\begin{equation}\nonumber\begin{split}
&\bm \tilde{\Omega}\\ \overline{\tilde{\Omega}}\endm_{dis}(t,.)=e^{it\calH}\bm \tilde{U}^{(t)}\\
\overline{\tilde{U}^{(t)}}\endm_{dis}(0,.)-i\int_{0}^{t}e^{i(t-s)\calH}[...]_{dis}ds\\&-2i\int_{0}^{t}e^{i(t-s)\calH}[\bm
0&-e^{+2i(\Psi-\Psi_{\infty})_{1}(s)-2i(\Psi-\Psi_{\infty})_{1}(t)}+1\\e^{-2i(\Psi-\Psi_{\infty})_{1}(s)+2i(\Psi-\Psi_{\infty})_{1}(t)}-1\endm
\\&\hspace{12cm}\times\phi_{0}^{4}\bm
\tilde{\Omega}^{(t)}(s,.)\\\overline{\tilde{\Omega}^{(t)}}(s,.)\endm]_{dis}ds\\
&-i\int_{0}^{t}e^{i(t-s)\calH}\widetilde{[...]}_{dis}(s)ds
\end{split}\end{equation}
where $[...]$ refers to the righthand side of \eqref{Omega}
translated by $+\lambda_{\infty}(\mu-\mu_{\infty})(t)$ in the
spatial variable and 'twisted' by $\bm
e^{i(\Psi_{\infty}-\Psi)_{1}(t)}&0\\0&e^{i(\Psi_{\infty}-\Psi)_{1}(t)}\endm$,
and we put
\begin{equation}\nonumber\begin{split}
&\widetilde{[...]}\\&=\bm
e^{i(\Psi_{\infty}-\Psi)_{1}(t)}&0\\0&e^{-i(\Psi_{\infty}-\Psi)_{1}(t)}\endm\\&\times\bm
-3\nu^{2}(s)\tilde{\phi}_{0}^{4}(.-\lambda_{\infty}(\mu-\mu_{\infty})(s))+3\tilde{\phi}_{0}^{4}&2\tilde{\phi}_{0}^{4}e^{2i(\Psi-\Psi_{\infty})_{1}(s)}(-e^{2i(\Psi-\Psi_{\infty})_{2}(t)}+1)\\
-2\tilde{\phi}_{0}^{4}e^{-2i(\Psi-\Psi_{\infty})_{1}(s)}(-e^{-2i(\Psi-\Psi_{\infty})_{2}(t)}+1)&3\nu^{2}(s)\tilde{\phi}_{0}^{4}(.-\lambda_{\infty}(\mu-\mu_{\infty})(s))-3\tilde{\phi}_{0}^{4}\endm\\&\hspace{10cm}\times\bm
\Omega\\
\bar{\Omega}\endm(s,.+\lambda_{\infty}(\mu-\mu_{\infty})(t))
\end{split}\end{equation}
and as before we put
$\tilde{\phi}_{0}(.)=\phi(.+\lambda_{\infty}(\mu-\mu_{\infty})(t))$.
From the iterative step and the fact that the modulation equations
are satisfied, we then deduce that
\begin{equation}\label{deltadis}\begin{split}
&\bm \tilde{\Omega}-\tilde{U}\\
\overline{\tilde{\Omega}-\tilde{U}}\endm_{dis}(t,.)\\&=-2i\int_{0}^{t}e^{i(t-s)\calH}[\bm
0&-e^{+2i(\Psi-\Psi_{\infty})_{1}(s)-2i(\Psi-\Psi_{\infty})_{1}(t)}+1\\e^{-2i(\Psi-\Psi_{\infty})_{1}(s)+2i(\Psi-\Psi_{\infty})_{1}(t)}-1\endm
\\&\hspace{10cm}\times\phi_{0}^{4}\bm
\tilde{\Omega}^{(t)}(s,.)-\tilde{U}^{(t)}(s,.)\\\overline{\tilde{\Omega}^{(t)}(s,.)-\tilde{U}^{(t)}(s,.)}\endm]_{dis}ds\\
&-i\int_{0}^{t}e^{i(t-s)\calH}\widetilde{\Delta [...]}_{dis}(s)ds,
\end{split}\end{equation}
where the expression $\widetilde{\Delta [...]}_{dis}(s)$ is
defined as above but with $\Omega$ replaced by the difference
$\Omega-U$. Next, we have
\begin{equation}\nonumber\begin{split}
&i\partial_{s}\la \bm \Omega\\ \bar{\Omega}\endm,
\tilde{\xi}_{i}\ra=\la (i\partial_{s}+\calH(s))\bm \Omega\\
\bar{\Omega}\endm,\tilde{\xi}_{i}\ra-\la  \bm \Omega\\
\bar{\Omega}\endm,
(i\partial_{s}+\calH^{*}(s))\tilde{\xi}_{i}\ra,\,i=2,\ldots,6
\end{split}\end{equation}
whence we obtain
\begin{equation}\nonumber
i\partial_{s}\la\bm \Omega\\ \bar{\Omega}\endm-\bm U\\
\bar{U}\endm, \tilde{\xi}_{i}\ra=-\la\bm \Omega\\ \bar{\Omega}\endm-\bm U\\
\bar{U}\endm,
(i\partial_{s}+\calH^{*}(s))\tilde{\xi}_{i}\ra,\,i=2,\ldots,6
\end{equation}
Thus we obtain from $(\bm \Omega\\ \bar{\Omega}\endm-\bm U\\
\bar{U}\endm)(0,.)=0$ the relation
\begin{equation}\label{deltai}\begin{split}
&i\la \bm \Omega\\ \bar{\Omega}\endm-\bm U\\
\bar{U}\endm, \tilde{\xi}_{i}\ra(t,.)=\la \bm \tilde{\Omega}\\ \overline{\tilde{\Omega}}\endm-\bm \tilde{U}\\
\overline{\tilde{U}}\endm, \Xi_{i}\ra(t,.)\\&=-\int_{0}^{t}\la\bm \Omega\\ \bar{\Omega}\endm-\bm U\\
\bar{U}\endm,
(i\partial_{s}+\calH^{*}(s))\tilde{\xi}_{i}\ra(s)ds,\,i=2,\ldots,6,
\end{split}\end{equation}
where the first equality can be used to define $\Xi_{i}$ in the
obvious fashion. We shall use the last relation to solve for the
coefficients $\mu_{i}$, $i=1,\ldots,5$. Finally, we consider the
coefficient $\mu_{6}$:
\begin{equation}\nonumber\begin{split}
&i\dot{\mu}_{6}(t)=\la i\partial_{t}\bm \Omega\\
\bar{\Omega}\endm, \xi_{1,\text{proper}}\ra\\
&=\la\bm \Omega\\
\bar{\Omega}\endm,\bm
\partial_{t}(\Psi-\Psi_{\infty})_{1}&0\\0&-\partial_{t}(\Psi-\Psi_{\infty})_{1}\endm\xi_{1,\text{proper}}\ra+\la
\partial_{t}(\lambda_{\infty}(\mu-\mu_{\infty})(t))\bm
\partial_{x}\tilde{\Omega}\\\overline{\partial_{x}\tilde{\Omega}}\endm,
\xi_{1,\text{proper}}\ra \\
&+\la \bm
e^{-i(\Psi-\Psi_{\infty})_{1}(t)}i\partial_{t}\Omega(t,y+\lambda_{\infty}(\mu-\mu_{\infty})(t))\\
\overline{-e^{-i(\Psi-\Psi_{\infty})_{1}(t)}i\partial_{t}\Omega(t,y+\lambda_{\infty}(\mu-\mu_{\infty})(t))}\endm,
\xi_{1,\text{proper}}\ra
\end{split}\end{equation}
Now we recall the corresponding identity in the derivation of
\eqref{lambda6}, namely \eqref{lambda6'}, take the difference of
the latter and the identity above, and proceed as the paragraphs
after \eqref{lambda6'}. Using the fact that
$\lambda_{6}(0)=\mu_{6}(0)$, we deduce the schematic identity
\begin{equation}\label{delta6}
(\lambda_{6}-\mu_{6})(t)=\int_{0}^{t}[\la \bm
\tilde{\Omega}\\\overline{\tilde{\Omega}}\endm(s,.)-\bm\tilde{U}\\
\overline{\tilde{U}}\endm(s,.), \phi(s)\ra+ \la \partial_{x}\bm
\tilde{\Omega}\\\overline{\tilde{\Omega}}\endm(s,.)-\partial_{x}\bm\tilde{U}\\
\overline{\tilde{U}}\endm(s,.), \psi(s)\ra] ds
\end{equation}
Now from \eqref{deltadis}, \eqref{deltai}, \eqref{delta6}, as well
as \eqref{Omegadecomp} and the linear estimate
Theorem~\ref{linearestimates}, we easily deduce the estimate
\begin{equation}\nonumber
\sup_{0\leq t\leq
T}[||[\tilde{\Omega}-\tilde{U}]_{dis}(t,.)||_{H^{1}}+\sum_{i=1}^{6}|\lambda_{i}-\mu_{i}|(t)]\lesssim
T\sup_{0\leq t\leq
T}[||[\tilde{\Omega}-\tilde{U}]_{dis}(t,.)||_{H^{1}}+\sum_{i=1}^{6}|\lambda_{i}-\mu_{i}|(t)]
\end{equation}
Now choose $T>0$ small enough to get the identity
$\Omega|_{[0,T]}= U|_{[0,T]}$. Continuing in this fashion implies
$U(.,.)=\Omega(.,.)$. Next, observe that the condition
\begin{equation}\label{600}
\bm \tilde{U}\\
\overline{\tilde{U}}\endm_{dis}(0,.)=P_{s}[\sum_{i}\eta_{i,\text{proper}}(.+\lambda_{\infty}(\mu-\mu_{\infty})(0))\la
\bm \tilde{U}\\
\overline{\tilde{U}}\endm(0,.),
\tilde{\xi}_{k(i),\text{proper}}\ra +\bm A\\
\bar{A}\endm(.+\lambda_{\infty}(\mu-\mu_{\infty})(0))]
\end{equation}
in addition to
$P_{root}\bm\tilde{U}\\\overline{\tilde{U}}\endm(0,.)=\sum_{i=1}^{6}\lambda_{i}(0)\eta_{i,\text{proper}}$
uniquely determines $\bm \tilde{U}\\
\overline{\tilde{U}}\endm_{dis}(0,.)$, hence in conjunction with
the values of the modulation parameters also $\bm U\\
\overline{U}\endm(0,.)$. Then one verifies that $\bm U\\
\bar{U}\endm(0,.)=\bm A\\ \bar{A}\endm
+\sum_{i}\alpha_{i}\eta_{i,\text{proper}}$ with the $\alpha_{i}$
defined as in \eqref{alphai} is consistent with \eqref{600}, as well as the root part of $\bm \tilde{U}\\ \overline{\tilde{U}}\endm(0,.)$.\\
Now reverse the algebraic manipulations that led to \eqref{Ueqn}.
We deduce that
\begin{equation}\nonumber
Z(t,x)=W(t,x)+\calT_{\infty}^{-1}[e^{is}U(s,.)](t,x)
\end{equation}
is indeed a non-generic blow-up solution of \eqref{TheEquation}.
One checks that
\begin{equation}\nonumber\begin{split}
&\calT_{\infty}^{-1}[e^{is}U(s,.)]\\&=(a_{\infty}^{-1}-b_{\infty}t)^{-\frac{1}{2}}e^{-i\frac{b_{\infty}x^{2}}{4(a_{\infty}^{-1}-b_{\infty}t)}}
e^{i(\gamma_{\infty}+v_{\infty}x-v_{\infty}^{2}\frac{a_{\infty}t}{a_{\infty}^{-1}-b_{\infty}t}-v_{\infty}y_{\infty})}
e^{i(\frac{a_{\infty}t}{a_{\infty}^{-1}-b_{\infty}t})}U(\frac{a_{\infty}t}{a_{\infty}^{-1}-b_{\infty}t},
\frac{x-2v_{\infty}a_{\infty}t}{a_{\infty}^{-1}-b_{\infty}t}-y_{\infty})
\end{split}\end{equation}
The assumptions in the Proposition imply that this remains bounded
with respect to $L^{\infty}$ as $t\rightarrow
t_{*}=\frac{1}{a_{\infty}b_{\infty}}$, while due to the asymptotic
relations \eqref{modulationasympto} the principal soliton part
$W(t,x)$ blows up according to the non-generic profile. Finally,
recall the decomposition \eqref{initialdecomp} in which we use
\eqref{alphai}. Our assumptions \eqref{modulationasympto} as well
as \eqref{rootasympto} imply the last statement of the
Proposition.
\end{proof}

\subsection{Deducing the fixed point from a priori estimates.}

We now need to demonstrate the existence of a fixed point for the
map $T_{A}$ on the set of tuples satisfying the above specified
inequalities. This will follow from an application of the
Schauder-Tychonoff fixed point Theorem, which we recall here:

\begin{theorem}\label{fixedpoint}(Schauder-Tychonoff) A non-empty compact convex
subset $S$ of a Banach space has the fixed point property, i. e.
for any continuous map $T: S\longrightarrow S$ there exists
$x_{T}\in S$ satisfying $T(x_{T})=x_{T}$.
\end{theorem}

We now need to locate such a set $S$. We construct this as
follows: first, for $M, N$ as before (\eqref{modulationasympto},
\eqref{rootasympto}) and very large\footnote{This parameter will
eventually depend on a time $T$.} $K>0$ introduce the
norm\footnote{We only include the parameters $N,K$ as superscripts
in the norm, since we shall only vary these.}
\begin{equation}\nonumber\begin{split}
&|||U|||_{S^{N,K}}\\&=\sum_{0\leq k\leq N}C_{k}^{-1}\sum_{2i+j\leq
k}[\sup_{s\geq 0}\la
s\ra^{\frac{1}{2}-25^{k}\delta_{2}}||\partial_{s}^{i}\partial_{y}^{j}U(s,y)||_{L_{s}^{M}L_{y}^{M}}+\sup_{s\geq
0}\la
s\ra^{-10^{k}\delta_{2}}||\partial_{s}^{i}\partial_{y}^{j}U(s,y)||_{L_{s}^{M}L_{y}^{2}}\\&+\sup_{\phi\in\calA}\sup_{s\geq
0}\la
s\ra^{1-20^{k}\delta_{2}}||\phi\partial_{s}^{i}\partial_{y}^{j}U(s,y)||_{L_{s}^{M}L_{y}^{M}}]+\sup_{s\geq
0}[\sup_{\phi\in\calA}\la s\ra^{\frac{3}{2}-\delta_{3}}||\phi
U(s,.)||_{L_{y}^{\infty}}+\sup_{s\geq
0}||CU(s,y)||_{L_{y}^{2}}\\&\hspace{8cm}+\sum_{1\leq k\leq
N-1}K^{-k}\sup_{2i+j=k}||C\partial_{s}^{i}\partial_{y}^{j}U(s,.)||_{L_{s}^{M}L_{y}^{2}}],
\end{split}\end{equation}
where $C$ is as in \eqref{disasympto}. The role of the last
summand is to ensure uniform spatial decay on finite time
intervals, again needed for compactness. Also, let
$|||U|||_{S^{N}}$ be defined as above, but with the last summand
replaced by $||C U||_{L_{y}^{2}}$, where we use
$C=x+2is\partial_{x}$. Define the Banach space $S^{N}$ as the
completion of ${\mathcal{S}}({\mathbf{R}}^{2})$ with
respect to these norms. Now for a tuple $\Gamma:=\{\bm \tilde{U}\\
\overline{\tilde{U}}\endm_{dis}, \ldots\}$ as before, define the
norm\footnote{We use the notation $\bm \tilde{U}\\
\overline{\tilde{U}}\endm_{dis}=\bm \tilde{U}_{dis}\\
\overline{\tilde{U}_{dis}}\endm$.} (as usual we let $\la s\ra
=|s|+1$)
\begin{equation}\nonumber\begin{split}
&|||\Gamma|||_{\tilde{S}^{N,K}}:=|||\tilde{U}_{dis}|||_{S^{N,K}}+
\delta^{-1}[\sup_{0\leq s<\infty}\la
s\ra^{\frac{1}{2}-\delta_{1}}|\nu_{1}(s)|+\sum_{1\leq k\leq
[\frac{N}{2}]}||\la
s\ra^{\frac{3}{2}-2\delta_{1}}\frac{d^{k}}{ds^{k}}\nu_{1}(s)||_{L^{M}}
+\sup_{0\leq s<\infty}\la
s\ra^{\frac{3}{2}-\delta_{1}}|\beta_{1}(s)|\\&+\sum_{1\leq k\leq
[\frac{N}{2}]}||\la
s\ra^{2-2\delta_{1}}\frac{d^{k}}{ds^{k}}\beta_{1}(s)||_{L^{M}}
+\sup_{0\leq s<\infty}|\la
s\ra^{\frac{3}{2}-\delta_{1}}\omega_{1}(s)|+\sum_{1\leq k\leq
[\frac{N}{2}]}||\la
s\ra^{\frac{3}{2}-2\delta_{1}}\frac{d^{k}}{ds^{k}}\omega_{1}(s)||_{L^{M}}\\&+
\sup_{0\leq s<\infty}|\la
s\ra^{\frac{1}{2}-\delta_{1}}\frac{d}{ds}\gamma_{1}(s)|+
\sum_{2\leq k\leq [\frac{N}{2}]}||\la
s\ra^{\frac{3}{2}-2\delta_{1}}\frac{d^{k}}{ds^{k}}\gamma_{1}(s)||_{L^{M}}+
\sup_{0\leq s<\infty}|\la
s\ra^{\frac{3}{2}-\delta_{1}}|\mu_{1}(s)|\\&+\sum_{1\leq k\leq
[\frac{N}{2}]}||\la
s\ra^{\frac{5}{2}-\delta_{1}}\frac{d^{k}}{ds^{k}}\mu_{1}(s)||_{L^{M}}
+\sum_{i=1}^{6}\sum_{0\leq k\leq [\frac{N}{2}]}||\la
t\ra^{2-2\delta_{1}}\frac{d^{k}}{dt^{k}}\lambda_{i}(t)||_{L^{M}}+|a_{\infty}^{(1)}|+|b_{\infty}^{(1)}|+|y_{\infty}^{(1)}|+|v_{\infty}^{(1)}|]
\end{split}\end{equation}
Also, let $|||\Gamma|||_{\tilde{S}^{N}}$ be defined as above but
with $||\tilde{U}_{dis}||_{S^{N,K}}$ replaced by
$||\tilde{U}_{dis}||_{S^{N}}$. Then we define the restrictions
$|||.|||_{S^{N}([0,T))}$ etc for any time interval $[0,T)$ in the
obvious fashion, and denote $S^{N}([0,T))$,
$\tilde{S}^{N,K}([0,T))$ as completion of
$\calS([0,T)\times{\mathbf{R}})$ and of
\begin{equation}\nonumber
\calS([0,T)\times{\mathbf{R}})\times
(C^{\infty}[0,T))^{11}\times[-100\delta^{2},100\delta^{2}]\times
[-100\delta^{2},100\delta^{2}]\times\prod_{i=1}^{3}[-10,10],
\end{equation}
respectively, with respect to the above norms localized to
$[0,T)\times{\mathbf{R}}$. Then we have
\begin{lemma} Fix $T$, $\infty>T>0$. For any $R\geq 0$, the set of
tuples\footnote{We omit the dependence of these sets on $R$ in the
notation, it being understood that $R$ below will be fixed
throughout.} $A^{(0)}_{[0,T)}:=\{\Gamma\,\text{on
$[0,T)\times{\mathbf{R}}$}\,|\,|||\Gamma|||_{\tilde{S}^{N,K}[0,T)}\leq
R\delta\}$ equipped with the norm
$|||.|||_{\tilde{S}^{N-4,K}[0,T)}$ is a compact convex subset of
\begin{equation}\nonumber
A_{[0,T)}:=\{\Gamma\,\text{defined on
$[0,T)\times{\mathbf{R}}$}\,|\,|||\Gamma|||_{\tilde{S}^{N-4,K}[0,T)}\leq
R\delta\}
\end{equation}
\end{lemma}

\begin{proof} We demonstrate the compactness assertion: thus
consider a sequence of tuples\\ $\{\Gamma_{i}\}_{i\geq 1}\subset
A^{(0)}_{[0,T)}$. Choose a subsequence for which the 'parameters
at infinity' $a_{\infty,i}^{(1)}$ etc converge. Now consider the
functions $\bm \tilde{U}_{i}\\
\overline{\tilde{U}_{i}}\endm_{dis}=\bm
\tilde{U}_{i,dis}\\\overline{\tilde{U}_{i,dis}}\endm$. By
assumption, letting $\phi_{\rho}(x):=\phi(\frac{x}{\rho})$, where
$\phi(.)$ smoothly localizes to $|x|>1$, we have that
\begin{equation}\nonumber
\lim_{\rho\rightarrow\infty}\sup_{l}\sum_{0\leq k\leq
N-4}\sum_{2i+j\leq k}||\phi_{\rho}(x)\partial_{t}^{i}\partial_{x}
^{j}\tilde{U}_{l,dis}||_{L_{t}^{M}L_{x}^{2}[0,T)}=0
\end{equation}
Indeed, this follows from uniform control over
$||C\partial_{t}^{i}\partial_{x}
^{j}\tilde{U}_{l,dis}||_{L_{t}^{M}L_{x}^{2}[0,T)}$. Combining this
with the fact that $\sum_{0\leq k\leq N}\sup_{2i+j\leq
k}||\partial_{s}^{i}\partial_{x}^{j}\tilde{U}_{i,dis}||_{L_{s}^{M}L_{x}^{2}[0,T)}$
is uniformly bounded and applying the Rellich-Kondrakhov
compactness Theorem as well as Sobolev embedding, we obtain a
subsequence (which we again label $\bm \tilde{U}_{i}\\
\overline{\tilde{U}_{i}}\endm_{dis}$) which converges with respect
to $\sum_{0\leq k\leq N-4}\sup_{2i+j\leq
k}||\partial_{s}^{i}\partial_{x}^{j}(.)||_{L_{s}^{M}L_{x}^{2}[0,T)}$
as well as the remaining norms in $||.||_{S^{N-4,K}}$ to a
limit\footnote{Clearly this limit satisfies the same orthogonality
relations, whence we may
apply the subscript $_{dis}$.} $\bm \tilde{U}\\
\overline{\tilde{U}}\endm_{dis}$. Passing to a further
subsequence, we may assume that all
$\partial_{s}^{i}\partial_{x}^{j}\tilde{U}_{i,dis}$, $N\geq
2i+j\geq N-4$, converge weakly\footnote{It is at this stage that
we need $L^{M}$ instead of $L^{\infty}$} with respect to
$L_{s}^{M}L_{x}^{2}[0,T)$, $L_{s}^{M}L_{x}^{M}[0,T)$, and one
checks that the corresponding limits necessarily equal
$\partial_{s}^{i}\partial_{x}^{j}\tilde{U}_{dis}$ in the
distributive sense, respectively. Also, $||| \bm \tilde{U}\\
\overline{\tilde{U}}\endm_{dis}|||_{S^{N}}\leq R\delta$.\\
Now consider the root part, i. e. the functions
$\lambda_{j,i}(t)$, $j=1,\ldots,6$. By assumption, we have a
uniform bound on $\sum_{0\leq k\leq [\frac{N}{2}]}||\la
t\ra^{2-2\delta_{1}}\frac{d^{k}}{dt^{k}}\lambda_{j,i}(t)||_{L^{M}}$.
By the Arzela-Ascoli Theorem, we can then choose a converging
subsequence with respect to $\sum_{0\leq k\leq
[\frac{N}{2}]-2}||\la
t\ra^{2-2\delta_{1}}\frac{d^{k}}{dt^{k}}(.)||_{L^{M}}$ whose limit
satisfies the desired estimates. The argument for the modulation
parameters $\nu_{1}(t)=\nu(t)-1$ etc. is identical.
\end{proof}

Now define the sets $A^{(n)}_{[0,T)},\,n\geq 1$ inductively as
follows: first, we can modify the inductive step $T_{A}$ to the
interval $[0,T)$, by simply replacing $\infty$ by $T$ in the
formulae for the modulation parameters and root parameters. By
abuse of notation refer to this by $T_{A}$ as well. Then put
\begin{equation}\nonumber
A^{(n)}_{[0,T]}:=A^{(0)}_{[0,T)}\cap
\text{convhull}\overline{(T_{A}A^{(n-1)}_{[0,T)})}
\end{equation}
The closure operation is always with respect to
$|||.|||_{\tilde{S}_{1}^{N-4,K}}$. Then clearly
$A^{(n)}_{[0,T)}\subset A^{(n-1)}_{[0,T)}$, and these are all
compact convex subsets of $A_{[0,T)}$.
\\

Everything now reduces to the following {\bf{core analytic
Theorem}}: first we make a definition:
\begin{defi}\label{Abounds} We call a function $A(x): {\mathbf{R}}\rightarrow {\mathbf{C}}$ admissible
provided the estimate
\begin{equation}\nonumber
|||A|||:=\sup_{0\leq k\leq N}||\la
x\ra^{100}\frac{d^{k}}{dx^{k}}A(x)||_{L^{1}\cap L^{2}}\leq \delta
\end{equation}
is satisfied.
\end{defi}

\begin{theorem}\label{core} Let $A$ be admissible. Let $N$ satisfy the specifications in \eqref{modulationasympto}, \eqref{disasympto}, and $\delta>0$ small enough. There exists $R>0$ sufficiently large such that with
the corresponding $A^{0}_{[0,T)}$ etc constructed as above we have
the following: For every $T>0$, there exists a number $K=K(N,T)$,
as well as an index $n_{0}(N,K)$ such that for $n\geq n_{0}$, we
have $T_{A}A^{(n)}_{[0,T)}\subset A^{(n)}_{[0,T)}$; moreover,
$T_{A}$ acts continuously on $A^{(n)}_{[0,T)}$.  The last
assertions are always non-vacuous if $\delta>0$ is sufficiently
small, since then $A^{(n)}_{[0,T)}\neq \emptyset$. Thus by
Theorem~\ref{fixedpoint} there exists a tuple $\Gamma_{T}\in
A^{(n)}_{[0,T)}$ with the property $T_{A}\Gamma_{T}=\Gamma_{T}$.
For $T<\tilde{T}$, one has the inequality
\begin{equation}\nonumber
|||\Gamma_{\tilde{T}}|_{[0,T)\times{\mathbf{R}}}|||_{\tilde{S}^{N,K(N,T)}}\lesssim
\delta
\end{equation}
Also, we get the uniform\footnote{where the implied constant is
independent of $T$} bounds
\begin{equation}\nonumber
\sup_{T>s\geq 0}||\tilde{U}_{T}(s,.)||_{L_{x}^{2}}\lesssim
\delta,\,\sup_{T>s\geq 0}\la
s\ra^{\frac{1}{2}}||\tilde{U}_{T}(s,.)||_{L_{x}^{\infty}}\lesssim
\delta
\end{equation}
\end{theorem}

Assuming this, we can now conclude the following:

\begin{theorem} There exists a fixed point $\Gamma$ for $T_{A}$
(acting on $[0,\infty)\times{\mathbf{R}}$) satisfying the
assumptions \eqref{modulationasympto}, \eqref{rootasympto},
\eqref{disasympto}, as well
$|a_{\infty}^{(1)}|+|b_{\infty}^{(1)}|+|y_{\infty}^{(1)}|+|v_{\infty}^{(1)}|\lesssim
\delta^{2}$. Thus the assumption of
Proposition~\ref{fixedpointreduction} is realizable.
\end{theorem}

\begin{proof} Let $T_{i}=i$, $i\geq 1$. Then construct fixed
points $\Gamma_{i}$ for the operation of $T_{A}|_{[0,T_{i}]}$ as
in the preceding Theorem. Thanks to the uniform bounds for
$|||\Gamma_{i}|_{[0,j)}|||_{\tilde{S}^{N,K(N,T_{j})}}$, $j\leq i$,
and invoking another compactness argument as before, we can select
a subsequence $\Gamma_{j,i}$ which converges on $[0, T_{j})$ with
respect to $|||.|||_{\tilde{S}^{N-4,K(N,j)}}$. Observe that we
only need a uniform bound on $K(N,T)$ for bounded $T$ here, as we
have arranged. Doing this for $j=1,2,\ldots$ and invoking the
Cantor diagonal argument, we then construct a subsequence, which
we again label $\Gamma_{i}$, which converges on every $[0,T_{j}]$,
$j\geq 1$ to a tuple in $\tilde{S}^{N,K(N,j)}$. The limits then
fit coherently to define a tuple $\Gamma$ on $[0,\infty)$ living
in $\tilde{S}^{N}$, which is the desired fixed point.
\end{proof}

\section{The proof of the core analytic estimates,
Theorem~\ref{core}.}

We shall show that $T_{A}(A^{(n)}_{[0,T)})\subset
A^{(0)}_{[0,T)}$, provided $n$ is large enough. The proof will
also reveal the continuity of the operation $T_{A}$. Also, we
shall show that if one iterates $T_{A}$ starting with the tuple
$\Gamma_{trivial}:=\{\bm 0\\0\endm,
0,\ldots,\lambda(0),\beta(0),0,0,0\}$, one always stays inside
$A^{(0)}$ provided $\delta, \delta_{i}$ etc are chosen suitably,
whence $A^{(i)}\neq \emptyset$ $\forall i\geq 0$. Then observe
that
\begin{equation}\nonumber
T_{A}(A^{(n+1)}_{[0,T)})\subset A^{(0)}_{[0,T)}
\end{equation}
as well as
\begin{equation}\nonumber
T_{A}(\text{convhull}(\overline{T_{A}A^{(n)}_{[0,T)}}))\subset
T_{A}(\text{convhull}(\overline{T_{A}A^{(n-1)}_{[0,T)}})\cap
A^{(0)}_{[0,T)})\subset
\text{convhull}(\overline{T_{A}(A^{(n)}_{[0,T)})})
\end{equation}
This then implies $T_{A}(A^{(n+1)})\subset A^{(n+1)}$, as desired.
Thus we need to prove

\begin{theorem}\label{coreredux} Under the assumptions of Theorem~\ref{core}, we
have $T_{A}(A^{(n)}_{[0,T)})\subset A^{(0)}_{[0,T)}$.
\end{theorem}

This shall occupy the rest of the paper
\begin{proof} In order to prove this, we iterate the map $T_{A}$,
i. e. we show that after applying sufficiently many iterations of
$T_{A}$ to a function in $A^{(0)}_{[0,T)}$, and assuming that each
iterate up to but excluding the last sits in $A^{(0)}_{[0,T)}$,
one again lands in $A_{0}$ for the last iterate. The proof will
also easily imply that applying any number of iterations
$T_{A}^{n}$ to the trivial tuple $\Gamma_{trivial}$ takes one into
$A^{(0)}_{[0,T)}$, and that the Theorem is true. Our procedure
shall simply be to apply $T_{A}$ to
a given tuple $\Gamma=\{\bm \tilde{U}\\
\overline{\tilde{U}}\endm_{dis},\ldots\}$, which is assumed to sit
in $A^{(0)}$, and if necessary assume that $\bm \tilde{U}\\
\overline{U}\endm_{dis}$ itself is defined by the righthand side
of \eqref{Udis} with respect to a different tuple $\Gamma'$, and
similarly for the other ingredients in $\Gamma$. In order to
simplify notation, we shall not even distinguish between these
tuples, i. e.  {\it{we proceed as for the derivation of a priori
estimates.}} We shall commence by retrieving the estimates
\eqref{disasympto} for the radiation part. However, before being
able to do so, we need to justify the assertion made earlier about
the phase $(\Psi-\Psi_{\infty})_{2}(t,y)$:

\begin{lemma}\label{Psi2} Assume that the relations \eqref{modulationasympto}
are satisfied. Then we have for $\phi\in \calA$ (recall the
definitions after \eqref{disasympto})
\begin{equation}\nonumber
\phi(y)|(\Psi-\Psi_{\infty})_{2}|(t,y)\lesssim \delta^{2}\la
t\ra^{-\frac{3}{2}+\delta_{1}}
\end{equation}
\end{lemma}
\begin{proof}
Recall the definition
\begin{equation}\nonumber\begin{split}
(\Psi-\Psi_{\infty})_{2}(t,y)=&y[\omega(t)\nu(t)-\frac{\beta(t)}{2}\nu(t)\lambda(t)(\mu_{\infty}-\mu)(t)-\frac{a_{\infty}v_{\infty}-\frac{b_{\infty}y_{\infty}}{2}}{a_{\infty}+b_{\infty}t}]
\\&+y^{2}[\frac{b_{\infty}}{4(a_{\infty}+b_{\infty}t)}-\frac{\beta}{4}\nu^{2}(t)-\frac{\beta}{4}\lambda^{2}(t)(\mu_{\infty}-\mu)^{2}(t)]
\end{split}\end{equation}
The claimed estimate now follows easily from the facts that
$c_{\infty}=a_{\infty}v_{\infty}-\frac{b_{\infty}y_{\infty}}{2}$
as well as
$|\omega(t)-\frac{c_{\infty}}{\lambda_{\infty}(t)}|\lesssim
\delta^{2}\la t\ra^{-\frac{3}{2}+\delta_{1}}$,
$|\beta(t)\nu(t)-\frac{b_{\infty}}{\lambda_{\infty}(t)}|\lesssim
\delta^{2}\la t\ra^{-\frac{3}{2}+\delta_{1}}$, $|\nu(t)-1|\lesssim
\delta^{2}\la t\ra^{-\frac{1}{2}+\delta_{1}}$.
\end{proof}
 We also need
 \begin{lemma}\label{Psi1}
 The following estimate holds under the same assumptions as in the
 preceding Lemma:
 \begin{equation}\nonumber
 |(\Psi-\Psi_{\infty})_{1}(t)|\lesssim \delta^{2} \la t\ra^{\frac{1}{2}+\delta_{1}}
\end{equation}
\end{lemma}
\begin{proof} This is along the same lines, although the algebra
is a bit more complicated. Observe that
\begin{equation}\nonumber
\frac{d}{ds}[\frac{v_{\infty}^{2}sa_{\infty}}{a_{\infty}+b_{\infty}s}-\frac{b_{\infty}v_{\infty}sy_{\infty}}{a_{\infty}+b_{\infty}s}+\gamma_{\infty}-\frac{b_{\infty}y_{\infty}^{2}}{4(a_{\infty}+b_{\infty}s)}]
=\frac{(a_{\infty}v_{\infty}-\frac{b_{\infty}y_{\infty}}{2})^{2}}{(a_{\infty}+b_{\infty}s)^{2}}=\frac{c_{\infty}^{2}}{\lambda_{\infty}^{2}(s)},
\end{equation}
while also
$|\partial_{s}[\gamma_{s}-s-c_{\infty}\lambda_{\infty}^{-2}(s)]|\lesssim
\delta^{2}s^{-\frac{1}{2}+\delta_{1}}$. The claim follows easily
from this and the definition of $(\Psi-\Psi_{\infty})_{1}$.
\end{proof}

We shall first dispose of the easier estimates in
\eqref{disasympto}, which happens to be everything with the
exception of the {\it{strong local dispersive estimate}} and the
{\it{pseudo-conformal almost conservation law}}, i. e. the third
and fourth inequality. Commence with the case $k=0$ and consider
$C_{0}^{-1}\la
s\ra^{\frac{1}{2}-\delta_{2}}||\tilde{U}_{dis}(s,y)||_{L_{y}^{M}}$.
We shall employ the customary bootstrap technique. Thus we assume
an estimate $C_{0}^{-1}\la
s\ra^{\frac{1}{2}-\delta_{2}}||\tilde{U}_{dis}(s,y)||_{L_{y}^{M}}\leq
\Lambda \delta$ for some sufficiently large\footnote{Here the size
depends on certain a priori constants independent of $\delta$}
$\Lambda$, and similarly for all the other norms in
\eqref{disasympto} as well as the modulation
parameters\footnote{More precisely, one replaces $\delta$ by
$\Lambda\delta$ and $\lesssim$ by $\leq$} etc., and then show that
choosing $\delta>0$ small enough implies the same inequalities
with $\frac{\Lambda}{2}$ instead. Then using the local solvability
and obvious continuous dependence of the norms on the
time-interval we infer the desired a priori bound. Now use
\eqref{Udis} as well as Duhamel's principle and
Theorem~\ref{linearestimates} to deduce that
\begin{equation}\nonumber\begin{split}
&C_{0}^{-1}\la
t\ra^{\frac{1}{2}-\delta_{2}}||\bm\tilde{U}\\\overline{\tilde{U}}\endm_{dis}(t,.)||_{L_{x}^{M}}
\\&\lesssim C_{0}^{-1}\la
t\ra^{\frac{1}{2}-\delta_{2}}||e^{it\calH}P_{s}[\bm
e^{-i(\Psi-\Psi_{\infty})_{1}(t)}&0\\0&e^{i(\Psi-\Psi_{\infty})_{1}(t)}\endm[\bm
A(.+\lambda_{\infty}(\mu-\mu_{\infty})(t))\\\bar{A}(.+\lambda_{\infty}(\mu-\mu_{\infty})(t))\endm_{dis}
\\&+\sum_{j=1}^{6}\alpha_{j}
\eta_{j,\text{proper}}(.+\lambda_{\infty}(\mu-\mu_{\infty})(t))]]||_{L_{x}^{M}}+C_{0}^{-1}\la
t\ra^{\frac{1}{2}-\delta_{2}}\int_{0}^{t}||e^{i(t-s)\calH}[...]_{dis}(s,.)||_{L_{x}^{M}}ds\\
\end{split}\end{equation}
The first two terms here are easy to estimate with respect to
$||.||_{L_{t}^{M}}$ on account of \eqref{modulationasympto} as
well as Theorem~\ref{linearestimates}, if one chooses $\Lambda$
large enough in relation to $||A||_{L_{x}^{1}\cap L_{x}^{2}}$; one
thus gets a bound of the form $\leq
\frac{\Lambda}{100}\delta+\Lambda^{2}\delta^{2}$ upon choosing
$\Lambda$ large enough for the contribution of these terms, and
this allows one to close if $\delta$ is chosen small enough in
relation to $\Lambda$. Thus we now focus on the last integral
term, in which we recall $[...]_{dis}$ stands for a sum of
expressions, which basically fall into two contributions, namely
local as well as non-local terms. Of the local ones, the most
difficult contribution is easily seen\footnote{All the other local
terms in \eqref{Udis} contain extra weights of at least the
strength of $\nu-1$, and can be handled similarly.} to come from
the following term\footnote{Recall the definition of
$\tilde{U}^{(t)}(s,.)$ in \eqref{Uts}.} in $[...]_{dis}$:
\begin{equation}\nonumber
[\bm
0&-e^{2i(\Psi-\Psi_{\infty})_{1}(s)-2i(\Psi-\Psi_{\infty})_{1}(t)}+1
\\ e^{-2i(\Psi-\Psi_{\infty})_{1}(s)+2i(\Psi-\Psi_{\infty})_{1}(t)}-1&0\endm
\phi_{0}^{4}\bm \tilde{U}^{(t)}(s)\\
\overline{\tilde{U}^{(t)}(s)}\endm]_{dis},
\end{equation}
which leads to the expression
\begin{equation}\nonumber\begin{split}
&C_{0}^{-1}\la
t\ra^{\frac{1}{2}-\delta_{2}}\int_{0}^{t}||e^{i(t-s)\calH}[\bm
0&-e^{2i(\Psi-\Psi_{\infty})_{1}(s)-2i(\Psi-\Psi_{\infty})_{1}(t)}+1
\\ e^{-2i(\Psi-\Psi_{\infty})_{1}(s)+2i(\Psi-\Psi_{\infty})_{1}(t)}-1&0\endm
\\&\hspace{12cm}\phi_{0}^{4}\bm \tilde{U}^{(t)}(s)\\
\overline{\tilde{U}^{(t)}(s)}\endm]_{dis}||_{L_{x}^{M}}ds\\
\end{split}\end{equation}
Observe from Lemma~\ref{Psi1} that we have
\begin{equation}\nonumber
\sup_{0\leq s\leq
\delta^{-\frac{1}{2}}}|(\Psi-\Psi_{\infty})_{1}|(s)\lesssim
\Lambda\delta^{\frac{7}{4}-\delta_{1}},
\end{equation}
whence we get if we restrict $t<\delta^{-\frac{1}{2}}$ and
$M>>\delta_{2}^{-1}$
\begin{equation}\nonumber\begin{split}
&C_{0}^{-1}||\chi_{<\delta^{-\frac{1}{2}}}(t)\la
t\ra^{\frac{1}{2}-\delta_{2}}\int_{0}^{t}||e^{i(t-s)\calH}[\bm
0&-e^{2i(\Psi-\Psi_{\infty})_{1}(s)-2i(\Psi-\Psi_{\infty})_{1}(t)}+1
\\ e^{-2i(\Psi-\Psi_{\infty})_{1}(s)+2i(\Psi-\Psi_{\infty})_{1}(t)}-1&0\endm
\\&\hspace{12cm}\phi_{0}^{4}\bm \tilde{U}^{(t)}(s)\\
\overline{\tilde{U}^{(t)}(s)}\endm]_{dis}||_{L_{x}^{M}}ds||_{L_{t}^{M}}\\
&\lesssim
C_{0}^{-1}\Lambda^{2}\delta^{\frac{11}{4}-\delta_{1}}||\la
t\ra^{\frac{1}{2}-\delta_{2}}\int_{0}^{t}(t-s)^{-\frac{1}{2}+\frac{1}{M}}\la
s\ra^{-\frac{3}{2}+\delta_{3}}ds||_{L_{t}^{M}}\lesssim
C_{0}^{-1}\Lambda^{2}\delta^{\frac{11}{4}-\delta_{1}}
\end{split}\end{equation}
upon invoking \eqref{disasympto}, \eqref{rootasympto}, which in
turn can be bounded by $\leq\frac{\Lambda}{100}\delta$ upon
choosing $\delta$ etc small enough. Now assume that $t\geq
\delta^{-\frac{1}{2}}$. Then we get
\begin{equation}\nonumber\begin{split}
&C_{0}^{-1}||\chi_{>\delta^{-\frac{1}{2}}}(t)\la
t\ra^{\frac{1}{2}-\delta_{2}}\int_{0}^{t}||e^{i(t-s)\calH}[\bm
0&-e^{2i(\Psi-\Psi_{\infty})_{1}(s)-2i(\Psi-\Psi_{\infty})_{1}(t)}+1
\\ e^{-2i(\Psi-\Psi_{\infty})_{1}(s)+2i(\Psi-\Psi_{\infty})_{1}(t)}-1&0\endm
\\&\hspace{12cm}\phi_{0}^{4}\bm \tilde{U}^{(t)}(s)\\
\overline{\tilde{U}^{(t)}(s)}\endm]_{dis}||_{L_{x}^{M}}ds||_{L_{t}^{M}}\\
&\lesssim C_{0}^{-1}\Lambda\delta
||\chi_{>\delta^{-\frac{1}{2}}}(t)\la
t\ra^{\frac{1}{2}-\delta_{2}}\int_{0}^{t}(t-s)^{-\frac{1}{2}+\frac{1}{M}}\la
s\ra^{-\frac{3}{2}+\delta_{3}}ds||_{L_{t}^{M}}\lesssim
C_{0}^{-1}\delta^{-\frac{1}{M}+\frac{\delta_{2}}{2}}\delta\Lambda,
\end{split}\end{equation}
which is also $\leq\frac{\Lambda}{100}\delta$ upon choosing
$\delta$ small enough. The remaining local terms in $[...]_{dis}$
can be handled analogously, so we now consider the contribution of
the non-local term, which is
\begin{equation}\nonumber
C_{0}^{-1}||\la
t\ra^{\frac{1}{2}-\delta_{2}}||\int_{0}^{t}e^{i(t-s)\calH}\bm
|\tilde{U}^{(t)}|^{4}\tilde{U}^{(t)}(s,.)\\-|\tilde{U}^{(t)}|^{4}\overline{\tilde{U}^{(t)}(s,.)}\endm
ds||_{L_{x}^{M}}||_{L_{t}^{M}}
\end{equation}
Again referring to \eqref{disasympto}, as well as
Theorem~\ref{linearestimates}, and using H\"{o}lder's inequality,
we can bound this by
\begin{equation}\nonumber
\lesssim C_{0}^{-1}||\la
t\ra^{\frac{1}{2}-\delta_{2}+\frac{5}{M}}(\Lambda
C_{0})^{5}\delta^{5}\int_{0}^{t}(t-s)^{-\frac{1}{2}+\frac{1}{M}}\la
s\ra^{-\frac{3}{2}+3\delta_{2}}ds||_{L_{t}^{M}} \lesssim
\frac{\Lambda}{100}\delta
\end{equation}
upon choosing $\delta>0$ small enough, as desired. The estimate
for $||\la
t\ra^{-\delta_{2}}||\tilde{U}(t,.)||_{L_{x}^{2}}||_{L_{t}^{M}}$ is
carried out similarly and omitted. Moreover, we postpone
retrieving the difficult strong local dispersive estimate, i. e.
the global bound for $||\la s\ra^{\frac{3}{2}-\delta_{3}}\phi
\tilde{U}||_{L_{s}^{\infty}L_{x}^{\infty}}$, as well as the
pseudo-conformal almost conservation law, until later. This then
completes the case $k=0$. \\
We move on to the case $k=1$. We start by retrieving control over
\begin{equation}\nonumber
\sup_{\phi\in\calA}||\la
s\ra^{1-20\delta_{2}}||\phi\partial_{x}\tilde{U}_{dis}(s,x)||_{L_{x}^{M}}||_{L_{s}^{M}},\,
||\la
s\ra^{-10\delta_{2}}||\partial_{x}\tilde{U}_{dis}(s,x)||_{L_{x}^{2}}||_{L_{s}^{M}},
\end{equation}
which we do in tandem. Thus assuming that the quantity
\begin{equation}\nonumber
B(t):=C_{1}^{-1}[\sup_{\phi\in \calA}||\la
s\ra^{1-20\delta_{2}}||\phi\partial_{x}\tilde{U}_{dis}(s,.)||_{L_{x}^{M}}||_{L_{s}^{M}([0,t])}+||\la
s\ra^{-10\delta_{2}}||\partial_{x}\tilde{U}_{dis}(s,.)||_{L_{x}^{2}}||_{L_{s}^{M}([0,t])}]\leq
\Lambda\delta,
\end{equation}
we shall boost this to the bound $\lesssim
\frac{\Lambda}{100}\delta$. Commence with the expression
\begin{equation}\nonumber
C_{1}^{-1}\sup_{\phi\in \calA}||\la
s\ra^{1-20\delta_{2}}||\phi\partial_{x}\tilde{U}_{dis}(s,x)||_{L_{x}^{M}}||_{L_{s}^{M}([0,t])}
\end{equation}
If we differentiate \eqref{Ueqn} with respect to $x$, we produce
additional local source terms of the schematic form $VU$ for
certain Schwartz functions $V$, in addition to terms involving
$\partial_{x}U$. The former terms are handled by using the already
established estimates for $k=0$ as well as the assumption
$C_{1}>>C_{0}$. As for the latter, the most difficult local term
leads to the Duhamel term
\begin{equation}\nonumber\begin{split}
&C_{1}^{-1}||\la s\ra^{1-20\delta_{2}}\int_{0}^{s}||\phi(x)
e^{i(s-\lambda)\calH}[\bm
0&-e^{2i(\Psi-\Psi_{\infty})_{1}(\lambda)-2i(\Psi-\Psi_{\infty})_{1}(s)}+1
\\ e^{-2i(\Psi-\Psi_{\infty})_{1}(\lambda)+2i(\Psi-\Psi_{\infty})_{1}(s)}-1&0\endm
\\&\hspace{11cm}\phi_{0}^{4}\bm \partial_{x}\tilde{U}^{(s)}(\lambda)\\
\overline{\partial_{x}\tilde{U}^{(s)}(\lambda)}\endm]_{dis}||_{L_{x}^{M}}d\lambda||_{L_{s}^{M}([0,t])}\\
\end{split}\end{equation}
First restrict integration to the interval
$[0,s-\delta^{-\frac{1}{2}}\la s\ra^{\frac{6}{M}}]$. This we can
estimate crudely by\footnote{One again uses H\"{o}lder's
inequality to handle the fact that we only control
$||\phi\partial_{x}\tilde{U}^{(s)}(\lambda,x)||_{L_{\lambda}^{M}L_{x}^{M}}$.}
\begin{equation}\nonumber
\lesssim \Lambda\delta C_{1}^{-1}||\la
s\ra^{1-20\delta_{2}+\frac{1}{M}}\int_{0}^{s-\delta^{-\frac{1}{2}}\la
s\ra^{\frac{6}{M}}}\la s\ra^{\frac{2}{M}}\la
s-\lambda\ra^{-\frac{3}{2}}\la
\lambda\ra^{-1+20\delta_{2}}d\lambda||_{L_{s}^{M}([0,t])} \lesssim
\frac{\Lambda}{100}\delta
\end{equation}
upon choosing $\delta$ small enough. On the interval
$[s-\delta^{-\frac{1}{2}}\la s\ra^{\frac{6}{M}},s]$, one proceeds
similarly, but exploits the fact that thanks to the proof of
Lemma~\ref{Psi1}, we have
\begin{equation}\nonumber
\sup_{\lambda\in [s-\delta^{-\frac{1}{2}}\la
s\ra^{\frac{6}{M}},s]}|-e^{2i(\Psi-\Psi_{\infty})_{1}(\lambda)-2i(\Psi-\Psi_{\infty})_{1}(s)}+1|\lesssim
\delta^{\frac{3}{2}} \la
s\ra^{-\frac{1}{2}+\delta_{1}+\frac{6}{M}}
\end{equation}
The remaining local terms involving $\partial_{x}\tilde{U}^{(t)}$
are easier and estimated similarly. Next, consider the non-local
term. We have to estimate
\begin{equation}\nonumber
C_{1}^{-1}||\la s\ra^{1-20\delta_{2}}||\phi\int_{0}^{s}e^{i(s-\lambda)\calH}\bm |\tilde{U}^{(s)}(\lambda,.)|^{4}\partial_{x}\tilde{U}^{(s)}(\lambda,.)\\
-|\tilde{U}^{(s)}(\lambda,.)|^{4}\overline{\partial_{x}\tilde{U}^{(s)}}(\lambda,.)\endm_{dis}
d\lambda||_{L_{x}^{M}}||_{L_{s}^{M}([0,t])}
\end{equation}
in addition to similar terms (Leibnitz rule). If we invoke the
assumption on $B(t)$, the weighted estimates in
Theorem~\ref{linearestimates} as well as the estimates for the
case $k=0$ we can bound this by
\begin{equation}\nonumber
\lesssim \Lambda^{5}C_{0}^{4}\delta^{5}||\la s\ra^{1-20\delta_{2}}
\int_{0}^{s}(s-\lambda)^{-1+20\delta_{2}-\frac{1}{M}}\la
\lambda\ra^{4\delta_{2}+10\delta_{2}-\frac{3}{2}}\la
\lambda\ra^{\frac{1}{2}-20\delta_{2}+\frac{6}{M}}d\lambda||_{L_{s}^{M}([0,t])}\lesssim
\frac{\Lambda}{100}\delta
\end{equation}
We move on to control $||\la
s\ra^{-10\delta_{2}}||\partial_{x}\tilde{U}_{dis}(s,x)||_{L_{x}^{2}}||_{L_{s}^{M}([0,t])}$.
Note that {\it{we cannot simply reiterate the Duhamel procedure
here, since this would lead to a loss for the local
terms\footnote{More precisely, it appears that the extra factor
$\la s\ra^{-10\delta_{2}}$ should allow one to gain a bit;
however, the loss from the weak local decay control, i. e.
$||\phi\partial_{x}U(s,.)||_{L_{x}^{\infty}}\lesssim \la
s\ra^{-1+20\delta_{2}}$, is too much. On the other hand, if one
strengthened the weight in the $L^{2}$-norm to $\la
s\ra^{-(20+)\delta_{2}}$, one would encounter difficulties in
retrieving the weak local estimate. The reader may ask why we
don't build in the strong local decay for all higher derivatives
to begin with. The problem with this is that we would have to gain
control over more and more derivatives that way, indeed forcing
control over infinitely many derivatives.}.}} Observe that
\eqref{Ueqn} implies the following schematic relation:
\begin{equation}\label{8}\begin{split}
&i\partial_{s}[\partial_{x}\tilde{U}^{(t)}_{dis}](s,.)+\partial_{x}^{2}[\partial_{x}U^{(t)}_{dis}](s,.)-\partial_{x}U^{(t)}_{dis}(s,.)\\&\hspace{2cm}=[V\partial_{x}U^{(t)}]_{dis}(s,.)+[V\overline{\partial_{x}U^{(t)}}]_{dis}(s,.)+\ldots+[|U^{(t)}|^{4}\partial_{x}U^{(t)}]_{dis}(s,.),
\end{split}\end{equation}
where the subscripts really refer to the top entry of the
dispersive part of the corresponding vector-valued function. As
usual $V$ refers to certain Schwartz functions which may depend on
$t$ and $s$. We deduce that
\begin{equation}\label{7}\begin{split}
&i\partial_{\lambda}\int_{-\infty}^{\infty}|\partial_{x}\tilde{U}^{(s)}|^{2}(\lambda,.)dx=2\Im\int_{-\infty}^{\infty}i\partial_{\lambda}[e^{i\lambda}\partial_{x}\tilde{U}^{(s)}](\lambda,.)\overline{e^{i\lambda}\partial_{x}\tilde{U}^{(s)}(\lambda,.)}dx
\\&=2\Im\int_{-\infty}^{\infty}[V\partial_{x}\tilde{U}^{(s)}(\lambda,.)\overline{\partial_{x}\tilde{U}^{(s)}(\lambda,.)}dx+\ldots+|U^{(s)}|^{4}\partial_{x}\tilde{U}^{(s)}\overline{\partial_{x}\tilde{U}^{(s)}}(\lambda,.)dx
\end{split}\end{equation}
Using the already improved local bound, we can estimate
\begin{equation}\nonumber
||\la
s\ra^{-20\delta_{2}}\int_{0}^{s}|\Im\int_{-\infty}^{\infty}[V\partial_{x}\tilde{U}^{(s)}(\lambda,.)\overline{\partial_{x}\tilde{U}^{(s)}(\lambda,.)}dx|d\lambda||_{L_{s}^{M}([0,t])}\lesssim
(\frac{C_{1}\Lambda}{100})^{2}\delta^{2}\int_{0}^{s}\la \lambda\ra
^{-2(1-20\delta_{2})}d\lambda,
\end{equation}
which is integrable in $\lambda$ upon choosing $\delta_{2}$ small
enough. The remaining expressions on the righthand side of
\eqref{7} can be estimated similarly, combining the preceding estimates we get the improved bound on $B(t)$, as desired.\\
Next, consider the norm $||\la
s\ra^{\frac{1}{2}-25^{k}\delta_{2}}||\partial_{x}\tilde{U}_{dis}(s,.)||_{L_{x}^{M}}||_{L_{s}^{M}([0,T))}$.
This we estimate by reverting to the usual Duhamel's formula; we
treat here the most difficult local term (schematically)
\[[-e^{-2i(\Psi-\Psi_{\infty})_{1}(\lambda)+2i(\Psi-\Psi_{\infty})_{1}(s)}+1]V\partial_{x}U(\lambda,.)\]
the others following in a similar vein\footnote{Indeed, one gains
extra weights like $\nu-1$ for the other local terms, and can
proceed analogously. For the non-local term, use that
$|||\tilde{U}^{(s)}|^{4}\partial_{x}\tilde{U}^{(s)}(\lambda,.)||_{L_{x}^{M'}}\lesssim
\la
\lambda\ra^{-\frac{3}{2}+4\delta_{2}}||\partial_{x}\tilde{U}^{(s)}(\lambda,.)||_{L_{x}^{2}}$.}:
we have
\begin{equation}\nonumber\begin{split}
&||\la
s\ra^{\frac{1}{2}-25\delta_{2}}||\int_{0}^{s}e^{i(s-\lambda)\calH}[-e^{-2i(\Psi-\Psi_{\infty})_{1}(\lambda)+2i(\Psi-\Psi_{\infty})_{1}(s)}+1]][V\partial_{x}U]_{dis}(\lambda,.)d\lambda||_{L_{x}^{M}}||_{L_{s}^{M}}
\\&\lesssim \frac{\Lambda C_{1}}{100}\delta||\la
s\ra^{\frac{1}{2}-25\delta_{2}}\int_{0}^{s}(s-\lambda)^{-\frac{1}{2}+\frac{1}{M}}\lambda^{-1+20\delta_{2}}d\lambda||_{L_{s}^{M}}
\end{split}\end{equation}
This is easily seen to also improve the bound $\Lambda\delta
C_{1}$, if necessary by improving the preceding estimates for
$B(t)$. In order to complete the case $k=1$, we still need to
retrieve control over $K^{-1}\sup_{t\in
[0,T]}||C\partial_{x}\tilde{U}^{(t)}_{dis}||_{L_{x}^{2}}$, which
we shall do later. \\ Proceeding to higher derivatives $k\geq 2$
is an elementary induction, recycling the same estimates. Observe
that if one differentiates the quintilinear non-local term $k$
times, one obtains schematically either
$|\tilde{U}^{(s)}|^{4}\partial_{x}^{k}(\tilde{U}^{(s)})$, or else
$\partial_{x}^{k-1}(\tilde{U}^{(s)})\partial_{x}(\tilde{U}^{(s)})(\tilde{U}^{(s)})^{2}\overline{\tilde{U}^{(s)}}$
or else terms of the form
$\partial_{x}^{\alpha_{1}}(\tilde{U}^{(s)})\partial_{x}^{\alpha_{2}}(\overline{\tilde{U}^{(s)}})\ldots\partial_{x}^{\alpha_{5}}(\tilde{U}^{(s)})$,
where all $\alpha_{i}<k-1$, or terms equivalent to these for all
intents and purposes. The first term in this list one can treat
just as before, using the estimates for $U$. For the 2nd, one
estimates
\begin{equation}\nonumber\begin{split}
&||(\partial_{x}^{k-1}(\tilde{U}^{(s)})\partial_{x}(\tilde{U}^{(s)})(\tilde{U}^{(s)})^{2}\overline{\tilde{U}^{(s)}})(\lambda,.)||_{L_{x}^{M'}}\leq
\Lambda^{5}C_{k-1}C_{1}C_{0}^{3}\delta^{5}\lambda^{10^{k-1}\delta_{2}+10\delta_{2}+3\delta_{2}}\la\lambda\ra^{-\frac{3}{2}}
\\&\times ||\la\lambda\ra^{\frac{1}{2}-\delta_{2}}\tilde{U}^{(s)}(\lambda,.)||_{L_{x}^{M}}||\la\lambda\ra^{-10^{k-1}\delta_{2}}\partial_{x}^{k-1}\tilde{U}^{(s)}(\lambda,.)||_{L_{x}^{2+}}||\la
\lambda\ra^{-\delta_{2}}\partial_{x}\tilde{U}^{(s)}(\lambda,.)||_{L_{x}^{2}},
\end{split}\end{equation}
and we have
\begin{equation}\nonumber
\lambda^{10^{k-1}\delta_{2}+10\delta_{2}+3\delta_{2}}\la\lambda\ra^{-\frac{3}{2}}\lesssim
\la \lambda\ra^{-\frac{3}{2}+15^{k}\delta_{2}}
\end{equation}
Thus one can comfortably absorb an extra weight
$\lambda^{\frac{1}{2}-20^{k}\delta_{2}+\frac{4}{M}}$ here, and
continues as before. The estimate for the last term in the above
list is similar. Time derivatives can now be handled upon turning them into spatial derivatives via \eqref{Ueqn}.\\
Finally, one retrieves control over
\begin{equation}\nonumber
\sum_{1\leq k\leq
N-1}K^{-k}\sup_{2i+j=k}||C\partial_{s}^{i}\partial_{y}^{j}U(s,.)||_{L_{s}^{M}L_{y}^{2}([0,T])}]
\end{equation}
by putting $K=\la T\ra^{100}$, say, using the fact that $C=x-2pt$
and $i\partial_{t}+\triangle$ commute, and using the already
improved estimates as well as crude bounds. We also observe that
the preceding estimates can easily be bootstrapped to yield the
final bounds in Theorem~\ref{core}.
\\

We now have to come to terms with the {\it{strong local dispersive
estimate}}, or SLDE,  as well as the {\it{pseudo-conformal almost
conservation}}, i. e. the expressions
\begin{equation}\nonumber
\sup_{\phi\in\calA}\sup_{t\geq 0}\la
t\ra^{\frac{3}{2}-\delta_{3}}||\phi
\tilde{U}_{dis}(s,.)||_{L_{x}^{\infty}},\,\sup_{t\geq 0}||C
\tilde{U}_{dis}(t,.)||_{L_{x}^{2}}
\end{equation}
We commence with SLDE in the following subsection.

\subsection{Retrieving the strong local
dispersion.}\footnote{The argument to follow is certainly not
optimal; however, it allows us to gain a refined understanding
which will play an important role for the bilinear estimates
needed to control $\lambda_{6}$. A large simplification would
result if one could improve the SLDE to not contain any losses.}

We now need to show that we can deduce the inequality
$\sup_{\phi\in \calA}\sup_{0\leq t}\la
t\ra^{\frac{3}{2}-\delta_{3}}||\phi
\tilde{U}(t,.)||_{L_{x}^{\infty}}\leq \frac{\Lambda}{100}\delta$.
For this we need to employ \eqref{Udis}, which forces us to
distinguish between the different kinds of expressions on the
right hand side. Clearly we may assume $t\geq 1$. We first observe
that Theorem~\ref{linearestimates} in conjunction with our
assumptions on $A(.)$ as well as \eqref{modulationasympto},
\eqref{rootasympto}, \eqref{disasympto} and Lemma~\ref{Psi1} imply
that the free contribution is acceptable, with a $t^{\delta_{3}}$
to spare:
\begin{equation}\nonumber\begin{split}
&\sup_{0\leq t}\la t\ra^{\frac{3}{2}}||e^{it\calH}P_{s}[\bm
e^{i(\Psi_{\infty}-\Psi)_{1}(t)}&0\\0&e^{-i(\Psi_{\infty}-\Psi)_{1}(t)}\endm[\bm
A(.+\lambda_{\infty}(\mu-\mu_{\infty})(t))\\\bar{A}(.+\lambda_{\infty}(\mu-\mu_{\infty})(t))\endm_{dis}
\\&\hspace{7cm}+\sum_{j=1}^{6}\alpha_{j}\eta_{j,\text{proper}}(.+\lambda_{\infty}(\mu-\mu_{\infty})(t))]||\leq
\frac{\Lambda}{100}\delta,
\end{split}\end{equation}
upon choosing $\Lambda>0$ large enough. We next subdivide
$[...]_{dis}$ in the integrand of the Duhamel term in \eqref{Udis}
into local and non-local contributions. As for the local
contributions, as usual the most difficult is (schematically)
\begin{equation}\nonumber
\int_{0}^{t}[-e^{-2i(\Psi-\Psi_{\infty})_{1}(s)+2i(\Psi-\Psi_{\infty})_{1}(t)}+1]e^{i(t-s)\calH}[V\tilde{U}^{(t)}]_{dis}(s,.)ds
\end{equation}
One estimates (using Theorem~\ref{linearestimates})
\begin{equation}\nonumber\begin{split}
&\la t\ra
^{\frac{3}{2}-\delta_{3}}||\phi(x)\int_{0}^{t-\delta^{-\frac{1}{2}}}[-e^{-2i(\Psi-\Psi_{\infty})_{1}(s)+2i(\Psi-\Psi_{\infty})_{1}(t)}+1]e^{i(t-s)\calH}[V\tilde{U}^{(t)}]_{dis}(s,.)ds||_{L_{x}^{\infty}}
\\&\lesssim \la t\ra
^{\frac{3}{2}-\delta_{3}}\Lambda\delta\int_{0}^{t-\delta^{-\frac{1}{2}}}\la
t-s\ra^{-\frac{3}{2}}\la s\ra^{-\frac{3}{2}+\delta_{3}}ds\lesssim
\Lambda\delta^{1+\frac{\delta_{3}}{2}},
\end{split}\end{equation}
which leads to the bound $\lesssim \frac{\Lambda}{100}\delta$ upon
choosing $\delta>0$ small enough. Moreover, using
Lemma~\ref{Psi1}, one gets the same estimate for the integral over
$[t-\delta^{-\frac{1}{2}},t]$, as desired. The remaining local
terms in $[...]_{dis}$ can be handled similarly, whence we now
turn to the real task, dealing with the non-local term, i. e. the
expression
\begin{equation}\nonumber
\la
t\ra^{\frac{3}{2}-\delta_{3}}||\phi(x)\int_{0}^{t}e^{i(t-s)\calH}\bm
|\tilde{U}^{(t)}|^{4}\tilde{U}^{(t)}(s,.)\\-|\tilde{U}^{(t)}|^{4}\overline{\tilde{U}^{(t)}(s,.)}\endm_{dis}
ds||_{L_{x}^{\infty}}
\end{equation}
We intend to turn this into an expression of the following form:
\begin{equation}\nonumber
\la t\ra^{\frac{3}{2}-\delta_{3}}\la\int_{0}^{t}e^{i(t-s)\calH}\bm
|\tilde{U}^{(t)}|^{4}\tilde{U}^{(t)}(s,.)\\-|\tilde{U}^{(t)}|^{4}\overline{\tilde{U}^{(t)}(s,.)}\endm_{dis}
ds, \bm \phi\\ \psi\endm\ra,
\end{equation}
for suitable Schwartz functions $\phi,\psi$. The device for
achieving this is the discrete Fourier transform. First, using a
partition of unity $\{\phi_{i}\}$ subordinate to intervals of
length $2\pi$, we reduce to estimating
\begin{equation}\nonumber
\phi_{j}(x)\phi(x)\int_{0}^{t}e^{i(t-s)\calH}\bm
|\tilde{U}^{(t)}|^{4}\tilde{U}^{(t)}(s,.)\\-|\tilde{U}^{(t)}|^{4}\overline{\tilde{U}^{(t)}(s,.)}\endm_{dis}
ds
\end{equation}
Write (we omit the superscripts $\tilde{}$ and $^{(t)}$ from now
on as they are irrelevant in this argument)
\begin{equation}\nonumber
i\phi\phi_{j}\int_{0}^{t}e^{i(t-s)\calH}\bm
|U|^{4}U(s)\\-|U|^{4}\bar{U}(s)\endm_{dis}
ds=\sum_{n\in{\mathbf{Z}}}\bm a_{nj}e^{in(x-x_{j})}\\
\overline{a_{nj}}e^{-in(x-x_{j})}\endm
\end{equation}
We have
\begin{equation}\nonumber\begin{split}
&\Re a_{nj}=\frac{i}{4\pi}\la \phi\phi_{j}
\int_{0}^{t}e^{i(t-s)\calH}\bm
|U|^{4}U(s)\\-|U|^{4}\bar{U}(s)\endm_{dis}ds, \bm e^{in(x-x_{j})}\\
e^{-in(x-x_{j})}\endm\ra,\\&\Im a_{nj}= \frac{1}{4\pi}\la
\phi\phi_{j} \int_{0}^{t}e^{i(t-s)\calH}\bm
|U|^{4}U(s)\\-|U|^{4}\bar{U}(s)\endm_{dis}ds, \bm e^{in(x-x_{j})}\\
-e^{-in(x-x_{j})}\endm\ra
\end{split}\end{equation}
Thus for example
\begin{equation}\nonumber
\Re a_{nj}=\frac{i}{4\pi}\la \int_{0}^{t}e^{i(t-s)\calH}\bm
|U|^{4}U(s)\\-|U|^{4}\bar{U}(s)\endm_{dis} ds, \bm \phi\phi_{j}
e^{in(x-x_{j})}\\ \phi\phi_{j}e^{-in(x-x_{j})}\endm\ra,
\end{equation}
and similarly for the imaginary part. In order to be able to carry
out the summation over $n$, we need to carry out an integration by
parts. Another way to approach this is to note that only
moderately small values of $n$, i. e. $|n|<t^{\epsilon}$,
contribute since we can independently control
$||\partial_{x}^{N}U||_{L_{x}^{2}}$, whence the large frequency
part of $U$ can be made arbitrarily small. Indeed, note that we
have
\begin{equation}\nonumber
\la \bm U\\ \bar{U}\endm, \phi e^{inx}\ra =\frac{1}{(in)^{N}}\la
\partial_{x}^{N}[\phi\bm U\\ \bar{U}\endm],  e^{inx}\ra
\end{equation}
Thus we may and shall assume that
$|n|<\delta^{-\epsilon}t^{\epsilon}$, for
$\epsilon>\epsilon_{0}(N)>0$, $N$ as in \eqref{disasympto}, such
that $\lim_{N\rightarrow\infty}\epsilon_{0}(N)=0$. We thus lose
$\delta^{-\epsilon}\la t\ra^{\epsilon}$ in the end, which we can
afford since we may arrange $\delta_{3}>>\epsilon_{0}(N)$. We now
need to estimate the following expression:
\begin{equation}\nonumber
\int_{0}^{t}\la \bm |U|^{4}U(s)\\ -|U|^{4}\bar{U}(s)\endm ,
e^{-i(t-s)\calH^{*}}\bm \phi\\\psi\endm_{dis}\ra ds
\end{equation}
We distinguish between the cases $s>\frac{t}{2}$, $s\leq
\frac{t}{2}$, which aside from a simple technicality are treated
by the same method. We shall use the notation $P_{\leq a}$,
$P_{a}$, $P_{>a}$, $a\in{\mathbf{R}}_{>0}$ dyadic\footnote{We
shall more generally mean $P_{\alpha}$ for
$\alpha\in{\mathbf{R}}_{>0}$ to denote $P_{j}$ if $2^{j-1}\leq
\alpha<2^{j}$}, for the standard Littlewood-Paley multipliers, see
e. g. [St]. We shall also assume that $\tilde{U}$ satisfies
pointwise-in-time estimates below in order to simplify the
exposition; thus we shall assume bounds of the form $\la
s\ra^{\frac{1}{2}-\delta_{2}}||\tilde{U}(s,.)||_{L_{x}^{\infty}}\leq
\delta\Lambda$ etc. It will be straightforward to adjust the
arguments below to the case of weighted-in-time norms, since we
assume $M>>\delta_{2}^{-1}$ and hence we can absorb losses of
order $\la s\ra^{O(\frac{1}{M})}$. Finally, as it is clear that we
gain lots of $\delta$'s below, we shall occasionally omit the
$\Lambda$.
\\

{\bf{Case A}}: $s<\frac{t}{2}$. The idea is to exploit the
pseudo-conformal operator to reduce at least one of the two
$|U|^{2}$'s to frequency $<s^{-\frac{3}{4}}$. In this case, one
exploits the fact that the distorted Fourier transform vanishes at
the origin. We start by chopping things apart: specialize to the
following two terms:
\begin{equation}\nonumber
\int_{0}^{\frac{t}{2}}\la \bm \chi_{>0}|U|^{4}U(s)\\
-\chi_{>0}|U|^{4}\bar{U}(s)\endm , e^{-i(t-s)\calH^{*}}\bm
\phi\\\psi\endm_{dis}\ra ds
\end{equation}
\begin{equation}\nonumber
\int_{0}^{\frac{t}{2}}\la \bm \chi_{<0}|U|^{4}U(s)\\
-\chi_{<0}|U|^{4}\bar{U}(s)\endm , e^{-i(t-s)\calH^{*}}\bm
\phi\\\psi\endm_{dis}\ra ds,
\end{equation}
where $\chi_{>0}$ is the Heaviside function localizing to $x>0$.
Both are treated the same way, so consider the first expression:
rewrite it as
\begin{equation}\nonumber
\int_{0}^{\frac{t}{2}}\la \bm \chi_{>0}|U|^{4}(s)\\
-\chi_{>0}|U|^{4}(s)\endm, \bm\bar{U}(s)\\U(s)\endm\times
e^{-i(t-s)\calH^{*}}\bm \phi\\\psi\endm_{dis}\ra ds,
\end{equation}
where $\times$ denotes componentwise multiplication\footnote{Also
observe that we use the subscript $_{dis}$ both with reference to
$\calH$ as well as $\calH^{*}$.}. We commence
by reducing $\bm \chi_{>0}|U|^{4}(s)\\
-\chi_{>0}|U|^{4}(s)\endm$ to its dispersive part. To achieve
this, note that
\begin{equation}\nonumber
\bm \chi_{>0}|U|^{4}(s)\\
-\chi_{>0}|U|^{4}(s)\endm=\bm \chi_{>0}|U|^{4}(s)\\
-\chi_{>0}|U|^{4}(s)\endm_{dis}+\sum_{j=1}^{6}a_{k(j)}\la\bm \chi_{>0}|U|^{4}(s)\\
-\chi_{>0}|U|^{4}(s)\endm,
\xi_{k(j),\text{proper}}\ra\eta_{j,\text{proper}}
\end{equation}
where $\xi_{j, \text{proper}}$ is the basis for the generalized
root space of $\calH^{*}$, while $\eta_{j,\text{proper}}$ is the
basis for the generalized root space of $\calH$, as explained
earlier. The $a_{k(j)}$ are suitable numerical coefficients. Then
observe that by the improved local dispersive estimate, we have
($\epsilon=\epsilon(\delta_{3})$)
\begin{equation}\nonumber
|\la\bm \chi_{>0}|U|^{4}(s)\\
-\chi_{>0}|U|^{4}(s)\endm, \xi_{j,\text{proper}}\ra|\lesssim
(\Lambda C_{0})^{4}\delta^{4} \la s\ra^{-6+\epsilon},
\end{equation}
whence we treat the contribution of this part to the above
integral expression by
\begin{equation}\nonumber
\lesssim (\Lambda
C_{0})^{4}\delta^{4}\int_{0}^{\frac{t}{2}}(t-s)^{-\frac{3}{2}}\la
s\ra^{-6+\epsilon}ds,
\end{equation}
which is better than what we need. Thus we now consider
\begin{equation}\nonumber
\int_{0}^{\frac{t}{2}}\la \bm \chi_{>0}|U|^{4}(s)\\
-\chi_{>0}|U|^{4}(s)\endm_{dis}, \bm\bar{U}(s)\\U(s)\endm\times
e^{-i(t-s)\calH^{*}}\bm \phi\\\psi\endm_{dis}\ra ds
\end{equation}
Using the distorted Plancherel's Theorem~\ref{DistortedPlancherel}
we can equate this with
\begin{equation}\nonumber
\sum_{\pm}\int_{-\infty}^{\infty}\int_{0}^{\frac{t}{2}}\calF_{\pm}\bm \chi_{>0}|U|^{4}(s)\\
-\chi_{>0}|U|^{4}(s)\endm(\xi)
\overline{\tilde{\calF}_{\pm}\bm\bar{U}(s)\\U(s)\endm\times
e^{-i(t-s)\calH^{*}}\bm \phi\\\psi\endm}(\xi) ds d\xi
\end{equation}
Observe that we have $\calF_{\pm}(\phi)(\xi)=\la \phi,
\sigma_{3}e_{\pm}(x,\xi)\ra$, $\tilde{\calF}(\phi)(\xi)=\la \phi,
e_{\pm}(x,\xi)\ra$. The cases $\pm$ are treated exactly
analogously, so we stick with the $+$ case. Break the
$\xi$-integral into two parts, one over $[0,\infty)$, the other
over $(-\infty,0]$. Commence with the case $\xi\in [0,\infty)$.
Write
\begin{equation}\nonumber\begin{split}
&\int_{0}^{\infty}\int_{0}^{\frac{t}{2}}\calF\bm \chi_{>0}|U|^{4}(s)\\
-\chi_{>0}|U|^{4}(s)\endm(\xi)\overline{\tilde{\calF}[\chi_{>0}\bm\bar{U}(s)\\U(s)\endm\times
e^{-i(t-s)\calH^{*}}\bm \phi\\\psi\endm_{dis}]}(\xi)ds d\xi\\
&=\int_{0}^{\infty}\int_{0}^{\frac{t}{2}}\la\bm\chi_{>0}|U|^{4}(s)\\
-\chi_{>0}|U|^{4}(s)\endm,
s(\xi)e^{ix\xi}\underline{e}+\phi(x,\xi)\ra\overline{\tilde{\calF}[\chi_{>0}\bm\bar{U}(s)\\U(s)\endm\times
e^{-i(t-s)\calH^{*}}\bm \phi\\ \psi\endm_{dis}]}(\xi)ds d\xi,\\
\end{split}\end{equation}
recalling Theorem~\ref{DistortedFourier}. We first treat the
simple contribution coming from the rapidly decaying function
$\phi(x,\xi)$. As before, observe that
\begin{equation}\nonumber
|\la\bm\chi_{>0}|U|^{4}(s)\\
\chi_{>0}|U|^{4}(s)\endm, \phi(x,\xi)\ra|\lesssim \la
s\ra^{-6+\epsilon},\,\epsilon=\epsilon(\delta_{3})
\end{equation}
Indeed, we can estimate the $L_{\xi}^{2}$-norm of the function on
the left in this fashion. Moreover, we have
\begin{equation}\nonumber
||\tilde{\calF}[\chi_{>0}\bm\bar{U}(s)\\U(s)\endm\times
e^{-i(t-s)\calH^{*}}\bm \phi\\
\psi\endm_{dis}](\xi)||_{L_{\xi}^{2}}\lesssim \la
t-s\ra^{-\frac{3}{2}}s^{1+\epsilon(\delta_{2})}
\end{equation}
We are using here the pseudo-conformal almost conservation, which
is part of our assumptions \eqref{disasympto}:
\begin{equation}\nonumber
||(x-2sp)U(s,.)||_{L_{x}^{2}}\lesssim
O(\delta),\,p=-i\frac{\partial}{\partial x}
\end{equation}
Thus feeding in the preceding two estimates easily implies the
desired bound. Now consider the difficult oscillatory part.
Decompose
\begin{equation}\nonumber\begin{split}
&\int_{0}^{\infty}\int_{0}^{\frac{t}{2}}\la\bm\chi_{>0}|U|^{4}(s)\\
-\chi_{>0}|U|^{4}(s)\endm,s(\xi)e^{ix\xi}\underline{e}\ra\overline{\tilde{\calF}[\chi_{>0}\bm\bar{U}(s)\\U(s)\endm\times
e^{-i(t-s)\calH^{*}}\bm \phi\\ \psi\endm_{dis}]}(\xi)ds d\xi=\\
&\int_{0}^{\infty}\int_{0}^{\frac{t}{2}}\la \bm P_{\geq a}[\chi_{>0}|U|^{2}(s)]P_{\geq a}[|U|^{2}]\\
-P_{\geq a}[\chi_{>0}|U|^{2}(s)]P_{\geq a}[|U^{2}|]\endm,
s(\xi)e^{ix\xi}\underline{e}\ra\overline{\tilde{\calF}[\chi_{>0}\bm\bar{U}(s)\\U(s)\endm\times
e^{-i(t-s)\calH^{*}}\bm \phi\\ \psi\endm_{dis}]}(\xi)ds d\xi\\
&+\int_{0}^{\infty}\int_{0}^{\frac{t}{2}}\la \bm P_{\geq a}[\chi_{>0}|U|^{2}(s)]P_{<a}[|U|^{2}]\\
-P_{\geq a}[\chi_{>0}|U|^{2}(s)]P_{< a}[|U^{2}|]\endm,
s(\xi)e^{ix\xi}\underline{e}\ra\overline{\tilde{\calF}[\chi_{>0}\bm\bar{U}(s)\\U(s)\endm\times
e^{-i(t-s)\calH^{*}}\bm \phi\\ \psi\endm_{dis}]}(\xi)ds d\xi\\
&+\int_{0}^{\infty}\int_{0}^{\frac{t}{2}}\la \bm P_{<a}[\chi_{>0}|U|^{2}(s)]P_{\geq a}[|U|^{2}]\\
-P_{<a}[\chi_{>0}|U|^{2}(s)]P_{\geq a}[|U^{2}|]\endm,
s(\xi)e^{ix\xi}\underline{e}\ra\overline{\tilde{\calF}[\chi_{>0}\bm\bar{U}(s)\\U(s)\endm\times
e^{-i(t-s)\calH^{*}}\bm \phi\\ \psi\endm_{dis}]}(\xi)ds d\xi\\
&+\int_{0}^{\infty}\int_{0}^{\frac{t}{2}}\la \bm P_{<a}[\chi_{>0}|U|^{2}(s)]P_{<a}[|U|^{2}]\\
-P_{<a}[\chi_{>0}|U|^{2}(s)]P_{< a}[|U^{2}|]\endm,
s(\xi)e^{ix\xi}\underline{e}\ra\overline{\tilde{\calF}[\chi_{>0}\bm\bar{U}(s)\\U(s)\endm\times
e^{-i(t-s)\calH^{*}}\bm \phi\\ \psi\endm_{dis}]}(\xi)ds d\xi\\
\end{split}\end{equation}
The cutoff $a$ here will be later chosen to be $\la
s\ra^{-\frac{3}{4}}$. We treat each of the above terms separately.
Start with the first, the {\bf{high-high case}}: note that we can
write
\begin{equation}\nonumber
\tilde{\calF}[\chi_{>0}\bm\bar{U}(s)\\U(s)\endm\times
e^{-i(t-s)\calH^{*}}\bm \phi\\ \psi\endm_{dis}](\xi)=\la
\chi_{>0}\bm\bar{U}(s)\\U(s)\endm\times e^{-i(t-s)\calH^{*}}\bm
\phi\\ \psi\endm_{dis},
s(\xi)e^{ix\xi}\underline{e}+\phi(x,\xi)\ra
\end{equation}
Thus we can estimate
\begin{equation}\nonumber\begin{split}
&|\int_{0}^{\infty}\int_{0}^{\frac{t}{2}}\la \bm P_{\geq a}[\chi_{>0}|U|^{2}(s)]P_{\geq a}[|U|^{2}]\\
-P_{\geq a}[\chi_{>0}|U|^{2}(s)]P_{\geq
a}[|U^{2}|]\endm,s(\xi)e^{ix\xi}\underline{e}\ra\overline{\la
\chi_{>0}\bm\bar{U}(s)\\U(s)\endm\times
e^{-i(t-s)\calH^{*}}\bm \phi\\ \psi\endm_{dis}, \phi(x,\xi)\ra}ds d\xi|\\
&\lesssim\int_{0}^{\frac{t}{2}}\la
s\ra^{-3+\epsilon(\delta_{3})}\la t-s\ra^{-\frac{3}{2}}ds\lesssim
\la t\ra^{-\frac{3}{2}}\\
\end{split}\end{equation}
Hence we reduce to estimating the expression
\begin{equation}\nonumber
\int_{0}^{\infty}\int_{0}^{\frac{t}{2}}\la\bm P_{\geq a}[\chi_{>0}|U|^{2}(s)]P_{\geq a}[|U|^{2}]\\
-P_{\geq a}[\chi_{>0}|U|^{2}(s)]P_{\geq a}[|U^{2}|]\endm,
s(\xi)e^{ix\xi}\underline{e}\ra
\overline{\la\chi_{>0}\bm\bar{U}(s)\\U(s)\endm\times
e^{-i(t-s)\calH^{*}}\bm \phi\\ \psi\endm_{dis},
s(\xi)e^{ix\xi}\underline{e}\ra} d\xi
\end{equation}
We intend to use the ordinary Plancherel's Theorem here. We break
this integral into two by including a multiplier
$\phi_{(t^{-1000},t^{1000})}(\xi)$ or
$\chi_{>0}(\xi)-\phi_{(t^{-1000},t^{1000})}(\xi)$, where
$\phi_{(t^{-1000},t^{1000})}(\xi)$ smoothly localizes to the
interval $(t^{-1000},t^{1000})$. It is easily seen that
contribution obtained upon including the latter is very small
(bounded by $\la t\ra ^{-500}$), whence we may focus on the
contribution of the former. By choosing
$\phi_{(t^{-1000},t^{1000})}(\xi)$ suitably, we may assume that
its Fourier transform has $L^{1}$-mass bounded by $\log t$. Now
denote the Fourier multiplier with symbol
$s(\xi)^{2}\phi_{(t^{-1000},t^{1000})}(\xi)$ by
$\Pi_{(t^{1000},t^{1000})}$. Using ordinary Plancherel, we now
reduce to estimating
\begin{equation}\nonumber \la P_{\geq
a}[\chi_{>0}|U|^{2}(s)]P_{\geq a}[|U|^{2}],
\Pi_{(t^{-1000},t^{1000})}\la
\underline{e},\chi_{>0}\bm\bar{U}(s)\\U(s)\endm\times
e^{-i(t-s)\calH^{*}}\bm \phi\\ \psi\endm_{dis}\ra
\end{equation}
We claim that we have
\begin{equation}\nonumber
|\Pi_{(t^{-1000},t^{1000})}\la
\underline{e},\chi_{>0}\bm\bar{U}(s)\\U(s)\endm\times
e^{-i(t-s)\calH^{*}}\bm \phi\\ \psi\endm_{dis}|\lesssim \la
t-s\ra^{-\frac{3}{2}}\la s\ra^{\frac{1}{2}+\epsilon(\delta_{2})}
\end{equation}
This follows from $||xU(s)||_{L_{x}^{\infty}}\lesssim \la
s\ra^{\frac{1}{2}+\epsilon(\delta_{2})}$, which in turn is a
consequence of
\begin{equation}\nonumber\begin{split}
&||xU||_{L_{x}^{\infty}}\leq
||(x+2is\partial_{x})U||_{L_{x}^{\infty}}+||2is\partial_{x}U||_{L_{x}^{\infty}}\\&\lesssim
||(x+2is\partial_{x})\partial_{x}U||_{L_{x}^{2}}+s||\partial_{x}U||_{L_{x}^{\infty}}+||(x+2is\partial_{x})U||_{L_{x}^{2}}
\end{split}\end{equation}
and the following bounds, the 2nd of which we establish later:
\begin{equation}\nonumber
||\partial_{x}U(s)||_{L_{x}^{\infty}}\lesssim \la
s\ra^{-\frac{1}{2}+\epsilon(\delta_{2})},\,||(x+2is\partial_{x})\nabla
U||_{L_{x}^{2}}\lesssim \la
s\ra^{\frac{1}{2}+\epsilon(\delta_{2})}
\end{equation}
Now consider
\begin{equation}\nonumber
P_{\geq a}[\chi_{>0}|U|^{2}(s)]P_{\geq a}[|U|^{2}]
\end{equation}
Observe that
\begin{equation}\nonumber
P_{\geq a}[\chi_{>0}|U|^{2}(s)]=P_{\geq
a}\partial_{x}\triangle^{-1}[\partial_{x}(\chi_{>0})|U|^{2}(s)]+P_{\geq
a}\partial_{x}\triangle^{-1}[\chi_{>0}\partial_{x}[|U|^{2}(s)]]
\end{equation}
Note that
\begin{equation}\nonumber
\partial_{x}(\chi_{>0})|U|^{2}(s)=\delta_{0}|U(0)|^{2},
\end{equation}
whence
\begin{equation}\nonumber
||P_{\geq
a}\partial_{x}\triangle^{-1}[\partial_{x}(\chi_{>0})|U|^{2}(s)]||_{L_{x}^{1}}\lesssim
a^{-1}\la s\ra^{-3+\epsilon(\delta_{3})}
\end{equation}
Next, use that
\begin{equation}\label{keyidentity}
is\partial_{x}[|U|^{2}]=is\partial_{x}U\bar{U}-U\overline{is\partial_{x}U}=(is\partial_{x}+\frac{x}{2})U\bar{U}-U\overline{(is\partial_{x}+\frac{x}{2})U},
\end{equation}
whence we get
\begin{equation}\nonumber
||P_{\geq
a}\partial_{x}\triangle^{-1}[\chi_{>0}\partial_{x}[|U|^{2}(s)]]||_{L_{x}^{2}}\lesssim
a^{-1}\la s\ra^{-1}\la s\ra^{-\frac{1}{2}}
\end{equation}
Arguing similarly for $P_{\geq a}[|U|^{2}]$, one gets
\begin{equation}\nonumber
||P_{\geq a}[\chi_{>0}|U|^{2}]P_{\geq
a}[|U|^{2}]||_{L_{x}^{1}}\lesssim a^{-2}\la s\ra^{-3}
\end{equation}
Combining with the bound on $\bm
\bar{U}\\U\endm\times e^{-i(t-s)\calH^{*}}\bm \phi\\
\psi\endm_{dis}$ given above, we can bound the whole expression by
\begin{equation}\nonumber\begin{split}
& |\la P_{\geq a}[\chi_{>0}|U|^{2}(s)]P_{\geq a}[|U|^{2}],
\Pi_{(t^{-1000}, t^{1000})}\la
\underline{e},\chi_{>0}\bm\bar{U}(s)\\U(s)\endm\times
e^{-i(t-s)\calH^{*}}\bm \phi\\ \psi\endm_{dis}\ra|\\
&\lesssim a^{-2}\la s\ra^{-3}\la
s\ra^{\frac{1}{2}+\epsilon(\delta_{2})}\la t-s\ra^{-\frac{3}{2}},
\end{split}\end{equation}
which yields something almost integrable in $s$
 upon omitting the factor $\la t-s\ra^{-\frac{3}{2}}$ provided we choose
$a=\la s\ra^{-\frac{3}{4}}$. This is good enough since by
assumption $\delta_{3}>>\delta_{2}$. We now consider the other
extreme, the case of {\bf{low-low frequency interactions}}, i. e.
the expression
\begin{equation}\nonumber
\int_{0}^{\infty}\int_{0}^{\frac{t}{2}}\la \bm P_{<a}[\chi_{>0}|U|^{2}(s)]P_{<a}[|U|^{2}]\\
-P_{<a}[\chi_{>0}|U|^{2}(s)]P_{< a}[|U^{2}|]\endm,
s(\xi)e^{ix\xi}\underline{e}\ra
\overline{\tilde{\calF}[\chi_{>0}\bm\bar{U}(s)\\U(s)\endm\times
e^{-i(t-s)\calH^{*}}\bm \phi\\ \psi\endm_{dis}]}(\xi)ds d\xi
\end{equation}
We can express $\tilde{\calF}...$ as before, and only the term
$s(\xi)e^{ix\xi}\underline{e}$ in the Fourier basis matters. On
account of the fact that
\begin{equation}\nonumber
\la\bm P_{<a}[\chi_{>0}|U|^{2}(s)]P_{<a}[|U|^{2}]\\
-P_{<a}[\chi_{>0}|U|^{2}(s)]P_{< a}[|U^{2}|]\endm,
s(\xi)e^{ix\xi}\underline{e}\ra
=\chi_{<a+O(1)}(\xi)\la \bm P_{<a}[\chi_{>0}|U|^{2}(s)]P_{<a}[|U|^{2}]\\
-P_{<a}[\chi_{>0}|U|^{2}(s)]P_{< a}[|U^{2}|]\endm,
s(\xi)e^{ix\xi}\underline{e}\ra,
\end{equation}
we can estimate
\begin{equation}\nonumber\begin{split}
&|\int_{0}^{\infty}\la\bm P_{<a}[\chi_{>0}|U|^{2}(s)]P_{<a}[|U|^{2}]\\
-P_{<a}[\chi_{>0}|U|^{2}(s)]P_{< a}[|U^{2}|]\endm,
s(\xi)e^{ix\xi}\underline{e}\ra\overline{\tilde{\calF}[\chi_{>0}\bm\bar{U}(s)\\U(s)\endm\times
e^{-i(t-s)\calH^{*}}\bm \phi\\ \psi\endm_{dis}]}(\xi)d\xi|\\
&\lesssim a^{2}\la s\ra^{-\frac{3}{2}}\la s\ra\la t-s\ra^{-\frac{3}{2}},\\
\end{split}\end{equation}
which for $a\sim \la s\ra^{-\frac{3}{4}}$ can be integrated in $s$
to yield the bound $\la t\ra^{-\frac{3}{2}}$. Finally, we consider
the {\bf{mixed case}}, i. e. the expression
\begin{equation}\nonumber
\int_{0}^{\infty}\int_{0}^{\frac{t}{2}}\la\bm P_{\geq a}[\chi_{>0}|U|^{2}(s)]P_{<a}[|U|^{2}]\\
- P_{\geq a}[\chi_{>0}|U|^{2}(s)]P_{< a}[|U^{2}|]\endm,
 s(\xi)e^{ix\xi}\underline{e}\ra\overline{\tilde{\calF}[\chi_{>0}\bm\bar{U}(s)\\U(s)\endm\times
e^{-i(t-s)\calH^{*}}\bm \phi\\ \psi\endm_{dis}]}(\xi)ds d\xi
\end{equation}
We proceed as before, simplifying $\tilde{\calF}(...)$ by
discarding the Schwartz term in the Fourier basis (as we may), and
using the ordinary Plancherel's Theorem to translate this to the
physical side. Arguing as before, we may do this by including a
multiplier $\Pi_{(t^{-1000},t^{1000})}$ which is given by a kernel
of $L^{1}$-mass $\lesssim \log t$. We indicate this by replacing
the functions $P_{\geq a}[\chi_{>0}|U|^{2}(s)]$, $P_{<a}[|U|^{2}]$
by translates, $T_{z}P_{\geq a}[\chi_{>0}|U|^{2}(s)]$ and
$T_{z}P_{<a}[|U|^{2}]$, where $(T_{z}f)(x):=f(x+z)$. The
integration over $z$ in the end will cost $\lesssim \log t$. Thus
we now need to consider the following expression:
\begin{equation}\nonumber
\la T_{z}P_{\geq a}[\chi_{>0}|U|^{2}(s)] T_{z}P_{<a}[|U|^{2}],
\la\underline{e}, \chi_{>0}\bm\bar{U}(s)\\U(s)\endm\times
e^{-i(t-s)\calH^{*}}\bm \phi\\ \psi\endm_{dis}\ra
\end{equation}
We re-arrange the terms here:
\begin{equation}\nonumber
\la T_{z}P_{<a}[|U|^{2}], T_{z}P_{\geq
a}[\chi_{>0}|U|^{2}(s)]\la\underline{e},
 \chi_{>0}\bm\bar{U}(s)\\U(s)\endm\times e^{-i(t-s)\calH^{*}}\bm
\phi\\ \psi\endm_{dis}\ra\ra
\end{equation}
Revert to vectorial notation:
\begin{equation}\nonumber
\la \bm T_{z}P_{<a}[|U|^{2}]\\ 0\endm, \bm T_{z}P_{\geq
a}[\chi_{>0}|U|^{2}(s)]\la\underline{e},
 \chi_{>0}\bm\bar{U}(s)\\U(s)\endm\times e^{-i(t-s)\calH^{*}}\bm
\phi\\ \psi\endm_{dis}\ra\\ 0\endm\ra
\end{equation}
We break this into two portions:
\begin{equation}\nonumber
\la \chi_{>0}(x)\bm T_{z}P_{<a}[|U|^{2}]\\ 0\endm, \bm
T_{z}P_{\geq a}[\chi_{>0}|U|^{2}(s)]\la\underline{e},
 \chi_{>0}\bm\bar{U}(s)\\U(s)\endm\times e^{-i(t-s)\calH^{*}}\bm
\phi\\ \psi\endm_{dis}\ra\\ 0\endm\ra
\end{equation}
\begin{equation}\nonumber
\la \chi_{<0}(x)\bm T_{z}P_{<a}[|U|^{2}]\\ 0\endm, \bm
T_{z}P_{\geq a}[\chi_{>0}|U|^{2}(s)]\la\underline{e},
 \chi_{>0}\bm\bar{U}(s)\\U(s)\endm\times e^{-i(t-s)\calH^{*}}\bm
\phi\\ \psi\endm_{dis}\ra\\ 0\endm\ra
\end{equation}
As these can be treated similarly, we consider only the first. Our
first step consists in reducing the factor $\chi_{>0}(x)\bm
T_{z}P_{<a}[|U|^{2}]\\ 0\endm$ to its dispersive part. Note that
if we substitute $\la \chi_{>0}(x)\bm T_{z}P_{<a}[|U|^{2}]\\
0\endm, \xi_{k(j)}\ra \eta_{j}$ for this expression instead, we
can estimate
\begin{equation}\nonumber\begin{split}
&\la\la \chi_{>0}(x)\bm T_{z}P_{<a}[|U|^{2}]\\
0\endm, \xi_{k(j)}\ra \eta_{j} , \bm T_{z}P_{\geq
a}[\chi_{>0}|U|^{2}(s)]\la\underline{e},
 \chi_{>0}\bm\bar{U}(s)\\U(s)\endm\times e^{-i(t-s)\calH^{*}}\bm
\phi\\ \psi\endm_{dis}\ra\\ 0\endm\ra\\
&\lesssim \la s\ra^{-\frac{3}{2}}\la t-s\ra^{-\frac{3}{2}}
\end{split}\end{equation}
This can be integrated in $s$ to yield the upper bound $\lesssim
\la t\ra^{-\frac{3}{2}}$. Now, with the left factor reduced to its
dispersive part, invoking the distorted Plancherel's
Theorem~\ref{DistortedPlancherel}, we need to estimate
\begin{equation}\nonumber\begin{split}
&\sum_{\pm}\int_{-\infty}^{\infty}\calF_{\pm}[\chi_{>0}(x)\bm
T_{z}P_{<a}[|U|^{2}]\\
0\endm](\xi)\\&\hspace{2cm}\overline{\tilde{\calF}_{\pm}[\bm
T_{z}P_{\geq a}[\chi_{>0}|U|^{2}(s)]\la\underline{e},
 \chi_{>0}\bm\bar{U}(s)\\U(s)\endm\times e^{-i(t-s)\calH^{*}}\bm
\phi\\ \psi\endm_{dis}\ra\\ 0\endm](\xi)} d\xi
\end{split}\end{equation}
We may and shall treat the case $+$, and omit the subscript for
simplicity. As before, we need to subdivide the $\xi$-integration
into two contributions, one from $(-\infty,0]$, the other from
$[0,\infty)$. We treat here the contribution from the latter, that
from the former being more complicated and treated below. We
decompose
\begin{equation}\nonumber
\calF[\chi_{>0}(x)\bm
T_{z}P_{<a}[|U|^{2}]\\
0\endm](\xi)=\calF[\chi_{>0}(x)\bm
T_{z}[|U|^{2}]\\
0\endm](\xi)-\calF[\chi_{>0}(x)\bm
T_{z}P_{\geq a}[|U|^{2}]\\
0\endm](\xi)
\end{equation}
Substituting the 2nd summand results in an expression which can be
treated like in the high-high case. Thus substitute the first
summand on the right, $\calF[\chi_{>0}(x)\bm
T_{z}[|U|^{2}]\\
0\endm](\xi)$. One explicitly writes out the Fourier transform,
and may discard the contribution from the local part $\phi(x,\xi)$
of the Fourier basis, reasoning as before. Then one obtains the
following expression:
\begin{equation}\nonumber\begin{split}
&\int_{0}^{\infty}\la\chi_{>0}(x)\bm
T_{z}[|U|^{2}]\\
0\endm,
s(\xi)e^{ix\xi}\underline{e}\ra\\&\hspace{2cm}\overline{\tilde{\calF}[\bm
T_{z}P_{\geq a}[\chi_{>0}|U|^{2}(s)]\la\underline{e},
 \chi_{>0}\bm\bar{U}(s)\\U(s)\endm\times e^{-i(t-s)\calH^{*}}\bm
\phi\\ \psi\endm_{dis}\ra\\ 0\endm](\xi)} d\xi
\end{split}\end{equation}
In this break $\chi_{>0}(x)\bm
T_{z}[|U|^{2}]\\
0\endm$ into two parts, a large frequency and a small frequency
part:
\begin{equation}\nonumber
\chi_{>0}(x)\bm
T_{z}[|U|^{2}]\\
0\endm=P_{\geq a}[\chi_{>0}(x)\bm
T_{z}[|U|^{2}]\\
0\endm] +P_{<a}[\chi_{>0}(x)\bm
T_{z}[|U|^{2}]\\
0\endm]
\end{equation}
Consider the first summand on the right: one may differentiate the
expression, replacing it by
\begin{equation}\nonumber
\partial_{x}\triangle^{-1}P_{\geq a}[\delta_{0}(x)\bm
T_{z}[|U|^{2}]\\
0\endm]+\partial_{x}\triangle^{-1}P_{\geq a}[\chi_{>0}(x)\bm
T_{z}\partial_{x}[|U|^{2}]\\
0\endm]
\end{equation}
Observe that the $a^{-1}$ from the operator
$\partial_{x}\triangle^{-1}P_{\geq a}$ is counteracted by the
factor $s(\xi)$ above. In order to treat the contribution from the
first summand, subdivide the interval $[a,\infty)$ into dyadic
intervals, and sum. Thus we need to estimate
\begin{equation}\nonumber\begin{split}
&\sum_{2^{j}\geq a}\int_{0}^{\infty}\la
\bm\partial_{x}\triangle^{-1}P_{2^{j}}[\delta_{0}(x)
T_{z}[|U|^{2}]\\0\endm,
s(\xi)e^{ix\xi}\underline{e}\ra\\&\hspace{2cm}\overline{\tilde{\calF}[\bm
T_{z}P_{\geq a}[\chi_{>0}|U|^{2}(s)]\la\underline{e},
 \chi_{>0}\bm\bar{U}(s)\\U(s)\endm\times e^{-i(t-s)\calH^{*}}\bm
\phi\\ \psi\endm_{dis}\ra\\ 0\endm](\xi)} d\xi
\end{split}\end{equation}
We have
\begin{equation}\nonumber
|\la\bm\partial_{x}\triangle^{-1}P_{2^{j}}[\delta_{0}(x)
T_{z}[|U|^{2}]\\0\endm, s(\xi)e^{ix\xi}\underline{e}\ra|\lesssim
\min\{2^{-j},1\}s^{-1}
\end{equation}
Moreover, we have
\begin{equation}\nonumber
||\tilde{\calF}[\bm T_{z}P_{\geq
a}[\chi_{>0}|U|^{2}(s)]\la\underline{e},
 \chi_{>0}\bm\bar{U}(s)\\U(s)\endm\times e^{-i(t-s)\calH^{*}}\bm
\phi\\ \psi\endm_{dis}\ra\\ 0\endm](\xi)||_{L_{\xi}^{2}}\lesssim
a^{-1}\la s\ra^{-\frac{3}{2}}\la
s\ra^{\frac{1}{2}+\epsilon(\delta_{2})}\la t-s\ra^{-\frac{3}{2}}
\end{equation}
Thus, using H\"{o}lder's inequality, we have
\begin{equation}\nonumber\begin{split}
&|\int_{0}^{\infty}\la\bm\partial_{x}\triangle^{-1}P_{2^{j}}[\delta_{0}(x)
T_{z}[|U|^{2}]\\0\endm,
s(\xi)e^{ix\xi}\underline{e}\ra\\&\hspace{2cm}\overline{\tilde{\calF}[\bm
T_{z}P_{\geq a}[\chi_{>0}|U|^{2}(s)]\la\underline{e},
 \chi_{>0}\bm\bar{U}(s)\\U(s)\endm\times e^{-i(t-s)\calH^{*}}\bm
\phi\\ \psi\endm_{dis}\ra\\ 0\endm](\xi)} d\xi|\\&\lesssim
\min\{2^{-\frac{j}{2}},1\}\la s\ra^{-1}a^{-1}\la
s\ra^{-\frac{3}{2}}\la s\ra^{\frac{1}{2}+\epsilon(\delta_{2})}\la
t-s\ra^{-\frac{3}{2}}
\end{split}\end{equation}
Summing over $j$ costs at most $\log s$, whence substituting
$a=\la s\ra^{-\frac{3}{4}}$ and integrating in $s$ yields the
upper bound $\lesssim \la t\ra^{-\frac{3}{2}+\delta_{3}}$, as
desired. If, on the other hand, we substitute
$\partial_{x}\triangle^{-1}P_{\geq a}[\chi_{>0}(x)
T_{z}\partial_{x}[|U|^{2}]$, we argue just as for the high-high
case. Now consider the expression
\begin{equation}\nonumber\begin{split}
&\int_{0}^{\infty}\la P_{<a}[\chi_{>0}(x)\bm
T_{z}[|U|^{2}]\\
0\endm],
s(\xi)e^{ix\xi}\underline{e}\ra\\&\hspace{2cm}\overline{\tilde{\calF}[\bm
T_{z}P_{\geq a}[\chi_{>0}|U|^{2}(s)]\la\underline{e},
 \chi_{>0}\bm\bar{U}(s)\\U(s)\endm\times e^{-i(t-s)\calH^{*}}\bm
\phi\\ \psi\endm_{dis}\ra\\ 0\endm](\xi)} d\xi
\end{split}\end{equation}
Estimate
\begin{equation}\nonumber
||\la P_{<a}[\chi_{>0}(x)\bm
T_{z}[|U|^{2}]\\
0\endm],s(\xi)e^{ix\xi}\underline{e}\ra||_{L_{\xi}^{2}}\lesssim a
\la s\ra^{-\frac{1}{2}}
\end{equation}
\begin{equation}\nonumber
||\calF[\bm T_{z}P_{\geq a}[\chi_{>0}|U|^{2}(s)]\la\underline{e},
 \chi_{>0}\bm\bar{U}(s)\\U(s)\endm\times e^{-i(t-s)\calH^{*}}\bm
\phi\\ \psi\endm_{dis}\ra\\ 0\endm](\xi)||_{L_{\xi}^{2}}\lesssim
a^{-1}\la s\ra^{-\frac{3}{2}}\la
s\ra^{\frac{1}{2}+\epsilon(\delta_{2})}\la t-s\ra^{-\frac{3}{2}}
\end{equation}
Putting these together results in the upper bound
\begin{equation}\nonumber\begin{split}
&|\int_{0}^{\infty}\la P_{<a}[\chi_{>0}(x)\bm
T_{z}[|U|^{2}]\\
0\endm],s(\xi)e^{ix\xi}\underline{e}\ra\\&\hspace{2cm}\overline{\tilde{\calF}[\bm
T_{z}P_{\geq a}[\chi_{>0}|U|^{2}(s)]\la\underline{e},
 \chi_{>0}\bm\bar{U}(s)\\U(s)\endm\times e^{-i(t-s)\calH^{*}}\bm
\phi\\ \psi\endm_{dis}\ra\\ 0\endm](\xi)} d\xi| \\& \lesssim a\la
s\ra^{-\frac{1}{2}}a^{-1}\la s\ra^{-\frac{3}{2}}\la
s\ra^{\frac{1}{2}+\epsilon(\delta_{2})}\la
t-s\ra^{-\frac{3}{2}}\lesssim \la
s\ra^{-\frac{3}{2}+\epsilon(\delta_{2})}\la t-s\ra^{-\frac{3}{2}},
\end{split}\end{equation}
which upon integration in $s$ again yields the desired upper bound
$\la t\ra ^{-\frac{3}{2}+\delta_{3}}$. The case when the
$\xi$-variable is restricted to $(-\infty,0]$ in the mixed case
will be treated further below.
\\

Now we consider
\begin{equation}\nonumber
\int_{-\infty}^{0}\int_{0}^{\frac{t}{2}}\calF\bm \chi_{>0}|U|^{4}(s)\\
-\chi_{>0}|U|^{4}(s)\endm(\xi)\overline{\tilde{\calF}[\chi_{>0}\bm\bar{U}(s)\\U(s)\endm\times
e^{-i(t-s)\calH^{*}}\bm \phi\\\psi\endm_{dis}]}(\xi)ds d\xi
\end{equation}
Reformulate this as
\begin{equation}\nonumber\begin{split}
&\int_{-\infty}^{0}\int_{0}^{\frac{t}{2}}\la\bm \chi_{>0}|U|^{4}(s)\\
-\chi_{>0}|U|^{4}(s)\endm,
(e^{ix\xi}-e^{-ix\xi})\underline{e}+(1+r(-\xi)e^{-ix\xi}\underline{e}+\phi(x,\xi)\ra\\&\hspace{6cm}\overline{\tilde{\calF}[\chi_{>0}\bm\bar{U}(s)\\U(s)\endm\times
e^{-i(t-s)\calH^{*}}\bm \phi\\\psi\endm_{dis}]}(\xi)ds d\xi
\end{split}\end{equation}
We can treat the contribution of $\phi(x,\xi)$ just as we did
before. Also, note that $|1+r(-\xi)|=O(|\xi|)$ around $\xi=0$, see
Theorem~\ref{DistortedFourier}, whence we can treat the
contribution of this part just like we did for the transmission
part before. The remaining part we break into a number of
contributions:
\begin{equation}\nonumber
\int_{-\infty}^{0}\int_{0}^{\frac{t}{2}}\la \bm P_{\geq
a}[\chi_{>0}|U|^{2}(s)]P_{\geq a}[|U|^{2}(s)]\\-P_{\geq
a}[\chi_{>0}|U|^{2}(s)]P_{\geq a}[|U|^{2}(s)]\endm,
(e^{ix\xi}-e^{-ix\xi})\underline{e}\ra
\overline{\tilde{\calF}[\chi_{>0}\bm\bar{U}(s)\\U(s)\endm\times
e^{-i(t-s)\calH^{*}}\bm \phi\\ \psi\endm]}(\xi)ds d\xi
\end{equation}
\begin{equation}\nonumber
\int_{-\infty}^{0}\int_{0}^{\frac{t}{2}}\la \bm
P_{<a}[\chi_{>0}|U|^{2}(s)]P_{\geq
a}[|U|^{2}(s)]\\-P_{<a}[\chi_{>0}|U|^{2}(s)]P_{\geq
a}[|U|^{2}(s)]\endm,(e^{ix\xi}-e^{-ix\xi})\underline{e}\ra
\overline{\tilde{\calF}[\chi_{>0}\bm\bar{U}(s)\\U(s)\endm\times
e^{-i(t-s)\calH^{*}}\bm \phi\\ \psi\endm]}(\xi)ds d\xi
\end{equation}
\begin{equation}\nonumber
\int_{-\infty}^{0}\int_{0}^{\frac{t}{2}}\la \bm P_{\geq
a}[\chi_{>0}|U|^{2}(s)]P_{<a}[|U|^{2}(s)]\\-P_{\geq
a}[\chi_{>0}|U|^{2}(s)]P_{<a}[|U|^{2}(s)]\endm,
(e^{ix\xi}-e^{-ix\xi})\underline{e}\ra
\overline{\tilde{\calF}[\chi_{>0}\bm\bar{U}(s)\\U(s)\endm\times
e^{-i(t-s)\calH^{*}}\bm \phi\\ \psi\endm]}(\xi)ds d\xi
\end{equation}
\begin{equation}\nonumber
\int_{-\infty}^{0}\int_{0}^{\frac{t}{2}}\la \bm P_{<
a}[\chi_{>0}|U|^{2}(s)]P_{<a}[|U|^{2}(s)]\\-P_{<
a}[\chi_{>0}|U|^{2}(s)]P_{<a}[|U|^{2}(s)]\endm,
(e^{ix\xi}-e^{-ix\xi})\underline{e}\ra
\overline{\tilde{\calF}[\chi_{>0}\bm\bar{U}(s)\\U(s)\endm\times
e^{-i(t-s)\calH^{*}}\bm \phi\\ \psi\endm]}(\xi)ds d\xi
\end{equation}
Start with the first term in this list: write
\begin{equation}\nonumber
\tilde{\calF}[\chi_{>0}\bm\bar{U}(s)\\U(s)\endm\times
e^{-i(t-s)\calH^{*}}\bm \phi\\
\psi\endm]=\la\chi_{>0}\bm\bar{U}(s)\\U(s)\endm\times
e^{-i(t-s)\calH^{*}}\bm \phi\\ \psi\endm,
[e^{ix\xi}+r(-\xi)e^{-ix\xi}]\underline{e}+\phi(x,\xi)\ra
\end{equation}
The contribution of $\phi(x,\xi)$ here is again straightforward,
and left out. Now one proceeds as for the transmission part
($\xi\geq 0$) treated before, using the ordinary
Plancherel's Theorem and introducing a multiplier $\Pi_{(t^{-1000},t^{1000})}$. \\
Next, we consider the low-low frequency interaction, i. e. the
expression
\begin{equation}\nonumber
\int_{-\infty}^{0}\int_{0}^{t}\la\bm P_{<
a}[\chi_{>0}|U|^{2}(s)]P_{<a}[|U|^{2}(s)]\\-P_{<
a}[\chi_{>0}|U|^{2}(s)]P_{<a}[|U|^{2}(s)]\endm,
(e^{ix\xi}-e^{-ix\xi})\underline{e}\ra
\overline{\tilde{\calF}[\chi_{>0}\bm\bar{U}(s)\\U(s)\endm\times
e^{-i(t-s)\calH^{*}}\bm \phi\\ \psi\endm]}(\xi)ds d\xi
\end{equation}
This calls for a different strategy than for the transmission
part, since the Fourier basis in this regime doesn't vanish
uniformly at $\xi=0$. First, we observe that
\begin{equation}\nonumber\begin{split}
&\la\bm P_{< a}[\chi_{>0}|U|^{2}(s)]P_{<a}[|U|^{2}(s)]\\-P_{<
a}[\chi_{>0}|U|^{2}(s)]P_{<a}[|U|^{2}(s)]\endm,
(e^{ix\xi}-e^{-ix\xi})\underline{e}\ra\\&=\chi_{<a+O(1)}(\xi)\la
\bm P_{< a}[\chi_{>0}|U|^{2}(s)]P_{<a}[|U|^{2}(s)]\\-P_{<
a}[\chi_{>0}|U|^{2}(s)]P_{<a}[|U|^{2}(s)]\endm,
(e^{ix\xi}-e^{-ix\xi})\underline{e}\ra
\end{split}\end{equation}
Hence we have
\begin{equation}\nonumber
||\la \bm P_{< a}[\chi_{>0}|U|^{2}(s)]P_{<a}[|U|^{2}(s)]\\-P_{<
a}[\chi_{>0}|U|^{2}(s)]P_{<a}[|U|^{2}(s)]\endm,
(e^{ix\xi}-e^{-ix\xi})\underline{e}\ra||_{L_{\xi}^{1}}\lesssim a
\la s\ra^{-1}
\end{equation}
Notice that putting $a=\la s\ra^{-\frac{3}{4}}$ is not quite good
enough yet to counterbalance the loss of $s$ arising when one
extracts the $(t-s)^{-\frac{3}{2}}$-gain. This extra gain of
$s^{-\frac{1}{4}}$ has to come from the 2nd factor
$\tilde{\calF}(...)$. Write for $\xi<0$
\begin{equation}\nonumber\begin{split}
&\tilde{\calF}[\chi_{>0}\bm \bar{U}\\U\endm\times
e^{-i(t-s)\calH^{*}}\bm\phi\\
\psi\endm_{dis}](\xi)\\&=\la\chi_{>0}\bm
\bar{U}\\U\endm\times e^{-i(t-s)\calH^{*}}\bm\phi\\
\psi\endm_{dis},
[e^{ix\xi}-e^{-ix\xi}]\underline{e}+(1+r(-\xi))e^{-ix\xi}\underline{e}+\phi(x,\xi)\ra
\end{split}\end{equation}
We first get rid of
$(1+r(-\xi))e^{-ix\xi}\underline{e}+\phi(x,\xi)$. Note that for
$\xi\lesssim a$, we have
\begin{equation}\nonumber
||\la \chi_{>0}\bm
\bar{U}\\U\endm\times e^{-i(t-s)\calH^{*}}\bm\phi\\
\psi\endm_{dis},
(1+r(-\xi))e^{-ix\xi}\underline{e}\ra||_{L_{\xi}^{2}}\lesssim a
\la s\ra^{1+\epsilon(\delta_{2})}\la t-s\ra^{-\frac{3}{2}}
\end{equation}
Combining this with
\begin{equation}\nonumber
||\la\bm P_{< a}[\chi_{>0}|U|^{2}(s)]P_{<a}[|U|^{2}(s)]\\-P_{<
a}[\chi_{>0}|U|^{2}(s)]P_{<a}[|U|^{2}(s)]\endm,
(e^{ix\xi}-e^{-ix\xi})\underline{e}\ra||_{L_{\xi}^{2}}\lesssim \la
s\ra^{-\frac{3}{2}}
\end{equation}
easily leads to the upper bound $\lesssim \la
s\ra^{-\frac{5}{4+}}\la t-s\ra^{-\frac{3}{2}}$ for this
contribution. Next, we have
\begin{equation}\nonumber
||\la\chi_{>0}\bm
\bar{U}\\U\endm\times e^{-i(t-s)\calH^{*}}\bm\phi\\
\psi\endm_{dis},  \phi(x,\xi)\ra||_{L_{\xi}^{\infty}}\lesssim \la
s\ra^{-\frac{3}{2}+\delta_{3}}\la t-s\ra^{-\frac{3}{2}},
\end{equation}
which similarly leads to an acceptable upper bound. We now reduce
to estimating the expression
\begin{equation}\nonumber\begin{split}
&\int_{-\infty}^{\infty}[e^{ix\xi}-e^{-ix\xi}]\la\underline{e},
\chi_{>0}\bm
\bar{U}\\U\endm\times e^{-i(t-s)\calH^{*}}\bm\phi\\
\psi\endm_{dis}\ra dx\\&=i\xi
\int_{0}^{\infty}[e^{ix\xi}+e^{-ix\xi}]\la\underline{e},
\int_{x}^{\infty}\chi_{>0}\bm
\bar{U}(s)\\U(s)\endm\times e^{-i(t-s)\calH^{*}}\bm\phi\\
\psi\endm_{dis}(y)dy\ra dx
\end{split}\end{equation}
In order to analyze the inner integral here, it appears useful to
express $U$ etc as Fresnel integrals, which makes the spatial
oscillations visible. First, using the distorted Fourier
transform, we write
\begin{equation}\nonumber
\chi_{>0}e^{-i(t-s)\calH^{*}}\bm\phi\\
\psi\endm_{dis}=\sum_{\pm}\int_{-\infty}^{\infty}e^{\pm i(t-s)(\xi^{2}+1)}\chi_{>0}\sigma_{3}e_{\pm}(x,\xi)\tilde{\calF}_{\pm}(\bm\phi\\
\psi\endm(\xi)d\xi
\end{equation}
Fix the $+$-sign here, the $-$-sign being treated accordingly; it
is important here that the oscillatory part of $e_{-}(x,\xi)$ only
has a lower component, i. e.
$e_{-}(x,\xi)=e^{ix\xi}\sigma_{1}\underline{e}+\ldots$, where
$\sigma_{1}=\bm 0&1\\1&0\endm$, $\underline{e}=\bm 0\\1\endm$,
while the oscillatory part of $e_{+}(x,\xi)$ only has an upper
component. Then we break the integral into two contributions:
\begin{equation}\label{01}
\int_{0}^{\infty}e^{i(t-s)(\xi^{2}+1)}\chi_{>0}\sigma_{3}e(x,\xi)\tilde{\calF}(\bm\phi\\
\psi\endm(\xi)d\xi
\end{equation}
\begin{equation}\label{02}
\int_{-\infty}^{0}e^{i(t-s)(\xi^{2}+1)}\chi_{>0}\sigma_{3}e(x,\xi)\tilde{\calF}(\bm\phi\\
\psi\endm(\xi)d\xi
\end{equation}
Write the first integral \eqref{01} as
\begin{equation}\nonumber
\int_{0}^{\infty}e^{i(t-s)(\xi^{2}+1)}\chi_{>0}\sigma_{3}[s(\xi)e^{ix\xi}\underline{e}+\phi(x,\xi)]\tilde{\calF}(\bm\phi\\
\psi\endm(\xi)d\xi
\end{equation}
Carry out an integration by parts in $\xi$, thereby replacing this
by
\begin{equation}\nonumber
\frac{1}{t-s}\int_{0}^{\infty}e^{i(t-s)(\xi^{2}+1)}\chi_{>0}(x)\sigma_{3}\partial_{\xi}([s(\xi)e^{ix\xi}\underline{e}+\phi(x,\xi)]\frac{\tilde{\calF}(\bm\phi\\
\psi\endm(\xi)}{\xi})d\xi
\end{equation}
The contribution of $\phi(x,\xi)$ here is again negligible, as is
easily seen. The worst case occurs when the derivative
$\partial_{\xi}$ falls on the phase $e^{ix\xi}$, costing a factor
$ix$. Explicitly, this is the following expression:
\begin{equation}\nonumber
\frac{ix}{t-s}\int_{0}^{\infty}e^{i(t-s)(\xi^{2}+1)}\chi_{>0}(x)\sigma_{3}s(\xi)e^{ix\xi}\underline{e}\frac{\tilde{\calF}(\bm\phi\\
\psi\endm(\xi)}{\xi}d\xi
\end{equation}
Break the $\xi$-integral into two, one over the interval
$[0,t^{1000}]$, the other over its complement on $[0,\infty)$. On
the latter, an additional integration by parts in $\xi$ easily
furnishes more than the needed gain in $t$. On the former
interval, observe that we may interpret the integral
\begin{equation}\nonumber
\int_{0}^{\infty}e^{i(t-s)\xi^{2}}\sigma_{3}[s(\xi)e^{ix\xi}\underline{e}\frac{\tilde{\calF}(\bm\phi\\
\psi\endm(\xi)}{\xi}\chi_{<t^{1000}}(\xi)d\xi
\end{equation}
as a solution for the free Schroedinger equation, evaluated at
time $t-s$, with initial data
\begin{equation}\nonumber
g(x)=\int_{0}^{\infty}\sigma_{3}s(\xi)e^{ix\xi}\underline{e}\frac{\tilde{\calF}(\bm\phi\\
\psi\endm(\xi)}{\xi}\chi_{<t^{1000}}(\xi)d\xi
\end{equation}
The definition of $\tilde{\calF}$ as well as further integrations
by parts in $\xi$ reveal that this decays like $x^{-2}$ for large
values of $x$, resulting in
\begin{equation}\nonumber
||\int_{0}^{\infty}\sigma_{3}s(\xi)e^{ix\xi}\underline{e}\frac{\tilde{\calF}(\bm\phi\\
\psi\endm(\xi)}{\xi}\chi_{<t^{1000}}(\xi)d\xi||_{L_{x}^{1}}\lesssim
\log t
\end{equation}
Thus we can now write
\begin{equation}\label{03}
\frac{ix}{t-s}\int_{0}^{\infty}e^{i(t-s)(\xi^{2}+1)}\chi_{<t^{1000}}(\xi)\sigma_{3}s(\xi)e^{ix\xi}\underline{e}\frac{\tilde{\calF}(\bm\phi\\
\psi\endm(\xi)}{\xi}d\xi
=\frac{ie^{i(t-s)}x}{t-s}\frac{1}{\sqrt{t-s}}\int_{-\infty}^{\infty}e^{-\frac{(x-y)^{2}}{i(t-s)}}g(y)dy
\end{equation}
Next, returning to \eqref{02}, we consider the integral
\begin{equation}\nonumber
\int_{-\infty}^{0}e^{i(t-s)(\xi^{2}+1)}\chi_{>0}\sigma_{3}e(x,\xi)\tilde{\calF}(\bm\phi\\
\psi\endm(\xi)d\xi
\end{equation}
In the regime under consideration we can write
$e(x,\xi)=[e^{ix\xi}-e^{-ix\xi}+(1+r(-\xi))e^{-ix\xi}]\underline{e}+\phi(x,\xi)$.
We proceed as before, arriving (up to error terms handled as
before) at the expression
\begin{equation}\label{04}\begin{split}
&\frac{ix}{t-s}\int_{-\infty}^{0}e^{i(t-s)(\xi^{2}+1)}\chi_{>0}(x)[e^{ix\xi}+e^{-ix\xi}]\underline{e}\frac{\tilde{\calF}\bm\phi\\\psi\endm(\xi)}{\xi}d\xi
\\&=\frac{ix}{t-s}\frac{1}{\sqrt{t-s}}\int_{-\infty}^{\infty}e^{-\frac{(x-y)^{2}}{i(t-s)}}\tilde{g}(y)dy,
\end{split}\end{equation}
where
\begin{equation}\nonumber
\tilde{g}(y)=\int_{-\infty}^{0}\chi_{>0}(x)[e^{ix\xi}+e^{-ix\xi}]\underline{e}\frac{\tilde{\calF}\bm\phi\\\psi\endm(\xi)}{\xi}d\xi
\end{equation}
Now substitute either \eqref{03} or \eqref{04} for the right hand factor in $\bm \bar{U}\\
U\endm\times e^{-i(t-s)\calH^{*}}\bm\phi\\ \psi\endm_{dis}$. We
replace the resulting $x\bar{U}$ by $s\partial_{x}\bar{U}$, upon
using control over $||CU||_{L_{x}^{2}}$. Thus, if we substitute
for example\footnote{The contribution of \eqref{03} is handled
similarly.} \eqref{04}, we need to estimate
\begin{equation}\nonumber
\chi_{>0}(x)s\partial_{x}\bar{U}\frac{e^{i(t-s)}}{t-s}\int_{-\infty}^{0}e^{i(t-s)\xi^{2}}[e^{ix\xi}+e^{-ix\xi}]\underline{e}\frac{\tilde{\calF}\bm\phi\\\psi\endm(\xi)}{\xi}d\xi
\end{equation}
We now schematically record the equation satisfied by
$\partial_{x}\bar{U}$ as follows:
\begin{equation}\nonumber
(-i\partial_{t}+\triangle)\partial_{x}\bar{U}= VU
+V\bar{U}+V\partial_{x}U+V\partial_{x}\bar{U}+\ldots+\partial_{x}[|U|^{4}\bar{U}]
\end{equation}
Here $V$ denotes certain Schwartz functions whose fine structure
is irrelevant. Thus we can write
\begin{equation}\nonumber
\partial_{x}\bar{U}(s,x)=\int_{-\infty}^{\infty}\int_{0}^{s}\frac{1}{\sqrt{s-\lambda}}e^{-\frac{(x-y)^{2}}{i(s-\lambda)}}[
VU(\lambda,y)+ \partial_{x}[|U|^{4}\bar{U}(\lambda,.)]] d\lambda
dy+\ldots
\end{equation}
We claim that we may replace the local terms $VU(\lambda,y)$ by\\
$\chi_{<s-s^{2\epsilon}}(\lambda)\chi_{<s^{\epsilon}}(y)VU(\lambda,y)$
for small $\epsilon>0$ (independent of $\delta_{i}$ etc), and the
non-local term $\partial_{x}[|U|^{4}\bar{U}(\lambda,.)]$ by
$\chi_{<s^{\frac{1}{2}}}(\lambda)\chi_{<s^{\frac{1}{2}}}(y)\partial_{x}[|U|^{4}\bar{U}(\lambda,.)]$.
To see this, note that
\begin{equation}\nonumber
||\int_{-\infty}^{\infty}\int_{0}^{s}\frac{1}{\sqrt{s-\lambda}}e^{-\frac{(x-y)^{2}}{i(s-\lambda)}}\chi_{\geq
s^{\epsilon}}(y)VU(\lambda,y)d\lambda dy||_{L_{x}^{2}}\lesssim
e^{-s^{\epsilon}}
\end{equation}
\begin{equation}\nonumber
||\int_{-\infty}^{\infty}\int_{0}^{s}\frac{1}{\sqrt{s-\lambda}}e^{-\frac{(x-y)^{2}}{i(s-\lambda)}}\chi_{<
s^{\epsilon}}(y)\chi_{\geq
s-s^{2\epsilon}}(\lambda)VU(\lambda,y)d\lambda
dy||_{L_{x}^{2}}\lesssim \la s\ra^{-\frac{3}{2}}\la
s\ra^{\epsilon+\delta_{3}}
\end{equation}
If one substitutes the corresponding terms in the Duhamel formula
for $U$, $\bar{U}$ directly for the fifth factors $U$, $\bar{U}$
in
\begin{equation}\nonumber
\int_{0}^{\frac{t}{2}}\la e^{i(t-s)\calH}\bm |U|^{4}U(s,.)\\
-|U|^{4}\bar{U}(s,.)\endm_{dis}, \phi\ra ds,
\end{equation}
one easily bounds this contribution by $\lesssim \la
t\ra^{-\frac{3}{2}}$. Similarly, we have
\begin{equation}\nonumber
||\int_{-\infty}^{\infty}\int_{0}^{s}\frac{1}{\sqrt{s-\lambda}}e^{-\frac{(x-y)^{2}}{i(s-\lambda)}}\chi_{>s^{\frac{1}{2}}}(\lambda)
|U|^{4}(\lambda,y)U(\lambda,y)d\lambda dy||_{L_{x}^{2}}\lesssim
\la s\ra^{-\frac{1}{2}},
\end{equation}
which leads to a similar conclusion upon substituting this
integral for the last factors $U$, $\bar{U}$. Finally, note that
on account of the pseudo-conformal conservation law, we have for
$x>\lambda$
\begin{equation}\nonumber
||U(\lambda,.)||_{L_{x}^{2}}\sim ||\frac{\lambda}{x}\nabla
U(\lambda,.)||_{L_{x}^{2}},
\end{equation}
whence we can estimate
\begin{equation}\nonumber
||\int_{-\infty}^{\infty}\int_{0}^{s^{\frac{1}{2}}}\frac{1}{\sqrt{s-\lambda}}e^{-\frac{(x-y)^{2}}{i(s-\lambda)}}|U(\lambda,y)|^{4}\chi_{>s^{\frac{1}{2}}}(y)U(y,\lambda)dy||_{L_{x}^{2}}
\lesssim
s^{-\frac{1}{2}}\int_{0}^{s^{\frac{1}{2}}}\lambda^{-2}\lambda
d\lambda\lesssim \la s\ra^{-\frac{1}{2}}\log \la s\ra,
\end{equation}
and the argument proceeds from here as before. This discussion
justifies us in substituting
\begin{equation}\nonumber
\partial_{x}\bar{U}(s,x)=\int_{-\infty}^{\infty}\int_{0}^{s}\frac{1}{\sqrt{s-\lambda}}e^{-\frac{(x-y)^{2}}{i(s-\lambda)}}[
\chi_{<s-s^{2\epsilon}}\chi_{<s^{\epsilon}}(y)VU(\lambda,y)+
\chi_{<s^{\frac{1}{2}}}(y)\chi_{<s^{\frac{1}{2}}}(\lambda)\partial_{y}[|U|^{4}\bar{U}(\lambda,y)]]
d\lambda dy+\ldots
\end{equation}
Next, write as before
\begin{equation}\nonumber
\int_{-\infty}^{0}e^{i(t-s)\xi^{2}}\chi_{<t^{1000}}(\xi)[e^{ix\xi}+e^{-ix\xi}]\underline{e}\frac{\tilde{\calF}\bm\phi\\\psi\endm(\xi)}{\xi}d\xi
=\frac{1}{\sqrt{t-s}}\int_{-\infty}^{\infty}e^{-\frac{(x-y')^{2}}{i(t-s)}}\tilde{g}(y')dy'
\end{equation}
where
\begin{equation}\nonumber
\tilde{g}(y')=\int_{-\infty}^{0}\chi_{<t^{1000}}(\xi)[e^{iy'\xi}+e^{-iy'\xi}]\underline{e}\frac{\tilde{\calF}\bm\phi\\\psi\endm(\xi)}{\xi}d\xi
\end{equation}
We observe that we may include the cutoff
$\chi_{<(t-s)^{\frac{1}{2}}}(y')$ in front of $\tilde{g}(y')$;
this is on account of the estimate $|\tilde{g}(y')|\lesssim
y'^{-2}\log t$. Finally, plugging these expressions into
\begin{equation}\nonumber
\chi_{>0}(x)s\partial_{x}\bar{U}\frac{e^{i(t-s)}}{t-s}\int_{-\infty}^{0}e^{i(t-s)\xi^{2}}[e^{ix\xi}+e^{-ix\xi}]\underline{e}\frac{\tilde{\calF}\bm\phi\\\psi\endm(\xi)}{\xi}d\xi,
\end{equation}
and keeping in mind that our point of departure was the expression
\begin{equation}\nonumber
\int_{x}^{\infty}\chi_{>0}\bm
\bar{U}(s)\\U(s)\endm\times e^{-i(t-s)\calH^{*}}\bm\phi\\
\psi\endm_{dis}(y)dy,
\end{equation}
we arrive at terms of the following form:
\begin{equation}\label{A}\begin{split}
&\frac{s}{(t-s)^{\frac{3}{2}}}\int_{x_{0}}^{\infty}\int_{-\infty}^{\infty}\int_{-\infty}^{\infty}\chi_{>0}(x)
\int_{0}^{s}\frac{1}{\sqrt{s-\lambda}}e^{-\frac{(x-y)^{2}}{i(s-\lambda)}}\chi_{<s^{\epsilon}}(y)\chi_{<s-s^{2\epsilon}}(\lambda)VU(\lambda,y)e^{-\frac{(x-y')^{2}}{i(t-s)}}\\&\hspace{8cm}\chi_{<(t-s)^{\frac{1}{2}}}(y')\tilde{g}(y')d\lambda
dy dy' dx,
\end{split}\end{equation}
\begin{equation}\label{B}\begin{split}
&\frac{s}{(t-s)^{\frac{3}{2}}}\int_{x_{0}}^{\infty}\int_{-\infty}^{\infty}\int_{-\infty}^{\infty}\chi_{>0}(x)
\int_{0}^{s}\frac{1}{\sqrt{s-\lambda}}e^{-\frac{(x-y)^{2}}{i(s-\lambda)}}|U(\lambda,y)|^{4}\chi_{<s^{\frac{1}{2}}}(y)\chi_{<s^{\frac{1}{2}}}(\lambda)U(y,\lambda)e^{-\frac{(x-y')^{2}}{i(t-s)}}\\&\hspace{8cm}\chi_{<(t-s)^{\frac{1}{2}}}(y')\tilde{g}(y')d\lambda
dy dy' dx,
\end{split}\end{equation}
where $x_{0}$ ranges over $[0,\infty]$, plus similar terms which
can be treated identically. Write
\begin{equation}\nonumber
e^{-\frac{(x-y)^{2}}{i(s-\lambda)}}e^{-\frac{(x-y')^{2}}{i(t-s)}}=e^{+i[(\frac{1}{s-\lambda}+\frac{1}{t-s})x^{2}-\frac{2xy}{s-\lambda}-\frac{2xy'}{t-s}]}
e^{+\frac{iy^{2}}{s-\lambda}+\frac{iy'^{2}}{t-s}}
\end{equation}
This can be rewritten as
\begin{equation}\nonumber
e^{i(x\sqrt{\frac{1}{s-\lambda}+\frac{1}{t-s}}-y_{1})^{2}}e^{iy_{2}}
\end{equation}
for certain functions $y_{1,2}(y,y',s,\lambda,t)$. Our
restrictions in either term \eqref{A} or \eqref{B} ensure that
$y_{1}=O(1)$. Carrying out the $x$-integration, we obtain
\begin{equation}\nonumber
\frac{1}{\sqrt{s-\lambda}}\int_{x_{0}}^{\infty}e^{i(x\sqrt{\frac{1}{s-\lambda}+\frac{1}{t-s}}-y_{1})^{2}}dx
=\frac{1}{\sqrt{s-\lambda}}(\frac{1}{s-\lambda}+\frac{1}{t-s})^{-\frac{1}{2}}S(x_{0}\sqrt{\frac{1}{s-\lambda}+\frac{1}{t-s}}-y_{1})
\end{equation}
where
$S(y)=\int_{y}^{\infty}e^{ix^{2}}dx=\frac{e^{iy^{2}}}{y}+O(y^{-2})$
as $y\rightarrow +\infty$. Finally, we need to estimate
\begin{equation}\nonumber
\xi\int_{0}^{\infty}[e^{ix_{0}\xi}+e^{-ix_{0}\xi}]\frac{1}{\sqrt{s-\lambda}}(\frac{1}{s-\lambda}+\frac{1}{t-s})^{-\frac{1}{2}}S(x_{0}\sqrt{\frac{1}{s-\lambda}+\frac{1}{t-s}}-y_{1})d
x_{0}
\end{equation}
Here it is important that $y_{1}$ be uniformly bounded. The
oscillatory nature of $S(y)$ allows us to bound this integral by
\begin{equation}\nonumber
\lesssim |\xi| \sqrt{\la s\ra}\lesssim \la s\ra^{-\frac{1}{4}},
\end{equation}
as desired. The remaining integrations over $y,y'$ are
straightforward to carry out on account of the integrability of
the functions $g(y)$, $g(y')$. One pays $(\log t)^{2}$, which is
irrelevant. This finally completes treating the low-low case.
\\

We proceed to the mixed frequency case. This is the expression
\begin{equation}\nonumber
\int_{-\infty}^{0}\int_{0}^{\frac{t}{2}}\la \bm P_{\geq
a}[\chi_{>0}|U|^{2}(s)]P_{<a}[|U|^{2}(s)]\\-P_{\geq
a}[\chi_{>0}|U|^{2}(s)]P_{<a}[|U|^{2}(s)]\endm,
(e^{ix\xi}-e^{-ix\xi})\underline{e}\ra
\overline{\tilde{\calF}[\chi_{>0}\bm\bar{U}(s)\\U(s)\endm\times
e^{-i(t-s)\calH^{*}}\bm \phi\\ \psi\endm_{dis}]}(\xi)ds d\xi
\end{equation}
We first employ the ordinary Plancherel's Theorem to replace this
(up to negligible error terms) by an expression
\begin{equation}\nonumber
\la T_{z}P_{<a}[|U|^{2}(s)], T_{z}P_{\geq
a}[\chi_{>0}|U|^{2}(s)]\la \underline{e},
\chi_{>0}\bm\bar{U}(s)\\U(s)\endm\times e^{-i(t-s)\calH^{*}}\bm
\phi\\ \psi\endm_{dis}\ra\ra,
\end{equation}
where we have to integrate over $z$ in the end which will cost
$\log t$, as before. Express this in vectorial form as
\begin{equation}\nonumber
\la \bm T_{z}P_{<a}[|U|^{2}(s)]\\0\endm, \bm T_{z}P_{\geq
a}[\chi_{>0}|U|^{2}(s)]\la \underline{e},
\chi_{>0}\bm\bar{U}(s)\\U(s)\endm\times e^{-i(t-s)\calH^{*}}\bm
\phi\\ \psi\endm\ra\\0\endm\ra,
\end{equation}
Decompose this into the following two terms:
\begin{equation}\nonumber
\la \chi_{>0}(x)\bm T_{z}P_{<a}[|U|^{2}(s)]\\0\endm, \bm
T_{z}P_{\geq a}[\chi_{>0}|U|^{2}(s)]\la \underline{e},
\chi_{>0}\bm\bar{U}(s)\\U(s)\endm\times e^{-i(t-s)\calH^{*}}\bm
\phi\\ \psi\endm_{dis}\ra\\0\endm\ra,
\end{equation}
\begin{equation}\nonumber
\la \chi_{<0}(x)\bm T_{z}P_{<a}[|U|^{2}(s)]\\0\endm, \bm
T_{z}P_{\geq a}[\chi_{>0}|U|^{2}(s)]\la \underline{e},
\chi_{>0}\bm\bar{U}(s)\\U(s)\endm\times e^{-i(t-s)\calH^{*}}\bm
\phi\\ \psi\endm_{dis}\ra\\0\endm\ra,
\end{equation}
These being treated similarly, we treat the first term: commence
by replacing $\chi_{>0}(x)\bm T_{z}P_{<a}[|U|^{2}(s)]\\0\endm$ by
its dispersive part. This is done as in the mixed frequency case
treated earlier. Then  use the distorted Plancherel's
Theorem~\ref{DistortedPlancherel}, which produces
\begin{equation}\nonumber\begin{split}
&\int_{-\infty}^{\infty}\calF[\chi_{>0}(x)\bm
T_{z}P_{<a}[|U|^{2}(s)]\\0\endm](\xi)\\&\hspace{2cm}\overline{\tilde{\calF}[\bm
T_{z}P_{\geq a}[\chi_{>0}|U|^{2}(s)]\la \underline{e},
\chi_{>0}\bm\bar{U}(s)\\U(s)\endm\times e^{-i(t-s)\calH^{*}}\bm
\phi\\ \psi\endm_{dis}\ra\\0\endm]}(\xi)d\xi
\end{split}\end{equation}
Divide this into the integral over $(-\infty,0]$ as well as the
integral over $[0,\infty)$. We treat the more difficult former
case, the latter already having been dealt with in the preceding.
We recast this as
\begin{equation}\nonumber\begin{split}
&\int_{-\infty}^{0}\la\chi_{>0}(x)\bm
T_{z}P_{<a}[|U|^{2}(s)]\\0\endm,
[e^{ix\xi}-e^{-ix\xi}+(1+r(-\xi))e^{-ix\xi}]\underline{e}+\phi(x,\xi)\ra\\&\hspace{2cm}\overline{\tilde{\calF}[\bm
T_{z}P_{\geq a}[\chi_{>0}|U|^{2}(s)]\la \underline{e},
\chi_{>0}\bm\bar{U}(s)\\U(s)\endm\times e^{-i(t-s)\calH^{*}}\bm
\phi\\ \psi\endm_{dis}\ra\\0\endm]}(\xi)d\xi
\end{split}\end{equation}
The contributions of $1+r(-\xi)$ and $\phi(x,\xi)$ are
straightforward, and handled as in the preceding. We then need to
estimate the following two contributions:
\begin{equation}\nonumber\begin{split}
&\int_{-\infty}^{0}\la P_{\geq a}[\chi_{>0}(x)\bm
T_{z}P_{<a}[|U|^{2}(s)]\\0\endm],
[e^{ix\xi}-e^{-ix\xi}]\underline{e}\ra\\&\hspace{2cm}\overline{\tilde{\calF}[\bm
T_{z}P_{\geq a}[\chi_{>0}|U|^{2}(s)]\la \underline{e},
\chi_{>0}\bm\bar{U}(s)\\U(s)\endm\times e^{-i(t-s)\calH^{*}}\bm
\phi\\ \psi\endm_{dis}\ra\\0\endm]}(\xi)d\xi
\end{split}\end{equation}
\begin{equation}\nonumber\begin{split}
&\int_{-\infty}^{0}\la P_{<a}[\chi_{>0}(x)\bm
T_{z}P_{<a}[|U|^{2}(s)]\\0\endm],[e^{ix\xi}-e^{-ix\xi}]\underline{e}\ra\\&\hspace{2cm}\overline{\tilde{\calF}[\bm
T_{z}P_{\geq a}[\chi_{>0}|U|^{2}(s)]\la \underline{e},
\chi_{>0}\bm\bar{U}(s)\\U(s)\endm\times e^{-i(t-s)\calH^{*}}\bm
\phi\\ \psi\endm_{dis}\ra\\0\endm]}(\xi)d\xi
\end{split}\end{equation}
Consider the first of these terms: it is straightforward to
replace $P_{\geq a}[\chi_{>0}(x)\bm
T_{z}P_{<a}[|U|^{2}(s)]\\0\endm]$ by \\ $P_{\geq
a}[\chi_{>0}(x)\bm T_{z}[|U|^{2}(s)]\\0\endm]$, by arguing as for
the high-high case. Then we replace this term by
\begin{equation}\nonumber\begin{split}
&\int_{-\infty}^{0}\la\partial_{x}\triangle^{-1}P_{\geq
a}\partial_{x}[\chi_{>0}(x)\bm T_{z}P_{<a}[|U|^{2}(s)]\\0\endm],
[e^{ix\xi}-e^{-ix\xi}]\underline{e}\ra\\&\hspace{2cm}\overline{\tilde{\calF}[\bm
T_{z}P_{\geq a}[\chi_{>0}|U|^{2}(s)]\la \underline{e},
\chi_{>0}\bm\bar{U}(s)\\U(s)\endm\times e^{-i(t-s)\calH^{*}}\bm
\phi\\ \psi\endm_{dis}\ra\\0\endm]}(\xi)d\xi
\end{split}\end{equation}
First let the inner derivative $\partial_{x}$ fall onto the factor
$\chi_{>0}(x)$. This results in
\begin{equation}\nonumber\begin{split}
&\int_{-\infty}^{0}\la\partial_{x}\triangle^{-1}P_{\geq
a}[\delta_{0}(x)\bm T_{z}P_{<a}[|U|^{2}(s)]\\0\endm],
[e^{ix\xi}-e^{-ix\xi}]\underline{e}\ra\\&\hspace{2cm}\overline{\tilde{\calF}[\bm
T_{z}P_{\geq a}[\chi_{>0}|U|^{2}(s)]\la \underline{e},
\chi_{>0}\bm\bar{U}(s)\\U(s)\endm\times e^{-i(t-s)\calH^{*}}\bm
\phi\\ \psi\endm_{dis}\ra\\0\endm]}(\xi)d\xi
\end{split}\end{equation}
In order to estimate this, we decompose it further into two
contributions:
\begin{equation}\nonumber\begin{split}
&\int_{-\infty}^{0}\la
\partial_{x}\triangle^{-1}P_{\la s\ra^{-\frac{1}{2}}\geq
.\geq a}[\delta_{0}(x)\bm T_{z}P_{<a}[|U|^{2}(s)]\\0\endm],
[e^{ix\xi}-e^{-ix\xi}]\underline{e}\ra\\&\hspace{2cm}\overline{\tilde{\calF}[\bm
T_{z}P_{\geq a}[\chi_{>0}|U|^{2}(s)]\la \underline{e},
\chi_{>0}\bm\bar{U}(s)\\U(s)\endm\times e^{-i(t-s)\calH^{*}}\bm
\phi\\ \psi\endm_{dis}\ra\\0\endm]}(\xi)d\xi
\end{split}\end{equation}
\begin{equation}\nonumber\begin{split}
&\int_{-\infty}^{0}\la
\partial_{x}\triangle^{-1}P_{>\la s\ra^{-\frac{1}{2}}}[\delta_{0}(x)\bm
T_{z}P_{<a}[|U|^{2}(s)]\\0\endm],
[e^{ix\xi}-e^{-ix\xi}]\underline{e}\ra\\&\hspace{2cm}\overline{\tilde{\calF}[\bm
T_{z}P_{\geq a}[\chi_{>0}|U|^{2}(s)]\la \underline{e},
\chi_{>0}\bm\bar{U}(s)\\U(s)\endm\times e^{-i(t-s)\calH^{*}}\bm
\phi\\ \psi\endm_{dis}\ra\\0\endm]}(\xi)d\xi
\end{split}\end{equation}
Freeze the frequency of $P_{\la s\ra^{-\frac{1}{2}}\geq .\geq
a}[\delta_{0}(x)\bm T_{z}P_{<a}[|U|^{2}(s)]\\0\endm]$ to dyadic
size $\sim b$. By Bernstein's inequality we get
\begin{equation}\label{07}
||\triangle^{-1}\partial_{x}P_{b}[\delta_{0}(x)\bm
T_{z}P_{<a}[|U|^{2}(s)]\\0\endm]||_{L_{x}^{2}}\lesssim
b^{-\frac{1}{2}}\la s\ra^{-1+\epsilon(\delta_{2})}
\end{equation}
Next, note that
\begin{equation}\nonumber\begin{split}
&\tilde{\calF}[\bm T_{z}P_{\geq a}[\chi_{>0}|U|^{2}(s)]\la
\underline{e}, \chi_{>0}\bm\bar{U}(s)\\U(s)\endm\times
e^{-i(t-s)\calH^{*}}\bm \phi\\ \psi\endm\ra\\0\endm](\xi)
\\&=\int_{0}^{\infty}\la\overline{[(e^{ix\xi}-e^{-ix\xi}+(1+r(-\xi)e^{-ix\xi}))\underline{e}+\phi(x,\xi)]},\\&\hspace{2cm}\bm
T_{z}P_{\geq a}[\chi_{>0}|U|^{2}(s)]\la \underline{e},
\chi_{>0}\bm\bar{U}(s)\\U(s)\endm\times e^{-i(t-s)\calH^{*}}\bm
\phi\\ \psi\endm_{dis}\ra\\0\endm\ra dx
\end{split}\end{equation}
We treat here the most difficult contribution which comes as usual
from $e^{ix\xi}-e^{-ix\xi}$. Carrying out an integration by parts,
we have to estimate the following terms:
\begin{equation}\nonumber
\int_{0}^{\infty}\xi[e^{ix\xi}+e^{-ix\xi}]\bm T_{z}P_{\geq
a}[\chi_{>0}|U|^{2}(s)]\int_{x}^{\infty}\la \underline{e},
\chi_{>0}\bm\bar{U}(s)\\U(s)\endm\times e^{-i(t-s)\calH^{*}}\bm
\phi\\ \psi\endm_{dis}(s,y)\ra\\0\endm dy dx
\end{equation}
\begin{equation}\nonumber
\int_{0}^{\infty}[e^{ix\xi}-e^{-ix\xi}]\bm T_{z}P_{\geq
a}\partial_{x}[\chi_{>0}|U|^{2}(s)]\int_{x}^{\infty}\la
\underline{e}, \chi_{>0}\bm\bar{U}(s)\\U(s)\endm\times
e^{-i(t-s)\calH^{*}}\bm \phi\\ \psi\endm_{dis}(s,y)\ra\\0\endm dy
dx
\end{equation}
Our calculations for the low-low case above have taught us that we
may assume\footnote{More precisely, we may write $\bm U\\
\bar{U}\endm$ as the sum of two functions, one of which leads to a
trivially estimable contribution, while the other satisfies the
above inequality.}
\begin{equation}\nonumber
|\int_{x}^{\infty}\la \underline{e},
\chi_{>0}\bm\bar{U}(s,y)\\U(s,y)\endm\times
e^{-i(t-s)\calH^{*}}\bm \phi\\ \psi\endm_{dis}\ra dy|\lesssim \la
s\ra\la t-s\ra^{-\frac{3}{2}}
\end{equation}
 Using that
$|\xi|\sim b$, the ordinary Plancherel's Theorem then implies that
we have
\begin{equation}\nonumber\begin{split}
&||\int_{0}^{\infty}\xi[e^{ix\xi}+e^{-ix\xi}]\bm T_{z}P_{\geq
a}[\chi_{>0}|U|^{2}(s)]\int_{x}^{\infty}\la \underline{e},
\chi_{>0}\bm\bar{U}(s)\\U(s)\endm\times e^{-i(t-s)\calH^{*}}\bm
\phi\\ \psi\endm_{dis}(s,y)\ra\\0\endm dy dx||_{L_{\xi}^{2}}
\\&\lesssim \la s\ra\la t-s\ra^{-\frac{3}{2}} b||T_{z}P_{\geq
a}[\chi_{>0}|U|^{2}(s)]||_{L_{x}^{2}}\lesssim b \la
s\ra^{-\frac{3}{4}}\la s\ra\la t-s\ra^{-\frac{3}{2}}
\end{split}\end{equation}
The contribution of the term
\begin{equation}\nonumber
\int_{-\infty}^{0}[e^{ix\xi}-e^{-ix\xi}]\bm T_{z}P_{\geq
a}\partial_{x}[\chi_{>0}|U|^{2}(s)]\int_{x}^{\infty}\la
\underline{e}, \chi_{>0}\bm\bar{U}(s)\\U(s)\endm\times
e^{-i(t-s)\calH^{*}}\bm \phi\\ \psi\endm_{dis}(s,y)\ra\\0\endm dy
dx
\end{equation}
is handled similarly, arguing as in the high-high case, using
\eqref{keyidentity}. Combining this with the bound \eqref{07} from
before, we estimate
\begin{equation}\nonumber\begin{split}
&|\int_{-\infty}^{0}\la\partial_{x}\triangle^{-1}P_{b}[\delta_{0}(x)\bm
T_{z}P_{<a}[|U|^{2}(s)]\\0\endm],
[e^{ix\xi}-e^{-ix\xi}]\underline{e}\ra\\&\hspace{2cm}\overline{\tilde{\calF}[\bm
T_{z}P_{\geq a}[\chi_{>0}|U|^{2}(s)]\la \underline{e},
\chi_{>0}\bm\bar{U}(s)\\U(s)\endm\times e^{-i(t-s)\calH^{*}}\bm
\phi\\ \psi\endm_{dis}\ra\\0\endm]}(\xi)d\xi|\\
&\lesssim b^{-\frac{1}{2}}\la s\ra^{-1+\epsilon(\delta_{2})} b\la
s\ra^{\frac{1}{4}}\la t-s\ra^{-\frac{3}{2}}
\end{split}\end{equation}
Summing over all dyadic $b$ with $a<b<\la s\ra^{-\frac{1}{2}}$ and
integrating over $s$ results in the bound $\lesssim \la
t\ra^{-\frac{3}{2}+\delta_{3}}$.  Next, we consider the
contribution of
\begin{equation}\nonumber\begin{split}
&\int_{-\infty}^{0}\la\partial_{x}\triangle^{-1}P_{>\la
s\ra^{-\frac{1}{2}}}[\delta_{0}(x)\bm
T_{z}P_{<a}[|U|^{2}(s)]\\0\endm],
[e^{ix\xi}-e^{-ix\xi}]\underline{e}\ra\\&\hspace{2cm}\overline{\tilde{\calF}[\bm
T_{z}P_{\geq a}[\chi_{>0}|U|^{2}(s)]\la \underline{e},
\chi_{>0}\bm\bar{U}(s)\\U(s)\endm\times e^{-i(t-s)\calH^{*}}\bm
\phi\\ \psi\endm_{dis}\ra\\0\endm]}(\xi)d\xi
\end{split}\end{equation}
We may again essentially replace $\tilde{\calF}$ by the ordinary
Fourier transform, and invoke the ordinary Plancherel's Theorem to
replace this by (up to negligible errors)
\begin{equation}\nonumber\begin{split}
&\la\partial_{x}\triangle^{-1}P_{>\la
s\ra^{-\frac{1}{2}}}[\delta_{0}(x)\bm
T_{z}P_{<a}[|U|^{2}(s)]\\0\endm],
\\&\hspace{2cm}\Pi_{(t^{-1000},t^{1000})}\bm T_{z}P_{\geq
a}[\chi_{>0}|U|^{2}(s)]\la \underline{e},
\chi_{>0}\bm\bar{U}(s)\\U(s)\endm\times e^{-i(t-s)\calH^{*}}\bm
\phi\\ \psi\endm_{dis}\ra\\0\endm\ra,
\end{split}\end{equation}
where $\Pi_{(t^{-1000},t^{1000})}$ is as in the discussion of the
high-high case. We bound this by
\begin{equation}\nonumber\begin{split}
&\lesssim ||\partial_{x}\triangle^{-1}P_{>\la
s\ra^{-\frac{1}{2}}}[\delta_{0}(x)\bm
T_{z}P_{<a}[|U|^{2}(s)]\\0\endm]||_{L_{x}^{1}}\\&\hspace{2cm}||\Pi_{(t^{-1000},t^{1000})}\bm
T_{z}P_{\geq a}[\chi_{>0}|U|^{2}(s)]\la \underline{e},
\chi_{>0}\bm\bar{U}(s)\\U(s)\endm\times e^{-i(t-s)\calH^{*}}\bm
\phi\\ \psi\endm_{dis}\ra\\0\endm||_{L_{x}^{\infty}}
\end{split}\end{equation}
We can bound the preceding expression by
\begin{equation}\nonumber
\lesssim \log t \la s\ra^{\epsilon(\delta_{2})}\la
s\ra^{\frac{1}{2}}\la s\ra^{-1}\la s\ra^{-1}\la
s\ra^{\frac{1}{2}}\la t-s\ra^{-\frac{3}{2}},
\end{equation}
which upon integration over $s$ leads to an acceptable bound. Thus
in order to complete the discussion for the case $s<\frac{t}{2}$,
we need to estimate the expression
\begin{equation}\nonumber\begin{split}
&\int_{-\infty}^{0}\la P_{<a}[\chi_{>0}(x)\bm
T_{z}P_{<a}[|U|^{2}(s)]\\0\endm],
[e^{ix\xi}-e^{-ix\xi}]\underline{e}\ra\\&\hspace{2cm}\overline{\tilde{\calF}[\bm
T_{z}P_{\geq a}[\chi_{>0}|U|^{2}(s)]\la \underline{e},
\chi_{>0}\bm\bar{U}(s)\\U(s)\endm\times e^{-i(t-s)\calH^{*}}\bm
\phi\\ \psi\endm_{dis}\ra\\0\endm]}(\xi)d\xi
\end{split}\end{equation}
Keep in mind that we put $a=\la s\ra^{-\frac{3}{4}}$. As usual we
simplify $\tilde{\calF}[...]$ and carry out an integration by
parts, replacing this by
\begin{equation}\nonumber
\xi\la [e^{ix\xi}+e^{-ix\xi}]\underline{e}, [\bm T_{z}P_{\geq
a}[\chi_{>0}|U|^{2}(s)]\la \underline{e},
\chi_{>0}\int_{x}^{\infty}\bm\bar{U}(s)\\U(s)\endm\times
e^{-i(t-s)\calH^{*}}\bm \phi\\ \psi\endm_{dis}\ra\\0\endm\ra
\end{equation}
\begin{equation}\nonumber
\la [e^{ix\xi}-e^{-ix\xi}]\underline{e}, [\bm T_{z}P_{\geq
a}\partial_{x}[\chi_{>0}|U|^{2}(s)]\la \underline{e},
\chi_{>0}\int_{x}^{\infty}\bm\bar{U}(s)\\U(s)\endm\times
e^{-i(t-s)\calH^{*}}\bm \phi\\ \psi\endm_{dis}\ra\\0\endm\ra
\end{equation}
Consider the first of these terms. The 2nd is treated similarly,
using \eqref{keyidentity}. We have
\begin{equation}\nonumber\begin{split}
&||\chi_{<\la s\ra^{-\frac{3}{4}}}(\xi)\xi\la
[e^{ix\xi}+e^{-ix\xi}]\underline{e}, [\bm T_{z}P_{\geq
a}[\chi_{>0}|U|^{2}(s)]\la \underline{e},
\chi_{>0}\int_{x}^{\infty}\bm\bar{U}(s)\\U(s)\endm\times
e^{-i(t-s)\calH^{*}}\bm
\phi\\
\psi\endm_{dis}\ra\\0\endm\ra||_{L_{\xi}^{\infty}}\\&\lesssim \la
s\ra^{-\frac{3}{4}} || T_{z}P_{\geq
a}[\chi_{>0}|U|^{2}(s)]||_{L_{x}^{2}}||\chi_{>0}(x)\int_{x}^{\infty}\bm\bar{U}(s)\\U(s)\endm\times
e^{-i(t-s)\calH^{*}}\bm \phi\\ \psi\endm_{dis}||_{L_{x}^{2}}
\end{split}\end{equation}
From our treatment of the low-low case we may assume that
\begin{equation}\nonumber
||\chi_{>0}(x)\int_{x}^{\infty}\bm\bar{U}(s)\\U(s)\endm\times
e^{-i(t-s)\calH^{*}}\bm \phi\\
\psi\endm_{dis}||_{L_{x}^{2}}\lesssim \la s\ra \la
s\ra^{\frac{1}{4}} \la t-s\ra^{-\frac{3}{2}},
\end{equation}
More precisely we may decompose $\bm U\\ \bar{U}\endm$ into two
constituents one of which upon substitution into the original
quintilinear expression immediately yields the desired estimate,
while the other constituent satisfies the above estimate, see the
discussion of the low-low case. We also have
\begin{equation}\nonumber
|| T_{z}P_{\geq a}[\chi_{>0}|U|^{2}(s)]||_{L_{x}^{2}}\lesssim
a^{-1}\la s\ra^{-\frac{3}{2}}\lesssim \la s\ra^{-\frac{3}{4}}
\end{equation}
Combining this with
\begin{equation}\nonumber
||\la P_{<a}[\chi_{>0}(x)\bm T_{z}P_{<a}[|U|^{2}(s)]\\0\endm],
[e^{ix\xi}-e^{-ix\xi}]\underline{e}\ra||_{L_{\xi}^{1}}\lesssim \la
s\ra^{-\frac{3}{4}},
\end{equation}
we can bound
\begin{equation}\nonumber\begin{split}
&|\int_{-\infty}^{0}\la P_{<a}[\chi_{>0}(x)\bm
T_{z}P_{<a}[|U|^{2}(s)]\\0\endm],
[e^{ix\xi}-e^{-ix\xi}]\underline{e}\ra\\&\hspace{2cm}\overline{\xi\la
e^{ix\xi}+e^{-ix\xi},[\bm T_{z}P_{\geq a}[\chi_{>0}|U|^{2}(s)]\la
\underline{e},
\chi_{>0}\int_{x}^{\infty}\bm\bar{U}(s)\\U(s)\endm\times
e^{-i(t-s)\calH^{*}}\bm
\phi\\ \psi\endm_{dis}\ra\\0\endm]\ra}d\xi|\\
&\lesssim \la s\ra^{\frac{5}{4}}\la s\ra^{-\frac{3}{4}}\la
s\ra^{-\frac{3}{4}}\la s\ra^{-\frac{3}{4}}\la
t-s\ra^{-\frac{3}{2}},
\end{split}\end{equation}
which is again as desired. \\

{\bf{Case B}:} $s\geq \frac{t}{2}$. The procedure here is
basically identical to the preceding case {\bf{A}}, so we shall be
relatively short here: one divides into the cases
\begin{equation}\nonumber
\int_{\frac{t}{2}}^{t}\la \bm \chi_{>0}|U|^{4}U(s,.)\\
-\chi_{>0}|U|^{4}\bar{U}(s,.)\endm, e^{-i(t-s)\calH^{*}}\bm \phi\\
\psi\endm_{dis}\ra ds
\end{equation}
\begin{equation}\nonumber
\int_{\frac{t}{2}}^{t}\la \bm \chi_{<0}|U|^{4}U(s,.)\\
-\chi_{<0}|U|^{4}\bar{U}(s,.)\endm, e^{-i(t-s)\calH^{*}}\bm \phi\\
\psi\endm_{dis}\ra ds
\end{equation}
Both being treated similarly, we shall only consider the first
term. We easily reduce $\bm \chi_{>0}|U|^{4}U(s,.)\\
-\chi_{>0}|U|^{4}\bar{U}(s,.)\endm$ to its dispersive part: note
that
\begin{equation}\nonumber
|\la \la\bm \chi_{>0}|U|^{4}U(s,.)\\
-\chi_{>0}|U|^{4}\bar{U}(s,.)\endm, \xi_{k(j)}\ra \eta_{j}, e^{-i(t-s)\calH^{*}}\bm \phi\\
\psi\endm_{dis}\ra|\lesssim s^{-6+\epsilon}(t-s)^{-\frac{3}{2}},
\end{equation}
which is significantly better than what we need. Now use the
distorted Plancherel's Theorem to rewrite what remains as
\begin{equation}\nonumber
\sum_{\pm}\int_{-\infty}^{\infty}\int_{\frac{t}{2}}^{t}\calF_{\pm}\bm \chi_{>0}|U|^{4}\\
-\chi_{>0}|U|^{4}\endm(\xi) \overline{\tilde{\calF}_{\pm}[\bm
\bar{U}\\U\endm\times e^{-i(t-s)\calH^{*}}\bm\phi\\
\psi\endm_{dis}](\xi)} ds d\xi
\end{equation}
We consider here the case $+$ and  $\xi\in[0,\infty)$ and how one
has to modify the argument in case {\bf{A}} to get the desired
estimate. Analogous modifications will then also give the result
for $\xi\in (-\infty,0]$. Write (leaving out the subscript)
\begin{equation}\nonumber\begin{split}
&\int_{0}^{\infty}\int_{\frac{t}{2}}^{t}\calF\bm \chi_{>0}|U|^{4}\\
-\chi_{>0}|U|^{4}\endm(\xi) \overline{\tilde{\calF}[\bm
\bar{U}\\U\endm\times e^{-i(t-s)\calH^{*}}\bm\phi\\
\psi\endm_{dis}](\xi)} ds d\xi\\
&=\int_{0}^{\infty}\int_{\frac{t}{2}}^{t}\la\bm \chi_{>0}|U|^{4}\\
-\chi_{>0}|U|^{4}\endm,  s(\xi)e^{ix\xi}+\sigma_{3}\phi(x,\xi)\ra
\overline{\tilde{\calF}[\bm
\bar{U}\\U\endm\times e^{-i(t-s)\calH^{*}}\bm\phi\\
\psi\endm_{dis}](\xi)} ds d\xi\\
\end{split}\end{equation}
The contribution of the local term $\phi(x,\xi)$ is again easy to
handle. As usual invoke the decomposition
\begin{equation}\nonumber\begin{split}
&\int_{0}^{\infty}\int_{\frac{t}{2}}^{t}\la\bm \chi_{>0}|U|^{4}\\
-\chi_{>0}|U|^{4}\endm, s(\xi)e^{ix\xi}\ra
\overline{\tilde{\calF}[\bm
\bar{U}\\U\endm\times e^{-i(t-s)\calH^{*}}\bm\phi\\
\psi\endm_{dis}](\xi)} ds d\xi\\
&=\int_{0}^{\infty}\int_{\frac{t}{2}}^{t}\la\bm P_{\geq a}[\chi_{>0}|U|^{2}]P_{\geq a}[|U|^{2}]\\
-P_{\geq a}[\chi_{>0}|U|^{2}]P_{\geq a}[|U|^{2}]\endm,
s(\xi)e^{ix\xi}\ra \overline{\tilde{\calF}[\bm
\bar{U}\\U\endm\times e^{-i(t-s)\calH^{*}}\bm\phi\\
\psi\endm_{dis}](\xi)} ds d\xi\\
&+\int_{0}^{\infty}\int_{\frac{t}{2}}^{t}\la \bm P_{\geq a}[\chi_{>0}|U|^{2}]P_{<a}[|U|^{2}]\\
-P_{\geq a}[\chi_{>0}|U|^{2}]P_{<a}[|U|^{2}]\endm,
s(\xi)e^{ix\xi}\ra \overline{\tilde{\calF}[\bm
\bar{U}\\U\endm\times e^{-i(t-s)\calH^{*}}\bm\phi\\
\psi\endm_{dis}](\xi)} ds d\xi\\
&+\int_{0}^{\infty}\int_{\frac{t}{2}}^{t}\la \bm P_{<a}[\chi_{>0}|U|^{2}]P_{\geq a}[|U|^{2}]\\
-P_{<a}[\chi_{>0}|U|^{2}]P_{\geq a}[|U|^{2}]\endm,
s(\xi)e^{ix\xi}\ra \overline{\tilde{\calF}[\bm
\bar{U}\\U\endm\times e^{-i(t-s)\calH^{*}}\bm\phi\\
\psi\endm_{dis}](\xi)} ds d\xi\\
&+\int_{0}^{\infty}\int_{\frac{t}{2}}^{t}\la \bm P_{<a}[\chi_{>0}|U|^{2}]P_{<a}[|U|^{2}]\\
-P_{<a}[\chi_{>0}|U|^{2}]P_{<a}[|U|^{2}]\endm, s(\xi)e^{ix\xi}\ra
\overline{\tilde{\calF}[\bm
\bar{U}\\U\endm\times e^{-i(t-s)\calH^{*}}\bm\phi\\
\psi\endm_{dis}](\xi)} ds d\xi\\
\end{split}\end{equation}
We consider here the first term. Use that for $\xi>0$ we have
\begin{equation}\nonumber
\tilde{\calF}[\chi_{>0}(x)\bm
\bar{U}\\U\endm\times e^{-i(t-s)\calH^{*}}\bm\phi\\
\psi\endm_{dis}](\xi)= \la\chi_{>0}(x)\bm
\bar{U}\\U\endm\times e^{-i(t-s)\calH^{*}}\bm\phi\\
\psi\endm_{dis}, s(\xi)e^{ix\xi}+\phi(x,\xi)\ra
\end{equation}
The contribution from $\phi(x,\xi)$ is easy to handle: note that
\begin{equation}\nonumber
||\la\sigma_{3}\phi(x,\xi), \chi_{>0}(x)\bm
\bar{U}\\U\endm\times e^{-i(t-s)\calH^{*}}\bm\phi\\
\psi\endm_{dis}\ra||_{L_{\xi}^{2}} \lesssim \la
t-s\ra^{-\frac{3}{2}}s^{-\frac{3}{2}+\epsilon},
\end{equation}
which is clearly good enough to close everything. Thus we now need
to consider
\begin{equation}\nonumber
\int_{0}^{\infty}\int_{\frac{t}{2}}^{t}\la\bm P_{\geq a}[\chi_{>0}|U|^{2}]P_{\geq a}[|U|^{2}]\\
-P_{\geq a}[\chi_{>0}|U|^{2}]P_{\geq a}[|U|^{2}]\endm,
s(\xi)e^{ix\xi}\underline{e}\ra \overline{\la \chi_{>0}(x)\bm
\bar{U}\\U\endm\times e^{-i(t-s)\calH^{*}}\bm\phi\\
\psi\endm_{dis}, s(\xi)e^{ix\xi}\underline{e}\ra }ds d\xi
\end{equation}
Using the ordinary Plancherel's Theorem, we replace this by (up to
negligible error terms)
\begin{equation}\nonumber
\int_{0}^{\infty}\int_{\frac{t}{2}}^{t}\la P_{\geq
a}[\chi_{>0}|U|^{2}]P_{\geq a}[|U|^{2}],
\Pi_{(t^{-1000},t^{1000})}\la \underline{e}, \chi_{>0}(x)\bm
\bar{U}\\U\endm\times e^{-i(t-s)\calH^{*}}\bm\phi\\
\psi\endm_{dis}\ra
\end{equation}
This is all like in the case $s<\frac{t}{2}$. At this point,
though, we don't pull down the full power of
$(t-s)^{-\frac{3}{2}}$, but only $(t-s)^{-\frac{1}{2}}$, which
costs nothing in terms of weights. In other words, we estimate
\begin{equation}\nonumber
|\Pi_{(t^{-1000}, t^{1000})}\la \underline{e}, \chi_{>0}(x)\bm
\bar{U}\\U\endm\times e^{-i(t-s)\calH^{*}}\bm\phi\\
\psi\endm_{dis}\ra|\lesssim (t-s)^{-\frac{1}{2}}s^{-\frac{1}{2}}
\end{equation}
Then  use again that
\begin{equation}\nonumber
||P_{\geq a}[\chi_{>0}|U|^{2}]P_{\geq
a}[|U|^{2}]||_{L_{x}^{1}}\lesssim a^{-2}s^{-3}
\end{equation}
Putting these together and integrating up over $s>\frac{t}{2}$
yields the upper bound $t^{-\frac{3}{2}}$ up to an arbitrarily
small exponential error independent of $\delta_{3}$. The remaining
terms above follow by similar modifications from the arguments for
the case $s<\frac{t}{2}$, and are omitted. This establishes the
strong local dispersive estimate up to demonstrating the bound
$||C\partial_{x}U(s,.)||_{L_{x}^{2}}\lesssim
s^{\frac{1}{2}+\epsilon(\delta_{2})}$, which we shall do in the
next subsection.

\subsection{Establishing the pseudo-conformal almost conservation law.}

We now demonstrate that\\ $\sup_{T>t\geq
0}||C\tilde{U}(t,.)||_{L_{x}^{2}}\leq \frac{\Lambda}{100}\delta$,
provided we have already improved all the estimates of
\eqref{disasympto} {\it{without the norm $\sup_{T\geq s\geq
0}||CU(s,y)||_{L_{y}^{2}}$}} to have righthand side
$\frac{\Lambda}{1000}\delta$; also, assume that we have already
improved the bound on \eqref{rootasympto}\footnote{These will
later be improved independently of this subsection, of course.} to
be $\leq \frac{\Lambda\delta^{2}}{1000}$, say. It is clear that
we may alternatively prove this estimate for $\bm U\\
\bar{U}\endm(t,.)$, in light of the already bootstrapped bounds.
Write\footnote{Thus by abuse of notation we use the same symbol
$C$ for this matrix-valued pseudo-conformal operator as for the
scalar operator $C=x-2tp$. This shouldn't cause confusion.}
\begin{equation}\nonumber
i\partial_{t}\la C\bm U\\ \bar{U}\endm, \bm U\\ \bar{U}\endm\ra
=\la i\dot{C}\bm U\\ \bar{U}\endm+iC\partial_{t} \bm U\\ \bar{U}\endm, \bm U\\
\bar{U}\endm\ra-\la C \bm U\\ \bar{U}\endm, i\partial_{t} \bm U\\
\bar{U}\endm\ra,
\end{equation}
where
\begin{equation}\nonumber
C=\bm (x-2tp)^{2}&0\\0&(x+2tp)^{2}\endm,\,p=-i\partial_{x}
\end{equation}
Then write the equation for $\bm U\\ \bar{U}\endm$ schematically
as follows: with $\calH_{0}=\bm \partial_{y}^{2}-1 & 0\\
0&1-\partial_{y}^{2}\endm$
\begin{equation}\nonumber
i\partial_{t}\bm U\\ \bar{U}\endm +\calH_{0} \bm U\\
\bar{U}\endm= -\bm 3\nu^{2}\tilde{\phi}_{0}^{4}&
2\nu^{2}\tilde{\phi}_{0}^{4}e^{2i(\Psi-\Psi_{\infty})}\\-2\nu^{2}\tilde{\phi}_{0}^{4}e^{-2i(\Psi-\Psi_{\infty})}&-3\nu^{2}\tilde{\phi}_{0}^{4}\endm\bm
U\\ \bar{U}\endm+\dot{\pi}\partial_{\pi}W+N(U,\pi)
\end{equation}
where
$\tilde{\phi}=\phi(\nu(t)(.-\lambda_{\infty}[\mu-\mu_{\infty}](t))$.
If we substitute this back into the preceding and expand
$\dot{C}$, we obtain the relation
\begin{equation}\nonumber\begin{split}
i\partial_{t}\la C\bm U\\ \bar{U}\endm, \bm U\\ \bar{U}\endm\ra=
&\la \bm 8itp^{2}-2-4ixp&0\\0&2+8itp^{2}+4ixp\endm\bm U\\
\bar{U}\endm\\&+C(-\calH_{0}\bm U\\ \bar{U}\endm-\bm
3\nu^{2}\tilde{\phi}_{0}^{4}&
2\nu^{2}\tilde{\phi}_{0}^{4}e^{2i(\Psi-\Psi_{\infty})}\\-2\nu^{2}\tilde{\phi}_{0}^{4}e^{-2i(\Psi-\Psi_{\infty})}&-3\nu^{2}\tilde{\phi}_{0}^{4}\endm\bm
U\\ \bar{U}\endm\\&+\dot{\pi}\partial_{\pi}W+N(U,\pi)),\bm U\\ \bar{U}\endm\ra\\&+\la C\bm U\\
\bar{U}\endm, \calH_{0}\bm U\\ \bar{U}\endm+\bm
3\nu^{2}\tilde{\phi}_{0}^{4}&
2\nu^{2}\tilde{\phi}_{0}^{4}e^{2i(\Psi-\Psi_{\infty})}\\-2\nu^{2}\tilde{\phi}_{0}^{4}e^{-2i(\Psi-\Psi_{\infty})}&-3\nu^{2}\tilde{\phi}_{0}^{4}\endm\bm
U\\ \bar{U}\endm\\&-\dot{\pi}\partial_{\pi}W-N(U,\pi)\ra
\end{split}\end{equation}
One can simplify the right hand side to the following:
\begin{equation}\nonumber\begin{split}
&\la -C\bm 3\nu^{2}\tilde{\phi}_{0}^{4}&
2\nu^{2}\tilde{\phi}_{0}^{4}e^{2i(\Psi-\Psi_{\infty})}\\-2\nu^{2}\tilde{\phi}_{0}^{4}e^{-2i(\Psi-\Psi_{\infty})}&-3\nu^{2}\tilde{\phi}_{0}^{4}\endm\bm
U\\ \bar{U}\endm\\&+\bm 3\nu^{2}\tilde{\phi}_{0}^{4}&
-2\nu^{2}\tilde{\phi}_{0}^{4}e^{2i(\Psi-\Psi_{\infty})}\\2\nu^{2}\tilde{\phi}_{0}^{4}e^{-2i(\Psi-\Psi_{\infty})}&-3\nu^{2}\tilde{\phi}_{0}^{4}\endm C\bm U\\ \bar{U}\endm, \bm U\\
\bar{U}\endm\ra\\&+\la C(\dot{\pi}\partial_{\pi}W+N(U,\pi)), \bm U\\
\bar{U}\endm\ra -\la C\bm U\\
\bar{U}\endm,\dot{\pi}\partial_{\pi}W+N(U,\pi)\ra
\end{split}\end{equation}
The last four terms here shall be fairly straightforward to
control. However, the first two appear to lead to a loss, as they
aren't absolutely integrable. Observe that we can rewrite the sum
of the first two terms as a commutator
\begin{equation}\nonumber
\la[\bm 3\nu^{2}\tilde{\phi}_{0}^{4}&
2\nu^{2}\tilde{\phi}_{0}^{4}e^{2i(\Psi-\Psi_{\infty})}\\2\nu^{2}\tilde{\phi}_{0}^{4}e^{-2i(\Psi-\Psi_{\infty})}&3\nu^{2}\tilde{\phi}_{0}^{4}\endm,
\bm 1&0\\0&-1\endm C]\bm U\\ \bar{U}\endm, \bm U\\ \bar{U}\endm\ra
\end{equation}
The trick here is to introduce a correction function
\begin{equation}\nonumber
t\rightarrow \theta(t):=t^{2}\la\bm 3\nu^{2}\tilde{\phi}_{0}^{4}&
2\nu^{2}\tilde{\phi}_{0}^{4}e^{2i(\Psi-\Psi_{\infty})}\\2\nu^{2}\tilde{\phi}_{0}^{4}e^{-2i(\Psi-\Psi_{\infty})}&3\nu^{2}\tilde{\phi}_{0}^{4}\endm\bm
U\\ \bar{U}\endm, \bm U\\ \bar{U}\endm\ra
\end{equation}
If one applies the time derivative to $\theta(t)$, the main
contribution comes from the terms when $\bm U\\ \bar{U}\endm$ gets
hit. Otherwise, one obtain at least an extra
$\partial_{t}[\Psi-\Psi_{\infty}]$, which makes the expression
absolutely integrable. Thus $\theta'(t)$ equals up to negligible
errors
\begin{equation}\nonumber\begin{split}
&i\theta'(t)\sim t^{2}\la\bm 3\nu^{2}\tilde{\phi}_{0}^{4}&
2e^{2i(\Psi-\Psi_{\infty})(t,y)}\nu^{2}\tilde{\phi}_{0}^{4}\\
2e^{-2i(\Psi-\Psi_{\infty})(t,y)}\nu^{2}\tilde{\phi}_{0}^{4}&
3\nu^{2}\tilde{\phi}_{0}^{4}\endm i\partial_{t}\bm U\\ \bar{U}\endm,\bm U\\
\bar{U}\endm\ra\\&-t^{2}\la\bm 3\nu^{2}\tilde{\phi}_{0}^{4}&
2e^{2i(\Psi-\Psi_{\infty})(t,y)}\nu^{2}\tilde{\phi}_{0}^{4}\\
2e^{-2i(\Psi-\Psi_{\infty})(t,y)}\nu^{2}\tilde{\phi}_{0}^{4}&
3\nu^{2}\tilde{\phi}_{0}^{4}\endm\bm
U\\ \bar{U}\endm,i\partial_{t}\bm U\\ \bar{U}\endm\ra \\
&=t^{2}\la\bm 3\nu^{2}\tilde{\phi}_{0}^{4}&
2e^{2i(\Psi-\Psi_{\infty})(t,y)}\nu^{2}\tilde{\phi}_{0}^{4}\\
2e^{-2i(\Psi-\Psi_{\infty})(t,y)}\nu^{2}\tilde{\phi}_{0}^{4}&
3\nu^{2}\tilde{\phi}_{0}^{4}\endm
\\&\hspace{2cm}[-\calH_{0}\bm U\\ \bar{U}\endm-\bm 3\nu^{2}\tilde{\phi}_{0}^{4}&
2e^{2i(\Psi-\Psi_{\infty})}\nu^{2}\tilde{\phi}_{0}^{4}\\
-2e^{-2i(\Psi-\Psi_{\infty})}\nu^{2}\tilde{\phi}_{0}^{4}& -3\nu^{2}\tilde{\phi}_{0}^{4}\endm\bm U\\ \bar{U}\endm+\ldots,\bm U\\
\bar{U}\endm\ra\\&-t^{2}\la\bm 3\nu^{2}\tilde{\phi}_{0}^{4}&
2e^{2i(\Psi-\Psi_{\infty})(t,y)}\nu^{2}\tilde{\phi}_{0}^{4}\\
2e^{-2i(\Psi-\Psi_{\infty})(t,y)}\nu^{2}\tilde{\phi}_{0}^{4}&
3\nu^{2}\tilde{\phi}_{0}^{4}\endm\bm
U\\ \bar{U}\endm,\\&\hspace{2cm}[-\calH_{0}\bm U\\
\bar{U}\endm-\bm 3\nu^{2}\tilde{\phi}_{0}^{4}&
2e^{2i(\Psi-\Psi_{\infty})}\nu^{2}\tilde{\phi}_{0}^{4}\\
-2e^{-2i(\Psi-\Psi_{\infty})}\nu^{2}\tilde{\phi}_{0}^{4}& -3\nu^{2}\tilde{\phi}_{0}^{4}\endm\bm U\\ \bar{U}\endm+\ldots, \bm U\\ \bar{U}\endm\ra \\
\end{split}\end{equation}
We shall see below that the terms denoted $\ldots$ lead to
absolutely integrable expressions. Observe the matrix identity
\begin{equation}\nonumber\begin{split}
&\bm 3\nu^{2}\tilde{\phi}_{0}^{4}&
2e^{2i(\Psi-\Psi_{\infty})}\nu^{2}\tilde{\phi}_{0}^{4}\\
2e^{-2i(\Psi-\Psi_{\infty})}\nu^{2}\tilde{\phi}_{0}^{4}&
3\nu^{2}\tilde{\phi}_{0}^{4}\endm \bm
3\nu^{2}\tilde{\phi}_{0}^{4}&
2e^{2i(\Psi-\Psi_{\infty})}\nu^{2}\tilde{\phi}_{0}^{4}\\
-2e^{-2i(\Psi-\Psi_{\infty})}\nu^{2}\tilde{\phi}_{0}^{4}&
-3\nu^{2}\tilde{\phi}_{0}^{4}\endm
\\&=\bm 3\nu^{2}\tilde{\phi}_{0}^{4}&
-2e^{2i(\Psi-\Psi_{\infty})}\nu^{2}\tilde{\phi}_{0}^{4}\\
2e^{-2i(\Psi-\Psi_{\infty})}\nu^{2}\tilde{\phi}_{0}^{4}&
-3\nu^{2}\tilde{\phi}_{0}^{4}\endm \bm
3\nu^{2}\tilde{\phi}_{0}^{4}&
2e^{2i(\Psi-\Psi_{\infty})}\nu^{2}\tilde{\phi}_{0}^{4}\\
2e^{-2i(\Psi-\Psi_{\infty})}\nu^{2}\tilde{\phi}_{0}^{4}&
3\nu^{2}\tilde{\phi}_{0}^{4}\endm
\end{split}\end{equation}
Thus the only contribution, up to smaller error terms, comes from
the commutator with $\calH_{0}$.  Observe that
\begin{equation}\nonumber\begin{split}
&t^{2}[\bm 3\nu^{2}\tilde{\phi}_{0}^{4}&
2e^{2i(\Psi-\Psi_{\infty})}\nu^{2}\tilde{\phi}_{0}^{4}\\
2e^{-2i(\Psi-\Psi_{\infty})}\nu^{2}\tilde{\phi}_{0}^{4}&
3\nu^{2}\tilde{\phi}_{0}^{4}\endm, \calH_{0}]\\&=[\bm
3\nu^{2}\tilde{\phi}_{0}^{4}&
2e^{2i(\Psi-\Psi_{\infty})}\nu^{2}\tilde{\phi}_{0}^{4}\\
2e^{-2i(\Psi-\Psi_{\infty})}\nu^{2}\tilde{\phi}_{0}^{4}&
3\nu^{2}\tilde{\phi}_{0}^{4}\endm,\tilde{C}\bm
1&0\\0&-1\endm]\\&+t^{2}[\bm 3\nu^{2}\tilde{\phi}_{0}^{4}&
2e^{2i(\Psi-\Psi_{\infty})}\nu^{2}\tilde{\phi}_{0}^{4}\\
2e^{-2i(\Psi-\Psi_{\infty})}\nu^{2}\tilde{\phi}_{0}^{4}&
3\nu^{2}\tilde{\phi}_{0}^{4}\endm, \bm 1&0\\0&-1\endm],
\end{split}\end{equation}
where we denote
\begin{equation}\nonumber
\tilde{C}=\bm
t^{2}\partial_{y}^{2}&0\\0&t^{2}\partial_{y}^{2}\endm
\end{equation}
Up to an error of order $\la(\Psi-\Psi_{\infty})_{2}U^{2},
\phi\ra$ and hence integrable, the last term in the preceding
expansion leads to the expression
\begin{equation}\label{symplectinullform}
t^{2}\la\tilde{U}^{2}(t,.)-\overline{\tilde{U}^{2}(t,.)},\phi\ra,
\end{equation}
for a suitable even Schwartz function\footnote{This function is
also time-dependent, but with uniform decay estimates in time.}
$\phi$, where we are reverting to the notation used in the
preceding section. Although this still isn't absolutely
integrable, it oscillates sufficiently (due to an inherent
symplectic cancellation structure) that one can integrate it over
$[0,T]$. Indeed, we shall later show that for $\tilde{T}\leq T$,
we have
$\int_{\tilde{T}}^{T}t^{2}\la\tilde{U}_{dis}^{2}(t,.)-\overline{\tilde{U}_{dis}^{2}(t,.)},\phi\ra
dt\lesssim \tilde{T}^{-(\frac{1}{2}-\delta_{1})}$, uniformly in
$T$. This in addition to the bounds on the
$|\lambda_{i}(t)|\lesssim \la t\ra^{-2+\delta_{1}}$ derived in the
next subsection will suffice. We now see that up to establishing
absolute integrability of the expressions
\begin{equation}\nonumber
t^{2}\la\dot\pi\partial_{\pi}W, \bm U\\
\bar{U}\endm\ra,\,t^{2}\la\dot\pi\partial_{\pi}W, \partial_{y}^{i}\bm U\\
\bar{U}\endm\ra,i=1,2,\,\la C\bm U\\\bar{U}\endm, \bm |U|^{4}U\\
-|U|^{4}\bar{U}\endm\ra
\end{equation}
plus similar contributions from the other local terms in
$N(U,\pi)$, we have that
\begin{equation}\nonumber
\int_{0}^{T}i\partial_{t}[\la C\bm U\\ \bar{U}\endm, \bm U\\
\bar{U}\endm\ra-\frac{\theta(t)}{4}]=O(1)
\end{equation}
This clearly allows us to retrieve pseudo-conformal
almost-conservation. Consider $t^{2}\la\dot\pi\partial_{\pi}W, \bm U\\
\bar{U}\endm\ra$. Using \eqref{modulation}, we see that these are
all equivalent to expressions of the form
\begin{equation}\nonumber
t^{2}(\nu-1)(t)\la U^{2},
\phi\ra,\,t^{2}\lambda_{6}^{2}(t),\,t^{2}\lambda_{6}(t)\la
U,\phi\ra
\end{equation}
as well as terms of higher order in $U$, $\lambda_{6}$. These are
all easily seen to be integrable in light of
\eqref{modulationasympto}, \eqref{rootasympto},
\eqref{disasympto}. One has to argue a bit differently when the
derivative $\partial_{y}$ in $C$ falls on $U$, since we haven't
built in a strong local dispersive estimate for $\nabla U$. In
this case, we need to leak a little extra: for a Schwartz function
$\phi$ write
\begin{equation}\nonumber
\phi\partial_{y}(U(s,.))=\phi
P_{<s^{\epsilon(\frac{1}{N})}}\partial_{y}(U(s,.))+\phi P_{\geq
s^{\epsilon(\frac{1}{N})}}\partial_{y}(U(s,.))
\end{equation}
Note that if we choose $\epsilon(\frac{1}{N})$ large enough, we
can ensure that
\begin{equation}\nonumber
||\phi P_{\geq
s^{\epsilon(\frac{1}{N})}}\partial_{y}(U(s,.))||_{L_{x}^{1}}\lesssim
\la s\ra ^{-N_{0}}
\end{equation}
for large $N_{0}=N_{0}(N)$.  Next, use a compactly supported
partition of unity $\{\phi_{j}\}$ with $\phi_{0}(y)$ centered at
$y=0$ to write
\begin{equation}\nonumber
\phi
P_{<s^{\epsilon(\frac{1}{N})}}\partial_{y}(U(s,.))=\sum_{j}\phi
P_{<s^{\epsilon(\frac{1}{N})}}\partial_{y}(\phi_{j}U(s,.))
\end{equation}
Then we have
\begin{equation}\nonumber
||\phi
P_{<s^{\epsilon(\frac{1}{N})}}\partial_{y}(\phi_{j}U(s,.))||_{L_{x}^{1}}\lesssim
j^{-\tilde{N}}
\end{equation}
for any $\tilde{N}$ and $j>s^{\epsilon_{1}(\frac{1}{N})}$, whence
we may restrict to $j\leq s^{\epsilon_{1}(\frac{1}{N})}$. In that
case, use the fact that the proof of the pseudo-conformal
conservation law only required control of finitely many weighted
estimates involving $\phi$ and its derivatives to conclude that
\begin{equation}\nonumber
||\phi
P_{<s^{\epsilon(\frac{1}{N})}}\partial_{y}(\phi_{j}U(s,.))||_{L_{x}^{1}}\lesssim
s^{\epsilon_{2}(\frac{1}{N})}s^{-\frac{3}{2}}
\end{equation}
Thus one obtains in summary (with a similar estimate for the 2nd
derivative)
\begin{equation}\nonumber
||\phi\partial_{y}U(s,.)||_{L_{x}^{1}}\lesssim \la
s\ra^{-\frac{3}{2}+\epsilon(\frac{1}{N})}
\end{equation}
Then one can proceed as before to estimate $t^{2}\la
\dot{\pi}\partial_{\pi}W, \partial_{y}^{1,2}U\ra$. The remaining
local terms in the nonlinearity shall be treatable along analogous
lines, hence we now turn to the contribution of the non-local
term, which is
\begin{equation}\label{pcquintilinear}\begin{split}
&\la C\bm |U|^{4}U\\
-|U|^{4}\bar{U}\endm(t,.), \bm U\\ \bar{U}\endm(t,.)\ra\\&=\la \bm x-2tp&0\\0&x+2tp\endm\bm |U|^{4}U\\
-|U|^{4}\bar{U}\endm(t,.), \bm x-2tp&0\\0&x+2tp\endm\bm U\\
\bar{U}\endm(t,.)\ra
\end{split}\end{equation}
We have
\begin{equation}\nonumber
(x-2tp)[|U|^{4}U]=-2tp[|U|^{4}]U+|U|^{4}(x-2tp)U
\end{equation}
Also, we have
\begin{equation}\nonumber
2tp[|U|^{4}]=4tp[|U|^{2}]|U|^{2}=2[-(x-2tp)U\bar{U}+U\overline{(x-2tp)U}]|U|^{2}
\end{equation}
Thus we can expand \eqref{pcquintilinear} as a sum of expressions
of the form
\begin{equation}\nonumber
[(x-2tp)U][(x+2tp)\bar{U}]|U|(t)^{2}U^{2}(t)
\end{equation}
plus similar terms. We can bound the $L_{x}^{1}$-norm of this by
$\lesssim \la t\ra^{-2+\epsilon(\delta_{2})}$, more than what we
need. We note that one may similarly deduce a global bound for
$||C\partial_{y}U(t,.)||_{L_{x}^{2}}$. But we only need the bound
$||C\partial_{y}U||_{L_{x}^{2}}\lesssim \la t\ra^{\frac{1}{2}}$,
anyways, see the last subsection. We are now done with
establishing the estimates for $U_{dis}$, up to bounding
\eqref{symplectinullform}, which we shall do later.

\subsection{Retrieving control over the root part.}

We now improve the bound on
\begin{equation}
\sup_{0\leq t<T}\sup_{0\leq k\leq [\frac{N}{2}]}||\la
t\ra^{2-2\delta_{1}}\frac{d^{k}}{dt^{k}}\lambda_{i}(t)||_{L^{M}}
\end{equation}
First, note from \eqref{lambda2}, \eqref{lambdai} as well as
\eqref{nu}, \eqref{beta} that we may immediately achieve this if
we have improved the estimates for $\lambda_{6}(t)$. Again we
shall suppress $\Lambda$ in the following, it being clear that
sufficiently many $\delta$'s will come up to improve the bound.
Now recall the ODE \eqref{lambda6}, as well as \eqref{lambda6''}.
We brutally expand $e^{\Lambda(s)}$ into a Taylor series. Then
note that the worst\footnote{This is of course a naive
qualification; this term actually oscillates a lot.} contribution
in $[...]$ in \eqref{lambda6''} comes from the
term $(\nu-1)(s)\la \bm \tilde{U}\\
\overline{\tilde{U}}\endm_{dis}, \phi\ra(s)$. Indeed, all other
expressions are easily seen to contribute a term decaying like
$\la s\ra^{-3+\delta_{1}+\delta_{3}}$ or faster. In order to treat
this bad term, we have to recycle the equation again(as we may).
Employing Duhamel expansion as usual, we have to bound the
following terms:
\begin{equation}\label{local01}\begin{split}
&\int_{t}^{\infty}(\nu-1)(s)\la
\int_{0}^{s}e^{i(s-\lambda)\calH}[\bm
0&-e^{-2i(\Psi-\Psi_{\infty})_{1}(s)+2i(\Psi-\Psi_{\infty})_{1}(\lambda)}+1\\e^{2i(\Psi-\Psi_{\infty})_{1}(s)-2i(\Psi-\Psi_{\infty})_{1}(\lambda)}-1&\endm\\&\hspace{10cm}\bm
\tilde{U}^{(s)}(\lambda,.)\\
\overline{\tilde{U}^{(s)}}(\lambda,.)\endm\phi]_{dis},\phi\ra
d\lambda ds
\end{split}\end{equation}
\begin{equation}\label{nonlocal01}
\int_{t}^{\infty}(\nu-1)(s)\la
\int_{0}^{s}e^{i(s-\lambda)\calH}\bm
|\tilde{U}^{(s)}|^{4}\tilde{U}^{(s)}(\lambda,.)\\
-|\tilde{U}^{(s)}|^{4}\overline{\tilde{U}^{(s)}}(\lambda,.)\endm_{dis},\phi\ra
d\lambda ds
\end{equation}
We commence with the first of these: introduce the vector valued
function
\begin{equation}\nonumber
\phi(s,\lambda,x):= \bm
0&-e^{-2i(\Psi-\Psi_{\infty})_{1}(s)+2i(\Psi-\Psi_{\infty})_{1}(\lambda)}+1\\e^{2i(\Psi-\Psi_{\infty})_{1}(s)-2i(\Psi-\Psi_{\infty})_{1}(\lambda)}-1&\endm\phi(x)
\end{equation}
Using the Plancherel's Theorem~\ref{DistortedPlancherel} for the
distorted Fourier transform, we can express the first term as
\begin{equation}\nonumber
\sum_{\pm}\int_{-\infty}^{\infty}\int_{t}^{\infty}(\nu-1)(s)\int_{0}^{s}e^{\pm i(s-\lambda)(\xi^{2}+1)}\calF_{\pm}[\phi(s,\lambda,x)\bm\tilde{U}^{(s)}\\
\overline{\tilde{U}^{(s)}}\endm(\lambda,.)
](\xi)\overline{\tilde{\calF}_{\pm}}\phi(\xi)d\lambda ds d\xi
\end{equation}
We perform an integration by parts in the $s$-variable, replacing
the above by the following terms:
\begin{equation}\label{09}
\sum_{\pm}\int_{-\infty}^{\infty}(\nu-1)(t)\int_{0}^{t}e^{\pm i(t-\lambda)(\xi^{2}+1)}\calF_{\pm}[\phi(t,\lambda,x)\bm\tilde{U}^{(t)}\\
\overline{\tilde{U}^{(t)}}\endm(\lambda)
](\xi)\frac{\overline{\tilde{\calF}_{\pm}\phi(\xi)}}{\xi^{2}+1}d\lambda
d\xi
\end{equation}
\begin{equation}\label{010}
\sum_{\pm}\int_{-\infty}^{\infty}\int_{t}^{\infty}(\nu-1)(s)\partial_{s}[\Psi-\Psi_{\infty}]_{1}(s)\int_{0}^{s}e^{\pm i(s-\lambda)(\xi^{2}+1)}\calF_{\pm}[\tilde{\phi}(s,\lambda,x)\bm\tilde{U}^{(s)}\\
\overline{\tilde{U}^{(s)}}\endm](\xi)\frac{\overline{\tilde{\calF}_{\pm}\phi(\xi)}}{\xi^{2}+1}d\lambda
ds d\xi
\end{equation}
\begin{equation}\begin{split}\label{011}
&\sum_{\pm}\int_{-\infty}^{\infty}\int_{t}^{\infty}(\nu-1)(s)\partial_{s}[\lambda_{\infty}(\mu-\mu_{\infty})](s)\\&\hspace{4cm}\int_{0}^{s}e^{\pm i(s-\lambda)(\xi^{2}+1)}\calF_{\pm}[\phi(s,\lambda,x)\partial_{y}\bm\tilde{U}^{(s)}\\
\overline{\tilde{U}^{(s)}}\endm](\xi)\frac{\overline{\tilde{\calF}_{\pm}\phi(\xi)}}{\xi^{2}+1}d\lambda
ds d\xi
\end{split}\end{equation}
\begin{equation}\label{012}
\sum_{\pm}\int_{-\infty}^{\infty}\int_{t}^{\infty}\dot{\nu}(s)\int_{0}^{s}e^{\pm i(s-\lambda)(\xi^{2}+1)}\calF_{\pm}[\phi(s,\lambda,x)\bm\tilde{U}^{(s)}\\
\overline{\tilde{U}^{(s)}}\endm(\lambda,.)
](\xi)\overline{\frac{\tilde{\calF}_{\pm}\phi(\xi)}{\xi^{2}+1}}d\lambda
ds d\xi
\end{equation}
Most of these are almost immediate to estimate. Note that upon
undoing the Fourier transform in the first expression, we obtain
\begin{equation}\nonumber
(\nu-1)(t)\int_{0}^{t}\la e^{i(t-\lambda)\calH}[\bm\tilde{U}(\lambda,.)\\
\overline{\tilde{U}(\lambda,.)}\endm \phi(t,\lambda,x)]_{dis},
(\calH^{*})^{-1}\phi(x)\ra d\lambda
\end{equation}
We claim that $||x^{2}\calH^{-1}\phi||_{L_{x}^{1+}}\lesssim O(1)$,
where $1+$ can be chosen arbitrarily close to $1$.  Write
\begin{equation}\label{019}
(\calH)^{-1}\phi=\sum_{\pm}\int_{-\infty}^{\infty}e_{\pm}(x,\xi)\frac{\tilde{\calF}_{\pm}(\phi)(\xi)}{\xi^{2}+1}d\xi
\end{equation}
We treat here the $+$ part, the other one being similar. We equate
the preceding for $x>0$ with
\begin{equation}\nonumber
\int_{0}^{\infty}[s(\xi)e^{ix\xi}\underline{e}+\phi(x,\xi)]\frac{\calF_{\pm}(\phi)(\xi)}{\xi^{2}+1}d\xi
+\int_{-\infty}^{0}([e^{ix\xi}+r(-\xi)e^{-ix\xi}]\underline{e}+\phi(x,\xi))\frac{\calF(\phi)(\xi)}{\xi^{2}+1}d\xi
\end{equation}
For the first integral, leaving out the trivial local part, we
have
\begin{equation}\nonumber\begin{split}
&\int_{0}^{\infty}s(\xi)e^{ix\xi}\underline{e}\frac{\calF_{\pm}(\phi)(\xi)}{\xi^{2}+1}d\xi=
-\frac{1}{ix}\int_{0}^{\infty}(\partial_{\xi}[s(\xi)]e^{ix\xi}\frac{\calF_{\pm}(\phi)(\xi)}{\xi^{2}+1}\underline{e}+s(\xi)e^{ix\xi}\underline{e}\partial_{\xi}[\frac{\calF_{\pm}(\phi)(\xi)}{\xi^{2}+1}])d\xi
\\&=+\frac{1}{-x^{2}}\int_{0}^{\infty}e^{ix\xi}\partial_{\xi}(\partial_{\xi}[s(\xi)]\frac{\calF_{\pm}(\phi)(\xi)}{\xi^{2}+1}\underline{e})d\xi
+\frac{1}{-x^{2}}\int_{0}^{\infty}e^{ix\xi}\partial_{\xi}[s(\xi)\underline{e}\partial_{\xi}[\frac{\calF_{\pm}(\phi)(\xi)}{\xi^{2}+1}]]d\xi
=O(\frac{1}{x^{3}})
\end{split}\end{equation}
Similarly (omitting the contributions from $1+r(-\xi)$,
$\phi(x,\xi)$), we have
\begin{equation}\nonumber\begin{split}
&\int_{-\infty}^{0}\frac{e^{ix\xi}-e^{-ix\xi}}{\xi^{2}+1}\underline{e}\calF_{\pm}(\phi)(\xi)d\xi
=-\frac{1}{ix}\int_{-\infty}^{0}[e^{ix\xi}+e^{-ix\xi}]\underline{e}\partial_{\xi}[\frac{\calF_{\pm}\phi(\xi)}{\xi^{2}+1}]d\xi
\\&=+\frac{1}{-x^{2}}\int_{-\infty}^{0}(e^{ix\xi}-e^{-ix\xi})\underline{e}\partial_{\xi}^{2}[\frac{\calF_{\pm}\phi(\xi)}{\xi^{2}+1}]d\xi=O(\frac{1}{x^{3}})
\end{split}\end{equation}
If we then repeat the steps in the proof of the strong local
dispersive estimate, we get
\begin{equation}\nonumber
|(\nu-1)(t)\int_{0}^{t}\la e^{\pm i(t-\lambda)\calH}[\bm\tilde{U}(\lambda,.)\\
\overline{\tilde{U}(\lambda,.)}\endm \phi(t,\lambda,x)]_{dis},
(\calH^{*})^{-1}\phi(x)\ra d\lambda|\lesssim
t^{-\frac{1}{2}+\delta_{1}}t^{-\frac{3}{2}+2\delta_{3}},
\end{equation}
as desired. In the expression \eqref{010}, carry out two
additional integrations by parts. This either produces additional
factors of at least the decay $\nu-1$, or else kills the integral
over $\lambda$, in which case one arrives at an
expression\footnote{Using the distorted Plancherel's
theorem~\ref{DistortedPlancherel}.}
\begin{equation}\nonumber
\int_{t}^{\infty}(\nu-1)(s)\frac{d}{ds}[\Psi-\Psi_{\infty}]_{1}(s)\la
\bm \overline{\tilde{U}}\\ \tilde{U}\endm(s,.)\phi(x),
(\calH^{*})^{-2}\phi_{dis}\ra ds
\end{equation}
Then note that $\la \bm \overline{\tilde{U}}\\
\tilde{U}\endm(s,.)\phi(x),
(\calH^{*})^{-2}\phi_{dis}\ra=\overline{\la \bm \tilde{U}\\
\overline{\tilde{U}}\endm(s,.)\phi(x),
(\calH^{*})^{-2}\overline{\phi_{dis}}\ra}$. Now use the customary
decomposition of $\bm \tilde{U}\\ \overline{\tilde{U}}\endm$ into
dispersive and root part. Recycling \eqref{lambdai}, we thus see
that up to an integral of the form
\begin{equation}\label{500}
\int_{t}^{\infty}(\nu-1)^{2}(t)\lambda_{6}(t)dt,
\end{equation}
we arrive at the expression we started out with but with an extra
weight of at least the strength
$(\nu-1)(s)\frac{d}{ds}[\Psi-\Psi_{\infty}]_{1}(s)\sim
(\nu(s)-1)^{2}$. Now iterate the procedure. The integral
\eqref{500} can be estimated using Proposition~\ref{Hard}.
The remaining terms are simpler; indeed, they can be integrated absolutely.\\
Now consider \eqref{nonlocal01}. Here we also pass to the Fourier
side, perform an integration by parts in $s$, and undo the Fourier
transform. This results in the extra weights
$\partial_{s}[\Psi-\Psi_{\infty}]_{1}$,
$\partial_{s}[\lambda_{\infty}(\mu-\mu_{\infty})]$, or else one
winds up with the expression
\begin{equation}\nonumber
\int_{t}^{\infty}\la \bm
|\tilde{U}|^{4}\tilde{U}(s)\\
-|\tilde{U}|^{4}\tilde{U}(s)\endm_{dis},\calH^{-1}\phi\ra ds
\end{equation}
The former cases are treated by repeating integration by parts in
$s$ if necessary and proceeding as in the proof of the strong
local dispersive estimate with $\phi$ replaced by
$\calH^{-1}\phi$, while in the latter case we simply estimate
(using the bound derived above on $||\la
x\ra^{2}\calH^{-1}\phi||_{L_{x}^{1+}}$)
\begin{equation}\nonumber
|\int_{t}^{\infty}\la \bm
|\tilde{U}|^{4}\tilde{U}(s)\\
-|\tilde{U}|^{4}\tilde{U}(s)\endm_{dis},\calH^{-1}\phi\ra
ds|\lesssim
\int_{t}^{\infty}s^{-\frac{9}{2+\epsilon(\delta_{2})}}ds\lesssim
t^{-2+2\delta_{1}}
\end{equation}
Finally, obtaining estimates for
$\frac{d^{k}}{dt^{k}}\lambda_{6}$, $1\leq k\leq [\frac{N}{2}]$, is
straightforward upon differentiating \eqref{lambda6}. We have used
a sleight of hand here, since suppressed the possible time
dependence of $\phi$ in $(\nu-1)\la \bm \tilde{U}\\
\overline{\tilde{U}}\endm_{dis}, \phi\ra$. However, as already
mentioned in the derivation of the equation for $\lambda_{6}$, the
derivative of $\phi$ with respect to time has at least the decay
of $\dot{\nu}$. Inspecting the above proof, one easily checks that
the additional terms generated upon integration by parts can be
handled by the same method. We are now done with the a priori
estimates for $U$.

\subsection{Interlude: deriving a refined estimate for
$\lambda_{6}(t)$.}

This is the most challenging subsection, and condenses all the
preceding considerations into one crucial estimate for
$\lambda_{6}(t)$ which appears indispensable to close the
estimates for the modulation parameters. Indeed, all our travails
in establishing the refined local decay estimates for $U$ are
really leading up to and flowing into this estimate, which we
state as the following Proposition:

\begin{proposition}\label{Hard} Let $\Gamma\in A^{(n)}[0,T)$, $0<T\leq \infty$, see Theorem~\ref{core}.
Then for $\tilde{T}\leq T$ and
\begin{equation}\nonumber
\Gamma=\{\bm
\tilde{U}\\\overline{U}\endm_{dis},\ldots,\lambda_{6}(t),\ldots\}
\end{equation}
we have the bound
\begin{equation}\nonumber
|\int_{\tilde{T}}^{T}t\lambda_{6}(t)dt|\lesssim \delta^{2}\la
\tilde{T}\ra^{-\frac{1}{2}+\delta_{1}}
\end{equation}
uniformly in $T$. Also, we have
\begin{equation}\nonumber
|\int_{\tilde{T}}^{T}(\nu-1)^{a}\lambda_{6}(t)dt|\lesssim
\delta^{2}\la
\tilde{T}\ra^{-\frac{3}{2}+\delta_{1}+a(-\frac{1}{2}+\delta_{1})},\,a\geq
0
\end{equation}
\end{proposition}

\begin{proof} This estimate is clearly significantly more
difficult than what we established in the previous
subsection\footnote{We did use the last conclusion of the
Proposition to deduce the point wise estimates for $\lambda_{6}$,
but this was certainly an overkill.}. We may put $T=\infty$, the
more general case being treated identically. Also, we shall prove
the first inequality, the 2nd following from the same proof. We
shall again recycle \eqref{lambda6''}, which amongst other
expressions will lead to $\int_{T}^{\infty}t\int_{t}^{\infty}\la
\tilde{U}^{2}(s,.)-\overline{\tilde{U}}^{2}(s,.),\phi\ra ds dt$.
The treatment of the latter shall also be applicable to
$\int_{T}^{\infty}t^{2}\la
\tilde{U}^{2}(t,.)-\overline{\tilde{U}}^{2}(t,.),\phi\ra dt$,
which will fill in the hole in retrieving the bound for
$\sup_{0\leq t<T}||CU(t,.)||_{L_{x}^{2}}$. The logic of the
argument below shall be that if an expression can't be integrated
absolutely, we integrate by parts until either it can be
integrated absolutely, or else we wind up essentially in the
position we started out with but with an extra gain. Thus in the
proof of the Lemmata below, it may be that we arrive at terms just
as the first one in Proposition~\ref{Hard}, but with
$t\lambda_{6}(t)$ replaced by $(\nu(t)-1)t\lambda_{6}(t)$. The
idea then is to reiterate the whole process again. Thus the
Lemmata below should be thought of as being proved in tandem. Now
from \eqref{lambda6} we need to estimate a number of expressions,
the first of which is
\begin{equation}\nonumber
\int_{T}^{\infty}t\int_{t}^{\infty}(\nu-1)(s)\la \bm \tilde{U}\\
\overline{\tilde{U}}\endm_{dis},\phi\ra ds
\end{equation}

\begin{lemma}\label{tedious}The following inequality holds under the assumptions
of Proposition~\ref{Hard}:
\begin{equation}\nonumber
|\int_{T}^{\infty}t\int_{t}^{\infty}(\nu-1)(s)\la \bm \tilde{U}\\
\overline{\tilde{U}}\endm_{dis}(s,.),\phi\ra ds|\lesssim
\delta^{3}\la T\ra^{-\frac{1}{2}+\delta_{1}}
\end{equation}
\end{lemma}
\begin{proof}
We Duhamel-expand $\bm \tilde{U}\\
\overline{\tilde{U}}\endm_{dis}$.  This leads to the expressions
\begin{equation}\label{013}
\int_{T}^{\infty}t\int_{t}^{\infty}(\nu-1)(s)\la
\int_{0}^{s}e^{i(s-\lambda)\calH}[\bm(1-e^{2i(\Psi-\Psi_{\infty})_{1}(\lambda)-2i(\Psi-\Psi_{\infty})_{1}(s)})\overline{\tilde{U}^{(s)}}(\lambda)\\
(-1+e^{-2i(\Psi-\Psi_{\infty})_{1}(\lambda)+2i(\Psi-\Psi_{\infty})_{1}(s)})\tilde{U}^{(s)}(\lambda)\endm\phi]_{dis}d\lambda,
\phi\ra ds dt
\end{equation}
\begin{equation}\label{014}
\int_{T}^{\infty}t\int_{t}^{\infty}(\nu-1)(s)\la
\int_{0}^{s}e^{i(s-\lambda)\calH}\bm
|\tilde{U}^{(s)}|^{4}\tilde{U}^{(s)}(\lambda,.)\\-|\tilde{U}^{(s)}|^{4}\overline{\tilde{U}}^{(s)}(\lambda,.)\endm_{dis}d\lambda,
\phi\ra ds dt
\end{equation}
as well as local terms with better decay behavior than the first
expression; these can be treated analogously. Start with the first
expression. We perform an integration by parts in $s$, and replace
it by the following list of terms upon going to the Fourier side.
We leave out $\pm$ for simplicity:
\begin{equation}\label{015}
\int_{-\infty}^{\infty}\int_{T}^{\infty}t(\nu-1)(t)
\int_{0}^{t}e^{i(t-\lambda)(\xi^{2}+1)}\calF[\bm(1-e^{2i(\Psi-\Psi_{\infty})_{1}(\lambda)-2i(\Psi-\Psi_{\infty})_{1}(t)})\overline{\tilde{U}^{(t)}}(\lambda,.)\\
(-1+e^{-2i(\Psi-\Psi_{\infty})_{1}(\lambda)+2i(\Psi-\Psi_{\infty})_{1}(t)})\tilde{U}^{(t)}(\lambda,.)\endm\phi](\xi)
\overline{\frac{\tilde{\calF}\phi}{\xi^{2}+1}} d\lambda dt d\xi
\end{equation}
\begin{equation}\label{016}
\int_{-\infty}^{\infty}\int_{T}^{\infty}t\int_{t}^{\infty}\dot{\nu}(s)
\int_{0}^{s}e^{i(s-\lambda)(\xi^{2}+1)}\calF[\bm(1-e^{2i(\Psi-\Psi_{\infty})_{1}(\lambda)-2i(\Psi-\Psi_{\infty})_{1}(s)})\overline{\tilde{U}^{(s)}}(\lambda,.)\\
(-1+e^{-2i(\Psi-\Psi_{\infty})_{1}(\lambda)+2i(\Psi-\Psi_{\infty})_{1}(s)})\tilde{U}^{(s)}(\lambda,.)\endm\phi](\xi)
\overline{\frac{\tilde{\calF}\phi}{\xi^{2}+1}} d\lambda dt d\xi
\end{equation}
\begin{equation}\label{017}\begin{split}
&\int_{-\infty}^{\infty}\int_{T}^{\infty}t\int_{t}^{\infty}(\nu-1)(s)\frac{d}{ds}[\Psi-\Psi_{\infty}]_{1}
\\&\hspace{2cm}\int_{0}^{s}e^{i(s-\lambda)(\xi^{2}+1)}\calF[\bm(1+e^{2i(\Psi-\Psi_{\infty})_{1}(\lambda)-2i(\Psi-\Psi_{\infty})_{1}(s)})\overline{\tilde{U}^{(s)}}(\lambda,.)\\
(1+e^{-2i(\Psi-\Psi_{\infty})_{1}(\lambda)+2i(\Psi-\Psi_{\infty})_{1}(s)})\tilde{U}^{(s)}(\lambda,.)\endm\phi](\xi)
\overline{\frac{\tilde{\calF}\phi}{\xi^{2}+1}} d\lambda dt d\xi
\end{split}\end{equation}
\begin{equation}\label{018}\begin{split}
&\int_{-\infty}^{\infty}\int_{T}^{\infty}t\int_{t}^{\infty}(\nu-1)(s)\frac{d}{ds}[\lambda_{\infty}(\mu-\mu_{\infty})]
\\&\hspace{2cm}\int_{0}^{s}e^{i(s-\lambda)(\xi^{2}+1)}\calF[\bm(1-e^{2i(\Psi-\Psi_{\infty})_{1}(\lambda)-2i(\Psi-\Psi_{\infty})_{1}(s)})\overline{\partial_{x}\tilde{U}^{(s)}}(\lambda,.)\\
(-1+e^{-2i(\Psi-\Psi_{\infty})_{1}(\lambda)+2i(\Psi-\Psi_{\infty})_{1}(s)})\partial_{x}\tilde{U}^{(s)}(\lambda,.)\endm\phi](\xi)
\overline{\frac{\tilde{\calF}\phi}{\xi^{2}+1}} d\lambda dt d\xi
\end{split}\end{equation}
A similar list of terms results from \eqref{014}. We commence with
the first term in this list. Perform an integration by parts in
$t$, thereby replacing it by
\begin{equation}\nonumber
\int_{-\infty}^{\infty}T(\nu-1)(T)
\int_{0}^{T}e^{i(T-\lambda)(\xi^{2}+1)}\calF[\bm(1-e^{2i(\Psi-\Psi_{\infty})_{1}(\lambda)-2i(\Psi-\Psi_{\infty})_{1}(T)})\overline{\tilde{U}^{(T)}}(\lambda,.)\\
(-1+e^{2i(\Psi-\Psi_{\infty})_{1}(\lambda)-2i(\Psi-\Psi_{\infty})_{1}(T)})\tilde{U}^{(T)}(\lambda,.)\endm\phi](\xi)
\overline{\frac{\tilde{\calF}\phi}{(\xi^{2}+1)^{2}}} d\lambda dT
d\xi
\end{equation}
\begin{equation}\nonumber
\int_{-\infty}^{\infty}\int_{T}^{\infty}(\nu-1)(t)
\int_{0}^{t}e^{i(t-\lambda)(\xi^{2}+1)}\calF[\bm(1-e^{2i(\Psi-\Psi_{\infty})_{1}(\lambda)-2i(\Psi-\Psi_{\infty})_{1}(t)})\overline{\tilde{U}^{(t)}}(\lambda,.)\\
(-1+e^{2i(\Psi-\Psi_{\infty})_{1}(\lambda)-2i(\Psi-\Psi_{\infty})_{1}(t)})\tilde{U}^{(t)}(\lambda,.)\endm\phi](\xi)
\overline{\frac{\tilde{\calF}\phi}{(\xi^{2}+1)^{2}}} d\lambda dt
d\xi
\end{equation}
\begin{equation}\nonumber\begin{split}
&\int_{-\infty}^{\infty}\int_{T}^{\infty}t(\nu-1)(t)\frac{d}{dt}[\Psi-\Psi_{\infty}]_{1}(t)
\\&\hspace{2cm}\int_{0}^{t}e^{i(t-\lambda)(\xi^{2}+1)}\calF[\bm(1+e^{2i(\Psi-\Psi_{\infty})_{1}(\lambda)-2i(\Psi-\Psi_{\infty})_{1}(t)})\overline{\tilde{U}^{(t)}}(\lambda,.)\\
(1+e^{-2i(\Psi-\Psi_{\infty})(\lambda)+2i(\Psi-\Psi_{\infty})(t)})\tilde{U}^{(t)}(\lambda,.)\endm\phi](\xi)
\overline{\frac{\tilde{\calF}\phi}{(\xi^{2}+1)^{2}}} d\lambda dt
d\xi
\end{split}\end{equation}
\begin{equation}\nonumber\begin{split}
&\int_{-\infty}^{\infty}\int_{T}^{\infty}t(\nu-1)(t)\frac{d}{dt}[\lambda_{\infty}(\mu-\mu_{\infty})](t)
\\&\hspace{2cm}\int_{0}^{t}e^{i(t-\lambda)(\xi^{2}+1)}\calF[\bm(1-e^{2i(\Psi-\Psi_{\infty})_{1}(\lambda)-2i(\Psi-\Psi_{\infty})_{1}(t)})\overline{\partial_{x}\tilde{U}^{(t)}}(\lambda,.)\\
(-1+e^{-2i(\Psi-\Psi_{\infty})_{1}(\lambda)+2i(\Psi-\Psi_{\infty})_{1}(t)})\partial_{x}\tilde{U}^{(t)}(\lambda,.)\endm\phi](\xi)
\overline{\frac{\tilde{\calF}\phi}{(\xi^{2}+1)^{2}}} d\lambda dt
d\xi
\end{split}\end{equation}
Each of these terms is straightforward to estimate: undo the
Fourier transform, using the distorted Plancherel's
Theorem~\ref{DistortedPlancherel}, thereby replacing $\phi$ by
$\calH^{-2}\phi$, which again satisfies
$||x^{2}\calH^{-2}\phi||_{L_{x}^{1+}}\lesssim O(1)$. Thus, for
example we get
\begin{equation}\nonumber\begin{split}
&|\int_{-\infty}^{\infty}T(\nu-1)(T)
\int_{0}^{T}e^{i(T-\lambda)(\xi^{2}+1)}\calF[\bm(1-e^{2i(\Psi-\Psi_{\infty})_{1}(\lambda)-2i(\Psi-\Psi_{\infty})_{1}(T)})\overline{\tilde{U}^{(T)}}(\lambda)\\
(-1+e^{-2i(\Psi-\Psi_{\infty})_{1}(\lambda)+2i(\Psi-\Psi_{\infty})_{1}(T)})\tilde{U}^{(T)}(\lambda)\endm\phi](\xi)
\overline{\frac{\tilde{\calF}\phi}{(\xi^{2}+1)^{2}}} d\lambda dT
d\xi|\\
&\lesssim T^{\frac{1}{2}+\delta_{1}}\la
T\ra^{-\frac{3}{2}}dT\lesssim T^{-1+\delta_{1}},
\end{split}\end{equation}
better than what we need. The 2nd and 4th term in the immediately
preceding list are estimated similarly. For the third term,
perform an additional integration by parts. Either one pulls down
an additional factor of at least the decay
$\frac{d}{dt}[\Psi-\Psi_{\infty}]_{1}(t)$, in which case one can
integrate absolutely to get the desired upper bound,
or\footnote{In the case when the inner integral gets abolished}
one obtains expressions of the
form\footnote{Use the decomposition $\bm \tilde{U}\\
\overline{\tilde{U}}\endm = \bm \tilde{U}\\
\overline{\tilde{U}}\endm_{dis}+\sum_{i=1}^{6}\lambda_{i}\eta_{i,\text{proper}}$;
we shall soon see that neither $\lambda_{2}$ nor $\lambda_{6}$
contribute anything to the expression in question, hence the
expressions, using \eqref{lambdai}.}
\begin{equation}\nonumber
\int_{T}^{\infty}t(\nu(t)-1)^{2+a}\la \bm \tilde{U}\\
\overline{\tilde{U}}\endm_{dis}(t,.), \phi\ra dt,\,a=0,1
\end{equation}
One can easily handle this by further Duhamel expansion, which we
leave out. We now turn to \eqref{016}. This is more difficult, as
the SLDE\footnote{Strong local dispersive estimate, i. e. $||\phi
U(t,.)||_{L_{x}^{\infty}}\lesssim \la
t\ra^{-\frac{3}{2}+\delta_{3}}$.} together with
\eqref{modulationasympto} do not suffice to obtain the needed $\la
T\ra^{-\frac{1}{2}+\delta_{1}}$ decay. Thus we need to exploit
that we may reiterate the equation, and use \eqref{nu}, which
implies
\begin{equation}\nonumber\begin{split}
\dot{\nu}(s)=&-b_{\infty}\int_{s}^{\infty}\lambda_{\infty}^{-1}(\sigma)[\nu(\sigma)(4i\kappa_{2})^{-1}E_{5}(\sigma)+\beta\nu(\nu-1)^{2}(\sigma)+(2\nu-1)(\beta\nu-b_{\infty}\lambda_{\infty}^{-1}(\sigma))]d\sigma\\
&+\nu(s)(4i\kappa_{2})^{-1}E_{5}(s)+\beta\nu(\nu-1)^{2}(s)+(2\nu-1)(\beta\nu-b_{\infty}\lambda_{\infty}^{-1})(s)
\end{split}\end{equation}
If we substitute the first line here for $\dot{\nu}(s)$ into
\eqref{016}, we perform further integrations by parts in $s$,
which either produces produces something which we can integrate
absolutely to get the desired upper bound (namely when we
differentiate the integral expression in the first row, or when we
pick up at least two factors
$\frac{d}{ds}[\Psi-\Psi_{\infty}]_{1}$), or else we arrive at an
expression just as in the statement of the Lemma but with an extra
weight decaying like $\dot{\nu}$, in which case we reiterate the
whole process. Also, substituting $\beta\nu(\nu-1)^{2}$ instead of
$\dot{\nu}$ is easily seen to yield more than what is required
upon performing further integrations by parts. The real difficulty
comes from substituting $\nu(4i\kappa_{2})^{-1}E_{5}(s)$, since a
priori it isn't clear whether the derivative of this term decays
faster than $\la s\ra^{-2+\delta_{1}+\delta_{3}}$. Now
recall\footnote{See \eqref{modulation}} that
\begin{equation}\nonumber
E_{5}=-\la N, \tilde{\xi}_{5}\ra+\la U,
(i\partial_{s}+\calH^{*}(s))\tilde{\xi}_{5}\ra
\end{equation}
The first term on the right decays at least as fast as $\la
s\ra^{-3+2\delta_{3}}$, and is easily seen to cause no problems.
The 2nd term on the right is essentially of the form
$(\nu-1)(s)\la \bm\tilde{U}\\
\overline{\tilde{U}}\endm_{dis},\phi\ra$, plus errors which are
again negligible.  Thus, in order to control \eqref{016}, we need
to substitute $(\nu-1)(s)\la \bm\tilde{U}\\
\overline{\tilde{U}}\endm_{dis},\phi\ra$ for $\dot{\nu}(s)$; our
only hope of succeeding here is to Duhamel-expand this 2nd
instance of $\bm \tilde{U}\\ \overline{\tilde{U}}\endm_{dis}$.
Pausing in our analysis here for a moment, observe that we run
into quite similar issues upon integrating by parts with respect
to $s$ in \eqref{014}. Thus what our problem really boils down to
is estimating the expression
\begin{equation}\label{generalbilinear}
\int_{T}^{\infty}t\int_{t}^{\infty}(\nu-1)(s)\la \bm
\tilde{U}\\\overline{\tilde{U}}\endm_{dis}(s,.),\phi\ra \la \bm
\tilde{U}\\\overline{\tilde{U}}\endm_{dis}(s,.),(\calH^{*})^{-1}\psi\ra
ds,
\end{equation}
where we use schematic notation. Let's assume for now that we can
bound this expression by $\la T\ra^{-\frac{1}{2}+\delta_{1}}$.
Finally, substituting
$\beta\nu-\frac{b_{\infty}}{\lambda_{\infty}}$ for $\dot{\nu}$ is
treated as in the preceding case: perform an additional
integration by parts in $s$; either one produces an extra factor
$\frac{d}{ds}[\Psi-\Psi_{\infty}]_{1}(s)$, in which case one
reiterates integration by parts to arrive either at a much
improved expression like in the statement of the Lemma, or to
arrive at an expression which can be integrated absolutely to
yield the bound $T^{-\frac{1}{2}+\delta_{1}}$. If on the other
hand one differentiates
$\beta\nu-\frac{b_{\infty}}{\lambda_{\infty}}$, use \eqref{beta}:
The worst contribution there comes from
$E_{2}(s)\sim\lambda_{6}(s)$. This one treats by another
integration by parts, which either produces an extra
$\frac{d}{ds}[\Psi-\Psi_{\infty}]_{1}(s)$ whence one can integrate
absolutely, or else one gets $\dot{\lambda_{6}}(s)$, which is
treated by recycling \eqref{lambda6}. There the only dangerous
contribution comes from $(\nu-1)\la \bm \tilde{U}\\
\overline{\tilde{U}}\endm_{dis}, \phi\ra$, which leads to an
expression as in \eqref{generalbilinear}. Let's move on to
\eqref{017}, reiterate integration by parts in $s$; either one
hits $(\nu-1)(s)\frac{d}{ds}[\Psi-\Psi_{\infty}]_{1}$, with
$\frac{d}{ds}$, or one abolishes the integration in $s$, or one
produces at least an extra factor
$\frac{d}{ds}[\Psi-\Psi_{\infty}]_{1}$. In the first case, we
obtain an expression like in \eqref{016} but with an extra weight
$\nu-1$, whence we can treat this case just as above (actually,
this time the contribution of $E_{5}$ can just be integrated
absolutely). In the third case, reiterate integration by parts,
which either takes one into the first two cases, or else produces
an additional $\frac{d}{ds}[\Psi-\Psi_{\infty}]_{1}$. In the last
case, keep integrating by parts, which either eventually produces
arbitrary gains in $s$, or else takes one into the first two
cases. The 2nd case is more tricky: Observe that then we obtain an
expression of the form\footnote{The function $\phi(x)$ here is
scalar- and real valued; indeed, one checks that
$\phi(x)=\phi_{0}^{4}(x)$.}
\begin{equation}\nonumber
\int_{T}^{\infty}t\int_{t}^{\infty}(\nu(s)-1)\frac{d}{ds}[\Psi-\Psi_{\infty}]_{1}(s)\la \bm\overline{\tilde{U}}\\
\tilde{U}\endm(s,x)\phi(x),\,(\calH^{*})^{-2}\phi_{dis}\ra ds
\end{equation}
Observe that a priori the factor $\bm \overline{\tilde{U}}\\
\tilde{U}\endm$ might imply the presence of a $\lambda_{6}(s)$,
which would lead to an extremely difficult term. However, close
inspection of the argument for the equation of $\lambda_{6}$ shows
that in every expression $(\nu-1)\la \bm\tilde{U}\\
\overline{\tilde{U}}\endm_{dis}, \phi\ra$ {\it{with exactly one
power}} of $\nu-1$, the (vector valued) function $\phi$ has the
form $\bm \alpha\\-\alpha\endm$ and is real valued.
This then implies that $\phi_{dis}=\bm \tilde{\alpha}\\
-\tilde{\alpha}\endm$, and then also
$(\calH^{*})^{-2}\phi_{dis}=\bm \beta\\-\beta\endm$, with real
valued scalar function $\beta$. Thus we get $\la
\eta_{6,\text{proper}}\phi(x),\bm \beta\\-\beta\endm=0$, and
similarly $\la \eta_{2,\text{proper}}\phi(x),\bm
\beta\\-\beta\endm=0$, whence
\begin{equation}\nonumber
\la \bm\overline{\tilde{U}}\\
\tilde{U}\endm(s,x)\phi(x),(\calH^{*})^{-2}\phi_{dis}\ra =\la\bm
\overline{\tilde{U}}\\
\tilde{U}\endm_{dis}\phi(x),(\calH^{*})^{-2}\phi_{dis}\ra+\sum_{i\neq
2, 6}\la
\eta_{i,\text{proper}}\phi(x),(\calH^{*})^{-2}\phi_{dis}\ra
\end{equation}
Thus using \eqref{lambdai} we have put ourselves into basically
the situation we started out with in the Lemma (up to negligible
error terms), but with a weight
$(\nu-1)(s)\frac{d}{ds}[\Psi-\Psi_{\infty}]_{1}\sim
(\nu-1)^{2}(s)$. Now reiterate the whole process. Observe that
performing the same reasoning for the expressions of the form
\begin{equation}\nonumber
\int_{T}^{\infty}t\int_{t}^{\infty}(\nu(s)-1)^{2}\la \bm
\tilde{U}\\ \overline{\tilde{U}}\endm_{dis}(s,.)\phi\ra ds dt
\end{equation}
which are also implied by \eqref{lambda6} we may very well arrive
at a term of the form
\begin{equation}\nonumber
\int_{T}^{\infty}t\int_{t}^{\infty}(\nu(s)-1)^{2}(\frac{d}{ds}[\Psi-\Psi_{\infty}](s))\lambda_{6}(s)ds
\end{equation}
Upon integration by parts\footnote{Also, use
$\frac{d}{ds}[\Psi-\Psi_{\infty}]_{1}\sim (\nu-1)(s)$.}, this
leads to
\begin{equation}\nonumber
\int_{T}^{\infty}t\int_{t}^{\infty}(\int_{s}^{\infty}(\nu(\sigma)-1)^{3}d\sigma)\dot{\lambda_{6}}(s)ds
\end{equation}
\begin{equation}\nonumber
\int_{T}^{\infty}t\int_{t}^{\infty}(\nu(\sigma)-1)^{3}d\sigma\lambda_{6}(t)dt
\end{equation}
The last term is what we started out with in
Proposition~\ref{Hard}, but with an extra factor
$\int_{t}^{\infty}(\nu(\sigma-1)^{3}d\sigma$. Then reiterate the
whole process. For the first expression, recycle \eqref{lambda6};
one winds up with terms just as in the Lemma, but with the extra
weight $\int_{s}^{\infty}(\nu(\sigma)-1)^{3}d\sigma$. Now
reiterate the process. Finally, the term \eqref{018} is also
treated by expanding
$\frac{d}{ds}[\lambda_{\infty}(\mu-\mu_{\infty})](s)$ using
\eqref{mu}. One obtains terms which can either be handled by
further integrations by parts in $s$, or else one is led to an
expression just as \eqref{generalbilinear}. As already mentioned,
the expression \eqref{014} is treated by exact analogy. One
reiterates the proof of the SLDE\footnote{Strong local dispersive
estimate} for each instance of $\la
\int_{0}^{t}\la e^{i(t-s)\calH}\bm |\tilde{U}|^{(t)}\tilde{U}^{(t)}\\
-|\tilde{U}|^{(t)}\overline{\tilde{U}^{(t)}}\endm_{dis}(s,.),(\calH^{*})^{-k}\phi\ra
ds$. This concludes the proof of the Lemma up to the assertion
concerning \eqref{generalbilinear}.
\end{proof}

We thus need

\begin{lemma}\label{general bilinear}\footnote{The functions $\phi,\psi$ are again generally time dependent in our applications,
with derivatives decaying at least like $\dot{\nu}$. This more
general case can be handled just as in the ensuing proof.} Under
the assumptions of the previous Lemma, we have
\begin{equation}\label{generalbilinear'}
|\int_{T}^{\infty}t\int_{t}^{\infty}(\nu-1)(s)\la \bm
\tilde{U}\\\overline{\tilde{U}}\endm_{dis}(s,.),\phi\ra \la \bm
\tilde{U}\\\overline{\tilde{U}}\endm_{dis}(s,.),(\calH^{*})^{-k}\psi\ra
ds|\lesssim \la T\ra^{-\frac{1}{2}+\delta_{1}}
\end{equation}
for any $k\geq 0$ and Schwartz functions $\phi, \psi$.
\end{lemma}

\begin{proof}
We Duhamel-expand each copy of $\bm
\tilde{U}\\\overline{\tilde{U}}\endm_{dis}(s,.)$:
\begin{equation}\nonumber\begin{split}
&\bm \tilde{U}\\ \bar{\tilde{U}}\endm_{dis}(s,.)=e^{is\calH}\bm
\tilde{U}^{(s)}\\
\overline{\tilde{U}^{(s)}}\endm_{dis}(0,.)\\&+\int_{0}^{s}e^{i(s-\lambda)\calH}[\bm
0&
1-e^{2i(\Psi-\Psi_{\infty})_{1}(\lambda)-2i(\Psi-\Psi_{\infty})_{1}(s)}\\1-e^{-2i(\Psi-\Psi_{\infty})_{1}(\lambda)+2i(\Psi-\Psi_{\infty})_{1}(s)}&0\endm
\bm \tilde{U}^{(s)}\\
\overline{\tilde{U}^{(s)}}\endm]_{dis}d\lambda+\ldots\\
&+\int_{0}^{s}e^{i(s-\lambda)\calH}\bm|\tilde{U}^{(s)}|^{4}\tilde{U}^{(s)}(\lambda,.)\\-|\tilde{U}^{(s)}|^{4}\overline{\tilde{U}^{(s)}(\lambda,.)}\endm_{dis}d\lambda:=
A+B+\ldots+C
\end{split}\end{equation}
The terms $\ldots$ are local terms with better decay properties
and can be treated in a simpler fashion, hence left out. We then
substitute either $A, B$ or $C$ for $\bm \tilde{U}\\
\bar{\tilde{U}}\endm_{dis}(s,.)$ in the right hand side of
\eqref{generalbilinear'} and check that the resulting expression
has the same decay as $(\nu-1)(T)$. The logic behind these
estimates is as follows: in an expression of the form
\begin{equation}\nonumber
\la e^{i(s-\lambda)\calH}E,\phi_{dis}\ra\la
e^{i(s-\lambda)\calH}F,\calH^{-k}\psi_{dis}\ra=\la E,
e^{-i(s-\lambda)\calH^{*}}\phi_{dis}\ra \la F,
e^{-i(s-\lambda)\calH^{*}}(\calH^{*})^{-k}\psi_{dis}\ra,
\end{equation}
distinguish between the case when $\phi_{dis}$ and $\psi_{dis}$
have separated Fourier support, respectively closely aligned
(correlated) Fourier support. In the former case, the product
oscillates strongly with respect to $s$, whence one can integrate
by parts and hope to gain. In the latter case, one should be able
to exploit some kind of 'diagonalization effect' in order to gain.
The argument proceeds by distinguishing between several
interactions:
\\

(AA): this is the expression
\begin{equation}\nonumber
\int_{T}^{\infty}t\int_{t}^{\infty}(\nu-1)(s)\la e^{is\calH}\bm
\tilde{U}^{(s)}\\
\overline{\tilde{U}^{(s)}}\endm_{dis}(0,.),\phi\ra\la
e^{is\calH}\bm
\tilde{U}^{(s)}\\
\overline{\tilde{U}^{(s)}}\endm_{dis}(0,.),(\calH^{*})^{-k}\psi\ra
ds
\end{equation}
This is straightforward to control upon invoking the spatial
localization on the initial data and
Theorem~\ref{linearestimates}: we get
\begin{equation}\nonumber\begin{split}
&|\int_{T}^{\infty}t\int_{t}^{\infty}(\nu-1)(s)\la e^{is\calH}\bm
\tilde{U}^{(s)}\\
\overline{\tilde{U}^{(s)}}\endm_{dis}(0,.),\phi\ra\la
e^{is\calH}\bm
\tilde{U}^{(s)}\\
\overline{\tilde{U}^{(s)}}\endm_{dis}(0,.),(\calH^{*})^{-k}\psi\ra ds|\\
&\lesssim
\int_{T}^{\infty}t\int_{t}^{\infty}s^{-\frac{1}{2}+\delta_{1}}s^{-3}ds
dt\lesssim T^{-\frac{1}{2}+\delta_{1}},
\end{split}\end{equation}
for $\delta_{1}$ small enough, where we use the estimate derived
after \eqref{019}, with a trivial modification.
\\

(BB): we employ schematic notation here; scalar quantities really
represent vectorial quantities: this is the expression
\begin{equation}\nonumber\begin{split}
&\int_{T}^{\infty}t\int_{t}^{\infty}(\nu-1)(s)\int_{0}^{s}\la
(1-e^{2i(\Psi-\Psi_{\infty})_{1}(\lambda)-2i(\Psi-\Psi_{\infty})_{1}(s)})U^{(s)}(\lambda,.)\phi,
e^{-i(s-\lambda)\calH^{*}}\phi_{dis}\ra d\lambda
\\&\int_{0}^{s}\la
(1-e^{2i(\Psi-\Psi_{\infty})_{1}(\lambda')-2i(\Psi-\Psi_{\infty})_{1}(s)})U^{(s)}(\lambda',.)\phi,
e^{-i(s-\lambda')\calH^{*}}(\calH^{*})^{-k}\psi_{dis}\ra d\lambda'
\end{split}\end{equation}
For either of the integrands on the right one obtains the bound
$\lambda^{-\frac{3}{2}+\delta_{3}}(s-\lambda)^{-\frac{3}{2}}$,
$\lambda'^{-\frac{3}{2}+\delta_{3}}(s-\lambda')^{-\frac{3}{2}}$.
This clearly suffices as long as $\lambda, \lambda'<\frac{s}{2}$.
If for example $\lambda>\frac{s}{2}$, we can close provided
$s-\lambda>s^{10\delta_{3}}$. In the opposite case, use that
\begin{equation}\nonumber
|1-e^{2i(\Psi-\Psi_{\infty})_{1}(\lambda)-2i(\Psi-\Psi_{\infty})_{1}(s)}|<s^{-\frac{1}{2}+\delta_{1}}s^{10\delta_{3}},
\end{equation}
which again allows us to close. The case (AB) is handled
similarly.
\\

(BC): This is the expression
\begin{equation}\label{BC}\begin{split}
&\int_{T}^{\infty}t\int_{t}^{\infty}(\nu-1)(s)\int_{0}^{s}\la
(1-e^{2i(\Psi-\Psi_{\infty})_{1}(\lambda)-2i(\Psi-\Psi_{\infty})_{1}(s)})U^{(s)}(\lambda,.)\phi,
e^{-i(s-\lambda)\calH^{*}}\phi_{dis}\ra d\lambda
\\&\int_{0}^{s}\la
\bm|\tilde{U}^{(s)}|^{4}\tilde{U}^{(s)}(\lambda',.)\\-|\tilde{U}^{(s)}|^{4}\overline{\tilde{U}^{(s)}(\lambda',.)}\endm_{dis},
e^{-i(s-\lambda')\calH^{*}}(\calH^{*})^{-k}\psi_{dis}\ra d\lambda'
\end{split}\end{equation}
This case is much more difficult. Our method here shall make heavy
use of microlocalization. The idea is to first reduce $\phi_{dis}$
in the local integrand (involving $\lambda$) to small
$\calH$-frequency\footnote{In the sense that the frequency is
close to either $1$ or $-1$.}. Then either $\psi_{dis}$ is at
small frequency, too, which case is handled by exploiting an extra
slack in the proof of the strong dispersive estimate, or else
there is a gap between the frequency supports of these functions
which forces sufficient oscillation in the $s$ variable to render
the full expression manageable. First, we observe that we may
reduce to $\lambda<\frac{s}{2}$. This follows from the preceding
calculation, since we gain a small power of $s$ in the case
$\lambda>\frac{s}{2}$. Now write
\begin{equation}\label{000}\begin{split}
&\chi_{>0}(x)e^{-i(s-\lambda)\calH^{*}}\phi_{dis}=\sum_{\pm}\chi_{>0}(x)\int_{0}^{\infty}e^{\pm
i(s-\lambda)(\xi^{2}+1)}(s(\xi)e^{ix\xi}\underline{e}_{\pm}+\phi(x,\xi))\tilde{\calF}_{\pm}(\phi_{dis})(\xi)d\xi
\\&+\sum_{\pm}\chi_{>0}(x)\int_{-\infty}^{0}e^{\pm
i(s-\lambda)(\xi^{2}+1)}[(e^{ix\xi}-e^{-ix\xi})\underline{e}_{\pm}+(1+r(-\xi))\underline{e}_{\pm}+\phi(x,\xi)]
\tilde{\calF}_{\pm}(\phi_{dis})(\xi)d\xi
\end{split}\end{equation}
We put $\underline{e}_{+}:=\underline{e}=\bm 1\\0\endm$,
$\underline{e}_{-}:=\sigma_{1}\underline{e}$. Recalling that
$\lambda<\frac{s}{2}$, we claim that we may build in a smooth
multiplier $\phi_{<s^{-\epsilon}}(\xi)$ localizing to a dilate of
the indicated region $|\xi|<s^{-\epsilon}$ in either integrand,
provided $\epsilon>0$ is small enough. Indeed, if on the flip side
we build in a multiplier of the form $\phi_{\geq
s^{-\epsilon}}(\xi)$, integration by parts in $\xi$ results in
arbitrary gains in $s$, at the cost of powers of $x$. These,
however, are absorbed by the local factor $U\phi$ above. Thus we
shall now replace
$\chi_{>0}(x)e^{-i(s-\lambda)\calH^{*}}\phi_{dis}$ by the sum of
the above two terms with an extra cutoff $\phi_{<
s^{-\epsilon}}(\xi)$ included. Denote this by
$\chi_{>0}(x)e^{-i(s-\lambda)\calH^{*}}\tilde{\phi}_{dis}$. Now we
consider the non-local integrand. As above write
\begin{equation}\nonumber\begin{split}
&\chi_{>0}(x)e^{-i(s-\lambda')\calH^{*}}(\calH^{*})^{-k}\psi_{dis}=\sum_{\pm}\chi_{>0}(x)\int_{0}^{\infty}e^{\pm
i(s-\lambda')(\xi^{2}+1)}(s(\xi)e^{ix\xi}\underline{e}_{\pm}+\phi(x,\xi))\frac{\tilde{\calF}_{\pm}(\psi_{dis})(\xi)}{(\xi^{2}+1)^{k}}d\xi
\\&+\sum_{\pm}\chi_{>0}(x)\int_{-\infty}^{0}e^{\pm
i(s-\lambda')(\xi^{2}+1)}[(e^{ix\xi}-e^{-ix\xi})\underline{e}_{\pm}+(1+r(-\xi))\underline{e}_{\pm}+\phi(x,\xi)]
\frac{\tilde{\calF}_{\pm}(\psi_{dis})(\xi)}{(\xi^{2}+1)^{k}}d\xi
\end{split}\end{equation}
We then distinguish between the cases $\lambda'<\frac{s}{2}$,
$\lambda'>\frac{s}{2}$.
\\

($\lambda'>\frac{s}{2}$). This case is simpler on account of the
fact that no extra integration by parts is required to produce the
gain of $(s-\lambda')^{-\frac{3}{2}}$, see the proof of strong
local dispersive estimate. The first step consists in reducing
$\lambda'$ to the range $[\frac{s}{2}, s-s^{\frac{1}{10}}]$. This
follows from the following simple calculation:
\begin{equation}\nonumber\begin{split}
&|\int_{s-s^{\frac{1}{10}}}^{s}\la \bm
|\widetilde{U^{(s)}}|^{4}\widetilde{U^{(s)}}(\lambda',.)\\
-|\widetilde{U^{(s)}}|^{4}\widetilde{U^{(s)}}(\lambda',.)\endm,
e^{-i(s-\lambda')\calH^{*}}(\calH^{*})^{-k}\psi_{dis}\ra
d\lambda'|\\&\lesssim
\int_{s-s^{\frac{1}{10}}}^{s}||\widetilde{U^{(s)}}(\lambda',.)||_{L_{x}^{\infty}}^{4}
||\widetilde{U^{(s)}}(\lambda',.)||_{L_{x}^{2}}||e^{-i(s-\lambda')\calH^{*}}(\calH^{*})^{-k}\psi_{dis}||_{L_{x}^{2}}d\lambda'\lesssim
s^{\frac{1}{10}}s^{-2+\epsilon(\delta_{2})},
\end{split}\end{equation}
whence we obtain
\begin{equation}\nonumber\begin{split}
&|\int_{T}^{\infty}t\int_{t}^{\infty}(\nu-1)(s)\int_{0}^{s}\la
(1-e^{2i(\Psi-\Psi_{\infty})_{1}(\lambda)-2i(\Psi-\Psi_{\infty})_{1}(s)})U^{(s)}(\lambda,.)\phi,
e^{-i(s-\lambda)\calH^{*}}\tilde{\phi}_{dis}\ra d\lambda
\\&\int_{s-s^{\frac{1}{10}}}^{s}\la
\bm|\tilde{U}^{(s)}|^{4}\tilde{U}^{(s)}(\lambda',.)\\-|\tilde{U}^{(s)}|^{4}\overline{\tilde{U}^{(s)}(\lambda',.)}\endm_{dis},
e^{-i(s-\lambda')\calH^{*}}(\calH^{*})^{-k}\psi_{dis}\ra d\lambda'|\\
&\lesssim
\int_{T}^{\infty}t\int_{t}^{\infty}s^{-\frac{1}{2}+\delta_{1}}s^{-\frac{3}{2}}s^{\frac{1}{10}}s^{-2+\epsilon(\delta_{2})}ds
dt\lesssim T^{-\frac{1}{2}+\delta_{1}}
\end{split}\end{equation}
Thus we now reduce to estimating
\begin{equation}\nonumber\begin{split}
&\int_{T}^{\infty}t\int_{t}^{\infty}(\nu-1)(s)\int_{0}^{\frac{s}{2}}\la
(1-e^{i(\Psi-\Psi_{\infty})_{1}(\lambda)-i(\Psi-\Psi_{\infty})_{1}(s)})U^{(s)}(\lambda,.)\phi,
e^{-i(s-\lambda)\calH^{*}}\tilde{\phi}_{dis}\ra d\lambda
\\&\int_{\frac{s}{2}}^{s-s^{\frac{1}{10}}}\la
\bm|\tilde{U}^{(s)}|^{4}\tilde{U}^{(s)}(\lambda',.)\\-|\tilde{U}^{(s)}|^{4}\overline{\tilde{U}^{(s)}(\lambda',.)}\endm_{dis},
e^{-i(s-\lambda')\calH^{*}}(\calH^{*})^{-k}\psi_{dis}\ra
d\lambda',
\end{split}\end{equation}
where it is to be kept in mind that
$e^{-i(s-\lambda)\calH^{*}}\tilde{\phi}_{dis}$ has modified
Fourier support as described above. We now mimic the proof of the
strong local dispersive estimate (SLDE) for the quintilinear
expression. Recall that we re-arrange the terms as follows:
\begin{equation}\nonumber\begin{split}
&\la
\bm|\tilde{U}^{(s)}|^{4}\tilde{U}^{(s)}(\lambda',.)\\-|\tilde{U}^{(s)}|^{4}\overline{\tilde{U}^{(s)}(\lambda',.)}\endm_{dis},
e^{-i(s-\lambda')\calH^{*}}(\calH^{*})^{-k}\psi_{dis}\ra\\&=\la
\bm|\tilde{U}^{(s)}|^{4}(\lambda',.)\\-|\tilde{U}^{(s)}|^{4}(\lambda',.)\endm,
\bm \overline{\tilde{U}^{(s)}(\lambda',.)}\\
\tilde{U}^{(s)}(\lambda',.)\endm \times
e^{-i(s-\lambda')\calH^{*}}(\calH^{*})^{-k}\psi_{dis}\ra
\end{split}\end{equation}
The first step consists in reducing both factors in the preceding
expression to their dispersive part: thus write
\begin{equation}\nonumber\begin{split}
&\la
\bm|\tilde{U}^{(s)}|^{4}(\lambda',.)\\-|\tilde{U}^{(s)}|^{4}(\lambda',.)\endm,
\bm \overline{\tilde{U}^{(s)}(\lambda',.)}\\
\tilde{U}^{(s)}(\lambda',.)\endm \times
e^{-i(s-\lambda')\calH^{*}}(\calH^{*})^{-k}\psi_{dis}\ra\\
&=\sum_{i}a_{i}\la
\la\bm|\tilde{U}^{(s)}|^{4}(\lambda',.)\\-|\tilde{U}^{(s)}|^{4}(\lambda',.)\endm,\xi_{k(i),proper}\ra\eta_{i,proper},
\bm \overline{\tilde{U}^{(s)}(\lambda',.)}\\
\tilde{U}^{(s)}(\lambda',.)\endm \times
e^{-i(s-\lambda')\calH^{*}}(\calH^{*})^{-k}\psi_{dis}\ra\\&+\la
\bm|\tilde{U}^{(s)}|^{4}(\lambda',.)\\-|\tilde{U}^{(s)}|^{4}(\lambda',.)\endm_{dis},
\bm \overline{\tilde{U}^{(s)}(\lambda',.)}\\
\tilde{U}^{(s)}(\lambda',.)\endm \times
e^{-i(s-\lambda')\calH^{*}}(\calH^{*})^{-k}\psi_{dis}\ra
\end{split}\end{equation}
Note that
\begin{equation}\nonumber\begin{split}
&|\sum_{i}a_{i}\la
\la\bm|\tilde{U}^{(s)}|^{4}(\lambda',.)\\-|\tilde{U}^{(s)}|^{4}(\lambda',.)\endm,\xi_{k(i),proper}\ra\eta_{i,proper},
\bm \overline{\tilde{U}^{(s)}(\lambda',.)}\\
\tilde{U}^{(s)}(\lambda',.)\endm \times
e^{-i(s-\lambda')\calH^{*}}(\calH^{*})^{-k}\psi_{dis}\ra|\\&\hspace{10cm}\lesssim
\lambda'^{-5.(\frac{3}{2}-\delta_{3})}(s-\lambda')^{-\frac{3}{2}},
\end{split}\end{equation}
and plugging this back into the above yields an acceptable bound.
Thus we now need to estimate
\begin{equation}\nonumber\begin{split}
&\int_{T}^{\infty}t\int_{t}^{\infty}(\nu-1)(s)\int_{0}^{\frac{s}{2}}\la
(1-e^{2i(\Psi-\Psi_{\infty})_{1}(\lambda)-2i(\Psi-\Psi_{\infty})_{1}(s)})U^{(s)}(\lambda,.)\phi,
e^{-i(s-\lambda)\calH^{*}}\tilde{\phi}_{dis}\ra d\lambda
\\&\int_{\frac{s}{2}}^{s-s^{\frac{1}{10}}}\la
\bm|\tilde{U}^{(s)}|^{4}(\lambda',.)\\-|\tilde{U}^{(s)}|^{4}(\lambda',.)\endm_{dis},
\bm \overline{\tilde{U}^{(s)}(\lambda',.)}\\
\tilde{U}^{(s)}(\lambda',.)\endm \times
e^{-i(s-\lambda')\calH^{*}}(\calH^{*})^{-k}\psi_{dis}\ra
d\lambda',
\end{split}\end{equation}
Recall from the proof of SLDE that we reformulate
\begin{equation}\nonumber\begin{split}
&\la
\bm|\tilde{U}^{(s)}|^{4}(\lambda',.)\\-|\tilde{U}^{(s)}|^{4}(\lambda',.)\endm_{dis},
\chi_{>0}(x)\bm \overline{\tilde{U}^{(s)}(\lambda',.)}\\
\tilde{U}^{(s)}(\lambda',.)\endm \times
e^{-i(s-\lambda')\calH^{*}}(\calH^{*})^{-k}\psi_{dis}\ra
\\&=\int_{-\infty}^{\infty}\calF\bm\chi_{>0}(x)|\tilde{U}^{(s)}|^{4}(\lambda',.)\\-\chi_{>0}(x)|\tilde{U}^{(s)}|^{4}(\lambda',.)\endm(\xi)
\overline{\tilde{\calF}[\bm \chi_{>0}(x)\overline{\tilde{U}^{(s)}(\lambda',.)}\\
\chi_{>0}(x)\tilde{U}^{(s)}(\lambda',.)\endm \times
e^{-i(s-\lambda')\calH^{*}}(\calH^{*})^{-k}\psi_{dis}](\xi)}d\xi
\end{split}\end{equation}
Break the integral into the contribution over $[0,\infty)$ and
$(-\infty,0]$. Consider for example the latter, the former being
treated similarly. We may write $\calF$, $\tilde{\calF}$ purely in
terms of the oscillatory part: for example, consider
\begin{equation}\label{101}
\int_{-\infty}^{0}\calF\bm\chi_{>0}(x)|\tilde{U}^{(s)}|^{4}(\lambda',.)\\-\chi_{>0}(x)|\tilde{U}^{(s)}|^{4}(\lambda',.)\endm(\xi)
\overline{\la [\bm \chi_{>0}(x)\overline{\tilde{U}^{(s)}(\lambda',.)}\\
\chi_{>0}(x)\tilde{U}^{(s)}(\lambda',.)\endm \times
e^{-i(s-\lambda')\calH^{*}}(\calH^{*})^{-k}\psi_{dis}],
\phi(x,\xi)\ra} d\xi
\end{equation}
Estimate this by
\begin{equation}\nonumber\begin{split}
&\lesssim
||\calF\bm\chi_{>0}(x)|\tilde{U}^{(s)}|^{4}(\lambda',.)\\-\chi_{>0}(x)|\tilde{U}^{(s)}|^{4}(\lambda',.)\endm||_{L_{\xi}^{2}}
||\la \phi(x,\xi), [\bm \chi_{>0}(x)\overline{\tilde{U}^{(s)}(\lambda',.)}\\
\chi_{>0}(x)\tilde{U}^{(s)}(\lambda',.)\endm \times
e^{-i(s-\lambda')\calH^{*}}(\calH^{*})^{-k}\psi_{dis}]\ra||_{L_{\xi}^{2}}\\&\lesssim
\lambda'^{-\frac{3}{2}}(s-\lambda')^{-\frac{3}{2}},
\end{split}\end{equation}
which then leads to an acceptable contribution. One argues
similarly for $\calF(...)$, and hence\footnote{As usual we only
consider $\calF_{+}$.} replaces \eqref{101} by
\begin{equation}\nonumber\begin{split}
&\int_{-\infty}^{0}\la\bm\chi_{>0}(x)|\tilde{U}^{(s)}|^{4}(\lambda',.)\\-\chi_{>0}(x)|\tilde{U}^{(s)}|^{4}(\lambda',.)\endm,
(e^{ix\xi}-e^{-ix\xi})\underline{e}+(1+r(-\xi))e^{-ix\xi}\underline{e}\ra
\\&\overline{\la \bm \chi_{>0}(x)\overline{\tilde{U}^{(s)}(\lambda',.)}\\
\chi_{>0}(x)\tilde{U}^{(s)}(\lambda',.)\endm \times
e^{-i(s-\lambda')\calH^{*}}(\calH^{*})^{-k}\psi_{dis},
(e^{ix\xi}-e^{-ix\xi})\underline{e}+(1+r(-\xi))e^{-ix\xi}\underline{e}\ra}
d\xi
\end{split}\end{equation}
We first deal with the contributions of the factors $(1+r(-\xi))$,
which are straightforward: note from the proof of SLDE that we
distinguish between the cases
$P_{><a}[\chi_{>0}(x)|\tilde{U}^{(s)}|^{2}]P_{><a}[|\tilde{U}^{(s)}|^{2}]$,
where we put $a=\lambda'^{-\frac{3}{4}}\sim s^{-\frac{3}{4}}$. We
then put $a=\lambda'^{-\frac{3}{4}+\epsilon}$ instead, which leads
to a small extra gain in $s$ for all cases except
$P_{<a}[\chi_{>0}|\tilde{U}^{(s)}|^{2}]P_{<a}[|\tilde{U}^{(s)}|^{2}]$.
Substituting this forces $|\xi|<s^{-\frac{3}{4}+\epsilon}$, and so
if either factor $(1+r(-\xi))$ occurs we again obtain a gain in
$s$ (if $\epsilon$ is small enough). Thus we now reduce to the
contribution of
\begin{equation}\nonumber\begin{split}
&\int_{-\infty}^{0}\la\bm\chi_{>0}(x)|\tilde{U}^{(s)}|^{4}(\lambda',.)\\-\chi_{>0}(x)|\tilde{U}^{(s)}|^{4}(\lambda',.)\endm,
(e^{ix\xi}-e^{-ix\xi})\underline{e}\ra
\\&\hspace{4cm}\la \bm \chi_{>0}(x)\overline{\tilde{U}^{(s)}(\lambda',.)}\\
\chi_{>0}(x)\tilde{U}^{(s)}(\lambda',.)\endm \times
e^{-i(s-\lambda')\calH^{*}}(\calH^{*})^{-k}\psi_{dis},
(e^{ix\xi}-e^{-ix\xi})\underline{e}\ra d\xi
\end{split}\end{equation}
Arguing as before, one may reduce here to the expression
\begin{equation}\label{102}\begin{split}
&\int_{-\infty}^{0}\la\bm
P_{<a}[\chi_{>0}(x)|\tilde{U}^{(s)}|^{2}](\lambda',.)P_{<a}[|\tilde{U}^{(s)}|^{2}](\lambda',.)\\-P_{<a}[\chi_{>0}(x)|\tilde{U}^{(s)}|^{2}](\lambda',.)P_{<a}[|\tilde{U}^{(s)}|^{2}](\lambda',.)\endm,
(e^{ix\xi}-e^{-ix\xi})\underline{e}\ra
\\&\hspace{4cm}\overline{\la \bm \chi_{>0}(x)\overline{\tilde{U}^{(s)}(\lambda',.)}\\
\chi_{>0}(x)\tilde{U}^{(s)}(\lambda',.)\endm \times
e^{-i(s-\lambda')\calH^{*}}(\calH^{*})^{-k}\psi_{dis},
(e^{ix\xi}-e^{-ix\xi})\underline{e}\ra }d\xi,
\end{split}\end{equation}
where $a=\lambda'^{-\frac{3}{4}+\epsilon}\sim
s^{-\frac{3}{4}+\epsilon}$. Still following the proof of SLDE, we
expand
\begin{equation}\nonumber\begin{split}
&\chi_{>0}(x)e^{-i(s-\lambda')\calH^{*}}(\calH^{*})^{-k}\psi_{dis}=\chi_{>0}(x)\sum_{\pm}\int_{0}^{\infty}e^{\pm
i(s-\lambda')(\xi^{2}+1)}(s(\xi)e^{ix\xi}\underline{e}_{\pm}+\phi_{\pm}(x,\xi))\frac{\calF_{\pm}(\psi_{dis})(\xi)}{(\xi^{2}+1)^{k}}d\xi
\\&\hspace{1cm}+\chi_{>0}(x)\sum_{\pm}\int_{-\infty}^{0}e^{\pm
i(s-\lambda')(\xi^{2}+1)}([e^{ix\xi}-e^{-ix\xi}]\underline{e}_{\pm}+(1+r(-\xi))e^{-ix\xi}\underline{e}_{\pm}+\phi_{\pm}(x,\xi))\frac{\calF_{\pm}(\psi_{dis})(\xi)}{(\xi^{2}+1)^{k}}d\xi
\end{split}\end{equation}
We shall localize $\xi$ here away from $0$. For example consider
the following expression:
\begin{equation}\nonumber
\int_{-\infty}^{0}e^{\pm
i(s-\lambda')(\xi^{2}+1)}[e^{ix\xi}-e^{-ix\xi}]\underline{e}_{\pm}\frac{\calF_{\pm}(\psi_{dis})(\xi)}{(\xi^{2}+1)^{k}}d\xi=\frac{e^{\pm
i(s-\lambda')}}{\sqrt{s-\lambda'}}\int_{-\infty}^{\infty}e^{\mp\frac{(x-y)^{2}}{i(s-\lambda')}}g(y)dy
\end{equation}
where
$g(y)=\int_{-\infty}^{0}[e^{iy\xi}-e^{-iy\xi}]\underline{e}_{\pm}\frac{\calF_{\pm}(\psi_{dis})(\xi)}{(\xi^{2}+1)^{k}}d\xi$.
It is straightforward to see that we may include a smooth
multiplier $\phi_{<s^{1000}}(\xi)$ here. With this modification we
then have $|g(y)|\lesssim \la y\ra ^{-2}\log s$,
whence\footnote{Following the proof of SLDE.} we may replace
$g(y)$ by
\begin{equation}\nonumber
\phi_{<s^{\frac{1}{10}}}(y)\int_{-\infty}^{0}\phi_{<s^{1000}}(\xi)[e^{iy\xi}-e^{-iy\xi}]\underline{e}_{\pm}\frac{\calF_{\pm}(\psi_{dis})(\xi)}{(\xi^{2}+1)^{k}}d\xi
\end{equation}
We now claim that we may further localize $\xi$ away from $0$. For
this include a sharp cutoff $\chi_{<s^{-\epsilon_{1}}}(\xi)$, i.
e. consider the contribution of
\begin{equation}\nonumber
\tilde{g}(y):=\phi_{<s^{\frac{1}{10}}}(y)\int_{-s^{-\epsilon_{1}}}^{0}\phi_{<s^{1000}}(\xi)[e^{iy\xi}-e^{-iy\xi}]\underline{e}_{\pm}\frac{\calF_{\pm}(\psi_{dis})(\xi)}{(\xi^{2}+1)^{k}}d\xi
\end{equation}
We integrate by parts here, and replace this expression by
\begin{equation}\nonumber
\phi_{<s^{\frac{1}{10}}}(y)\frac{e^{iy s^{-\epsilon_{1}}}+e^{-iy
s^{-\epsilon_{1}}}}{iy}\underline{e}_{\pm}\frac{\calF_{\pm}(\psi_{dis})(-s^{-\epsilon_{1}})}{(s^{-2\epsilon_{1}}+1)^{k}}-\phi_{<s^{\frac{1}{10}}}(y)\int_{-s^{-\epsilon_{1}}}^{0}\frac{[e^{iy\xi}+e^{-iy\xi}]}{iy}\underline{e}_{\pm}\partial_{\xi}[\phi_{<s^{1000}}(\xi)\frac{\calF_{\pm}(\psi_{dis})(\xi)}{(\xi^{2}+1)^{k}}]d\xi
\end{equation}
One now easily verifies that $||\tilde{g}(y)||_{L_{y}^{1}}\lesssim
\log s s^{-\epsilon_{1}}$, which extra gain suffices to close all
the other estimates upon reiterating the proof of SLDE. We now see
that we may replace $g(y)$ by the expression
\begin{equation}\nonumber
h(y):=\phi_{<s^{\frac{1}{10}}}(y)\int_{-\infty}^{-s^{-\epsilon_{1}}}\phi_{<s^{1000}}(\xi)[e^{iy\xi}-e^{-iy\xi}]\underline{e}_{\pm}\frac{\calF_{\pm}(\psi_{dis})(\xi)}{(\xi^{2}+1)^{k}}d\xi
\end{equation}
The Heisenberg uncertainty principle implies that this function
has frequency $\gtrsim s^{-\epsilon_{1}}$ (for
$\epsilon_{1}<<\frac{1}{10}$) up to errors of size $s^{-N}$, and
hence negligible, and the same comment applies to the function
\begin{equation}\nonumber
\frac{e^{\pm
i(s-\lambda')}}{\sqrt{s-\lambda'}}\int_{-\infty}^{\infty}e^{\mp\frac{(x-y)^{2}}{i(s-\lambda')}}h(y)dy=\int_{-\infty}^{\infty}e^{\pm
i(s-\lambda')(\xi'^{2}+1)}e^{ix\xi'}\hat{h}(\xi')d\xi'
\end{equation}
Thus we may replace\footnote{Generating negligible error terms}
this function by the following, where $\phi_{\geq
s^{-\epsilon_{1}}}$ is a smooth cutoff localizing to
$|\xi'|\gtrsim s^{-\epsilon_{1}}$:
\begin{equation}\nonumber
\int_{-\infty}^{\infty}e^{\pm
i(s-\lambda')(\xi'^{2}+1)}e^{ix\xi'}\phi_{\geq
s^{-\epsilon_{1}}}(\xi')\hat{h}(\xi')d\xi'
\end{equation}
Recalling \eqref{102} we now need to estimate the contribution of
\begin{equation}\nonumber\begin{split}
&\int_{-\infty}^{0}\la\bm
P_{<a}[\chi_{>0}(x)|\tilde{U}^{(s)}|^{2}](\lambda',.)P_{<a}[|\tilde{U}^{(s)}|^{2}](\lambda',.)\\-P_{<a}[\chi_{>0}(x)|\tilde{U}^{(s)}|^{2}](\lambda',.)P_{<a}[|\tilde{U}^{(s)}|^{2}](\lambda',.)\endm,
(e^{ix\xi}-e^{-ix\xi})\underline{e}\ra
\\&\hspace{2cm}\overline{\la \bm \chi_{>0}(x)\overline{\tilde{U}^{(s)}(\lambda',.)}\\
\chi_{>0}(x)\tilde{U}^{(s)}(\lambda',.)\endm \times
\int_{-\infty}^{\infty}e^{\pm
i(s-\lambda')(\xi'^{2}+1)}e^{ix\xi'}\phi_{\geq
s^{-\epsilon_{1}}}(\xi')\hat{h}(\xi')d\xi',(e^{ix\xi}-e^{-ix\xi})\underline{e}\ra}
d\xi,
\end{split}\end{equation}
This then gets substituted for
\begin{equation}\nonumber
\la
\bm|\tilde{U}^{(s)}|^{4}(\lambda',.)\\-|\tilde{U}^{(s)}|^{4}(\lambda',.)\endm_{dis},
\bm \overline{\tilde{U}^{(s)}(\lambda',.)}\\
\tilde{U}^{(s)}(\lambda',.)\endm \times
e^{-i(s-\lambda')\calH^{*}}(\calH^{*})^{-k}\psi_{dis}\ra
\end{equation}
in \eqref{BC}. It is now important to recall that in \eqref{BC} we
already reduced $e^{-i(s-\lambda)\calH^{*}}\tilde{\phi}_{dis}$ to
frequency $<s^{-\epsilon}$, see the paragraph before \eqref{BC}.
Thus writing $e^{-i(s-\lambda)\calH^{*}}\tilde{\phi}_{dis}$ in
terms of its Fourier expansion with variable $\xi$ and
re-arranging exponentials we get the phase
$e^{is(\xi^{2}-\xi'^{2})}$ in the case of destructive resonance
(and $e^{is(\xi^{2}+\xi'^{2}+2)}$for constructive resonance). If
we arrange $\epsilon_{1}>>\epsilon$, as we may, we have
$|\xi^{2}-\xi'^{2}|\gtrsim s^{-\epsilon_{1}}$. Then switch orders
of integration in \eqref{BC}, first performing an the integration
with respect to $s$. This costs $s^{\epsilon_{1}}$ but either
demolishes the integral in $\lambda'$ or else produces at least
extra factors
$\partial_{s}[\lambda_{\infty}[\mu-\mu_{\infty}]](s)$,
$\frac{d}{ds}(\Psi-\Psi_{\infty})_{1}\sim(\nu-1)(s)$. In order to
decouple the variables $\xi,\xi'$, notice that for $|\xi'|\lesssim
s^{\epsilon}$ we have\footnote{We may easily reduce to
$|\xi'|\lesssim s^{\epsilon}$, see e.g. the argument for
Lemma~\ref{nullform}}
\begin{equation}\nonumber
\phi_{<s^{-\epsilon}}(\xi)\phi_{\geq
s^{-\epsilon_{1}}}(\xi')\frac{e^{i(\xi^{2}-\xi'^{2})s}}{\xi^{2}-\xi'^{2}}=\sum_{n,m\in
s^{\epsilon}{\mathbf{Z}}^{2}}e^{i(\xi^{2}-\xi'^{2})s}a_{n,m}e^{in\xi+im\xi'},
\end{equation}
where one has $\sum_{n,m}[|n|+|m|]^{C}|a_{n,m}|\lesssim
s^{C\epsilon_{1}}$. If one substitutes this back into \eqref{BC}
and proceeds as in the proof of SLDE, one gets an extra gain in
$s$ upon choosing $\epsilon_{1}$ small enough, as desired. In
detail, consider
\begin{equation}\nonumber\begin{split}
&\int_{T}^{\infty}t\int_{t}^{\infty}(\nu-1)(s)\int_{0}^{\frac{s}{2}}\la
(1-e^{2i(\Psi-\Psi_{\infty})_{1}(\lambda)-2i(\Psi-\Psi_{\infty})_{1}(s)})U^{(s)}(\lambda,.)\phi,
e^{-i(s-\lambda)\calH^{*}}\tilde{\phi}_{dis}\ra d\lambda
\\&\times\int_{\frac{s}{2}}^{s-s^{\frac{1}{10}}}\int_{-\infty}^{0}\la\bm
P_{<a}[\chi_{>0}(x)|\tilde{U}^{(s)}|^{2}](\lambda',.)P_{<a}[|\tilde{U}^{(s)}|^{2}](\lambda',.)\\-P_{<a}[\chi_{>0}(x)|\tilde{U}^{(s)}|^{2}](\lambda',.)P_{<a}[|\tilde{U}^{(s)}|^{2}](\lambda',.)\endm,
(e^{ix\xi}-e^{-ix\xi})\underline{e}\ra
\\&\times\overline{\la \bm \chi_{>0}(x)\overline{\tilde{U}^{(s)}(\lambda',.)}\\
\chi_{>0}(x)\tilde{U}^{(s)}(\lambda',.)\endm \times
\int_{-\infty}^{\infty}e^{\pm
i(s-\lambda')(\xi'^{2}+1)}e^{ix\xi'}\phi_{\geq
s^{-\epsilon_{1}}}(\xi')\hat{h}(\xi')d\xi',
(e^{ix\xi}-e^{-ix\xi})\underline{e}\ra} d\xi d\lambda',
\end{split}\end{equation}
Then substitute \eqref{000} with a suitable frequency cutoff as
discussed there, which amongst similar terms results in
\begin{equation}\nonumber\begin{split}
&\int_{T}^{\infty}t\int_{t}^{\infty}(\nu-1)(s)\int_{0}^{\frac{s}{2}}\la
(1-e^{2i(\Psi-\Psi_{\infty})_{1}(\lambda)-2i(\Psi-\Psi_{\infty})_{1}(s)})U^{(s)}(\lambda,.)\phi,\\&\hspace{7cm}
\int_{-\infty}^{0}e^{\pm
i(s-\lambda)(\eta^{2}+1)}\phi_{<s^{-\epsilon}}(\eta)(e^{ix\eta}-e^{-ix\eta})\calF(\phi_{dis})\eta\ra
d\lambda
\\&\times\int_{\frac{s}{2}}^{s-s^{\frac{1}{10}}}\int_{-\infty}^{0}\la\bm
P_{<a}[\chi_{>0}(x)|\tilde{U}^{(s)}|^{2}](\lambda',.)P_{<a}[|\tilde{U}^{(s)}|^{2}](\lambda',.)\\-P_{<a}[\chi_{>0}(x)|\tilde{U}^{(s)}|^{2}](\lambda',.)P_{<a}[|\tilde{U}^{(s)}|^{2}](\lambda',.)\endm,
(e^{ix\xi}-e^{-ix\xi})\underline{e}\ra
\\&\times\overline{\la \bm \chi_{>0}(x)\overline{\tilde{U}^{(s)}(\lambda',.)}\\
\chi_{>0}(x)\tilde{U}^{(s)}(\lambda',.)\endm \times
\int_{-\infty}^{\infty}e^{\pm
i(s-\lambda')(\xi'^{2}+1)}e^{ix\xi'}\phi_{\geq
s^{-\epsilon_{1}}}(\xi')\hat{h}(\xi')d\xi',
(e^{ix\xi}-e^{-ix\xi})\underline{e}\ra} d\xi d\lambda',
\end{split}\end{equation}
Rewrite this as
\begin{equation}\nonumber
\int_{T}^{\infty}t\int_{t}^{\infty}(\nu-1)(s)\int_{-\infty}^{0}\int_{-\infty}^{0}\int_{-\infty}^{\infty}\int_{0}^{\infty}\int_{0}^{\infty}e^{is(\pm
\eta^{2}\pm\xi'^{2})}\Psi(s,\lambda,\lambda',\eta,\xi,\xi')
d\lambda d\lambda' d\xi' d\eta d\xi ds dt
\end{equation}
Now perform an integration by parts in $s$, after switching the
orders of integration, then restore the original order of
integration:
\begin{equation}\nonumber\begin{split}
&\int_{T}^{\infty}t\int_{t}^{\infty}(\nu-1)(s)\int_{-\infty}^{0}\int_{-\infty}^{0}\int_{-\infty}^{\infty}\int_{0}^{\infty}\int_{0}^{\infty}e^{is(\pm
\eta^{2}\pm\xi'^{2})}\Psi(s,\lambda,\lambda',\eta,\xi,\xi')d\eta
d\xi d\xi' d\lambda d\lambda'ds dt\\
&=-\int_{T}^{\infty}t(\nu-1)(t)\int_{-\infty}^{0}\int_{-\infty}^{0}\int_{-\infty}^{\infty}\int_{0}^{\infty}\int_{0}^{\infty}\frac{e^{it(\pm
\eta^{2}\pm\xi'^{2})}}{i(\pm
\eta^{2}\pm\xi'^{2})}\Psi(t,\lambda,\lambda',\eta,\xi,\xi')d\eta
d\xi d\xi' d\lambda d\lambda'ds dt\\
&+\int_{T}^{\infty}t\int_{t}^{\infty}\dot{\nu}(s)\int_{-\infty}^{0}\int_{-\infty}^{0}\int_{-\infty}^{\infty}\int_{0}^{\infty}\int_{0}^{\infty}\frac{e^{is(\pm
\eta^{2}\pm\xi'^{2})}}{i(\pm \eta^{2}\pm
\xi'^{2})}\Psi(s,\lambda,\lambda',\eta,\xi,\xi')d\eta
d\xi d\xi' d\lambda d\lambda'ds dt\\
&+\int_{T}^{\infty}t\int_{t}^{\infty}(\nu-1)(s)\int_{-\infty}^{0}\int_{-\infty}^{0}\int_{-\infty}^{\infty}\int_{0}^{\infty}\int_{0}^{\infty}\frac{e^{is(\pm
\eta^{2}\pm\xi'^{2})}}{i(\pm
\eta^{2}\pm\xi'^{2})}\partial_{s}\Psi(s,\lambda,\lambda',\eta,\xi,\xi')d\eta
d\xi d\xi' d\lambda d\lambda'ds dt
\end{split}\end{equation}
We have $|\xi'^{2}\pm\eta^{2}|\gtrsim s^{-\epsilon_{1}}$ on the
support of each integrand. If one then decouples the variables
$\xi'$, $\eta$ as outlined above and then proceeds as in the proof
of SLDE, one checks that each of these terms can be bounded by
$\lesssim T^{-\frac{1}{2}+\delta_{1}-\mu(\epsilon_{1})}$ upon
choosing $\epsilon_{1}$ small enough, which suffices. This
concludes the case $\lambda'>\frac{s}{2}$.
\\

($\lambda'<\frac{s}{2}$), still in case (BC). This is the
expression
\begin{equation}\label{BC'}\begin{split}
&\int_{T}^{\infty}t\int_{t}^{\infty}(\nu-1)(s)\int_{0}^{\frac{s}{2}}\la
(1-e^{2i(\Psi-\Psi_{\infty})_{1}(\lambda)-2i(\Psi-\Psi_{\infty})_{1}(s)})U^{(s)}(\lambda,.)\phi,
e^{-i(s-\lambda)\calH^{*}}\tilde{\phi}_{dis}\ra d\lambda
\\&\int_{0}^{\frac{s}{2}}\la
\bm|\tilde{U}^{(s)}|^{4}(\lambda',.)\\-|\tilde{U}^{(s)}|^{4}(\lambda',.)\endm_{dis},
\bm \overline{\tilde{U}^{(s)}(\lambda',.)}\\
\tilde{U}^{(s)}(\lambda',.)\endm \times
e^{-i(s-\lambda')\calH^{*}}(\calH^{*})^{-k}\psi_{dis}\ra
d\lambda',
\end{split}\end{equation}
We start again by reducing
\begin{equation}
\la
\chi_{>0}\bm|\tilde{U}^{(s)}|^{4}(\lambda',.)\\-|\tilde{U}^{(s)}|^{4}(\lambda',.)\endm,
\bm \overline{\tilde{U}^{(s)}(\lambda',.)}\\
\tilde{U}^{(s)}(\lambda',.)\endm \times
e^{-i(s-\lambda')\calH^{*}}(\calH^{*})^{-k}\psi_{dis}\ra,
\end{equation}
to
\begin{equation}\label{103}
\la
\chi_{>0}\bm|\tilde{U}^{(s)}|^{4}(\lambda',.)\\-|\tilde{U}^{(s)}|^{4}(\lambda,.)\endm_{dis},
\bm \overline{\tilde{U}^{(s)}(\lambda',.)}\\
\tilde{U}^{(s)}(\lambda',.)\endm \times
e^{-i(s-\lambda')\calH^{*}}(\calH^{*})^{-k}\psi_{dis}\ra,
\end{equation}
This follows from
\begin{equation}
\sum_{i}|\la
\chi_{>0}\la\bm|\tilde{U}^{(s)}|^{4}\\-|\tilde{U}^{(s)}|^{4}\endm,\xi_{i}\ra\eta_{i},
\bm \overline{\tilde{U}^{(s)}(\lambda',.)}\\
\tilde{U}^{(s)}(\lambda',.)\endm \times
e^{-i(s-\lambda')\calH^{*}}(\calH^{*})^{-k}\psi_{dis}\ra|\lesssim
\lambda'^{-5(\frac{3}{2}-\delta_{3})}(s-\lambda')^{-\frac{3}{2}},
\end{equation}
 Now
replace \eqref{103} by
\begin{equation}
\int_{-\infty}^{\infty}\calF[\chi_{>0}\bm|\tilde{U}^{(s)}|^{4}(\lambda',.)\\-|\tilde{U}^{(s)}|^{4}(\lambda,.)\endm](\xi)
\overline{\tilde{\calF}[\chi_{>0}(x)\bm \overline{\tilde{U}^{(s)}(\lambda',.)}\\
\tilde{U}^{(s)}(\lambda',.)\endm \times
e^{-i(s-\lambda')\calH^{*}}\psi_{dis}](\xi)}d\xi
\end{equation}
We shall again simplify the Fourier transform here: for example,
consider the contribution of
\begin{equation}
\int_{-\infty}^{\infty}\la\chi_{>0}\bm|\tilde{U}^{(s)}|^{4}(\lambda',.)\\-|\tilde{U}^{(s)}|^{4}(\lambda,.)\endm,\phi(x,\xi)\ra
\overline{\tilde{\calF}[\chi_{>0}(x)\bm \overline{\tilde{U}^{(s)}(\lambda',.)}\\
\tilde{U}^{(s)}(\lambda',.)\endm \times
e^{-i(s-\lambda')\calH^{*}}(\calH^{*})^{-k}\psi_{dis}](\xi)}d\xi
\end{equation}
Proceeding as in the proof of SLDE one bounds this by $\lesssim
\lambda'^{-4(\frac{3}{2}-\delta_{3})}\lambda'(s-\lambda')^{-\frac{3}{2}}$,
which is more than enough. Further, for example the contribution
of
\begin{equation}
\int_{-\infty}^{0}\la\chi_{>0}\bm|\tilde{U}^{(s)}|^{4}(\lambda',.)\\-|\tilde{U}^{(s)}|^{4}(\lambda,.)\endm,(1+r(-\xi))e^{-ix\xi}\underline{e}\ra
\overline{\tilde{\calF}[\chi_{>0}(x)\bm \overline{\tilde{U}^{(s)}(\lambda',.)}\\
\tilde{U}^{(s)}(\lambda',.)\endm \times
e^{-i(s-\lambda')\calH^{*}}(\calH^{*})^{-k}\psi_{dis}](\xi)}d\xi
\end{equation}
is treated like in the case $\lambda'>\frac{s}{2}$ (one doesn't
gain in $s$ but in $\lambda'$). Thus focusing on the more
difficult reflection part, we need to estimate the contribution of
\begin{equation}\label{104}\begin{split}
&\int_{-\infty}^{0}\la\chi_{>0}\bm|\tilde{U}^{(s)}|^{4}(\lambda',.)\\-|\tilde{U}^{(s)}|^{4}(\lambda,.)\endm,(e^{ix\xi}-e^{-ix\xi})\underline{e}\ra
\\&\times\overline{\la[\chi_{>0}(x)\bm \overline{\tilde{U}^{(s)}(\lambda',.)}\\
\tilde{U}^{(s)}(\lambda',.)\endm \times
e^{-i(s-\lambda')\calH^{*}}(\calH^{*})^{-k}\psi_{dis}],
(e^{ix\xi}-e^{-ix\xi})\underline{e}\ra} d\xi
\end{split}\end{equation}
Arguing as in the case $\lambda'>\frac{s}{2}$, we may reduce this
expression further to
\begin{equation}\nonumber\begin{split}
&\int_{-\infty}^{0}\la\chi_{>0}\bm
P_{<a}[\chi_{>0}|\tilde{U}^{(s)}|^{2}](\lambda',.)P_{<a}[|\tilde{U}^{(s)}|^{2}](\lambda',.)\\-P_{<a}[\chi_{>0}|\tilde{U}^{(s)}|^{2}](\lambda',.)P_{<a}[|\tilde{U}^{(s)}|^{2}](\lambda',.)\endm,(e^{ix\xi}-e^{-ix\xi})\underline{e}\ra
\\&\times\overline{\la (e^{ix\xi}-e^{-ix\xi})\underline{e}, [\chi_{>0}(x)\bm \overline{\tilde{U}^{(s)}(\lambda',.)}\\
\tilde{U}^{(s)}(\lambda',.)\endm \times
e^{-i(s-\lambda')\calH^{*}}(\calH^{*})^{-k}\psi_{dis}]\ra} d\xi
\end{split}\end{equation}
where $a=\lambda'^{-\frac{3}{4}+\epsilon}$. We shall next show
that we may localize the Fourier support of\\
$e^{-i(s-\lambda')\calH^{*}}(\calH^{*})^{-k}\psi_{dis}$ away from
zero, in which case we can conclude as in the case
$\lambda'>\frac{s}{2}$, exploiting the frequency separation in
order to perform an integration by parts in $s$. Recall from the
proof of the SLDE that in this case, we need to perform an extra
integration by parts in the frequency variable in order to obtain
the gain $(s-\lambda')^{-\frac{3}{2}}$. More precisely, in the
expression
\begin{equation}\nonumber
\la
\chi_{>0}\bm|\tilde{U}^{(s)}|^{4}\\-|\tilde{U}^{(s)}|^{4}\endm_{dis},
\bm \overline{\tilde{U}^{(s)}(\lambda',.)}\\
\tilde{U}^{(s)}(\lambda',.)\endm \times
e^{-i(s-\lambda')\calH^{*}}(\calH^{*})^{-k}\psi_{dis}\ra,
\end{equation}
we replace the expression
\begin{equation}\nonumber
\chi_{>0}(x)\bm \overline{\tilde{U}^{(s)}(\lambda',.)}\\
\tilde{U}^{(s)}(\lambda',.)\endm \times
e^{-i(s-\lambda')\calH^{*}}(\calH^{*})^{-k}\psi_{dis}
\end{equation}
by
\begin{equation}\label{105}
\chi_{>0}(x)\bm x\overline{\tilde{U}^{(s)}(\lambda',.)}\\
x\tilde{U}^{(s)}(\lambda',.)\endm \times
\frac{1}{s-\lambda'}\int_{-\infty}^{0}e^{-i(s-\lambda')(\xi^{2}+1)}(e^{ix\xi}+e^{-ix\xi})\underline{e}[\frac{\calF(\psi_{dis})(\xi)}{(\xi^{2}+1)^{k}\xi}]d\xi
\end{equation}
\begin{equation}\nonumber
\chi_{>0}(x)\bm \overline{\tilde{U}^{(s)}(\lambda',.)}\\
\tilde{U}^{(s)}(\lambda',.)\endm \times
\frac{1}{s-\lambda'}\int_{-\infty}^{0}e^{-i(s-\lambda')(\xi^{2}+1)}(e^{ix\xi}-e^{-ix\xi})\underline{e}\partial_{\xi}[\frac{\calF(\psi_{dis})(\xi)}{(\xi^{2}+1)^{k}\xi}]d\xi
\end{equation}
\begin{equation}\nonumber
\chi_{>0}(x)\bm x\overline{\tilde{U}^{(s)}(\lambda',.)}\\
x\tilde{U}^{(s)}(\lambda',.)\endm
\times\frac{1}{s-\lambda'}\int_{0}^{\infty}e^{-i(s-\lambda')(\xi^{2}+1)}s(\xi)e^{ix\xi}\underline{e}\frac{\calF\psi_{dis}(\xi)}{(\xi^{2}+1)^{k}\xi}d\xi
\end{equation}
\begin{equation}\nonumber
\chi_{>0}(x)\bm \overline{\tilde{U}^{(s)}(\lambda',.)}\\
\tilde{U}^{(s)}(\lambda',.)\endm
\times\frac{1}{s-\lambda'}\int_{0}^{\infty}e^{-i(s-\lambda')(\xi^{2}+1)}\partial_{\xi}s(\xi)e^{ix\xi}\underline{e}\frac{\calF\psi_{dis}(\xi)}{(\xi^{2}+1)^{k}\xi}d\xi
\end{equation}
as well as the expressions
\begin{equation}\nonumber
\chi_{>0}(x)\bm \overline{\tilde{U}^{(s)}(\lambda',.)}\\
\tilde{U}^{(s)}(\lambda',.)\endm \times
\frac{1}{s-\lambda'}\int_{-\infty}^{0}e^{-i(s-\lambda')(\xi^{2}+1)}e^{-ix\xi}\partial_{\xi}[1+r(-\xi)]\underline{e}[\frac{\calF(\psi_{dis})(\xi)}{(\xi^{2}+1)^{k}\xi}]d\xi
\end{equation}
\begin{equation}\nonumber
\chi_{>0}(x)\bm \overline{\tilde{U}^{(s)}(\lambda',.)}\\
\tilde{U}^{(s)}(\lambda',.)\endm \times
\frac{1}{s-\lambda'}\int_{-\infty}^{0}e^{-i(s-\lambda')(\xi^{2}+1)}\partial_{\xi}[\phi(x,\xi)[\frac{\calF(\psi_{dis})(\xi)}{(\xi^{2}+1)^{k}\xi}]]d\xi
\end{equation}
as well as similar terms which can be treated identically. The
last term but one here is equivalent to the last term but two for
all intents and purposes. Moreover, the last term can be treated
by the same argument as for the last term but two, so we shall now
consider the four terms after and including \eqref{105}. Start
with \eqref{105}: we have
\begin{equation}\nonumber
\int_{-\infty}^{0}e^{-i(s-\lambda')(\xi^{2}+1)}(e^{ix\xi}+e^{-ix\xi})\underline{e}[\frac{\calF(\psi_{dis})(\xi)}{(\xi^{2}+1)^{k}\xi}]d\xi
=\frac{1}{\sqrt{s-\lambda'}}\int_{-\infty}^{\infty}e^{\frac{(x-y)^{2}}{i(s-\lambda')}}g(y)dy,
\end{equation}
where
$g(y)=\int_{-\infty}^{0}(e^{iy\xi}+e^{-iy\xi})\underline{e}[\frac{\calF(\phi_{dis})(\xi)}{(\xi^{2}+1)^{k}\xi}]d\xi$
satisfies $|g(y)|\lesssim \la y\ra ^{-2}$, whence we may replace
it by $\tilde{g}(y)=\phi_{<\lambda'^{\frac{1}{10}}}(y)g(y)$. Now
we specialize this further and consider the contribution of
\begin{equation}\nonumber
h(y):=\phi_{<\lambda'^{\frac{1}{10}}}(y)\int_{-\infty}^{0}\phi_{<\lambda'^{-\epsilon_{1}}}(\xi)(e^{iy\xi}+e^{-iy\xi})\underline{e}[\frac{\calF(\phi_{dis})(\xi)}{(\xi^{2}+1)^{k}\xi}]d\xi
\end{equation}
By the Heisenberg principle, it has frequency in the interval
$[0,\lambda'^{-\frac{\epsilon_{1}}{2}}]$ up to errors of size
$\lambda'^{-N}$, which we may neglect. Now consider the bracket
\begin{equation}\nonumber
\la (e^{ix\xi}-e^{-ix\xi})\underline{e}, [\chi_{>0}(x)\bm x\overline{\tilde{U}^{(s)}(\lambda',.)}\\
x\tilde{U}^{(s)}(\lambda',.)\endm \times
\frac{1}{(\sqrt{s-\lambda'})^{3}}\int_{-\infty}^{\infty}e^{\frac{(x-y)^{2}}{i(s-\lambda')}}h(y)dy]\ra,
\end{equation}
where we have the restriction $|\xi|\lesssim
\lambda'^{-\frac{3}{4}+\epsilon}$. We can then rewrite this as
\begin{equation}\nonumber
\chi_{\lesssim\lambda'^{-\frac{3}{4}+\epsilon}}(\xi)\la (e^{ix\xi}-e^{-ix\xi})\underline{e}, [P_{<\lambda'^{-\frac{\epsilon_{1}}{2}}}[\bm \chi_{>0}(x)x\overline{\tilde{U}^{(s)}(\lambda',.)}\\
\chi_{>0}(x)x\tilde{U}^{(s)}(\lambda',.)\endm] \times
\frac{1}{(\sqrt{s-\lambda'})^{3}}\int_{-\infty}^{\infty}e^{\frac{(x-y)^{2}}{i(s-\lambda')}}h(y)dy]\ra,
\end{equation}
Then we use Littlewood-Paley dichotomy in order to get
\begin{equation}\nonumber\begin{split}
&P_{<\lambda'^{-\frac{\epsilon_{1}}{2}}}[\bm \chi_{>0}(x)x\overline{\tilde{U}^{(s)}(\lambda',.)}\\
\chi_{>0}(x)x\tilde{U}^{(s)}(\lambda',.)\endm]=P_{<\lambda'^{-\frac{\epsilon_{1}}{2}}}[\bm \chi_{>0}(x)P_{<\lambda'^{-\frac{\epsilon_{1}}{2}}+10}[x\overline{\tilde{U}^{(s)}}](\lambda',.)\\
\chi_{>0}(x)P_{<\lambda'^{-\frac{\epsilon_{1}}{2}}+10}[x\tilde{U}^{(s)}](\lambda',.)\endm]
\\&\hspace{5cm}+P_{<\lambda'^{-\frac{\epsilon_{1}}{2}}}[\bm P_{\geq\lambda'^{-\frac{\epsilon_{1}}{2}}}[\chi_{>0}(x)][P_{\geq\lambda'^{-\frac{\epsilon_{1}}{2}}+10}][x\overline{\tilde{U}^{(s)}}](\lambda',.)\\
P_{\geq\lambda'^{-\frac{\epsilon_{1}}{2}}}[\chi_{>0}(x)]P_{\geq\lambda'^{-\frac{\epsilon_{1}}{2}}+10}[x\tilde{U}^{(s)}](\lambda',.)\endm]
\end{split}\end{equation}
We first consider the contribution from the 2nd term on the right.
We substitute this back into \eqref{104} undo the Fourier
transform using Plancherel's Theorem, and estimate this by
\begin{equation}\nonumber\begin{split}
&\lesssim ||\chi_{>0}\bm
P_{<a}[\chi_{>0}|\tilde{U}^{(s)}|^{2}](\lambda',.)P_{<a}[|\tilde{U}^{(s)}|^{2}](\lambda',.)\\-P_{<a}[\chi_{>0}|\tilde{U}^{(s)}|^{2}](\lambda',.)P_{<a}[|\tilde{U}^{(s)}|^{2}](\lambda',.)\endm||_{L_{x}^{\infty}}
||P_{\geq\lambda'^{-\frac{\epsilon_{1}}{2}}}[\chi_{>0}(x)]||_{L_{x}^{1}}\\&||x\tilde{U}^{(s)}(\lambda',.)||_{L_{x}^{\infty}}||\frac{1}{(\sqrt{s-\lambda'})^{3}}\int_{-\infty}^{\infty}e^{\frac{(x-y)^{2}}{i(s-\lambda')}}h(y)dy||_{L_{x}^{\infty}}
\lesssim
\lambda'^{-\frac{3}{2}+\epsilon_{1}}(s-\lambda')^{-\frac{3}{2}},
\end{split}\end{equation}
which is then seen to lead to an acceptable contribution upon
substitution into \eqref{BC'}. Thus we may now focus on the
contribution of
\begin{equation}\nonumber
\la (e^{ix\xi}-e^{-ix\xi})\underline{e}, [\chi_{>0}(x)\bm P_{<\lambda'^{-\frac{\epsilon_{1}}{2}+10}}[x\overline{\tilde{U}^{(s)}}](\lambda',.)\\
P_{<\lambda'^{-\frac{\epsilon_{1}}{2}+10}}[x\tilde{U}^{(s)}](\lambda',.)\endm
\times
\frac{1}{(\sqrt{s-\lambda'})^{3}}\int_{-\infty}^{\infty}e^{\frac{(x-y)^{2}}{i(s-\lambda')}}h(y)dy]\ra,
\end{equation}
always keeping in mind that
$|\xi|<\lambda'^{-\frac{3}{4}+\epsilon}$. We now replicate the
proof of SLDE for the low-low case (keep in mind that the full
expression we estimate is \eqref{104}). Thus we write
($p=-i\partial_{x}$)
\begin{equation}\nonumber
P_{<\lambda'^{-\frac{\epsilon_{1}}{2}+10}}[x\tilde{U}^{(s)}](t,.)=P_{<\lambda'^{-\frac{\epsilon_{1}}{2}+10}}[(x+2ip
t)\tilde{U}^{(s)}](t,.)-P_{<\lambda'^{-\frac{\epsilon_{1}}{2}+10}}[2ip
t\tilde{U}^{(s)}](t,.)
\end{equation}
Then we have
\begin{equation}\nonumber
[i\partial_{t}+\triangle]P_{<\lambda'^{-\frac{\epsilon_{1}}{2}+10}}\nabla
U(t,.)=P_{<\lambda'^{-\frac{\epsilon_{1}}{2}+10}}\nabla[
VU+\ldots+|U|^{4}U]
\end{equation}
just as in the proof of the strong local dispersive estimate. Now
one further manipulates the expressions on the right just as in
the proof of SLDE. Note that the operator
$P_{<\lambda'^{-\frac{\epsilon_{1}}{2}+10}}\nabla$ will smear out
the supports a bit, but this is easily seen to be harmless. Of
course one gains $\lambda'^{-\frac{\epsilon_{1}}{2}}$ in the
process, which overcomes any small losses in the proof of SLDE.
One can now restrict to $|\xi|>\lambda'^{-\epsilon_{1}}$, i. e.
include a multiplier $\phi_{\geq \lambda'^{-\epsilon_{1}}}$ in the
definition of $h(y)$, and then finish the argument just as in the
case $\lambda'>\frac{s}{2}$. This concludes estimating the
contribution from the term \eqref{105}. The contribution of the
third term in that list is treated analogously. We now turn to the
contribution of the 2nd term there, i. e. the expression
\begin{equation}\nonumber\begin{split}
&\la
\chi_{>0}\bm|\tilde{U}^{(s)}|^{4}\\-|\tilde{U}^{(s)}|^{4}\endm_{dis}(\lambda',.),\\
&\chi_{>0}(x)\bm \overline{\tilde{U}^{(s)}}(\lambda',.)\\
\tilde{U}^{(s)}(\lambda',.)\endm \times
\frac{1}{s-\lambda'}\int_{-\infty}^{0}e^{-i(s-\lambda')(\xi^{2}+1)}(e^{ix\xi}-e^{-ix\xi})\underline{e}\partial_{\xi}[\frac{\calF(\psi_{dis})(\xi)}{(\xi^{2}+1)^{k}\xi}]d\xi\ra
\end{split}\end{equation}
But this is easily seen to be estimable by
\begin{equation}\nonumber
\lesssim \lambda'^{-\frac{3}{2}}(s-\lambda')^{-\frac{3}{2}},
\end{equation}
which upon substitution into \eqref{BC'} yields an acceptable
contribution. The fourth term after \eqref{105} is handled
analogously. We are done with the case (BC). Clearly the case (AC)
can be handled analogously.
\\

(CC): the most difficult case. This is the expression
\begin{equation}\label{CC}\begin{split}
&\int_{T}^{\infty}t\int_{t}^{\infty}(\nu-1)(s)\int_{0}^{s}\la
\bm|\tilde{U}^{(s)}|^{4}\tilde{U}^{(s)}(\lambda,.)\\-|\tilde{U}^{(s)}|^{4}\overline{\tilde{U}^{(s)}(\lambda,.)}\endm_{dis},
e^{-i(s-\lambda)\calH^{*}}\phi_{dis}\ra d\lambda
\\&\hspace{4cm}\int_{0}^{s}\la
\bm|\tilde{U}^{(s)}|^{4}\tilde{U}^{(s)}(\lambda',.)\\-|\tilde{U}^{(s)}|^{4}\overline{\tilde{U}^{(s)}(\lambda',.)}\endm_{dis},
e^{-i(s-\lambda')\calH^{*}}(\calH^{*})^{-k}\psi_{dis}\ra d\lambda'
\end{split}\end{equation}
Start with the case $\max\{\lambda, \lambda'\}<\frac{s}{2}$. We
may restrict integration to the range $\lambda>\lambda'$.
Rearrange either of the factors in the integral as
\begin{equation}\nonumber
\la
\chi_{>0}\bm|\tilde{U}^{(s)}|^{4}(\lambda',.)\\-|\tilde{U}^{(s)}|^{4}(\lambda',.)\endm,
\bm \overline{\tilde{U}^{(s)}(\lambda',.)}\\
\tilde{U}^{(s)}(\lambda',.)\endm \times
e^{-i(s-\lambda')\calH^{*}}(\calH^{*})^{-k}\psi_{dis}\ra,
\end{equation}
\begin{equation}\nonumber
\la
\chi_{>0}\bm|\tilde{U}^{(s)}|^{4}(\lambda,.)\\-|\tilde{U}^{(s)}|^{4}(\lambda,.)\endm,
\bm \overline{\tilde{U}^{(s)}(\lambda,.)}\\
\tilde{U}^{(s)}(\lambda,.)\endm \times
e^{-i(s-\lambda)\calH^{*}}\phi_{dis}\ra,
\end{equation}
As usual we first need to reduce both factors in either bracket to
their dispersive part. This time, though, we have to analyze each
constituent more carefully, since they all interact with each
other. Thus we now write
\begin{equation}\nonumber
\chi_{>0}\bm|\tilde{U}^{(s)}|^{4}(\lambda',.)\\-|\tilde{U}^{(s)}|^{4}(\lambda',.)\endm
=\sum_{i}a_{i}\la\chi_{>0}\bm|\tilde{U}^{(s)}|^{4}(\lambda',.)\\-|\tilde{U}^{(s)}|^{4}(\lambda',.)\endm,
\xi_{k(i)}\ra\eta_{i}+\bm|\chi_{>0}\tilde{U}^{(s)}|^{4}(\lambda',.)\\-\chi_{>0}|\tilde{U}^{(s)}|^{4}(\lambda',.)\endm_{dis}:=\alpha+\beta
\end{equation}
First consider the contribution from $\alpha(\lambda',.)$, i. e.
the expression
\begin{equation}\nonumber
\la \alpha(\lambda',.), \bm \overline{\tilde{U}^{(s)}(\lambda',.)}\\
\tilde{U}^{(s)}(\lambda',.)\endm \times
e^{-i(s-\lambda')\calH^{*}}(\calH^{*})^{-k}\psi_{dis}\ra
\end{equation}
As usual we expand
\begin{equation}\nonumber
e^{-i(s-\lambda')\calH^{*}}(\calH^{*})^{-k}\psi_{dis}=\sum_{\pm}\int_{-\infty}^{\infty}e^{\pm
i(s-\lambda')(\xi^{2}+1)}e_{\pm}(x,\xi)\frac{\calF_{\pm}(\psi_{dis})}{(\xi^{2}+1)^{k}}d\xi
\end{equation}
We claim that we may sneak in a smooth cutoff
$\phi_{<s^{-\epsilon}}(\xi)$ into this integrand, which we then
denote as
$e^{-i(s-\lambda')\calH^{*}}(\calH^{*})^{-k}\tilde{\psi}_{dis}$.
This is because integration by parts in $\xi$ costs in addition to
$s^{\epsilon}$ at most $\max\{|x|,\xi\}$, and for $\alpha$ we may
assume $|x|$ to be bounded, whence choosing $\epsilon$ small
enough results in a gain in $s$. Of course we use
$\lambda'<\frac{s}{2}$ here. The same comment applies to
$\alpha(\lambda,.)$. Our strategy shall be to achieve a
localization away from zero for the frequencies of
$e^{-i(s-\lambda)\calH^{*}}\phi_{dis}$,
$e^{-i(s-\lambda)\calH^{*}}(\calH^{*})^{-k}\psi_{dis}$, occuring
in the contribution from $\beta(\lambda,.)$, $\beta(\lambda',.)$.
This ensures that $\alpha(\lambda',.)$ and $\beta(\lambda,.)$ etc
interact weakly. We now distinguish between the following cases:
\\

($\alpha\alpha$): this is the expression
\begin{equation}\label{alphaalpha}\begin{split}
&\int_{T}^{\infty}t\int_{t}^{\infty}(\nu-1)(s)\int_{0}^{\frac{s}{2}}\la
\la\chi_{>0}\bm|\tilde{U}^{(s)}|^{4}(\lambda,.)\\-|\tilde{U}^{(s)}|^{4}(\lambda,.)\endm,
\xi_{i}\ra\eta_{i} , \bm \overline{\tilde{U}^{(s)}}(\lambda,.)\\
\tilde{U}^{(s)}(\lambda,.)\endm\times
e^{-i(s-\lambda)\calH^{*}}\tilde{\phi}_{dis}\ra d\lambda
\\&\hspace{4cm}\int_{0}^{\frac{s}{2}}\la
\la\chi_{>0}\bm|\tilde{U}^{(s)}|^{4}(\lambda',.)\\-|\tilde{U}^{(s)}|^{4}(\lambda',.)\endm,
\xi_{i}\ra\eta_{i} , \bm \overline{\tilde{U}^{(s)}}(\lambda',.)\\
\tilde{U}^{(s)}(\lambda',.)\endm\times
e^{-i(s-\lambda')\calH^{*}}(\calH^{*})^{-k}\tilde{\psi}_{dis}\ra
d\lambda'
\end{split}\end{equation}
We can easily estimate this by
\begin{equation}\nonumber
\lesssim\int_{T}^{\infty}t\int_{t}^{\infty}s^{-\frac{1}{2}+\delta_{1}}\int_{0}^{\frac{s}{2}}\lambda^{-5(\frac{3}{2}-\delta_{3})}(s-\lambda)^{-\frac{3}{2}}d\lambda
\int_{0}^{\frac{s}{2}}\lambda'^{-5(\frac{3}{2}-\delta_{3})}(s-\lambda')^{-\frac{3}{2}}d\lambda'\lesssim
T^{-\frac{1}{2}+\delta_{1}}
\end{equation}
\\

($\alpha\beta$): the expression
\begin{equation}\label{alphabeta}\begin{split}
&\int_{T}^{\infty}t\int_{t}^{\infty}(\nu-1)(s)\int_{0}^{\frac{s}{2}}\la
\la\chi_{>0}\bm|\tilde{U}^{(s)}|^{4}(\lambda,.)\\-|\tilde{U}^{(s)}|^{4}(\lambda,.)\endm,
\xi_{i}\ra\eta_{i} , \bm \overline{\tilde{U}^{(s)}}(\lambda,.)\\
\tilde{U}^{(s)}(\lambda,.)\endm\times
e^{-i(s-\lambda)\calH^{*}}\tilde{\phi}_{dis}\ra d\lambda
\\&\hspace{4cm}\int_{0}^{\frac{s}{2}}\la \bm\chi_{>0}|\tilde{U}^{(s)}|^{4}(\lambda',.)\\-\chi_{>0}|\tilde{U}^{(s)}|^{4}(\lambda',.)\endm_{dis}
, \bm \overline{\tilde{U}^{(s)}}(\lambda',.)\\
\tilde{U}^{(s)}(\lambda',.)\endm\times
e^{-i(s-\lambda')\calH^{*}}(\calH^{*})^{-k}\psi_{dis}\ra d\lambda'
\end{split}\end{equation}
To proceed, we restrict the frequency of
$e^{-i(s-\lambda')\calH^{*}}\psi_{dis}$ away from $0$. The
procedure for this is identical to the one outlined in case (BC).
Having achieved frequency separation, we have of course achieved
rapid oscillation in $s$, whence we can close this case like at
the end of case (BC), by integration by parts in $s$. The case
$(\beta\alpha)$ is handled analogously.
\\

($\beta\beta$): the expression
\begin{equation}\label{betabeta}\begin{split}
&\int_{T}^{\infty}t\int_{t}^{\infty}(\nu-1)(s)\int_{0}^{\frac{s}{2}}\la
\la \bm\chi_{>0}|\tilde{U}^{(s)}|^{4}(\lambda,.)\\-\chi_{>0}|\tilde{U}^{(s)}|^{4}(\lambda,.)\endm_{dis}, \bm \overline{\tilde{U}^{(s)}}(\lambda,.)\\
\tilde{U}^{(s)}(\lambda,.)\endm\times
e^{-i(s-\lambda)\calH^{*}}\phi_{dis}\ra d\lambda
\\&\hspace{3cm}\int_{0}^{\frac{s}{2}}\la \bm\chi_{>0}|\tilde{U}^{(s)}|^{4}(\lambda',.)\\-\chi_{>0}|\tilde{U}^{(s)}|^{4}(\lambda',.)\endm_{dis}
, \bm \overline{\tilde{U}^{(s)}}(\lambda',.)\\
\tilde{U}^{(s)}(\lambda',.)\endm\times
e^{-i(s-\lambda')\calH^{*}}(\calH^{*})^{-k}\psi_{dis}\ra d\lambda'
\end{split}\end{equation}
As before we mimic the proof of SLDE. Thus we perform an
integration by parts in the Fourier representation for
$e^{-i(s-\lambda)\calH^{*}}\phi_{dis}$ etc and produce the
following list of terms provided the integration by parts results
in a loss of $x$. Call this list $\beta_{1}$:
\begin{equation}\label{beta1}
\chi_{>0}(x)\bm x\overline{\tilde{U}^{(s)}(\lambda',.)}\\
x\tilde{U}^{(s)}(\lambda',.)\endm \times
\frac{1}{s-\lambda'}\int_{-\infty}^{0}e^{-i(s-\lambda')(\xi^{2}+1)}(e^{ix\xi}+e^{-ix\xi})\underline{e}[\frac{\calF(\psi_{dis})(\xi)}{(\xi^{2}+1)^{k}\xi}]d\xi
\end{equation}
\begin{equation}\nonumber
\chi_{>0}(x)\bm x\overline{\tilde{U}^{(s)}(\lambda',.)}\\
x\tilde{U}^{(s)}(\lambda',.)\endm
\times\frac{1}{s-\lambda'}\int_{0}^{\infty}e^{-i(s-\lambda')(\xi^{2}+1)}s(\xi)e^{ix\xi}\underline{e}\frac{\calF\psi_{dis}(\xi)}{(\xi^{2}+1)^{k}\xi}d\xi
\end{equation}
\begin{equation}
\chi_{>0}(x)\bm x\overline{\tilde{U}^{(s)}(\lambda,.)}\\
x\tilde{U}^{(s)}(\lambda,.)\endm \times
\frac{1}{s-\lambda}\int_{-\infty}^{0}e^{-i(s-\lambda)(\xi^{2}+1)}(e^{ix\xi}+e^{-ix\xi})\underline{e}[\frac{\calF(\phi_{dis})(\xi)}{\xi}]d\xi
\end{equation}
\begin{equation}\nonumber
\chi_{>0}(x)\bm x\overline{\tilde{U}^{(s)}(\lambda,.)}\\
x\tilde{U}^{(s)}(\lambda,.)\endm
\times\frac{1}{s-\lambda}\int_{0}^{\infty}e^{-i(s-\lambda)(\xi^{2}+1)}s(\xi)e^{ix\xi}\underline{e}\frac{\calF\phi_{dis}(\xi)}{\xi}d\xi
\end{equation}
These get complemented by the following terms, which we refer to
as $\beta_{2}$:
\begin{equation}\label{beta2}
\chi_{>0}(x)\bm \overline{\tilde{U}^{(s)}(\lambda',.)}\\
\tilde{U}^{(s)}(\lambda',.)\endm \times
\frac{1}{s-\lambda'}\int_{-\infty}^{0}e^{-i(s-\lambda')(\xi^{2}+1)}(e^{ix\xi}-e^{-ix\xi})\underline{e}\partial_{\xi}[\frac{\calF(\psi_{dis})(\xi)}{(\xi^{2}+1)^{k}\xi}]d\xi
\end{equation}
\begin{equation}\nonumber
\chi_{>0}(x)\bm \overline{\tilde{U}^{(s)}(\lambda',.)}\\
\tilde{U}^{(s)}(\lambda',.)\endm
\times\frac{1}{s-\lambda'}\int_{0}^{\infty}e^{-i(s-\lambda')(\xi^{2}+1)}\partial_{\xi}s(\xi)e^{ix\xi}\underline{e}\frac{\calF\psi_{dis}(\xi)}{(\xi^{2}+1)^{k}\xi}d\xi
\end{equation}
\begin{equation}\nonumber
\chi_{>0}(x)\bm \overline{\tilde{U}^{(s)}(\lambda,.)}\\
\tilde{U}^{(s)}(\lambda,.)\endm \times
\frac{1}{s-\lambda}\int_{-\infty}^{0}e^{-i(s-\lambda)(\xi^{2}+1)}(e^{ix\xi}-e^{-ix\xi})\underline{e}\partial_{\xi}[\frac{\calF(\phi_{dis})(\xi)}{\xi}]d\xi
\end{equation}
\begin{equation}\nonumber
\chi_{>0}(x)\bm \overline{\tilde{U}^{(s)}(\lambda,.)}\\
\tilde{U}^{(s)}(\lambda,.)\endm
\times\frac{1}{s-\lambda}\int_{0}^{\infty}e^{-i(s-\lambda)(\xi^{2}+1)}\partial_{\xi}s(\xi)e^{ix\xi}\underline{e}\frac{\calF\phi_{dis}(\xi)}{\xi}d\xi
\end{equation}
\\

($\beta_{2}\beta_{2}$): this is an expression of the form
\begin{equation}\nonumber
\int_{T}^{\infty}t\int_{t}^{\infty}(\nu-1)(s)\int_{0}^{\frac{s}{2}}\la
\bm\chi_{>0}|\tilde{U}^{(s)}|^{4}(\lambda',.)\\-\chi_{>0}|\tilde{U}^{(s)}|^{4}(\lambda',.)\endm_{dis},
\beta_{2}(\lambda',.)\ra d\lambda' \int_{0}^{\frac{s}{2}}\la
\bm\chi_{>0}|\tilde{U}^{(s)}|^{4}(\lambda,.)\\-\chi_{>0}|\tilde{U}^{(s)}|^{4}(\lambda,.)\endm_{dis},\beta_{2}(\lambda,.)\ra
d\lambda
\end{equation}
where $\beta_{2}(\lambda',.), \beta_{2}(\lambda,.)$ stand for
certain terms of the 2nd list. This type of interaction is easy to
control: one bounds this by
\begin{equation}\nonumber
\lesssim
\int_{T}^{\infty}t\int_{t}^{\infty}s^{-\frac{1}{2}+\delta_{1}}\int_{0}^{\frac{s}{2}}\lambda'^{-\frac{3}{2}}(s-\lambda')^{-\frac{3}{2}}d\lambda'
\int_{0}^{\lambda}\lambda^{-\frac{3}{2}}(s-\lambda)^{-\frac{3}{2}}d\lambda
\lesssim T^{-\frac{1}{2}+\delta_{1}}
\end{equation}
\\

($\beta_{1}\beta_{2}$): this is an expression of the form
\begin{equation}\nonumber
\int_{T}^{\infty}t\int_{t}^{\infty}(\nu-1)(s)\int_{0}^{\frac{s}{2}}\la
\bm\chi_{>0}|\tilde{U}^{(s)}|^{4}(\lambda,.)\\-\chi_{>0}|\tilde{U}^{(s)}|^{4}(\lambda,.)\endm_{dis},
\beta_{1}(\lambda,.)\ra d\lambda \int_{0}^{\frac{s}{2}}\la
\bm\chi_{>0}|\tilde{U}^{(s)}|^{4}(\lambda',.)\\-\chi_{>0}|\tilde{U}^{(s)}|^{4}(\lambda',.)\endm_{dis},\beta_{2}(\lambda',.)\ra
d\lambda'
\end{equation}
Assume for example(the other cases being treated by exact analogy)
that $\beta_{2}(\lambda',.)$ has the following form:
\begin{equation}\nonumber
\beta_{2}(\lambda',.)=\chi_{>0}(x)\bm \overline{\tilde{U}^{(s)}(\lambda',.)}\\
\tilde{U}^{(s)}(\lambda',.)\endm
\times\frac{1}{s-\lambda'}\int_{0}^{\infty}e^{-i(s-\lambda')(\xi^{2}+1)}\partial_{\xi}s(\xi)e^{ix\xi}\underline{e}\frac{\calF\phi_{dis}(\xi)}{\xi}d\xi
\end{equation}
Note that on account of
\begin{equation}\nonumber
|\la\bm\chi_{>0}|\tilde{U}^{(s)}|^{4}(\lambda',.)\\-\chi_{>0}|\tilde{U}^{(s)}|^{4}(\lambda',.)\endm_{dis},\beta_{2}(\lambda',.)\ra|\lesssim
\lambda'^{-\frac{3}{2}}(s-\lambda')^{-\frac{3}{2}},
\end{equation}
we may assume that $\lambda'<s^{\epsilon}$ for a small
$\epsilon>0$. But then on account of the pseudo-conformal almost
conservation we may apply a localizer $\phi_{<s^{2\epsilon}}(x)$
to the quadrilinear term: indeed, we have
\begin{equation}\nonumber\begin{split}
&|\la\bm\chi_{>0}\phi_{\geq
s^{2\epsilon}}(x)|\tilde{U}^{(s)}|^{4}(\lambda',.)\\-\chi_{>0}\phi_{\geq
s^{2\epsilon}}(x)|\tilde{U}^{(s)}|^{4}(\lambda',.)\endm_{dis},\beta_{2}(\lambda',.)\ra|\\&\lesssim
||\frac{1}{|x|}\phi_{\geq
s^{2\epsilon}}(|x|)[(x+2ip\lambda')U-2ip\lambda'
U]||_{L_{x}^{2}}||U(\lambda,.)||_{L_{x}^{2}}||U(\lambda',.)||_{L_{x}^{\infty}}^{3}(s-\lambda')^{-\frac{3}{2}}
\\&\lesssim
s^{-\epsilon+\delta_{2}}\lambda'^{-\frac{3}{2}}(s-\lambda')^{-\frac{3}{2}},\\
\end{split}\end{equation}
which leads to an acceptable contribution above.  Finally, we may
reduce the frequency $\xi$ in the relation defining
$\beta_{2}(\lambda',.)$ above to absolute size $<s^{-\epsilon}$ by
inclusion of a suitable smooth cutoff
$\phi_{<s^{-\epsilon}}(\xi)$. This is since upon including a
smooth cutoff $\phi_{\geq s^{-\epsilon}}(\xi)$ for suitably small
$\epsilon$ results in an expression which can be integrated by
parts in $\xi$, resulting in losses of at most
$\max\{|x|,s^{\epsilon}\}s^{\epsilon}$ for each integration while
resulting in a gain of $s-\lambda'$. Choosing $\epsilon$ small
enough results in arbitrary gains in $s$. Next, we consider
$\beta_{1}(\lambda,.)$. Using the same argument as in case (BC),
we reduce the frequency to size $>\lambda^{-\epsilon_{1}}$. But
then we have again achieved frequency separation and can integrate
by parts in $s$. The case $(\beta_{2}\beta_{1})$ is simpler, as
one gains $\lambda^{-\frac{1}{2}}$ which suffices (since
$\lambda>\lambda'$).
\\

($\beta_{1}\beta_{1}$): This is an expression of the form
\begin{equation}\nonumber
\int_{T}^{\infty}t\int_{t}^{\infty}(\nu-1)(s)\int_{0}^{\frac{s}{2}}\la
\bm\chi_{>0}|\tilde{U}^{(s)}|^{4}(\lambda',.)\\-\chi_{>0}|\tilde{U}^{(s)}|^{4}(\lambda',.)\endm_{dis},
\beta_{1}(\lambda',.)\ra d\lambda' \int_{0}^{\frac{s}{2}}\la
\bm\chi_{>0}|\tilde{U}^{(s)}|^{4}(\lambda,.)\\-\chi_{>0}|\tilde{U}^{(s)}|^{4}(\lambda,.)\endm_{dis},\beta_{1}(\lambda,.)\ra
d\lambda
\end{equation}
Keep in mind that we assume $\lambda>\lambda'$. Use the distorted
Plancherel's Theorem~\ref{DistortedPlancherel} to rewrite this as
\begin{equation}\nonumber\begin{split}
&\sum_{\pm,\pm}\int_{T}^{\infty}t\int_{t}^{\infty}(\nu-1)(s)\int_{0}^{\frac{s}{2}}\int_{-\infty}^{\infty}
\calF_{\pm}\bm\chi_{>0}|\tilde{U}^{(s)}|^{4}(\lambda',.)\\-\chi_{>0}|\tilde{U}^{(s)}|^{4}(\lambda',.)\endm(\xi)
\overline{\tilde{\calF}_{\pm}[\beta_{1}(\lambda',.)]}d\xi
d\lambda'\\&
\hspace{2cm}\int_{0}^{\frac{s}{2}}\int_{-\infty}^{\infty}
\calF_{\pm}\bm\chi_{>0}|\tilde{U}^{(s)}|^{4}(\lambda,.)\\-\chi_{>0}|\tilde{U}^{(s)}|^{4}(\lambda,.)\endm(\xi')\overline{\tilde{\calF}_{\pm}[\beta_{1}(\lambda,.)](\xi')}
d\xi' d\lambda
\end{split}\end{equation}
We may and shall restrict to the $+$ case and omit the subscript,
and restrict both the $\xi$ and $\xi'$ integral to the range
$(-\infty,0]$, the other case being similar but simpler. We then
need to decompose each of the Fourier transforms $\calF(...)$ etc
into various constituents, i. e. write
\begin{equation}\nonumber\begin{split}
&\calF\bm\chi_{>0}|\tilde{U}^{(s)}|^{4}(\lambda',.)\\-\chi_{>0}|\tilde{U}^{(s)}|^{4}(\lambda',.)\endm(\xi)
=\la
\bm\chi_{>0}|\tilde{U}^{(s)}|^{4}(\lambda',.)\\-\chi_{>0}|\tilde{U}^{(s)}|^{4}(\lambda',.)\endm,
\phi(x,\xi)\ra
\\&+\la\bm\chi_{>0}|\tilde{U}^{(s)}|^{4}(\lambda',.)\\-\chi_{>0}|\tilde{U}^{(s)}|^{4}(\lambda',.)\endm,
(1+r(-\xi))\underline{e}e^{-ix\xi}\ra +\la
\bm\chi_{>0}|\tilde{U}^{(s)}|^{4}(\lambda',.)\\-\chi_{>0}|\tilde{U}^{(s)}|^{4}(\lambda',.)\endm,
(e^{ix\xi}-e^{-ix\xi})\underline{e}\ra
\end{split}\end{equation}
We shall consider the contribution from the first and third term,
the 2nd being treated similarly to the third. Moreover, performing
the same decomposition for $\tilde{\calF}[\beta_{1}(\lambda',.)]$
as well as
$\calF\bm\chi_{>0}|\tilde{U}^{(s)}|^{4}(\lambda,.)\\-\chi_{>0}|\tilde{U}^{(s)}|^{4}(\lambda,.)\endm(\xi')$,
it is easy to see that we may restrict to the contribution from
the third term, as the others are simpler. We commence with the
following expression:
\begin{equation}\label{108}\begin{split}
&\sum_{\pm,\pm}\int_{T}^{\infty}t\int_{t}^{\infty}(\nu-1)(s)\\&\int_{0}^{\frac{s}{2}}\int_{-\infty}^{0}
\la\bm\chi_{>0}|\tilde{U}^{(s)}|^{4}(\lambda',.)\\-\chi_{>0}|\tilde{U}^{(s)}|^{4}(\lambda',.)\endm,
(e^{ix\xi}-e^{-ix\xi})\underline{e}\ra \overline{\la
\beta_{1}(\lambda',.), (e^{ix\xi}-e^{-ix\xi})\underline{e}\ra}
d\xi d\lambda'\\&
\hspace{2cm}\int_{0}^{\frac{s}{2}}\int_{-\infty}^{0}
\la\bm\chi_{>0}|\tilde{U}^{(s)}|^{4}(\lambda,.)\\-\chi_{>0}|\tilde{U}^{(s)}|^{4}(\lambda,.)\endm,
(e^{ix\xi'}-e^{-ix\xi'})\underline{e}\ra\overline{\la\beta_{1}(\lambda,.),
(e^{ix\xi'}-e^{-ix\xi'})\underline{e}\ra} d\xi'd\lambda
\end{split}\end{equation}
If we recapitulate the proof of SLDE for both bracket factors, we
see that we may reduce to estimating
\begin{equation}\nonumber\begin{split}
&\sum_{\pm,\pm}\int_{T}^{\infty}t\int_{t}^{\infty}(\nu-1)(s)\\&\int_{0}^{\frac{s}{2}}\int_{-\infty}^{0}
\la\bm
P_{<a'}[\chi_{>0}|\tilde{U}^{(s)}|^{2}]P_{<a'}[|\tilde{U}^{(s)}|^{2}](\lambda',.)\\-P_{<a'}[\chi_{>0}|\tilde{U}^{(s)}|^{2}]P_{<a'}[|\tilde{U}^{(s)}|^{2}](\lambda',.)\endm,
(e^{ix\xi}-e^{-ix\xi})\underline{e}\ra \overline{\la
\beta_{1}(\lambda',.), (e^{ix\xi}-e^{-ix\xi})\underline{e}\ra}
d\xi d\lambda'\\& \int_{0}^{\frac{s}{2}}\int_{-\infty}^{0} \la \bm
P_{<a}[\chi_{>0}|\tilde{U}^{(s)}|^{2}]P_{<a}[|\tilde{U}^{(s)}|^{2}](\lambda,.)\\-P_{<a}[\chi_{>0}|\tilde{U}^{(s)}|^{2}]P_{<a}[|\tilde{U}^{(s)}|^{2}](\lambda,.)\endm,
(e^{ix\xi'}-e^{-ix\xi'})\underline{e}\ra\overline{\la\beta_{1}(\lambda,.),
(e^{ix\xi'}-e^{-ix\xi'})\underline{e}\ra} d\xi'd\lambda
\end{split}\end{equation}
where $a=\lambda^{-\frac{3}{4}+\epsilon_{2}}$,
$a'=\lambda'^{-\frac{3}{4}}\lambda^{\epsilon_{2}}$,
$\epsilon_{2}=\epsilon_{2}(\delta_{2})$. Recalling the product
representation of $\beta_{1}(\lambda',.)$ as in \eqref{beta1}, we
first reduce the frequency of the right hand integral factor of
both $\beta_{1}(\lambda,.),\beta_{1}(\lambda',.)$ to size
$>\lambda^{-\epsilon}$. This is achieved as in case (BC). For
technical reasons we shall effect this by means of a sharp cutoff
$\chi_{>\lambda^{-\epsilon}}(\xi)$ etc. Thus for
example\footnote{The same argument applies to all terms of the
list $\beta_{1}$.} we shall put
\begin{equation}\nonumber
\beta_{1}(\lambda,x)=\bm x\overline{\tilde{U}^{(s)}(\lambda,.)}\\
x\tilde{U}^{(s)}(\lambda,.)\endm \times\chi_{>0}(x)\frac{e^{\pm
i(s-\lambda)}}{(s-\lambda)^{\frac{3}{2}}}\int_{-\infty}^{\infty}e^{\frac{\pm(x-y)^{2}}{i(s-\lambda)}}g(y)dy,
\end{equation}
where
\begin{equation}\nonumber
g(y)=\phi_{<\lambda^{\frac{1}{10}}}(y)\int_{-\infty}^{0}\chi_{>\lambda^{-\epsilon}}(\xi)(e^{iy\xi}+e^{-iy\xi})\underline{e}[\frac{\calF(\phi_{dis})(\xi)}{\xi}]d\xi
\end{equation}
where $\chi_{>\lambda^{-\epsilon}}(\xi)$ is a Heavyside function.
Of course we have
\begin{equation}\nonumber
\frac{1}{\sqrt{s-\lambda}}\int_{-\infty}^{\infty}e^{\frac{\pm(x-y)^{2}}{i(s-\lambda)}}g(y)dy
=\int_{-\infty}^{\infty}e^{\mp
i(s-\lambda)\xi^{2}}e^{ix\xi}\hat{g}(\xi)d\xi
\end{equation}
Similar observations apply to $\beta_{1}(\lambda',.)$, for example
\begin{equation}\nonumber
\beta_{1}(\lambda',x)=\bm x\overline{\tilde{U}^{(s)}(\lambda',.)}\\
x\tilde{U}^{(s)}(\lambda',.)\endm \times\chi_{>0}(x)\frac{e^{\mp
i(s-\lambda')}}{(s-\lambda')^{\frac{3}{2}}}\int_{-\infty}^{\infty}e^{\pm\frac{(x-y)^{2}}{i(s-\lambda)}}\tilde{g}(y)dy,
\end{equation}
where
\begin{equation}\nonumber
\tilde{g}(y)=\phi_{<\lambda^{\frac{1}{10}}}(y)\int_{0}^{\infty}\chi_{>\lambda^{-\epsilon}}(\xi')s(\xi')e^{iy\xi'}\underline{e}\frac{\calF\phi_{dis}(\xi')}{(\xi'^{2}+1)^{k}\xi'}d\xi'
\end{equation}
We now further specialize the frequency support of $g(y)$,
$\tilde{g}(y)$, by including cutoffs
$\chi_{I_{i}}(\xi),\,\chi_{I_{j}}(\xi')$ corresponding to
intervals $I_{i,j}$ of length $\lambda^{-\epsilon}$, i. e.
introduce
\begin{equation}\nonumber
g_{i}(y)=\phi_{<\lambda^{\frac{1}{10}}}(y)\int_{-\infty}^{0}\chi_{>\lambda^{-\epsilon}}(\xi)\chi_{I_{i}}(\xi)(e^{iy\xi}+e^{-iy\xi})\underline{e}[\frac{\calF(\phi_{dis})(\xi)}{\xi}]d\xi
\end{equation}
\begin{equation}\nonumber
\tilde{g}_{j}(y)=\phi_{<\lambda^{\frac{1}{10}}}(y)\int_{0}^{\infty}\chi_{>\lambda^{-\epsilon}}(\xi')s(\xi')e^{iy\xi'}\chi_{j}(\xi')\underline{e}\frac{\calF\phi_{dis}(\xi')}{(\xi'^{2}+1)^{k}\xi'}d\xi'
\end{equation}
Clearly if $|i-j|>>1$ these functions have separated Fourier
supports (of distance $\gtrsim \lambda^{-\epsilon}$) up to errors
of order $\lambda^{-N}$, hence negligible. Now introduce
$\beta_{1,i}(\lambda',.)$, $\beta_{1,j}(\lambda,.)$ exactly as
above with $g(y)$, $\tilde{g}(y)$ replaced by $g_{i}(y)$,
$\tilde{g}_{j}(y)$. It is easy to see that we can restrict both
$|\xi|, |\xi'|$ to size $<\lambda^{\epsilon_{3}}$ for
$\epsilon_{3}=\epsilon_{3}(\delta_{2})$, since otherwise one gains
enough to overcome any losses in the proof of SLDE.
\\

Case 1: $|i-j|>>1$. Here we exploit integration by parts in $s$.
Write
\begin{equation}\nonumber
\beta_{1,i}(\lambda,x)=\bm x\overline{\tilde{U}^{(s)}(\lambda,.)}\\
x\tilde{U}^{(s)}(\lambda,.)\endm \times\chi_{>0}(x)\frac{e^{\pm
i(s-\lambda)}}{(s-\lambda)^{\frac{3}{2}}}\int_{-\infty}^{\infty}e^{\mp
i(s-\lambda)\xi^{2}}e^{ix\xi}\widehat{g_{i}}(\xi)d\xi
\end{equation}
and similarly for $\beta_{1,j}(\lambda',.)$. Then re-write
\begin{equation}\label{109}\begin{split}
&\sum_{\pm,\pm}\int_{T}^{\infty}t\int_{t}^{\infty}(\nu-1)(s)\\&\int_{0}^{\frac{s}{2}}\int_{-\infty}^{0}
\la\bm
P_{<a'}[\chi_{>0}|\tilde{U}^{(s)}|^{2}]P_{<a'}[|\tilde{U}^{(s)}|^{2}](\lambda',.)\\-P_{<a'}[\chi_{>0}|\tilde{U}^{(s)}|^{2}]P_{<a'}[|\tilde{U}^{(s)}|^{2}](\lambda',.)\endm,
(e^{ix\xi}-e^{-ix\xi})\underline{e}\ra \overline{\la
\beta_{1,j}(\lambda',.), (e^{ix\xi}-e^{-ix\xi})\underline{e}\ra}
d\xi d\lambda'\\& \int_{0}^{\frac{s}{2}}\int_{-\infty}^{0} \la \bm
P_{<a}[\chi_{>0}|\tilde{U}^{(s)}|^{2}]P_{<a}[|\tilde{U}^{(s)}|^{2}](\lambda,.)\\-P_{<a}[\chi_{>0}|\tilde{U}^{(s)}|^{2}]P_{<a}[|\tilde{U}^{(s)}|^{2}](\lambda,.)\endm,
(e^{ix\xi'}-e^{-ix\xi'})\underline{e}\ra\overline{\la\beta_{1,i}(\lambda,.),
(e^{ix\xi'}-e^{-ix\xi'})\underline{e}\ra} d\xi'd\lambda
\end{split}\end{equation}
as
\begin{equation}\nonumber
\int_{T}^{\infty}t\int_{t}^{\infty}(\nu-1)(s)\int_{-\infty}^{0}\int_{-\infty}^{0}\int_{0}^{\infty}\int_{0}^{\infty}e^{is(\pm\xi^{2}\pm\xi'^{2})}\Psi_{ij}(s,\lambda,\lambda',\xi,\xi')d\xi
d\xi'd\lambda d\lambda'ds dt
\end{equation}
switch the order of integration, integrate by parts in $s$,
decouple the variables $\xi,\xi'$ by means of discrete Fourier
transform and proceed as in the proof of the SLDE. Choosing
$\epsilon>0$ small enough results in a gain in $\lambda$, even
upon summing over $i,j$. This concludes Case 1.
\\

Case 2: $i=j+O(1)$. First write
\begin{equation}\nonumber
\beta_{1,i}=2\sqrt{-1}\beta^{a}_{1,i}+\beta^{b}_{1,i},
\end{equation}
where
\begin{equation}\nonumber
\beta^{a}_{1,i}=\lambda\bm \overline{\partial_{x}\tilde{U}^{(s)}(\lambda,.)}\\
-\partial_{x}\tilde{U}^{(s)}(\lambda,.)\endm
\times\chi_{>0}(x)\frac{e^{\pm
i(s-\lambda)}}{(s-\lambda)^{\frac{3}{2}}}\int_{-\infty}^{\infty}e^{\frac{\pm(x-y)^{2}}{i(s-\lambda)}}g_{i}(y)dy
\end{equation}
\begin{equation}\nonumber
\beta^{b}_{1,i}=\bm \overline{(x+2\lambda p)\tilde{U}^{(s)}(\lambda,.)}\\
(x+2\lambda p)\tilde{U}^{(s)}(\lambda,.)\endm
\times\chi_{>0}(x)\frac{e^{\pm
i(s-\lambda)}}{(s-\lambda)^{\frac{3}{2}}}\int_{-\infty}^{\infty}e^{\frac{\pm(x-y)^{2}}{i(s-\lambda)}}g_{i}(y)dy
\end{equation}
From the proof of SLDE, recall that we use integration by parts to
write
\begin{equation}\label{110}\begin{split}
&\la \beta^{(a)}_{1,i}(\lambda,.),
(e^{ix\xi}-e^{-ix\xi})\underline{e}\ra\\&=i\lambda\xi\int_{0}^{\infty}(e^{ix\xi}+e^{-ix\xi})\la\underline{e},
\int_{x}^{\infty}\chi_{>0}(y)\bm
\overline{\partial_{y}\tilde{U}^{(s)}}(\lambda,y)\\
-\partial_{y}\tilde{U}^{(s)}(\lambda,y)\endm\times\frac{e^{\pm
i(s-\lambda)}}{(s-\lambda)^{\frac{3}{2}}}\int_{-\infty}^{\infty}e^{\frac{\pm(y-z)^{2}}{i(s-\lambda)}}g_{i}(z)dz
dy
\end{split}\end{equation}
Still following the proof of SLDE in the low-low case, we then use
the free parametrix to write schematically
\begin{equation}\nonumber
\overline{\tilde{U}^{(s)}}(\lambda,y)=\int_{0}^{\lambda}\frac{1}{\sqrt{\lambda-\mu}}\int_{-\infty}^{\infty}e^{\frac{(y-z')^{2}}{i(\lambda-\mu)}}\partial_{z'}[|U|^{4}U(\mu,z')+VU(\mu,z')]dz'
d\mu
\end{equation}
We then break this into the contributions from
$\chi_{><\lambda^{\frac{1}{2}+\epsilon_{3}}}(\mu)\chi_{><(\lambda-\mu)^{\frac{1}{2}+\epsilon_{3}}}(|y|)|U|^{4}U(\mu,z')$,\\
$\chi_{><\lambda-\lambda^{2\epsilon_{3}}}(\mu)\chi_{><\lambda^{\epsilon_{3}}}(z')VU(\mu,z')$.
Indeed, from the argument in the proof of SLDE, it follows that if
one $>$ sign is chosen in these cutoffs, the corresponding
contribution leads to a small extra gain in $\lambda$, which then
suffices to close, provided
$\epsilon_{3}=\epsilon_{3}(\delta_{2})$. Thus we now choose
everywhere the $<$ sign, and substitute this into $\eqref{110}$.
Collecting the exponentials, we encounter the following phase
function, just as in the proof of SLDE:
\begin{equation}\nonumber\begin{split}
&e^{\frac{y^{2}}{i}[\frac{1}{\lambda-\mu}+\frac{1}{s-\lambda}]-\frac{2x}{i}[\frac{z'}{\lambda-\mu}+\frac{z}{s-\lambda}]}e^{\frac{z'^{2}}{i(\lambda-\mu)}+\frac{z^{2}}{i(s-\lambda)}}\\
&=e^{-i(y\sqrt{\frac{1}{\lambda-\mu}+\frac{1}{s-\lambda}}-[\frac{z'}{\lambda-\mu}+\frac{z}{s-\lambda}]\sqrt{\frac{1}{\lambda-\mu}+\frac{1}{s-\lambda}}^{-1})^{2}}
e^{\frac{z'^{2}}{i(\lambda-\mu)}+\frac{z^{2}}{i(s-\lambda)}+i[\frac{z'}{\lambda-\mu}+\frac{z}{s-\lambda}]^{2}[\frac{1}{\lambda-\mu}+\frac{1}{s-\lambda}]^{-1}}\\
&:=e^{-i(y\sqrt{\frac{1}{\lambda-\mu}+\frac{1}{s-\lambda}}-y_{1})^{2}}e^{iy_{2}}
\end{split}\end{equation}
where we have
\begin{equation}\nonumber
|y_{1,2}|=|y_{1,2}(z,z',\lambda,\mu,s)|\lesssim
\lambda^{\epsilon_{3}}
\end{equation}
on the support of the integrand in \eqref{110}. Thus plugging this
into \eqref{110} and omitting the integration in $\mu',\mu$ for
now, we obtain the expression
\begin{equation}\nonumber\begin{split}
&\frac{1}{\sqrt{(s-\lambda)(\lambda-\mu)}}\int_{-\infty}^{\infty}\int_{-\infty}^{\infty}\int_{0}^{\infty}
(e^{ix\xi}+e^{-ix\xi})\int_{x}^{\infty}e^{-i(y\sqrt{\frac{1}{\lambda-\mu}+\frac{1}{s-\lambda}}-y_{1})^{2}}e^{iy_{2}}
g_{i}(z)g_{1}(\mu,z')dy dx dz dz'\\
&=\frac{[\frac{1}{s-\lambda}+\frac{1}{\lambda-\mu}]^{-1}}{\sqrt{(s-\lambda)(\lambda-\mu)}}\int_{-\infty}^{\infty}\int_{-\infty}^{\infty}e^{iy_{2}}
\int_{-y_{1}}^{\infty}[e^{i\xi[\tilde{x}+y_{1}]\sqrt{\frac{1}{s-\lambda}+\frac{1}{\lambda-\mu}}^{-1}}+e^{-i\xi[\tilde{x}+y_{1}]\sqrt{\frac{1}{s-\lambda}+\frac{1}{\lambda-\mu}}^{-1}}]\int_{\tilde{x}}^{\infty}e^{i\rho^{2}}d\rho
d\tilde{x}\\
&\hspace{4cm}\phi_{<\lambda^{\frac{1}{10}}}(z)\int_{-\infty}^{0}\chi_{>\lambda^{-\epsilon}}(\eta)\chi_{I_{i}}(\eta)(e^{iz\eta}+e^{-iz\eta})\underline{e}[\frac{\calF(\phi_{dis})(\eta)}{\eta}]d\eta
g_{1}(z',\mu) dz dz'
\end{split}\end{equation}
Now assume $I_{i}=[a_{i},b_{i}]$ where $\min\{|a_{i}|,
|b_{i}|\}\geq\lambda^{-\epsilon}$. Then integrate by parts in
\begin{equation}\nonumber\begin{split}
&\int_{-\infty}^{0}\chi_{>\lambda^{-\epsilon}}(\eta)\chi_{I_{i}}(\eta)(e^{iz\eta}+e^{-iz\eta})\underline{e}[\frac{\calF(\phi_{dis})(\eta)}{\eta}]d\eta\\
&=\frac{e^{iza_{i}}\frac{\calF(\phi_{dis})(a_{i})}{a_{i}}-e^{izb_{i}}\frac{\calF(\phi_{dis})(b_{i})}{b_{i}}}{iz}
-\frac{1}{iz}\int_{-\infty}^{0}\chi_{>\lambda^{-\epsilon}}(\eta)\chi_{I_{i}}(\eta)(e^{iz\eta}-e^{-iz\eta})\underline{e}\partial_{\eta}[\frac{\calF(\phi_{dis})(\eta)}{\eta}]d\eta
\end{split}\end{equation}
Observe that
\begin{equation}\nonumber
||\phi_{<\lambda^{\frac{1}{10}}}(z)\frac{1}{iz}\int_{-\infty}^{0}\chi_{>\lambda^{-\epsilon}}(\eta)\chi_{I_{i}}(\eta)(e^{iz\eta}-e^{-iz\eta})\underline{e}\partial_{\eta}[\frac{\calF(\phi_{dis})(\eta)}{\eta}]d\eta||_{L_{z}^{1}}
\lesssim \log\lambda\lambda^{-\epsilon},
\end{equation}
whence this expression has a negligible contribution upon
continuing the proof of SLDE and choosing
$\epsilon_{3}<<\epsilon$. We further observe that the restriction
$|\xi|<\lambda^{-\frac{3}{4}+\epsilon}$ as well as the
restrictions on $|z|$, $|z'|$, $\lambda$ and $\mu$ specified
further above imply that
\begin{equation}\nonumber
|\partial_{z}[e^{iy_{2}}
\int_{-y_{1}}^{\infty}[e^{i\xi[\tilde{x}+y_{1}]\sqrt{\frac{1}{s-\lambda}+\frac{1}{\lambda-\mu}}^{-1}}+e^{-i\xi[\tilde{x}+y_{1}]\sqrt{\frac{1}{s-\lambda}+\frac{1}{\lambda-\mu}}^{-1}}]\int_{\tilde{x}}^{\infty}e^{i\rho^{2}}d\rho
d\tilde{x}]|\lesssim \lambda^{-\frac{1}{2}+\epsilon_{3}}
\end{equation}
\begin{equation}\nonumber
|[e^{iy_{2}}
\int_{-y_{1}}^{\infty}[e^{i\xi[\tilde{x}+y_{1}]\sqrt{\frac{1}{s-\lambda}+\frac{1}{\lambda-\mu}}^{-1}}+e^{-i\xi[\tilde{x}+y_{1}]\sqrt{\frac{1}{s-\lambda}+\frac{1}{\lambda-\mu}}^{-1}}]\int_{\tilde{x}}^{\infty}e^{i\rho^{2}}d\rho
d\tilde{x}]|\lesssim \lambda^{\epsilon_{3}}
\end{equation}
Thus for all intents and purposes we can replace the latter
function by a constant as far as its dependence of $z$ on
$[-\lambda^{\frac{1}{10}}, \lambda^{\frac{1}{10}}]$ is concerned.
But then one calculates
\begin{equation}\nonumber
|\int_{-\infty}^{\infty}\phi_{<\lambda^{\frac{1}{10}}}(z)\frac{e^{iza_{i}}\frac{\calF(\phi_{dis})(a_{i})}{a_{i}}-e^{izb_{i}}\frac{\calF(\phi_{dis})(b_{i})}{b_{i}}}{iz}dz|
\lesssim \lambda^{-\epsilon}
\end{equation}
Putting everything after \eqref{110} together, we see that
\begin{equation}\nonumber
|\la \beta^{(a)}_{1,i}(\lambda,.),
(e^{ix\xi}-e^{-ix\xi})\underline{e}\ra|\lesssim
\lambda\lambda^{-\frac{3}{4}+\epsilon_{2}}\lambda^{\frac{1}{2}}\log\lambda\lambda^{-\epsilon}(s-\lambda)^{-\frac{3}{2}},
\end{equation}
which then yields
\begin{equation}\nonumber\begin{split}
&|\int_{-\infty}^{0} \la \bm
P_{<a}[\chi_{>0}|\tilde{U}^{(s)}|^{2}]P_{<a}[|\tilde{U}^{(s)}|^{2}](\lambda,.)\\-P_{<a}[\chi_{>0}|\tilde{U}^{(s)}|^{2}]P_{<a}[|\tilde{U}^{(s)}|^{2}](\lambda,.)\endm,
(e^{ix\xi'}-e^{-ix\xi'})\underline{e}\ra\overline{\la\beta^{(a)}_{1,i}(\lambda,.),
(e^{ix\xi'}-e^{-ix\xi'})\underline{e}\ra} d\xi'|\\&\lesssim
\lambda^{-1}\lambda^{-\frac{3}{4}+\epsilon_{2}}\lambda\lambda^{-\frac{3}{4}+\epsilon_{2}}\lambda^{\frac{1}{2}}\log\lambda\lambda^{-\epsilon}(s-\lambda)^{-\frac{3}{2}}
\lesssim
\log\lambda\lambda^{-1+2\epsilon_{2}-\epsilon}(s-\lambda)^{-\frac{3}{2}}
\end{split}\end{equation}
We shall choose $0<\epsilon_{2}<<\epsilon$. Analogously to
\eqref{110}, we also need to estimate the contribution of
\begin{equation}\nonumber
\la \beta^{(b)}_{1,i}(\lambda,.),
(e^{ix\xi}-e^{-ix\xi})\underline{e}\ra
\end{equation}
This, however, is more elementary, as we can estimate
\begin{equation}\nonumber\begin{split}
&||\la (e^{ix\xi}-e^{-ix\xi})\underline{e}, \bm \overline{(x+2\lambda p)\tilde{U}^{(s)}(\lambda,.)}\\
(x+2\lambda p)\tilde{U}^{(s)}(\lambda,.)\endm
\times\chi_{>0}(x)\frac{e^{\pm
i(s-\lambda)}}{(s-\lambda)^{\frac{3}{2}}}\int_{-\infty}^{\infty}e^{\frac{\pm(x-y)^{2}}{i(s-\lambda)}}g_{i}(y)dy||_{L_{\xi}^{2}}
\\&\lesssim ||(x+2\lambda
p)\tilde{U}^{(s)}(\lambda,.)||_{L_{x}^{2}}(s-\lambda)^{-\frac{3}{2}}\lesssim
(s-\lambda)^{-\frac{3}{2}},
\end{split}\end{equation} whence we get
\begin{equation}\nonumber\begin{split}
&|\int_{-\infty}^{0} \la \bm
P_{<a}[\chi_{>0}|\tilde{U}^{(s)}|^{2}]P_{<a}[|\tilde{U}^{(s)}|^{2}](\lambda',.)\\-P_{<a}[\chi_{>0}|\tilde{U}^{(s)}|^{2}]P_{<a}[|\tilde{U}^{(s)}|^{2}](\lambda',.)\endm,
(e^{ix\xi'}-e^{-ix\xi'})\underline{e}\ra\overline{\la\beta^{(a)}_{1,i}(\lambda,.),
(e^{ix\xi'}-e^{-ix\xi'})\underline{e}\ra} d\xi'|\\
&\lesssim \lambda^{-\frac{3}{2}}(s-\lambda)^{-\frac{3}{2}},
\end{split}\end{equation}
which then leads to an acceptable contribution. Now of course we
eventually need to estimate the expression \eqref{109} under our
current assumption $i=j+O(1)$, and then sum over $i$. One can
replicate the preceding arguments for $\beta_{1,j}(\lambda',.)$ as
long as $\lambda'>\lambda^{10\epsilon}$, say. But we can exclude
the opposite case, since if $\lambda'\leq \lambda^{10\epsilon}$,
we can restrict $|U(\lambda',x)$ to the range
$|x|<\lambda^{20\epsilon}$, say, using pseudo-conformal almost
conservation, and then we can restrict
$e^{-i(s-\lambda')\calH^{*}}(\calH^{*})^{-k}\phi_{dis}$ to
frequency $<s^{-10\epsilon}$, say, provided $\epsilon$ is small
enough. Thus we can reduce this situation to the separated
frequency case. Finally, if $\lambda'>\lambda^{10\epsilon}$, we
get
\begin{equation}\nonumber\begin{split}
&|\int_{-\infty}^{0} \la \bm
P_{<a}[\chi_{>0}|\tilde{U}^{(s)}|^{2}]P_{<a}[|\tilde{U}^{(s)}|^{2}](\lambda,.)\\-P_{<a}[\chi_{>0}|\tilde{U}^{(s)}|^{2}]P_{<a}[|\tilde{U}^{(s)}|^{2}](\lambda,.)\endm,
(e^{ix\xi'}-e^{-ix\xi'})\underline{e}\ra\overline{\la\beta^{(a)}_{1,i}(\lambda,.),
(e^{ix\xi'}-e^{-ix\xi'})\underline{e}\ra} d\xi'|\\& \times
|\int_{-\infty}^{0} \la \bm
P_{<a}[\chi_{>0}|\tilde{U}^{(s)}|^{2}]P_{<a}[|\tilde{U}^{(s)}|^{2}](\lambda',.)\\-P_{<a}[\chi_{>0}|\tilde{U}^{(s)}|^{2}]P_{<a}[|\tilde{U}^{(s)}|^{2}](\lambda',.)\endm,
(e^{ix\xi'}-e^{-ix\xi'})\underline{e}\ra\overline{\la\beta^{(b)}_{1,j}(\lambda,.),
(e^{ix\xi'}-e^{-ix\xi'})\underline{e}\ra} d\xi'|\\&\lesssim
\lambda^{-1}\lambda^{-\frac{3}{4}+\epsilon_{2}}\lambda\lambda^{-\frac{3}{4}+\epsilon_{2}}\lambda^{\frac{1}{2}}\log\lambda\lambda^{-\epsilon}(s-\lambda)^{-\frac{3}{2}}
\lambda'^{-1}\lambda'^{-\frac{3}{4}}\lambda^{\epsilon_{2}}\lambda'\lambda'^{-\frac{3}{4}}\lambda^{\epsilon_{2}}\lambda'^{\frac{1}{2}}\log\lambda\lambda^{-\epsilon}(s-\lambda')^{-\frac{3}{2}}
\\&\lesssim
(\log\lambda)^{2}\lambda^{-1+4\epsilon_{2}-2\epsilon}\lambda'^{-1}(s-\lambda)^{-\frac{3}{2}}(s-\lambda')^{-\frac{3}{2}}
\end{split}\end{equation}
Upon summing over $O(\lambda^{\epsilon+\epsilon_{3}})$ indices
$i=j+O(1)$, this leads to a small gain provided
$\epsilon_{2}+\epsilon_{3}<<\epsilon$, which we may arrange. This
concludes the treatment of \eqref{108}.
\\

We next consider the contribution of the term
\begin{equation}\label{111}\begin{split}
&\sum_{\pm,\pm}\int_{T}^{\infty}t\int_{t}^{\infty}(\nu-1)(s)\\&\int_{0}^{\frac{s}{2}}\int_{-\infty}^{0}
\la\bm\chi_{>0}|\tilde{U}^{(s)}|^{4}(\lambda',.)\\-\chi_{>0}|\tilde{U}^{(s)}|^{4}(\lambda',.)\endm,
\phi(x,\xi)\ra \overline{\la \beta_{1}(\lambda',.),
(e^{ix\xi}-e^{-ix\xi})\underline{e}\ra} d\xi d\lambda'\\&
\hspace{2cm}\int_{0}^{\frac{s}{2}}\int_{-\infty}^{0}
\la\bm\chi_{>0}|\tilde{U}^{(s)}|^{4}(\lambda,.)\\-\chi_{>0}|\tilde{U}^{(s)}|^{4}(\lambda,.)\endm,
(e^{ix\xi'}-e^{-ix\xi'})\underline{e}\ra\overline{\la\beta_{1}(\lambda,.),
(e^{ix\xi'}-e^{-ix\xi'})\underline{e}\ra} d\xi'd\lambda
\end{split}\end{equation}
We observe that we may innocuously include cutoffs
$\phi_{<>\lambda'\lambda^{\epsilon}}(|x|)$
simultaneously\footnote{I. e. either both have the $<$ or the $>$
subscript.} in front of
\begin{equation}\nonumber
\bm\chi_{>0}|\tilde{U}^{(s)}|^{4}(\lambda',.)\\-\chi_{>0}|\tilde{U}^{(s)}|^{4}(\lambda',.)\endm,\,
\beta_{1}(\lambda',.)
\end{equation}
Including the cutoff $\phi_{>\lambda'\lambda^{\epsilon}}(|x|)$
clearly leads to the desired extra gain, while including the
cutoff $\phi_{<\lambda'\lambda^{\epsilon}}(|x|)$ allows us to
restrict $e^{-i(s-\lambda')\calH^{*}}(\calH^{*})^{-k}\psi_{dis}$
in $\beta_{1,\lambda'}(\lambda',.)$ to small frequency. As we we
can always restrict the frequency of
$e^{-i(s-\lambda)\calH}\phi_{dis}$ in $\beta_{1}(\lambda,.)$ away
from zero, we have then achieved frequency separation and can
argue as before in case (BC). This concludes the case (CC)
provided we have $\max\{\lambda,\lambda'\}<\frac{s}{2}$. The case
$\max\{\lambda,\lambda'\}\geq \frac{s}{2}$ is more elementary and
omitted. We are now done with the proof of Lemma~\ref{general
bilinear}, whence also the Lemma before it.
\end{proof}

We now continue with the proof of Proposition~\ref{Hard}; note
that the expressions
\begin{equation}\nonumber
\int_{T}^{\infty}t\int_{t}^{\infty}(\nu-1)^{a}(s)\la \bm
\tilde{U}\\\overline{\tilde{U}}\endm_{dis}(s,.),\phi\ra ds
dt,\,\int_{T}^{\infty}t\int_{t}^{\infty}\dot{\nu}(s)\la \bm
\tilde{U}\\\overline{\tilde{U}}\endm_{dis}(s,.),\phi\ra ds
dt,\,a\geq 1
\end{equation}
can be treated exactly like above, using Lemma~\ref{general
bilinear} and the relation \eqref{nu}. Thus up to terms which can
either be estimated using Lemma~\ref{general bilinear} or else can
even be absolutely integrated, there is only one potentially
troublesome expression, namely
\begin{equation}\nonumber
\int_{T}^{\infty}t\int_{t}^{\infty}\la
\tilde{U}^{2}-\overline{\tilde{U}}^{2}, \phi\ra ds dt
\end{equation}
where $\phi$ is an even time-independent Schwartz function. We
handle this by the following Lemma, which hinges on a
{\it{symplectic structure}}:

\begin{lemma}\label{nullform} Let $\Gamma\in A^{(n)}_{[0,T)}$, $\Gamma=\{\bm \tilde{U}\\ \overline{\tilde{U}}\endm_{dis},\ldots\}$. Then, for $\tilde{T}\leq T$ and $\phi$ an even Schwartz function, we have
\begin{equation}\nonumber
\int_{\tilde{T}}^{T}t\int_{t}^{T}\la
\tilde{U}_{dis}^{2}(s,.)-\overline{\tilde{U}_{dis}}^{2}(s,.),
\phi\ra ds dt \lesssim \tilde{T}^{-\frac{1}{2}+\delta_{1}},
\end{equation}
provided $\delta_{1}$ is large enough in relation to $\delta_{2},
\delta_{3}$. Similarly, we have the bound
\begin{equation}\nonumber
\int_{\tilde{T}}^{T}t^{2}\la
\tilde{U}_{dis}^{2}(t,.)-\overline{\tilde{U}_{dis}}^{2}(t,.),
\phi\ra dt \lesssim \tilde{T}^{-\frac{1}{2}+\delta_{1}},
\end{equation}
Both bounds are uniform in $T$.
\end{lemma}
\begin{proof} It relies on identifying a special cancellation in this expression. We treat the first expression
in detail, the 2nd following the same reasoning.  Also, we may put
$T=\infty$. To begin with, write
\begin{equation}\label{Fourier}
\bm \tilde{U}^{(s)}\\ \overline{
\tilde{U}^{(s)}}\endm_{dis}=\sum_{\pm}\int_{-\infty}^{\infty}e_{\pm}(x,\xi)\calF_{\pm}\bm
\tilde{U}^{(s)}\\ \overline{ \tilde{U}^{(s)}}\endm(\xi) d\xi
\end{equation}
This gets substituted into
\begin{equation}\label{null-form}
\la
(\tilde{U}^{(s)}_{dis})^{2}(s,.)-(\overline{\tilde{U}^{(s)}_{dis}})^{2}(s,.),
\phi\ra=\la\bm \tilde{U}^{(s)}\\ \overline{
\tilde{U}^{(s)}}\endm_{dis}^{t}\bm 1&0\\0&-1\endm \bm
\tilde{U}^{(s)}\\ \overline{ \tilde{U}^{(s)}}\endm_{dis}, \phi\ra
\end{equation}
The key here is that
\begin{equation}\nonumber
\la e_{+}(x,\xi)\bm 1&0\\0&-1\endm e_{-}(x,\xi'), \phi(x)\ra=\la
e_{+}(x,\xi)\bm 0&1\\-1& 0\endm e_{+}(x,\xi'), \phi(x)\ra,
\end{equation}
which is easily seen to vanish for $\xi=\xi'$. Moreover, this also
vanishes for $\xi=-\xi'$, since then it is the inner product of an
odd and an even function (use that
$e_{\pm}(x,-\xi)=e_{\pm}(-x,\xi)$). To proceed, we now substitute
the Duhamel expression for $\bm \tilde{U}^{(s)}\\ \overline{
\tilde{U}^{(s)}}\endm_{dis}$ into \eqref{null-form}. We shall then
treat the most difficult term which results when we substitute the
non-local source term for both factors, i. e. the expression
\begin{equation}\nonumber\begin{split}
&\sum_{\pm,\pm}\int_{-\infty}^{\infty}\int_{-\infty}^{\infty}\int_{T}^{\infty}t\int_{t}^{\infty}\la
e_{\pm}(x,\xi)\bm 1&0\\0&-1\endm e_{\pm}(x,\xi'), \phi(x)\ra \\
&\int_{0}^{s}e^{\mp i(s-\lambda)(\xi^{2}+1)}\calF_{\pm}\bm
|\tilde{U}^{(s)}|^{4}\tilde{U}^{(s)}(\lambda,.)\\-|\tilde{U}^{(s)}|^{4}\overline{\tilde{U}^{(s)}}(\lambda,.)\endm(\xi)
d\lambda \\&\hspace{4cm}\int_{0}^{s}e^{\mp
i(s-\lambda')(\xi'^{2}+1)}\calF_{\pm}\bm
|\tilde{U}^{(s)}|^{4}\tilde{U}^{(s)}(\lambda',.)\\-|\tilde{U}^{(s)}|^{4}\overline{\tilde{U}^{(s)}}(\lambda',.)\endm(\xi')
d\lambda' d\xi d\xi'ds dt
\end{split}\end{equation}
The remaining (local) source terms are handled by the exact same
method but much easier, hence omitted. Then we focus on the most
difficult case when the $s$-phases cancel each other, i. e. when
there is a $+$ and a $-$ sign. In order to render the Fourier
transforms explicit, we write as usual
\begin{equation}\nonumber
\bm
|\tilde{U}^{(s)}|^{4}\tilde{U}^{(s)}(\lambda,.)\\-|\tilde{U}^{(s)}|^{4}\overline{\tilde{U}^{(s)}}(\lambda,.)\endm
=\chi_{>0}(x)\bm
|\tilde{U}^{(s)}|^{4}\tilde{U}^{(s)}(\lambda,.)\\-|\tilde{U}^{(s)}|^{4}\overline{\tilde{U}^{(s)}}(\lambda,.)\endm+\chi_{\leq
0}(x)\bm
|\tilde{U}^{(s)}|^{4}\tilde{U}^{(s)}(\lambda,.)\\-|\tilde{U}^{(s)}|^{4}\overline{\tilde{U}^{(s)}}(\lambda,.)\endm,
\end{equation}
with a similar expression for the 2nd factor in $\lambda'$. W. l.
o. g. we shall then include the $\chi_{>0}(x)$-cutoff in both
cases. Then we subdivide the $(\xi,\xi')$-plane into the four
standard quadrants. If $(\xi,\xi')$ is in the first or third
quadrant, observe that
\begin{equation}\nonumber
\frac{\la \frac{e_{+}(x,\xi)}{\xi}\bm
0&1\\-1&0\endm\frac{e_{+}(x,\xi')}{\xi'},
\phi(x)\ra}{\xi-\xi'}=\frac{\la
[\frac{e_{+}(x,\xi)}{\xi}-\frac{e_{+}(x,\xi')}{\xi'}]\bm
0&1\\-1&0\endm\frac{e_{+}(x,\xi')}{\xi'}, \phi(x)\ra}{\xi-\xi'},
\end{equation}
whence this is smooth and bounded with bounded derivatives in the
interior of these quadrants, and continuous up to the boundary. If
$(\xi,\xi')$ is in one of the other quadrants, we have
\begin{equation}\nonumber
\frac{\la \frac{e_{+}(x,\xi)}{\xi}\bm
0&1\\-1&0\endm\frac{e_{+}(x,\xi')}{\xi'},
\phi(x)\ra}{\xi+\xi'}=\frac{\la
[\frac{e_{+}(x,\xi)}{\xi}-\frac{e_{+}(x,-\xi')}{-\xi'}]\bm
0&1\\-1&0\endm\frac{e_{+}(x,\xi')}{\xi'}, \phi(x)\ra}{\xi+\xi'},
\end{equation}
whence the same comment applies. Now let $(\xi,\xi')$ be in the
third quadrant. Then we can write
\begin{equation}\nonumber\begin{split}
&\calF_{+}\bm
\chi_{>0}(x)|\tilde{U}^{(s)}|^{4}\tilde{U}^{(s)}(\lambda,.)\\-\chi_{>0}(x)|\tilde{U}^{(s)}|^{4}\overline{\tilde{U}^{(s)}}(\lambda,.)\endm(\xi)
\\&=\la\bm
\chi_{>0}(x)|\tilde{U}^{(s)}|^{4}\tilde{U}^{(s)}(\lambda,.)\\-\chi_{>0}(x)|\tilde{U}^{(s)}|^{4}\overline{\tilde{U}^{(s)}}(\lambda,.)\endm,
(e^{ix\xi}-e^{-ix\xi})\underline{e}+(1+r(-\xi))e^{-ix\xi}\underline{e}+\phi(x,\xi)\ra,
\end{split}\end{equation}
with a similar expression for $\calF_{-}\bm
\chi_{>0}(x)|\tilde{U}^{(s)}|^{4}\tilde{U}^{(s)}(\lambda',.)\\-\chi_{>0}(x)|\tilde{U}^{(s)}|^{4}\overline{\tilde{U}^{(s)}}(\lambda',.)\endm(\xi')$,
where $\underline{e}$ gets replaced by $\sigma_{1}\underline{e}$.
We shall treat in detail the contribution of
$(e^{ix\xi}-e^{-ix\xi})\underline{e}$, the other contributions
being treated similarly. As usual we need to distinguish between
different frequency ranges: first assume
$\max\{|\xi|,|\xi'|\}>s^{\epsilon(\delta_{2})}$, for suitable
$\epsilon(\delta_{2})$.  For example assume
$|\xi|>s^{\epsilon(\delta_{2})}$, effected by means of a smooth
cutoff $\phi_{>s^{\epsilon}}(\xi)$. Observe that then
\begin{equation}\nonumber\begin{split}
&\la \chi_{>0}(x)\bm
|\tilde{U}^{(s)}|^{4}\tilde{U}^{(s)}(\lambda,.)\\-|\tilde{U}^{(s)}|^{4}\overline{\tilde{U}^{(s)}}(\lambda,.)\endm,
(e^{ix\xi}-e^{-ix\xi})\underline{e}\ra+=O(\frac{1}{s^{\epsilon}})[\la\delta_{0}(x)\bm
|\tilde{U}^{(s)}|^{4}\tilde{U}^{(s)}(\lambda,.)\\-|\tilde{U}^{(s)}|^{4}\overline{\tilde{U}^{(s)}}(\lambda,.)\endm,
(e^{ix\xi}+e^{-ix\xi})\underline{e}\ra\\
&\hspace{8cm}+\la \chi_{>0}(x)\partial_{x}\bm
|\tilde{U}^{(s)}|^{4}\tilde{U}^{(s)}(\lambda,.)\\-|\tilde{U}^{(s)}|^{4}\overline{\tilde{U}^{(s)}}(\lambda,.)\endm,
(e^{ix\xi}+e^{-ix\xi})\underline{e}\ra]
\end{split}\end{equation}
For the boundary term, integrate by parts in $\xi$ in order to
score arbitrary gains in $s$. For the 2nd term, keep integrating
by parts in $x$ until a boundary term results or else enough
powers
of $s$ are gained. \\
Thus we may now include smooth cutoffs
$\phi_{<s^{\epsilon}}(\xi)$, $\phi_{<s^{\epsilon}}(\xi')$.
Commence with the case $\max\{\lambda,\lambda'\}<\frac{s}{2}$. We
perform integrations by parts in $\xi$, $\xi'$, and obtain the
following list of integrals
\begin{equation}\nonumber\begin{split}
&A=\int_{-\infty}^{0}\int_{-\infty}^{0}\int_{T}^{\infty}t\int_{t}^{\infty}\phi_{<s^{\epsilon}}(\xi)\phi_{<s^{\epsilon}}(\xi')
\la \frac{e(x,\xi)}{\xi}\bm 0&1\\-1&0\endm\frac{e(x,\xi')}{\xi'},
\phi(x)\ra\\
&\int_{0}^{\frac{s}{2}}\frac{1}{s-\lambda}e^{-i(s-\lambda)(\xi^{2}+1)}\partial_{\xi}\calF_{+}\bm
\chi_{>0}(x)|\tilde{U}^{(s)}|^{4}\tilde{U}^{(s)}(\lambda,.)\\-\chi_{>0}(x)|\tilde{U}^{(s)}|^{4}\overline{\tilde{U}^{(s)}}(\lambda,.)\endm(\xi)
d\lambda\\
&\hspace{2cm}\int_{0}^{\frac{s}{2}}\frac{1}{s-\lambda'}e^{+i(s-\lambda)(\xi'^{2}+1)}\partial_{\xi'}\calF_{-}\bm
\chi_{>0}(x)|\tilde{U}^{(s)}|^{4}\tilde{U}^{(s)}(\lambda',.)\\-\chi_{>0}(x)|\tilde{U}^{(s)}|^{4}\overline{\tilde{U}^{(s)}}(\lambda',.)\endm(\xi')
d\lambda' ds dt d\xi d\xi'\\
\end{split}\end{equation}
\begin{equation}\nonumber\begin{split}
&B=\int_{-\infty}^{0}\int_{-\infty}^{0}\int_{T}^{\infty}t\int_{t}^{\infty}\phi_{<s^{\epsilon}}(\xi)\phi_{<s^{\epsilon}}(\xi')
\la \partial_{\xi}[\frac{e(x,\xi)}{\xi}]\bm
0&1\\-1&0\endm\frac{e(x,\xi')}{\xi'},
\phi(x)\ra\\
&\int_{0}^{\frac{s}{2}}\frac{1}{s-\lambda}e^{-i(s-\lambda)(\xi^{2}+1)}\calF_{+}\bm
\chi_{>0}(x)|\tilde{U}^{(s)}|^{4}\tilde{U}^{(s)}(\lambda,.)\\-\chi_{>0}(x)|\tilde{U}^{(s)}|^{4}\overline{\tilde{U}^{(s)}}(\lambda,.)\endm(\xi)
d\lambda\\
&\hspace{2cm}\int_{0}^{\frac{s}{2}}\frac{1}{s-\lambda'}e^{+i(s-\lambda')(\xi'^{2}+1)}\partial_{\xi'}\calF_{-}\bm
\chi_{>0}(x)|\tilde{U}^{(s)}|^{4}\tilde{U}^{(s)}(\lambda',.)\\-\chi_{>0}(x)|\tilde{U}^{(s)}|^{4}\overline{\tilde{U}^{(s)}}(\lambda',.)\endm(\xi')
d\lambda' ds dt d\xi d\xi'\\
\end{split}\end{equation}
\begin{equation}\nonumber\begin{split}
&C=\int_{-\infty}^{0}\int_{-\infty}^{0}\int_{T}^{\infty}t\int_{t}^{\infty}\phi_{<s^{\epsilon}}(\xi)\phi_{<s^{\epsilon}}(\xi')
\la \frac{e(x,\xi)}{\xi}\bm
0&1\\-1&0\endm\partial_{\xi'}[\frac{e(x,\xi')}{\xi'}],
\phi(x)\ra\\
&\int_{0}^{\frac{s}{2}}\frac{1}{s-\lambda}e^{-i(s-\lambda)(\xi^{2}+1)}\partial_{\xi}\calF_{+}\bm
\chi_{>0}(x)|\tilde{U}^{(s)}|^{4}\tilde{U}^{(s)}(\lambda,.)\\-\chi_{>0}(x)|\tilde{U}^{(s)}|^{4}\overline{\tilde{U}^{(s)}}(\lambda,.)\endm(\xi)
d\lambda\\
&\hspace{2cm}\int_{0}^{\frac{s}{2}}\frac{1}{s-\lambda'}e^{+i(s-\lambda')(\xi'^{2}+1)}\calF_{-}\bm
\chi_{>0}(x)|\tilde{U}^{(s)}|^{4}\tilde{U}^{(s)}(\lambda',.)\\-\chi_{>0}(x)|\tilde{U}^{(s)}|^{4}\overline{\tilde{U}^{(s)}}(\lambda',.)\endm(\xi')
d\lambda' ds dt d\xi d\xi'\\
\end{split}\end{equation}
\begin{equation}\nonumber\begin{split}
&D=\int_{-\infty}^{0}\int_{-\infty}^{0}\int_{T}^{\infty}t\int_{t}^{\infty}\phi_{<s^{\epsilon}}(\xi)\phi_{<s^{\epsilon}}(\xi')
\la \partial_{\xi}[\frac{e(x,\xi)}{\xi}]\bm
0&1\\-1&0\endm\partial_{\xi'}[\frac{e(x,\xi')}{\xi'}],
\phi(x)\ra\\
&\int_{0}^{\frac{s}{2}}\frac{1}{s-\lambda}e^{-i(s-\lambda)(\xi^{2}+1)}\calF_{+}\bm
\chi_{>0}(x)|\tilde{U}^{(s)}|^{4}\tilde{U}^{(s)}(\lambda,.)\\-\chi_{>0}(x)|\tilde{U}^{(s)}|^{4}\overline{\tilde{U}^{(s)}}(\lambda,.)\endm(\xi)
d\lambda\\
&\hspace{2cm}\int_{0}^{\frac{s}{2}}\frac{1}{s-\lambda'}e^{+i(s-\lambda')(\xi'^{2}+1)}\calF_{-}\bm
\chi_{>0}(x)|\tilde{U}^{(s)}|^{4}\tilde{U}^{(s)}(\lambda',.)\\-\chi_{>0}(x)|\tilde{U}^{(s)}|^{4}\overline{\tilde{U}^{(s)}}(\lambda',.)\endm(\xi')
d\lambda' ds dt d\xi d\xi'\\
\end{split}\end{equation}

We don't worry about the case when $\partial_{\xi}$ or
$\partial_{\xi'}$ falls on one of the cutoffs since then we have
either $|\xi|\sim s^{\epsilon}$ or $|\xi'|\sim s^{\epsilon}$,
which case is treated as before. Observe that we can write
\begin{equation}\nonumber
\calF_{+}\bm
|\tilde{U}^{(s)}|^{4}\tilde{U}^{(s)}(\lambda,.)\\-|\tilde{U}^{(s)}|^{4}\overline{\tilde{U}^{(s)}}(\lambda,.)\endm(\xi)
=\xi\int_{0}^{1}\partial_{\xi}\calF_{+}\bm
|\tilde{U}^{(s)}|^{4}\tilde{U}^{(s)}(\lambda,.)\\-|\tilde{U}^{(s)}|^{4}\overline{\tilde{U}^{(s)}}(\lambda,.)\endm(\alpha\xi)d\alpha
\end{equation}
Of course a similar identity applies to $\calF_{-}(...)$, whence
we reduce to estimating expressions of the form
\begin{equation}\label{201}\begin{split}
&\int_{-\infty}^{0}\int_{-\infty}^{0}\int_{T}^{\infty}t\int_{t}^{\infty}\phi_{<s^{\epsilon}}(\xi)\phi_{<s^{\epsilon}}(\xi')\xi
g(\xi,\xi')\\
&\hspace{2cm}\times\int_{0}^{\frac{s}{2}}\frac{1}{s-\lambda}e^{-i(s-\lambda)(\xi^{2}+1)}\partial_{\xi}\calF_{+}\bm
|\tilde{U}^{(s)}|^{4}\tilde{U}^{(s)}(\lambda,.)\\-|\tilde{U}^{(s)}|^{4}\overline{\tilde{U}^{(s)}}(\lambda,.)\endm(\xi)
d\lambda\\
&\hspace{2cm}\times\int_{0}^{\frac{s}{2}}\frac{1}{s-\lambda'}e^{+i(s-\lambda')(\xi'^{2}+1)}\partial_{\xi'}\calF_{-}\bm
|\tilde{U}^{(s)}|^{4}\tilde{U}^{(s)}(\lambda,.)\\-|\tilde{U}^{(s)}|^{4}\overline{\tilde{U}^{(s)}}(\lambda',.)\endm(\xi')
d\lambda' d\xi d\xi',
\end{split}\end{equation}
where the function $g(\xi,\xi')$ is smooth and bounded in the
interior of the third quadrant as well as continuous up to the
boundary. Expand as usual
\begin{equation}\nonumber
\partial_{\xi}\calF_{+}\bm
\chi_{>0}(x)|\tilde{U}^{(s)}|^{4}\tilde{U}^{(s)}(\lambda,.)\\-\chi_{>0}(x)|\tilde{U}^{(s)}|^{4}\overline{\tilde{U}^{(s)}}(\lambda,.)\endm(\xi)
=\la\bm
x\chi_{>0}(x)|\tilde{U}^{(s)}|^{4}\tilde{U}^{(s)}(\lambda,.)\\-x\chi_{>0}(x)|\tilde{U}^{(s)}|^{4}\overline{\tilde{U}^{(s)}}(\lambda,.)\endm,
(e^{ix\xi}+e^{-ix\xi})\underline{e}\ra+\ldots
\end{equation}
where $\ldots$ represent terms that can be treated similarly. Now
assume that we localize $x$ to dyadic range $|x|\sim 2^{k}$,
$k\geq 0$. If then we have $|\xi|>\max\{s^{-\frac{1}{2+}},
s^{-1-}2^{k}\}$, effected by means of a smooth cutoff, we obtain
arbitrary gains in $s$ by integration by parts in $\xi$. Thus we
shall now include a localizer $\phi_{<\max\{s^{-\frac{1}{2+}},
s^{-1-}2^{k}\}}(\xi)$ upon localizing $x$ to dyadic range $|x|\sim
2^{k}$, i. e. we reduce to considering
\begin{equation}\label{202}\begin{split}
&\sum_{k}\int_{-\infty}^{0}\int_{-\infty}^{0}\int_{T}^{\infty}t\int_{t}^{\infty}\phi_{<\max\{s^{-\frac{1}{2+}},
s^{-1-}2^{k}\}}(\xi)\phi_{<s^{\epsilon}}(\xi)\phi_{<s^{\epsilon}}(\xi')\xi
g(\xi,\xi')\\
&\hspace{2cm}\times\int_{0}^{\frac{s}{2}}\frac{1}{s-\lambda}e^{-i(s-\lambda)(\xi^{2}+1)}\la\bm
\phi_{\sim
2^{k}}(x)x|\tilde{U}^{(s)}|^{4}\tilde{U}^{(s)}(\lambda,.)\\-\phi_{\sim
2^{k}}(x)x|\tilde{U}^{(s)}|^{4}\overline{\tilde{U}^{(s)}}(\lambda,.)\endm,
(e^{ix\xi}+e^{-ix\xi})\underline{e}\ra
d\lambda\\
&\hspace{2cm}\times\int_{0}^{\frac{s}{2}}\frac{1}{s-\lambda'}e^{+i(s-\lambda')(\xi'^{2}+1)}\la\bm
x|\tilde{U}^{(s)}|^{4}\tilde{U}^{(s)}(\lambda,.)\\-x|\tilde{U}^{(s)}|^{4}\overline{\tilde{U}^{(s)}}(\lambda,.)\endm,
(e^{ix\xi'}+e^{-ix\xi'})\sigma_{1}\underline{e}\ra d\lambda' d\xi
d\xi',
\end{split}\end{equation}
Note that summing over $k$ will amount to an extra $\log s$ at
most, whence we shall safely discard this summation. Our strategy
shall be to perform an integration by parts in $\la\bm \phi_{\sim
2^{k}}(x)|\tilde{U}^{(s)}|^{4}\tilde{U}^{(s)}(\lambda,.)\\-\phi_{\sim
2^{k}}(x)|\tilde{U}^{(s)}|^{4}\overline{\tilde{U}^{(s)}}(\lambda,.)\endm,
(e^{ix\xi}+e^{-ix\xi})\underline{e}\ra$. For this to be useful,
though, we need to achieve some preliminary reductions in the last
factor $\tilde{U}^{(s)}$, just as in the proof of the SLDE. Recall
that we can write
\begin{equation}\nonumber
\tilde{U}^{(s)}(\lambda,.)=\int_{0}^{\lambda}\sqrt{\lambda-\mu}^{-1}\int_{-\infty}^{\infty}e^{\frac{(y-z)^{2}}{i(\lambda-\mu)}}g(\mu,z)dz
d\mu
\end{equation}
where
$g(\mu,z)=|\tilde{U}^{(s)}|^{4}\tilde{U}^{(s)}(\mu,z)+\ldots$.
Decompose
\begin{equation}\nonumber
|\tilde{U}^{(s)}|^{4}\tilde{U}^{(s)}(\mu,z)=\chi_{<2^{k}}(\mu)|\tilde{U}^{(s)}|^{4}\tilde{U}^{(s)}(\mu,z)
+\chi_{\geq 2^{k}}(\mu)|\tilde{U}^{(s)}|^{4}\tilde{U}^{(s)}(\mu,z)
\end{equation}
Then observe that
\begin{equation}\nonumber
||\int_{0}^{\lambda}\sqrt{\lambda-\mu}^{-1}\int_{-\infty}^{\infty}e^{\frac{(y-z)^{2}}{i(\lambda-\mu)}}\chi_{\geq
2^{k}}(\mu)|\tilde{U}^{(s)}|^{4}\tilde{U}^{(s)}(\mu,z)dz||_{L_{y}^{2}}\lesssim
\int_{2^{k}}^{\lambda}\mu^{-2}d\mu\lesssim 2^{-k}
\end{equation}
Similarly, we have
\begin{equation}\nonumber\begin{split}
&||\int_{0}^{\lambda}\sqrt{\lambda-\mu}^{-1}\int_{-\infty}^{\infty}e^{\frac{(y-z)^{2}}{i(\lambda-\mu)}}\chi_{<
2^{k}}(\mu)\chi_{\geq
2^{k-10}}(z)|\tilde{U}^{(s)}|^{4}\tilde{U}^{(s)}(\mu,z)dz||_{L_{y}^{2}}
\\&\lesssim
\int_{0}^{\lambda}||\chi_{\geq
s^{k}}(z)\frac{2i\mu\partial_{z}\tilde{U}^{(s)}(\mu,z)-C\tilde{U}^{(s)}(\mu,z)}{z}||_{L_{z}^{2}}\mu^{-2}d\mu\lesssim
\log\lambda 2^{-k}
\end{split}\end{equation}
One obtains similar estimates if one substitutes the remaining
local terms in $g(\mu,z)$, localized to $|z|>2^{k-10}$, in the
preceding integral. Thus if we substitute
\begin{equation}\nonumber\begin{split}
V^{(s)}(\lambda,y):=&\int_{0}^{\lambda}\sqrt{\lambda-\mu}^{-1}\int_{-\infty}^{\infty}e^{\frac{(y-z)^{2}}{i(\lambda-\mu)}}[\chi_{\geq
2^{k}}(\mu)|\tilde{U}^{(s)}|^{4}\tilde{U}^{(s)}(\mu,z)\\&+\chi_{<
2^{k}}(\mu)\chi_{\geq
2^{k-10}}(z)|\tilde{U}^{(s)}|^{4}\tilde{U}^{(s)}(\mu,z)+\chi_{\geq
2^{k-10}}(z)...]dz d\mu
\end{split}\end{equation}
instead of
$\tilde{U}^{(s)}$ for the last factor in
$|\tilde{U}^{(s)}|^{4}\tilde{U}^{(s)}(\lambda,y)$, we get
\begin{equation}\nonumber
||\phi_{\sim 2^{k}}(x)x
|\tilde{U}^{(s)}|^{4}V^{(s)}(\lambda,x)||_{L_{x}^{1}}\lesssim
\lambda^{-\frac{3}{2}}
\end{equation}
We now show how this suffices to control
\begin{equation}\nonumber\begin{split}
&\int_{-\infty}^{0}\int_{-\infty}^{0}\int_{T}^{\infty}t\int_{t}^{\infty}\phi_{<s^{\epsilon}}(\xi)\phi_{<s^{\epsilon}}(\xi')\phi_{<\max\{s^{-\frac{1}{2+}},
s^{-1-}2^{k}\}}(\xi)\xi
g(\xi,\xi')\\
&\hspace{2cm}\times\int_{0}^{\frac{s}{2}}\frac{1}{s-\lambda}e^{-i(s-\lambda)(\xi^{2}+1)}\la
\bm x\phi_{\sim
2^{k}}(x)|\tilde{U}^{(s)}|^{4}V^{(s)}\lambda,x)\\-x\phi_{\sim
2^{k}}(x)|\tilde{U}^{(s)}|^{4}\overline{V^{(s)}}(\lambda,.)\endm,
(e^{ix\xi}+e^{-ix\xi})\underline{e}\ra
d\lambda\\
&\hspace{2cm}\times\int_{0}^{\frac{s}{2}}\frac{1}{s-\lambda'}e^{+i(s-\lambda)(\xi'^{2}+1)}\la
\bm
x|\tilde{U}^{(s)}|^{4}\tilde{V}^{(s)}(\lambda',.)\\-x|\tilde{U}^{(s)}|^{4}\overline{\tilde{V}^{(s)}}(\lambda',.)\endm,
(e^{ix\xi'}+e^{-ix\xi'})\sigma_{1}\underline{e}\ra d\lambda' d\xi
d\xi',
\end{split}\end{equation}
where $\tilde{V}^{(s)}(\lambda',.)$ is defined analogously. Thus
we can estimate this by
\begin{equation}\nonumber\begin{split}
&\int_{T}^{\infty}t\int_{t}^{\infty}\max\{s^{-\frac{1}{2+}},
s^{-1-}2^{k}\}\int_{0}^{\frac{s}{2}}\frac{1}{(s-\lambda)^{\frac{3}{2}}}||\bm
x\phi_{\sim
2^{k}}(x)|\tilde{U}^{(s)}|^{4}V^{(s)}\lambda,x)\\-x\phi_{\sim
2^{k}}(x)|\tilde{U}^{(s)}|^{4}\overline{V^{(s)}}(\lambda,.)\endm||_{L_{x}^{1}}
d\lambda\\
&\int_{0}^{\frac{s}{2}}\frac{1}{(s-\lambda')^{\frac{3}{2}}}||\bm
x|\tilde{U}^{(s)}|^{4}\tilde{V}^{(s)}(\lambda',.)\\-x|\tilde{U}^{(s)}|^{4}\overline{\tilde{V}^{(s)}}(\lambda',.)\endm||_{L_{x}^{1}}d\lambda'\lesssim
\int_{T}^{\infty}t\int_{t}^{\infty}s^{-\frac{1}{2+}}s^{-3}ds
dt\lesssim T^{-\frac{1}{2+}}
\end{split}\end{equation}
We are exploiting here the pseudo-conformal almost conservation,
which implies that
\begin{equation}\nonumber
||x\phi_{\sim
2^{k}}(x)|\tilde{U}^{(s)}|^{4}V^{(s)}\lambda,x)||_{L_{x}^{1}}\lesssim
\frac{\lambda}{2^{k}}\lambda^{-\frac{3}{2}}
\end{equation}
Thus we now replace at least one of $\tilde{U}^{(s)}(\lambda,.)$,
$\tilde{U}^{(s)}(\lambda',.)$ by
\begin{equation}\nonumber
\tilde{W}^{(s)}(\lambda,.):=\int_{-\infty}^{\infty}\int_{0}^{\lambda}\frac{1}{\sqrt{\lambda-\mu}}e^{\frac{(y-z)^{2}}{i(\lambda-\mu)}}[\chi_{<2^{k}}(\mu)\chi_{<2^{k-10}}(z)|\tilde{U}^{(s)}|^{4}\tilde{U}^{(s)}(\mu,z)+\ldots]d\mu
dz\,\text{etc}
\end{equation}
Then we decouple the $\xi,\xi'$ variables in \eqref{202}, which
can be achieved by means of discrete Fourier expansion:
\begin{equation}\nonumber
\phi_{<s^{\epsilon}}(|\xi|)\phi_{<s^{\epsilon}}(|\xi'|)\xi
g(\xi,\xi')=\xi\sum_{n,m\in
s^{-\epsilon}{\mathbf{Z}}}a_{nm}e^{in\xi+im\xi'},\,|a_{nm}|\lesssim
[s^{\epsilon}|n|+s^{\epsilon}|m|]^{-N}
\end{equation}
Consider the case when we replace the fifth
$\tilde{U}^{(s)}(\lambda,.)$ by $\tilde{W}^{(s)}(\lambda,.)$. We
are thus led to estimating contributions of the form
\begin{equation}\nonumber
\int_{-\infty}^{0}\xi\phi_{<s^{\epsilon}}(\xi)e^{i(s-\lambda)(\xi^{2}+1)}\la
\chi_{>0}(x)\phi_{\sim
2^{k}}(x)x|\tilde{U}^{(s)}|^{4}\tilde{W}^{(s)}(\lambda,.),
(e^{i(x+n)\xi}+e^{-i(x-n)\xi})\underline{e}\ra d\xi
\end{equation}
We shall put $n=0$ since the other cases are dealt with similarly.
We treat here the contribution of $e^{-ix\xi}$, the one of
$e^{+ix\xi}$ being treated similarly. Switch the order of
integration in this, and introduce the new variable
$\tilde{\xi}:=\sqrt{s-\lambda}\xi+\frac{x}{2\sqrt{s-\lambda}}$.
Then we can rewrite the preceding expression as
\begin{equation}\nonumber
\frac{1}{\sqrt{s-\lambda}}\int_{0}^{\infty}\int_{-\infty}^{\frac{x}{2\sqrt{s-\lambda}}}\phi_{<s^{\epsilon}}(\frac{\tilde{\xi}-\frac{x}{2\sqrt{s-\lambda}}}{\sqrt{s-\lambda}})
\frac{\tilde{\xi}-\frac{x}{2\sqrt{s-\lambda}}}{\sqrt{s-\lambda}}e^{i\tilde{\xi}^{2}}e^{\frac{x^{2}}{4i(s-\lambda)}}x\phi_{\sim
2^{k}}(x)|\tilde{U}^{(s)}|^{4}\tilde{W}^{(s)}(\lambda,.)d\tilde{\xi}dx
\end{equation}
Now perform an integration by parts in the $x$-variable, and
replace the preceding by the sum of multiples of the following
expressions(as well as equivalent terms):
\begin{equation}\label{203}
\frac{1}{\sqrt{s-\lambda}}\int_{0}^{\infty}\int_{-\infty}^{\frac{x}{2\sqrt{s-\lambda}}}\partial_{x}[\phi_{<s^{\epsilon}}(\frac{\tilde{\xi}-\frac{x}{2\sqrt{s-\lambda}}}{\sqrt{s-\lambda}})
\frac{\tilde{\xi}-\frac{x}{2\sqrt{s-\lambda}}}{\sqrt{s-\lambda}}]e^{i\tilde{\xi}^{2}}d\tilde{\xi}x\phi_{\sim
2^{k}}(x)|\tilde{U}^{(s)}|^{4}[\int_{x}^{\infty}e^{\frac{y^{2}}{4i(s-\lambda)}}\tilde{W}^{(s)}(\lambda,y)dy]dx
\end{equation}
\begin{equation}\label{204}
\frac{1}{\sqrt{s-\lambda}}\int_{0}^{\infty}\int_{-\infty}^{\frac{x}{2\sqrt{s-\lambda}}}[\phi_{<s^{\epsilon}}(\frac{\tilde{\xi}-\frac{x}{2\sqrt{s-\lambda}}}{\sqrt{s-\lambda}})
\frac{\tilde{\xi}-\frac{x}{2\sqrt{s-\lambda}}}{\sqrt{s-\lambda}}]e^{i\tilde{\xi}^{2}}d\tilde{\xi}\phi_{\sim
2^{k}}(x)|\tilde{U}^{(s)}|^{4}[\int_{x}^{\infty}e^{\frac{y^{2}}{4i(s-\lambda)}}\tilde{W}^{(s)}(\lambda,y)dy]dx
\end{equation}
\begin{equation}\label{205}\begin{split}
&\frac{1}{\sqrt{s-\lambda}}\int_{0}^{\infty}\int_{-\infty}^{\frac{x}{2\sqrt{s-\lambda}}}[\phi_{<s^{\epsilon}}(\frac{\tilde{\xi}-\frac{x}{2\sqrt{s-\lambda}}}{\sqrt{s-\lambda}})
\frac{\tilde{\xi}-\frac{x}{2\sqrt{s-\lambda}}}{\sqrt{s-\lambda}}]e^{i\tilde{\xi}^{2}}d\tilde{\xi}x\phi_{\sim
2^{k}}(x)\partial_{x}[|\tilde{U}^{(s)}|^{4}(\lambda,x)]\\&\hspace{10cm}[\int_{x}^{\infty}e^{\frac{y^{2}}{4i(s-\lambda)}}\tilde{W}^{(s)}(\lambda,y)dy]dx
\end{split}\end{equation}
Now write
\begin{equation}\nonumber
\tilde{W}^{(s)}(\lambda,y)=\int_{-\infty}^{\infty}\int_{0}^{\lambda}\frac{1}{\sqrt{\lambda-\mu}}e^{\frac{(y-z)^{2}}{i(\lambda-\mu)}}g(\mu,z)d\mu
dz,
\end{equation}
where we have $|z|<2^{k-10}$ on the support of $g(\mu,z)$.  Thus
we obtain
\begin{equation}\nonumber\begin{split}
&\int_{x}^{\infty}e^{\frac{y^{2}}{4i(s-\lambda)}}\tilde{W}^{(s)}(\lambda,y)dy
=\int_{-\infty}^{\infty}\int_{x}^{\infty}\int_{0}^{\lambda}\frac{1}{\sqrt{\lambda-\mu}}e^{\frac{y^{2}}{4i(s-\lambda)}}e^{\frac{(y-z)^{2}}{i(\lambda-\mu)}}g(\mu,z)d\mu
dydz\\
&=\int_{-\infty}^{\infty}\int_{0}^{\lambda}O(\frac{1}{\sqrt{\lambda-\mu}\sqrt{\frac{1}{4(s-\lambda)}+\frac{1}{\lambda-\mu}}}\frac{1}{x\sqrt{\frac{1}{4(s-\lambda)}+\frac{1}{\lambda-\mu}}
-\frac{z}{(\lambda-\mu)\sqrt{\frac{1}{4(s-\lambda)}+\frac{1}{\lambda-\mu}}}})
g(\mu, z) d\mu \phi_{<2^{k-10}}(|z|) dz,
\end{split}\end{equation}
which is seen for $x\sim 2^{k}$ to be of order
$(\frac{x}{\sqrt{\lambda}})^{-1}$, upon using the definition of
$g(\mu,z)$. Now plug this into \eqref{203}. For example, we get
\begin{equation}\nonumber\begin{split}
&|\frac{1}{\sqrt{s-\lambda}}\int_{0}^{\infty}\int_{-\infty}^{\frac{x}{2\sqrt{s-\lambda}}}\frac{1}{s-\lambda}[\phi_{<s^{\epsilon}}'(\frac{\tilde{\xi}-\frac{x}{2\sqrt{s-\lambda}}}{\sqrt{s-\lambda}})
\frac{x}{2(s-\lambda)}e^{i\tilde{\xi}^{2}}d\tilde{\xi}x\phi_{\sim
2^{k}}(x)|\tilde{U}^{(s)}|^{4}\\&\hspace{10cm}[\int_{x}^{\infty}e^{\frac{y^{2}}{4i(s-\lambda)}}\tilde{W}^{(s)}(\lambda,y)dy]dx|
\\&\lesssim\frac{1}{\sqrt{s-\lambda}}\frac{1}{(s-\lambda)^{2}}\sqrt{\lambda}\lambda^{\epsilon(\delta_{2})}
\end{split}\end{equation}
Observe that we are using pseudo-conformal almost conservation
here. Integrating over $\lambda<\frac{s}{2}$ results in the upper
bound $\lesssim s^{-1}$.  The remaining contributions to
\eqref{203} (Leibnitz rule) are treated similarly, as is the
contribution of \eqref{204}. Now consider \eqref{205}. Here we
invoke the same trick as in the proof of SLDE:
\begin{equation}\nonumber
\partial_{x}[|\tilde{U}^{(s)}|^{2}(\lambda,x)]=\frac{1}{i\lambda}[C\tilde{U}^{(s)}
\overline{\tilde{U}^{(s)}}(\lambda,x)-\tilde{U}^{(s)}\overline{C\tilde{U}^{(s)}}(\lambda,x)]
\end{equation}
For example, we can estimate
\begin{equation}\nonumber\begin{split}
&|\frac{1}{\sqrt{s-\lambda}}\int_{0}^{\infty}\int_{-\infty}^{\frac{x}{2\sqrt{s-\lambda}}}[\phi_{<s^{\epsilon}}(\frac{\tilde{\xi}-\frac{x}{2\sqrt{s-\lambda}}}{\sqrt{s-\lambda}})
\frac{x}{(s-\lambda)}e^{i\tilde{\xi}^{2}}d\tilde{\xi}x\phi_{\sim
2^{k}}(x)\partial_{x}[|\tilde{U}^{(s)}|^{4}(\lambda,x)]\\&\hspace{10cm}[\int_{x}^{\infty}e^{\frac{y^{2}}{4i(s-\lambda)}}\tilde{W}^{(s)}(\lambda,y)dy]dx|\\
&\lesssim
|\frac{1}{\sqrt{s-\lambda}}\int_{0}^{\infty}\int_{-\infty}^{\frac{x}{2\sqrt{s-\lambda}}}[\phi_{<s^{\epsilon}}(\frac{\tilde{\xi}-\frac{x}{2\sqrt{s-\lambda}}}{\sqrt{s-\lambda}})
\frac{1}{(s-\lambda)}e^{i\tilde{\xi}^{2}}d\tilde{\xi}\phi_{\sim
2^{k}}(x)\frac{1}{\lambda}|C\tilde{U}^{(s)}||x\tilde{U}^{(s)}||\tilde{U}^{(s)}|^{2}(\lambda,x)]\sqrt{\lambda}
d\tilde{\xi} dx\\
&\lesssim
(s-\lambda)^{-\frac{3}{2}}\sqrt{\lambda}^{-1}\lambda^{\epsilon(\delta_{2})},
\end{split}\end{equation}
which upon integration over $\lambda<\frac{s}{2}$ again results in
the estimate $s^{-1+}$. The 2nd contribution to \eqref{205} is
treated similarly.
\\

Keeping in mind that we need to eventually estimate \eqref{202},
we next consider the expression
\begin{equation}\nonumber
\int_{-\infty}^{0}\phi_{<s^{\epsilon}}(\xi')e^{i(s-\lambda)(\xi'^{2}+1)}\la\chi_{>0}(x)\phi_{\sim
2^{k}}(x)x|\tilde{U}^{(s)}|^{4}\tilde{W}^{(s)}(\lambda',x),
(e^{ix\xi'}+e^{-ix\xi'})\underline{e}\ra d\xi'
\end{equation}
One proceeds analogously and obtains expressions as in
\eqref{203}, \eqref{204}, \eqref{205} but without the factor
$\frac{\tilde{\xi'}-\frac{x}{2\sqrt{s-\lambda'}}}{\sqrt{s-\lambda'}}$.
One then has to argue somewhat differently for the expression
\begin{equation}\nonumber
\frac{1}{\sqrt{s-\lambda'}}\int_{0}^{\infty}\int_{-\infty}^{\frac{x}{2\sqrt{s-\lambda'}}}\phi_{<s^{\epsilon}}(\frac{\tilde{\xi}-\frac{x}{2\sqrt{s-\lambda'}}}{\sqrt{s-\lambda'}})
e^{i\tilde{\xi}^{2}}d\tilde{\xi}\phi_{\sim
2^{k}}(x)|\tilde{U}^{(s)}|^{4}(\lambda',x)[\int_{x}^{\infty}e^{\frac{y^{2}}{4i(s-\lambda')}}\tilde{W}^{(s)}(\lambda',y)dy]dx
\end{equation}
Here we use that
\begin{equation}\nonumber
|\tilde{U}^{(s)}|(\lambda',x)[\frac{\la
x\ra^{\frac{1}{2}}}{\sqrt{\lambda'}}]^{-1}\lesssim
\sqrt{\lambda'}\lambda'^{-1+\epsilon(\delta_{2})},
\end{equation}
whence we can bound the above expression by
\begin{equation}\nonumber
\frac{1}{\sqrt{s-\lambda'}}\lambda'^{-\frac{3}{2}+\epsilon(\delta_{2})}
\end{equation}
Integrating over $\lambda'<\frac{s}{2}$ results in the upper bound
$\sqrt{s}^{-1+\epsilon(\delta_{2})}$.  Combining all these
estimates, we now obtain
\begin{equation}\nonumber\begin{split}
&|\sum_{k}\int_{-\infty}^{0}\int_{-\infty}^{0}\int_{T}^{\infty}t\int_{t}^{\infty}\phi_{<s^{\epsilon}}(\xi)\phi_{<s^{\epsilon}}(\xi')\xi
g(\xi,\xi')\\
&\hspace{2cm}\times\int_{0}^{\frac{s}{2}}\frac{1}{s-\lambda}e^{-i(s-\lambda)(\xi^{2}+1)}\la
\bm \phi_{\sim
2^{k}}(x)|\tilde{U}^{(s)}|^{4}\tilde{W}^{(s)}(\lambda,.)\\-\phi_{\sim
2^{k}}(x)|\tilde{U}^{(s)}|^{4}\overline{\tilde{W}^{(s)}}(\lambda,.)\endm,(e^{ix\xi}+e^{-ix\xi})\underline{e}\ra
d\lambda\\
&\hspace{2cm}\times\int_{0}^{\frac{s}{2}}\frac{1}{s-\lambda'}e^{+i(s-\lambda)(\xi'^{2}+1)}\la\bm
x|\tilde{U}^{(s)}|^{4}\tilde{W}^{(s)}(\lambda,.)\\-x|\tilde{U}^{(s)}|^{4}\overline{\tilde{W}^{(s)}}(\lambda,.)\endm,
(e^{ix\xi'}+e^{-ix\xi'})\sigma_{1}\underline{e}\ra d\lambda' d\xi
d\xi'|\lesssim T^{-\frac{1}{2+}}
\end{split}\end{equation}
\begin{equation}\nonumber\begin{split}
&|\sum_{k}\int_{-\infty}^{0}\int_{-\infty}^{0}\int_{T}^{\infty}t\int_{t}^{\infty}\phi_{<s^{\epsilon}}(\xi)\phi_{<s^{\epsilon}}(\xi')\xi
g(\xi,\xi')\\
&\hspace{2cm}\times\int_{0}^{\frac{s}{2}}\frac{1}{s-\lambda}e^{-i(s-\lambda)(\xi^{2}+1)}\la\bm
\phi_{\sim
2^{k}}(x)|\tilde{U}^{(s)}|^{4}\tilde{W}^{(s)}(\lambda,.)\\-\phi_{\sim
2^{k}}(x)|\tilde{U}^{(s)}|^{4}\overline{\tilde{W}^{(s)}}(\lambda,.)\endm,
(e^{ix\xi}+e^{-ix\xi})\underline{e}\ra
d\lambda\\
&\hspace{2cm}\times\int_{0}^{\frac{s}{2}}\frac{1}{s-\lambda'}e^{+i(s-\lambda)(\xi'^{2}+1)}\la\bm
x|\tilde{U}^{(s)}|^{4}\tilde{V}^{(s)}(\lambda,.)\\-x|\tilde{U}^{(s)}|^{4}\overline{\tilde{V}^{(s)}}(\lambda,.)\endm,
(e^{ix\xi'}+e^{-ix\xi'})\sigma_{1}\underline{e}\ra d\lambda' d\xi
d\xi'|\lesssim T^{-\frac{1}{2+}}
\end{split}\end{equation}
\begin{equation}\nonumber\begin{split}
&|\sum_{k}\int_{-\infty}^{0}\int_{-\infty}^{0}\int_{T}^{\infty}t\int_{t}^{\infty}\phi_{<s^{\epsilon}}(\xi)\phi_{<s^{\epsilon}}(\xi')\xi
g(\xi,\xi')\\
&\hspace{2cm}\times\int_{0}^{\frac{s}{2}}\frac{1}{s-\lambda}e^{-i(s-\lambda)(\xi^{2}+1)}\la\bm
\phi_{\sim
2^{k}}(x)|\tilde{U}^{(s)}|^{4}\tilde{V}^{(s)}(\lambda,.)\\-\phi_{\sim
2^{k}}(x)|\tilde{U}^{(s)}|^{4}\overline{\tilde{V}^{(s)}}(\lambda,.)\endm,
(e^{ix\xi}+e^{-ix\xi})\underline{e}\ra
d\lambda\\
&\hspace{2cm}\times\int_{0}^{\frac{s}{2}}\frac{1}{s-\lambda'}e^{+i(s-\lambda)(\xi'^{2}+1)}\la\bm
x|\tilde{U}^{(s)}|^{4}\tilde{W}^{(s)}(\lambda,.)\\-x|\tilde{U}^{(s)}|^{4}\overline{\tilde{W}^{(s)}}(\lambda,.)\endm,
(e^{ix\xi'}+e^{-ix\xi'})\sigma_{1}\underline{e}\ra d\lambda' d\xi
d\xi'|\lesssim T^{-\frac{1}{2+}},
\end{split}\end{equation}
which together with \eqref{201} implies
\begin{equation}\nonumber\begin{split}
&|\sum_{k}\int_{-\infty}^{0}\int_{-\infty}^{0}\int_{T}^{\infty}t\int_{t}^{\infty}\phi_{<s^{\epsilon}}(\xi)\phi_{<s^{\epsilon}}(\xi')\xi
g(\xi,\xi')\\
&\hspace{2cm}\times\int_{0}^{\frac{s}{2}}\frac{1}{s-\lambda}e^{-i(s-\lambda)(\xi^{2}+1)}\la
\bm \phi_{\sim
2^{k}}(x)|\tilde{U}^{(s)}|^{4}\tilde{U}^{(s)}(\lambda,.)\\-\phi_{\sim
2^{k}}(x)|\tilde{U}^{(s)}|^{4}\overline{\tilde{U}^{(s)}}(\lambda,.)\endm,(e^{ix\xi}+e^{-ix\xi})\underline{e}\ra
d\lambda\\
&\hspace{2cm}\times\int_{0}^{\frac{s}{2}}\frac{1}{s-\lambda'}e^{+i(s-\lambda)(\xi'^{2}+1)}\la\bm
x|\tilde{U}^{(s)}|^{4}\tilde{U}^{(s)}(\lambda,.)\\-x|\tilde{U}^{(s)}|^{4}\overline{\tilde{U}^{(s)}}(\lambda,.)\endm,(e^{ix\xi'}+e^{-ix\xi'})\sigma_{1}\underline{e}\ra
d\lambda' d\xi d\xi'|\lesssim T^{-\frac{1}{2+}},
\end{split}\end{equation}
as desired. \\
The case $\max\{\lambda,\lambda'\}>\frac{s}{2}$ is more of the
same. This concludes the proof of the estimate for the bilinear
symplectic form. We have now also filled the gap in retrieving
control over $||CU||_{L_{x}^{2}}$: while the $\phi$ in the
expression $\la \tilde{U}^{2}-\overline{\tilde{U}}^{2}, \phi\ra$
encountered there was time dependent, one easily checks that up to
an error which leads to an absolutely integrable expression, one
may replace this by a constant function.
\end{proof}
We are now also done with the proof of Proposition~\ref{Hard},
since all terms arising upon substituting \eqref{lambda6''},
\eqref{lambda6} into $\int_{T}^{\infty}t\lambda_{6}(t)dt$ are
controlled either by Lemma~\ref{general bilinear},
Lemma~\ref{nullform}, or else can be absolutely integrated.
\end{proof}

\subsection{Retrieving control over the modulation parameters, and parameters at infinity.}

We commence with
$\beta_{1}=\beta\nu-b_{\infty}\lambda_{\infty}^{-1}$, which is
given by the righthand side of \eqref{beta}. Observe that
schematically we have
\begin{equation}\nonumber
E_{2}(\sigma)=-\la N, \tilde{\xi}_{2}\ra
+\sum_{a=0,1}(\nu-1)^{a}\lambda_{6}(\sigma)+(\nu-1)(\sigma)\la \bm \tilde{U}\\
\overline{\tilde{U}}\endm_{dis}, \phi\ra
\end{equation}
Thus we recover the desired estimate for $\beta_{1}(s)$ upon using
Proposition~\ref{Hard}, if we also show that
\begin{equation}\nonumber
|\lambda_{\infty}^{-1}(s)\int_{s}^{T}\lambda_{\infty}(\sigma)(\nu-1)(\sigma)\la \bm \tilde{U}\\
\overline{\tilde{U}}\endm_{dis}, \phi\ra|\lesssim \la
s\ra^{-\frac{3}{2}+\delta_{1}}
\end{equation}
We state the
\begin{lemma}\label{Technical} Let $\Gamma\in A^{(n)}_{[0,T)}$, as usual $n$
sufficiently large. Then we have
\begin{equation}\nonumber
\int_{s}^{T}\lambda_{\infty}(\sigma)(\nu-1)(\sigma)\la \bm \tilde{U}\\
\overline{\tilde{U}}\endm_{dis}(\sigma), \phi\ra|\lesssim \la
s\ra^{-\frac{1}{2}+\delta_{1}}
\end{equation}
\end{lemma}
\begin{proof}
This is proved as usual by integration by parts in $t$, and
Duhamel-expanding $ \bm \tilde{U}\\
\overline{\tilde{U}}\endm_{dis}(\sigma)$. Thus we rewrite the
expression as the sum of the terms
\begin{equation}\nonumber
\int_{-\infty}^{\infty}\int_{T}^{\infty}t(\nu-1)(t)\int_{0}^{t}e^{i(t-\lambda)(\xi^{2}+1)}\calF\bm
(1-e^{2i(\Psi-\Psi_{\infty})_{1}(t)-2i(\Psi-\Psi_{\infty})_{1}(\lambda)})\overline{\tilde{U}^{(t)}(\lambda,.)}\phi_{0}^{4}\\
(-1+e^{2i(\Psi-\Psi_{\infty})_{1}(t)-2i(\Psi-\Psi_{\infty})_{1}(\lambda)})\tilde{U}^{(t)}(\lambda,.)\phi_{0}^{4}\endm(\xi)
\overline{\tilde{\calF}\phi(\xi)}d\lambda dt d\xi
\end{equation}
\begin{equation}\nonumber
\int_{-\infty}^{\infty}\int_{T}^{\infty}t(\nu-1)(t)\int_{0}^{t}e^{i(t-\lambda)(\xi^{2}+1)}\calF\bm|\tilde{U}^{(t)}|^{4}(\lambda,.)
\tilde{U}^{(t)}(\lambda,.)\\-|\tilde{U}^{(t)}|^{4}(\lambda,.)
\tilde{U}^{(t)}(\lambda,.)\endm(\xi)\overline{\tilde{\calF}\phi(\xi)}d\lambda
dt d\xi,
\end{equation}
as well as faster decaying local terms which can be handled
similarly to the first term. Let's look at the 2nd term here, the
first being treated as in the proof of lemma~\ref{tedious} by
integrations by parts and further Duhamel expansion. Perform an
integration by parts, replacing this by the sum of suitable
multiples of
\begin{equation}\nonumber
A=\int_{-\infty}^{\infty}T(\nu-1)(T)\int_{0}^{T}e^{i(T-\lambda)(\xi^{2}+1)}\calF\bm|\tilde{U}^{(T)}|^{4}(\lambda,.)
\tilde{U}^{(T)}(\lambda,.)\\-|\tilde{U}^{(T)}|^{4}(\lambda,.)
\tilde{U}^{(T)}(\lambda,.)\endm(\xi)\overline{\frac{\tilde{\calF}\phi(\xi)}{\xi^{2}+1}}d\lambda
dt d\xi,
\end{equation}
\begin{equation}\nonumber
B=\int_{-\infty}^{\infty}\int_{T}^{\infty}\frac{d}{dt}[t(\nu-1)(t)]\int_{0}^{t}e^{i(t-\lambda)(\xi^{2}+1)}\calF\bm|\tilde{U}^{(t)}|^{4}(\lambda,.)
\tilde{U}^{(t)}(\lambda,.)\\-|\tilde{U}^{(t)}|^{4}(\lambda,.)
\tilde{U}^{(t)}(\lambda,.)\endm(\xi)\overline{\frac{\tilde{\calF}\phi(\xi)}{\xi^{2}+1}}d\lambda
dt d\xi,
\end{equation}
\begin{equation}\nonumber\begin{split}
&C=\int_{-\infty}^{\infty}\int_{T}^{\infty}t(\nu-1)(t)\frac{d}{dt}[\lambda_{\infty}(\mu-\mu_{\infty})(t)]\\&\hspace{4cm}\int_{0}^{t}e^{i(t-\lambda)(\xi^{2}+1)}\calF\bm|\tilde{U}^{(t)}|^{4}(\lambda,.)
\partial_{x}\tilde{U}^{(t)}(\lambda,.)\\-|\tilde{U}^{(t)}|^{4}(\lambda,.)
\partial_{x}\tilde{U}^{(t)}(\lambda,.)\endm(\xi)\overline{\frac{\tilde{\calF}\phi(\xi)}{\xi^{2}+1}}d\lambda
dt d\xi,
\end{split}\end{equation}
\begin{equation}\nonumber\begin{split}
&D=\int_{-\infty}^{\infty}\int_{T}^{\infty}t(\nu-1)(t)\frac{d}{dt}[\lambda_{\infty}(\mu-\mu_{\infty})(t)]\\&\hspace{4cm}\int_{0}^{t}e^{i(t-\lambda)(\xi^{2}+1)}\calF\bm\partial_{x}[|\tilde{U}^{(t)}|^{4}](\lambda,.)
\tilde{U}^{(t)}(\lambda,.)\\-\partial_{x}[|\tilde{U}^{(t)}|^{4}](\lambda,.)
\tilde{U}^{(t)}(\lambda,.)\endm(\xi)\overline{\frac{\tilde{\calF}\phi(\xi)}{\xi^{2}+1}}d\lambda
dt d\xi,
\end{split}\end{equation}
\begin{equation}\nonumber
E=\int_{-\infty}^{\infty}\int_{T}^{\infty}t(\nu-1)(t)\frac{d}{dt}[\Psi-\Psi_{\infty}]_{1}(t)\int_{0}^{t}e^{i(t-\lambda)(\xi^{2}+1)}\calF\bm|\tilde{U}^{(t)}|^{4}(\lambda,.)
\tilde{U}^{(t)}(\lambda,.)\\|\tilde{U}^{(t)}|^{4}(\lambda,.)
\tilde{U}^{(t)}(\lambda,.)\endm(\xi)\overline{\frac{\tilde{\calF}\phi(\xi)}{\xi^{2}+1}}d\lambda
dt d\xi,
\end{equation}
\begin{equation}\nonumber
F=\int_{-\infty}^{\infty}\int_{T}^{\infty}t(\nu-1)(t)\calF\bm|\tilde{U}|^{4}(t,.)
\tilde{U}(t,.)\\-|\tilde{U}|^{4}(t,.)
\tilde{U}(t,.)\endm(\xi)\overline{\frac{\tilde{\calF}\phi(\xi)}{\xi^{2}+1}}
dt d\xi,
\end{equation}
Now, for $A$, repeat the proof of SLDE\footnote{SLDE} to bound it
by $\lesssim
T^{\frac{1}{2}+\delta_{1}}T^{-\frac{3}{2}+\delta_{3}}$, better
than what is needed. For $B$, use that
$|\frac{d}{dt}[t(\nu-1)(t)]|\lesssim
t^{-\frac{1}{2}+2\delta_{1}}$, see \eqref{modulationasympto}. $C$
is handled analogously. For $D$, observe that\footnote{Of course,
one should include the $\pm$ subscripts for $\calF$,
$\tilde{\calF}$.}
\begin{equation}\nonumber\begin{split}
&\int_{-\infty}^{\infty}\int_{0}^{t}e^{i(t-\lambda)(\xi^{2}+1)}\calF\bm\partial_{x}[|\tilde{U}^{(t)}|^{4}](\lambda,.)
\tilde{U}^{(t)}(\lambda,.)\\-\partial_{x}[|\tilde{U}^{(t)}|^{4}](\lambda,.)
\tilde{U}^{(t)}(\lambda,.)\endm(\xi)\overline{\frac{\tilde{\calF}\phi(\xi)}{\xi^{2}+1}}d\lambda
d\xi\\&= \la
\int_{0}^{t}e^{i(t-\lambda)\calH}\bm\partial_{x}[|\tilde{U}^{(t)}|^{4}](\lambda,.)
\tilde{U}^{(t)}(\lambda,.)\\-\partial_{x}[|\tilde{U}^{(t)}|^{4}](\lambda,.)
\tilde{U}^{(t)}(\lambda,.)\endm_{dis} d\lambda,\calH^{-1}\phi\ra
\end{split}\end{equation}
Proceeding as in the proof of SLDE, i. e. using
\eqref{keyidentity} as well as pseudo-conformal almost
conservation, we can bound the preceding expression by $\lesssim
\la t\ra^{-\frac{3}{2}+\delta_{3}}$. From here one proceeds as for
$C$ etc. Finally, $F$ is more elementary, as we have
\begin{equation}\nonumber
|\int_{-\infty}^{\infty}\calF\bm|\tilde{U}|^{4}(t,.)
\tilde{U}(t,.)\\-|\tilde{U}|^{4}(t,.)
\tilde{U}(t,.)\endm(\xi)\overline{\frac{\tilde{\calF}\phi(\xi)}{\xi^{2}+1}}d\xi|
=|\la \bm|\tilde{U}|^{4}(t,.)
\tilde{U}(t,.)\\-|\tilde{U}|^{4}(t,.) \tilde{U}(t,.)\endm_{dis},
\calH^{-1}\phi\ra|\lesssim t^{-\frac{9}{2}+\epsilon(\delta_{2})}
\end{equation}
\end{proof}
It is now straightforward to retrieve the desired bound for
$\beta_{1}(s)$ via \eqref{beta}, and from here one easily infers
the desired bound for $\nu_{1}(s)=\nu(s)-1$, via \eqref{nu}. Using
\eqref{ainfty}, \eqref{binfty}, one easily infers the bounds for
$a_{\infty},\,b_{\infty}$ outlined in \eqref{modulationasympto}.
\\

Next, consider $\omega$ satisfying \eqref{omega}. From what we
have established so far we can retrieve the bound
\begin{equation}\nonumber
B(s)^{-1}=c\lambda_{\infty}^{-1}(s)+O(s^{-\frac{3}{2}+\delta_{1}})
\end{equation}
From here, \eqref{omega} and Lemma~\ref{Technical} one infers the
existence of $c_{\infty}$ such that
\begin{equation}\nonumber
|\omega-c_{\infty}\lambda_{\infty}^{-1}(s)|\lesssim
s^{-\frac{3}{2}+\delta_{1}}
\end{equation}
Moving on to $\mu$ satisfying \eqref{mu}, we deduce the existence
of $v_{\infty}, y_{\infty}$ such that
\begin{equation}\nonumber
|\mu(s)-\frac{2v_{\infty}s+y_{\infty}}{a_{\infty}+b_{\infty}s}|\lesssim
s^{-\frac{3}{2}+\delta_{1}},
\end{equation}
and furthermore we have
$c_{\infty}=v_{\infty}a_{\infty}-\frac{b_{\infty}y_{\infty}}{2}$.
Finally, one easily deduces the bound on $\gamma(s)-s$ from
\eqref{gamma}. \\
To get the estimates specified in \eqref{modulationasympto} on the
derivatives $\dot{\nu}(s)$ etc., one simply differentiates the
relations \eqref{beta} etc, and uses the assumptions. This is
elementary and hence omitted. Finally we have completed the proof
of Theorem~\ref{coreredux}, up to the continuity assertion. To
show the latter, {\it{note that we are working on the finite time
interval $[0,T)$}}, whence continuity can be derived by using very
crude bounds from linear theory.
\end{proof}

\bibliographystyle{amsplain}

\noindent

\end{document}